\documentclass[12pt,openany,oneside]{book}
\usepackage[pdftex]{graphicx}
\usepackage[a4paper,top=2.5cm,bottom=2.5cm,left=3cm,right=3cm]{geometry} 
\usepackage[T1]{fontenc}
\usepackage{psfrag}
\usepackage[english]{babel}
\usepackage{graphicx}
\usepackage{float}
\usepackage{caption}
\usepackage{mathtools}
\usepackage{ae}
\usepackage{longtable, lscape}
\usepackage{subcaption}
\usepackage{amsmath}
\usepackage{amssymb}
\usepackage{adjustbox}
\usepackage{mathdots}
\usepackage{bm}
\usepackage{siunitx}
\usepackage{placeins}
\usepackage{titlesec}
\usepackage{sectsty}
\usepackage{titletoc,  fmtcount}
\usepackage{titlecaps}
\usepackage{hyperref}
\usepackage{booktabs}
\usepackage{multicol}

\addto\captionsenglish{} 
\titlecontents{chapter}
[0pt]
{}
{\bfseries\chaptername\ \thecontentslabel\quad}
{}
{\bfseries\hfill\contentspage}
\titleformat{\chapter}[display]
{\normalfont\huge\bfseries}{\chaptertitlename\ \thechapter}{10pt}{\LARGE}
\titlespacing*{\chapter}{0pt}{0pt}{20pt}
\chaptertitlefont{\fontsize{22pt}{30pt}\selectfont}

\hypersetup{ 
	colorlinks,
	citecolor=black,
	filecolor=black,
	linkcolor=black,
	urlcolor=black
}
\newlength{\defbaselineskip}
\newcommand{\setlinespacing}[1]%
{\setlength{\baselineskip}{#1 \defbaselineskip}}

\newcommand{\rk}{\operatorname{rank}}
\newcommand{\iu}{\mathrm{i}}
\newcommand{\glt}{\mathrm{GLT}}
\DeclarePairedDelimiter{\norm}{\lVert}{\rVert}
\DeclarePairedDelimiter{\abs}{\lvert}{\rvert}

\setlinespacing{1.5}

\newtheorem{definition}{Definition}
\newtheorem{lemma}{Lemma}[chapter]
\newtheorem{proposition}{Proposition}[chapter]
\newtheorem{theorem}{Theorem}[chapter]
\newtheorem{remark}{Remark}[chapter]

\newtheorem{conjecture}{Conjecture}[section]
\newtheorem{corollary}{Corollary}[chapter]
\newtheorem{gltax}{GLT}
\def\abs#1{\left|#1\right|} 

\def\softd{{\leavevmode\setbox1=\hbox{d}%
		\hbox to 1.05\wd1{d\kern-0.4ex{\char039}\hss}}}
\def \R{\mathbb R}
\def \N{\mathbb N}
\def \C{\mathbb C}

\newenvironment{proof}%
{\medskip\par\noindent{\bf Proof:}}%
{\hspace{\fill}$\Box$\medskip\par}

\numberwithin{equation}{chapter}
\numberwithin{theorem}{chapter}
\numberwithin{lemma}{chapter}
\numberwithin{proposition}{chapter}
\numberwithin{corollary}{chapter}
\numberwithin{example}{chapter}
\numberwithin{remark}{chapter}
\numberwithin{definition}{chapter}

\begin{document}
	\pagestyle{plain} 
	\pagenumbering{gobble} 
	\begin{titlepage}
\begin{figure}[!htb]
    \centering
    \includegraphics[width=5cm]{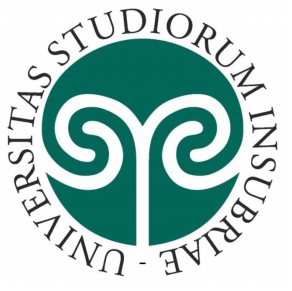}
\end{figure}

\begin{center}
    \Large{\textbf{UNIVERSITÀ DEGLI STUDI DELL'INSUBRIA}}
    \vspace{3mm}
    \\ \normalsize{Dipartimento di Scienza e Alta Tecnologia (DiSAT)}
\end{center}

\vspace{10mm}
\begin{center}
    \LARGE{\textbf{Inverse problems with fractional operators: an analytical study and its numerical linear algebra treatment via the GLT theory}}
\end{center}

\vspace*{\fill}

\begin{minipage}[t]{1\textwidth}
    \begin{multicols}{2}
    	{\normalsize{\textbf{Advisors:}}{\normalsize\vspace{1mm}
    			\\ \normalsize{ Prof. Sean (Yu Sing) Hon}
        \\ \normalsize{Prof. Stefano Serra-Capizzano}}} 

        
         \columnbreak

         \begin{flushright}
            {\normalsize{\textbf{Presented by:}}{\normalsize\vspace{1mm}
            \\ \normalsize{Asim Ilyas}\\
            \normalsize{755725}}} \\
         \end{flushright}
    \end{multicols}
\end{minipage}

\begin{center}
    {\normalsize{\textbf{Academic year: 2025}}}
\end{center}

\end{titlepage}
	\chapter*{Dedication}
To my beloved parents, for their unconditional love, sacrifices, and prayers.
 \vspace*{0.45cm}
 
To my precious son, Muhammad Salim Asim, the light of my life. 
\vspace*{0.45cm}

To my devoted wife, for her patience, love, and constant encouragement.
 \vspace*{0.45cm}
 
To my dear brother and sister, for their unwavering support. 
\vspace*{0.45cm}

To my esteemed supervisors, Prof. Stefano Serra-Capizzano and Prof. Sean (Yu Sing) Hon, for 
their invaluable guidance and mentorship. 
	\chapter*{Acknowledgments}
In the name of {\textbf{ALLAH}}, the Most Gracious, the Most Merciful. All praise is due to Allah Almighty 
for blessing me with strength, perseverance, and wisdom to accomplish this work. 
\smallskip

\noindent I extend my deepest gratitude to my supervisors, Prof. Stefano Serra-Capizzano and Prof. Sean 
(Yu Sing) Hon, for their exceptional guidance, patience, and invaluable feedback. Their 
expertise and encouragement were instrumental in shaping this research. 
\smallskip

\noindent I am immensely thankful to my beloved parents for their unconditional love and prayers, which 
have been my bedrock. I deeply appreciate my wife for her patience and unwavering support, 
and my precious son, Muhammad Salim Asim, for being my constant motivation. I also thank 
my brother, sister, and all family members for their steadfast belief in me. 
\smallskip

\noindent I am grateful to my friends and colleagues for their encouragement and camaraderie, and to the 
reviewers for their thorough evaluation and constructive feedback, which enhanced this thesis. 
Finally, I thank everyone who contributed to this work. May Allah bless them all abundantly. 

	\chapter*{Abstract}
Fractional derivatives have emerged as a powerful mathematical tool for modeling complex 
physical phenomena exhibiting memory effects and nonlocal interactions, which cannot be 
adequately captured by classical integer-order models. Fractional differential equations (FDEs) 
incorporating these derivatives provide a more accurate framework for describing anomalous 
diffusion, transport in heterogeneous media, and other processes with hereditary dynamics. 
However, the nonlocal nature that enhances their modeling capability also introduces 
substantial challenges, particularly in inverse problems where unknown sources, coefficients, 
or initial conditions must be recovered from indirect and often noisy measurements. 
\smallskip

\noindent Inverse problems for FDEs are typically ill-posed in the Hadamard sense, necessitating 
regularization techniques such as the quasi-boundary value method to restore stability and 
ensure meaningful reconstructions. The numerical solution of these problems relies on 
discretization techniques that produce large-scale linear systems whose coefficient matrices 
inherit the nonlocal structure of the fractional operators. Understanding the asymptotic spectral 
properties of these matrix sequences is essential for designing efficient solvers and 
preconditioning strategies.
\smallskip

\noindent The theory of Generalized Locally Toeplitz (GLT) sequences provides a powerful and unified 
framework for this purpose, characterizing the spectral distribution of large structured matrices 
arising from discretizations and offering rigorous guidelines for preconditioner design. This 
thesis leverages GLT theory to analyze matrix sequences from inverse problems of FDEs, 
develop structure-preserving preconditioners, and investigate fundamental questions 
concerning nonlinear operations within the GLT algebra. In particular, it establishes maximal 
distribution results for geometric means of GLT sequences, proving that when the associated 
symbols commute, the geometric mean sequence remains GLT without requiring symbol 
invertibility, thereby resolving a long-standing conjecture while demonstrating that the result 
is maximal in the noncommuting degenerate case. 
\smallskip

\noindent Collectively, this thesis contributes to both the theory and numerical solution of inverse 
problems for FDEs, integrating regularization theory, spectral analysis, and preconditioning 
into a coherent computational framework. The theoretical advancements concerning GLT 
sequences and their geometric means further enrich the mathematical foundations of structured 
matrix analysis, providing practical strategies for solving inverse problems in complex physical 
systems governed by nonlocal dynamics.

	\tableofcontents  
	\newpage 
	\pagenumbering{arabic}
	\setcounter{chapter}{0} 
	\listoffigures
	\addcontentsline{toc}{chapter}{List of Figures}
	\listoftables
	\addcontentsline{toc}{chapter}{List of Tables}
	\chapter*{General notations and acronyms}

\begin{itemize}
	\item $\mathbb{N}$, $\mathbb{R}$, $\mathbb{C}$: Sets of natural, real, and complex numbers, respectively.
	\item $\mathbb{R}^+$, $\mathbb{N}^+$: Sets of positive real and natural numbers.
	\item $\mathbb{Z}$: Set of integers.
	\item $\mathbb{K}$: A field, usually $\mathbb{R}$ or $\mathbb{C}$.
	\item $\mathbb{K}^r$: Space of vectors of length $r$ with entries in $\mathbb{K}$.
	\item $\mathbb{K}^{r \times s}$: Space of $r \times s$ matrices with entries in $\mathbb{K}$.
	\item  $\mathbb{K}^{n\times n}$: Space of square \(n\times n\) matrices over \(\mathbb{K}\).
	\item $\bm{x}, \bm{y}$: Vectors in $\mathbb{K}^r$ (bold letters usually denote vectors).
	\item $A, B$: Matrices in $\mathbb{K}^{n\times n}$ unless specified otherwise.
	\item $\mathbb{I}_n$ or $I_n$: $n \times n$ identity matrix.
	\item $\mathbf{0}_n$ or $0_n$: $n \times n$ zero matrix.
	\item $(\cdot)^{\top}$: denotes the transpose of a vector or matrix.
	\item $(\cdot)^*$: denotes the conjugate transpose of a vector or matrix.
	\item $(\cdot)^\dagger$: denotes the Moore-Penrose pseudoinverse of a square matrix.
	\item $\lambda_j(A)$: denotes the $j$-th eigenvalue of $A \in \mathbb{K}^{n\times n}$, for $j = 1, \dots, n$.
	\item $\sigma_j(A)$: denotes the $j$-th singular value of $A \in \mathbb{K}^{n\times n}$, for $j = 1, \dots, n$.
	\item $\mathcal{H}$: Hilbert space.
	\item $\sigma_x, \sigma_y, \sigma_z$: Pauli matrices.
	\item $\sigma^{(i)}$: Pauli operator acting on the $i$-th spin.
	\item $H$: Hamiltonian operator of the system.
	\item $\rho(H)$: empirical spectral measure of $H$.
	\item $\chi$: Magnetic susceptibility.
	\item $\beta$: Inverse temperature parameter.
	\item $Z$: Partition function.
	\item $\langle \cdot \rangle$: Expectation with respect to the Gibbs measure.
	\item $A \geq B$, where $A, B \in \mathbb{K}^{n\times n}$, means that $A - B$ is positive semidefinite. Likewise, $A \leq B$ means that $A - B$ is negative semidefinite.
	\item $A > B$, where $A, B \in\mathbb{K}^{n\times n}$, means that $A - B$ is positive definite. Likewise, $A < B$ means that $A - B$ is negative definite.
	\item $A \sim B$, where $A, B \in \mathbb{K}^{n\times n}$, denotes that $A$ and $B$ are similar matrices (i.e., there exists an invertible $P$ such that $A = PBP^{-1}$).
	
	\item $\| \cdot \|_p$: denotes the standard $p$-norm of a vector or matrix, $1 \leq p \leq \infty$.
	\item $\| \cdot \|$: denotes the spectral norm (largest singular value), coincides with $\| \cdot \|_{\infty}$. If the matrix is normal, the spectral norm coincides with the spectral radius.
	
	\item $\otimes$: Kronecker (tensor) product. For $A \in \mathbb{C}^{p \times q}$, $B \in \mathbb{C}^{r \times s}$,
	\[
	A \otimes B =
	\begin{bmatrix}
		a_{11} B & \cdots & a_{1q} B \\
		\vdots   &        & \vdots   \\
		a_{p1} B & \cdots & a_{pq} B
	\end{bmatrix}
	\in \mathbb{C}^{pr \times qs}
	\]
	\item $L^p(D)$: Space of $p$-integrable functions on domain $D$.
	\item $C([a,b])$: Space of continuous functions on $[a,b]$.
	\item $\mathcal{O}(\cdot)$, $o(\cdot)$: Big-$O$ and little-$o$ notation.
	\item $\mathcal{C}_c(\mathbb{R})$, $\mathcal{C}_c(\mathbb{C})$: The set of continuous, complex-valued functions with compact support on $\mathbb{R}$ or $\mathbb{C}$.
	\item \(\partial D\): boundary of a set \(D\subseteq\mathbb{R}^d\);
	in particular, \(\partial\Omega\) denotes the boundary of the spatial domain \(\Omega\).
\end{itemize}
\textbf{Multi-index notation:}
\begin{itemize}
	\item $\bm{i}$ denotes a multi–index $(i_1, \ldots, i_d) \in \mathbb{Z}^d$. Its size will be clear from context.
	\item $\bm{0}, \bm{1}, \bm{2}, \ldots$ are vectors of all zeros, ones, twos, etc.
	\item $\bm{h} + \bm{k}$ is the multi–index whose components are $h_r + k_r$ for $r = 1, \ldots, d$. More generally, operations between multi-indices such as subtraction, multiplication, division, and so on, are performed componentwise.
	\item $\bm{h} \leq \bm{k}$ means that $h_r \leq k_r$ for all $r = 1, \ldots, d$. More generally, relations between multi-indices are evaluated componentwise.
	\item $\bm{h}, \ldots, \bm{k}$ denotes the multi-index range $\{\bm{\ell} \in \mathbb{Z}^d : \bm{h} \leq \bm{\ell} \leq \bm{k}\}$. With $\bm{\ell} = \bm{h}, \ldots, \bm{k}$ we mean that $\bm{\ell}$ varies in this range from $\bm{h}$ to $\bm{k}$ according to the standard lexicographic ordering:
	\[
	\left[
	\cdots
	\left[
	\left[
	(\ell_{1},\ldots,\ell_{d})\right]_{\ell_d = h_d, \ldots, k_d}
	\right]_{\ell_{d-1} = h_{d-1}, \ldots, k_{d-1}}
	\cdots
	\right]_{\ell_1 = h_1, \ldots, k_1}
	.\]
	\item $\bm{m} \to \infty$ means that $\min(\bm{m}) = \min_{r=1, \ldots, d} \bm{m}_r \to \infty$, or equivalently, that all the components tend to infinity.
	\item $N(\bm{m})$ is the product of all the components of $\bm{m}$, i.e., $N(\bm{m}) := \prod_{r=1}^d m_r$.
\end{itemize}
\begin{itemize}
	\item $\{A_n\}_n$: Sequence of matrices indexed by size/discretization parameter $n$.
	\item $T_n(f)$: Toeplitz matrix of size $n$ generated by symbol $f$.
	\item $C_n(f)$: Circulant matrix of size $n$ generated by symbol $f$.
	\item $H_n(f)$: Hankel matrix of size $n$ generated by symbol $f$.
	\item Symbol $(\kappa)$: Measurable function associated with a matrix-sequence, e.g., $\{A_n\}_n \sim_{\mathrm{GLT}} \kappa$.
	\item Symbol $(\xi)$: Measurable function associated with a matrix-sequence, e.g., $\{B_n\}_n \sim_{\mathrm{GLT}} \xi$.
	\item $D = [0,1]^d \times [-\pi, \pi]^d$: Typical domain for the symbol of a $d$-level GLT sequence.
	\item $\sim_\sigma$: Singular value distribution.
	\item $\sim_\lambda$: Eigenvalue distribution.
	\item $\sim_{\mathrm{GLT}}$: Matrix-sequence is a GLT sequence with symbol.
	\item $\operatorname{ran} (M)$: The column space (range) of a matrix $M$.
	\item \(\operatorname{ess\,Ran}(g)\): the essential range of a measurable
	function \(g:D\subseteq\mathbb{R}^{t}\to Y\), where \(Y\) is \(\mathbb{C}\),
	\(\mathbb{C}^{r}\), \(\mathbb{C}^{r\times r}\), or the Grassmannian of
	subspaces of \(\mathbb{C}^{r}\), i.e.
	\[
	\operatorname{ess\,Ran}(g)
	= \bigl\{y\in Y : \mu_{t}\bigl(\{x\in D : g(x)\in B_{Y}(y,\epsilon)\}\bigr) > 0
	\ \ \forall \epsilon > 0\bigr\},
	\]
	where \(\mu_{t}\) is the Lebesgue measure on \(D\) and \(B_{Y}(y,\epsilon)\)
	is the ball of radius \(\epsilon\) centered at \(y\) in \(Y\) (with respect
	to the natural metric on \(Y\)).
\end{itemize}

\textbf{Acronyms and Abbreviations:}
\begin{itemize}
	\item GLT: Generalized Locally Toeplitz
	\item LT: Locally Toeplitz
	\item HPD: Hermitian positive definite
	\item SPD: Symmetric positive definite
	\item a.c.s.: Approximating class of sequences
	\item FC: Fractional Calculus
	\item R-LFI: Riemann–Liouville Fractional Integral
	\item R-LFD: Riemann–Liouville Fractional Derivative
	\item CFD: Caputo Fractional Derivative
	\item FDE : Fractional Differential/Diffusion Equation
	\item SDE: Standard Diffusion Equation
	\item IP: Inverse Problem
	\item ISP: Inverse Source Problem
	\item MLF: Mittag-Leffler Function
	\item WMT: Weierstrass M-test
	\item a.e.: almost everywhere
	\item  LPO : Linear positive operator 
	\item FFT : Fast Fourier Transform
\end{itemize}

	\addcontentsline{toc}{chapter}{General notations and acronyms}
	\chapter*{Introduction}
Fractional calculus (FC) provides a natural and powerful extension of classical
differential calculus by allowing derivatives and integrals of non-integer
order \cite{Oldham1974,Podlubny1998}. Although fractional derivatives already appeared in the early development of calculus, for a long time they remained mainly of theoretical interest. In recent decades, however, fractional differential equations (FDEs) have attracted growing attention from both pure and applied mathematics due to their remarkable ability to model complex phenomena that cannot be adequately described by classical integer-order models \cite{Kilbas2006}. The key feature underlying their success is the intrinsic nonlocality of fractional operators, which introduces memory effects in time and long-range interactions in space, allowing a more faithful representation of many real-world processes. FDEs play a central role in the modeling of
anomalous diffusion, a class of transport phenomena characterized by
deviations from standard Brownian motion \cite{Metzler2000}. Such deviations arise naturally in heterogeneous and porous media, plasma physics, viscoelastic materials, image reconstruction, finance, and wave propagation in biological tissues \cite{Hilfer2000,Magin2010}. Depending on the physical mechanism involved, fractional derivatives may act in time, space, or both variables. Time-fractional models typically describe non-Markovian dynamics with memory effects, while space-fractional models
capture heavy-tailed jump distributions and long-range spatial interactions \cite{Klafter2011,Mainardi2022}. Widely used operators include the Riemann--Liouville, Caputo, and Gr\"unwald--Letnikov derivatives, as well as Riesz operators and the fractional Laplacian, together with their tempered variants, which introduce exponential damping to account for finite-range effects \cite{Meerschaert2014,Samko1993}. 

Despite their modeling flexibility and expressive power, FDEs pose significant analytical and computational challenges. Closed-form solutions are rarely available and, when they exist, they are often unsuitable for practical computations. As a consequence,
numerical methods play a fundamental role in the study and simulation of
fractional models. Many numerical methods have been developed to efficiently
solve FDEs, including finite difference methods
\cite{Thomas2013}, finite element methods
\cite{Dhatt2012}, finite volume methods
\cite{Eymard2000}, and spectral methods
\cite{Gottlieb1977}. Although these approaches are capable of producing accurate approximations, achieving high accuracy typically requires very fine spatial and temporal discretizations. As a consequence, the resulting algebraic systems become
large and computationally demanding, especially in multi-dimensional
settings. A major source of difficulty lies in the intrinsic nonlocal nature
of (fractional) differential operators. Unlike classical local operators, their
discretization leads to dense coefficient matrices, resulting in high memory
requirements and increased computational cost.

In many applications, fractional models arise not only in direct problems,
where the governing equations and parameters are known, but also in inverse
problems (IPs). Typical IPs include the identification of unknown
sources, coefficients, initial conditions, boundary data, or even the
fractional order itself from indirect and often noisy measurements \cite{Jin2015}.
IPs are particularly challenging because they combine the
intrinsic ill-posedness of IPs with the nonlocal behavior of
fractional operators. From a mathematical standpoint, these problems are
often ill-posed in the sense of Hadamard, meaning that solutions may not
exist, may not be unique, or may depend discontinuously on the data. In
practice, even small perturbations in the measurements can produce large
errors in the reconstructions \cite{Engl1996}. To address the challenges arising from ill-posed IPs in FDEs, various regularization techniques have
been developed to restore stability and ensure meaningful reconstructions.
Motivated by their wide-ranging practical applications, IPs for FDEs have been extensively studied using diverse analytical and computational approaches \cite{IlyasMalik2022,IlyasSaif2023,SuhaibIlyas2023,IlyasMalik2024CAM,IlyasIqbal2024,IlyasKhalid2024,IlyasMalik2025,Jin2015,Samreen2025}. In general, regularization methods incorporate additional information or constraints on the solution in order to control instability and suppress the amplification of noise. Among the most widely used approaches, Tikhonov regularization and its variants, including simplified and generalized formulations, stabilize the solution by introducing penalty terms that control smoothness or enforce \emph{a priori} assumptions on the unknown solution \cite{Li2018,Nguyen2016}. Another effective approach is the quasi-boundary value method, also known as the nonlocal boundary value method, which has been successfully applied to a variety of ill-posed IPs, particularly for parabolic equations and fractional diffusion models \cite{Denche2005,Hao2008}. In FDEs, regularization methods must ensure stability and accuracy while accounting for the nonlocal nature of the operators.

The numerical treatment of both direct and IPs for FDEs ultimately leads to the solution of large linear systems
of the form
\[
A_n x_n = b_n,
\]
where $A_n \in \mathbb{C}^{d_n \times d_n}$ is a square matrix and the index
$n \in \mathbb{N}$ is related to the fineness of the discretization mesh. As the
grid is refined, the dimension $d_n$ increases and typically satisfies
$d_n \to +\infty$ as $n \to +\infty$, so that one naturally deals with a
\emph{sequence} of matrices $\{A_n\}_n$ of growing size rather than a single
matrix. The collective spectral behavior of $\{A_n\}_n$ is therefore a key
object of interest, since it reveals important numerical features and
potential difficulties, such as severe ill-conditioning and slow convergence
of iterative methods \cite{extreme3,extreme2,extreme4}. In this matrix-sequence
setting, the relevant spectral information is not contained in any individual
matrix, but rather in the way the spectra of $\{A_n\}_n$ evolve as the
discretization is refined. This naturally leads to the notion of asymptotic
eigenvalue and singular value distributions in the Weyl sense, which capture
the limiting spectral content of the entire sequence as $n \to +\infty$. For a
wide class of discretizations, this limiting behavior can be described by a
measurable function, referred to as the \emph{symbol}, which reflects the
combined influence of the underlying (fractional) differential operator, its
coefficients, the computational domain, and the discretization scheme. This
symbol-based spectral description plays a central role in numerical analysis,
as it explains the observed conditioning of the linear systems and provides a
rigorous framework for the design and analysis of effective preconditioners
and scalable iterative methods.

Structured matrices play a central role in this analysis, because in many
relevant settings the discretization of FDEs produces coefficient matrices
with Toeplitz or Toeplitz-like patterns that can be exploited to build fast
solvers and preconditioners. In particular, when FDEs have constant
coefficients and are discretized on uniform grids, the resulting matrices
often exhibit a Toeplitz or multilevel Toeplitz structure
\cite{BiniPan1994,BoettcherSilbermann2006,ChanNg1996,Ng2004}. Such matrices
arise from the shift-invariant nature of the underlying operators and possess
a rich spectral theory developed over the last century. Moreover, Toeplitz
matrices admit fast matrix--vector products and efficient algorithms based on
the fast Fourier transform, and their eigenvalue and singular value
distributions are directly linked to their generating functions
\cite{BoettcherSilbermann2006}. In more general settings, such as variable
coefficients, nonuniform grids, or non-Cartesian domains, the classical
Toeplitz structure is typically lost. Nevertheless, many discretization
matrices arising in practice belong to the broader class of GLT matrix sequences
\cite{SerraCapizzano2003,garoni2017,garoni2018}. The GLT theory provides a
unifying framework for describing the asymptotic spectral behavior of large
structured and partially structured matrices. GLT sequences form a rich
$*$-algebra that includes Toeplitz sequences, diagonal sampling matrices, and
zero-distributed sequences as fundamental components. Each GLT sequence is
uniquely associated with a symbol, and algebraic operations on matrix
sequences correspond, almost everywhere, to pointwise operations on their
symbols \cite{garoni2017,garoni2018}. This framework supports the design and
analysis of numerical solvers, preconditioners, and multigrid methods, since
the convergence of Krylov subspace and multigrid schemes depends strongly on
conditioning and eigenvalue clustering, which can be predicted and optimized
via GLT symbols \cite{garoni2017,ChanNg1996,Ng2004}. Recent developments in symbol-based preconditioning for evolutionary PDEs include the block Toeplitz preconditioner for wave equations introduced by Hon and Serra-Capizzano \cite{HonSerra2023}, and the sine transform based preconditioned MINRES method for all-at-once systems developed by Hon, Fung, Dong, and Serra-Capizzano \cite{HonFungDong2024}. These works highlight the effectiveness of combining structural properties with fast transforms, and the present thesis extends these ideas to IPs and variable-coefficient FDEs.

A particularly rich connection within the theory of GLT matrix sequences
concerns the notion of geometric means of matrices. The geometric mean of
Hermitian positive definite matrices was rigorously formalized by Ando, Li,
and Mathias and has since become a fundamental concept in matrix analysis and
operator theory. Given two Hermitian positive definite matrices $A$ and $B$,
their geometric mean is defined by
\[
G(A,B)
=
A^{1/2}\bigl(A^{-1/2} B A^{-1/2}\bigr)^{1/2} A^{1/2}
=
G(B,A),
\]
a definition that preserves symmetry, positive definiteness, and key
variational properties \cite{Ando2004,bhatia2007}. Geometric means naturally
arise in preconditioning strategies, operator interpolation, and the modeling
of composite physical systems, where nonlinear matrix operations are
unavoidable. When applied to sequences of Hermitian positive definite
matrices, the geometric mean interacts in a subtle way with asymptotic
spectral properties. Within the GLT framework, this interaction raises
fundamental questions concerning symbol calculus, commutativity of
matrix-valued symbols, degeneracy, and maximality of distribution results.

\section*{Structure of the thesis}
This thesis is structured as follows.

\noindent
\textbf{Chapter \ref{chapter:Fractional Calculus}} introduces the analytical foundations required throughout the
thesis. We review the main definitions of fractional calculus employed in our
work, including the Riemann--Liouville, Caputo, Gr\"unwald--Letnikov, and Riesz
fractional operators, together with their relationships under suitable
regularity assumptions. Generalized formulations in the form of tempered
fractional derivatives are also discussed, due to their relevance in modeling
finite-range diffusion processes. The chapter further presents Mittag--Leffler
and generalized Mittag--Leffler type functions, along with their key analytical properties used in the study of FDEs and IPs.

In \textbf{Chapter \ref{chapter:Structured Matrix Sequences}}, we introduce the foundational concepts related to matrix
sequences that are essential for the subsequent spectral analysis. The chapter
focuses on singular value and eigenvalue distributions in the asymptotic sense,
approximation classes of matrix sequences, and multi-level structures. These
topics provide the mathematical background required for understanding the
spectral behavior of large structured matrices arising from discretized
differential operators.

In \textbf{Chapter \ref{chap:Structured matrices}}, we review several classes of classical structured matrices that
play a central role in numerical linear algebra. In particular, we discuss
Toeplitz, circulant, $\omega$-circulant, $\tau$, and Hankel matrices, presenting
their definitions together with their main algebraic and spectral properties.
The importance of these matrix classes in both numerical computations and
asymptotic linear algebra is highlighted, with particular emphasis on their
spectral features and computational advantages.

\textbf{Chapter \ref{GLT}} is devoted to the development of the theory of GLT matrix sequences, which constitutes the core theoretical framework for the spectral analysis carried out in this thesis. We introduce the notion of GLT sequences and present the structure of GLT algebras together with their fundamental properties. The chapter establishes the axioms and algebraic rules governing GLT sequences, distinguishing between matrix
sequences with explicit and hidden asymptotic structures. Key concepts such as
zero-distributed sequences, block Toeplitz and block diagonal sampling matrices,
and the construction of $d$-level, $r$-block GLT $\ast$-algebras are also
discussed in detail.

\textbf{Chapter \ref{chapter:TFDE}} is devoted to the study of IPs for a two-dimensional generalized diffusion equation involving memory-dependent time-fractional derivatives. The chapter considers the identification of both a space-dependent source term and a time-dependent coefficient from additional measurement data. The IPs are shown to be ill-posed in the sense of Hadamard, and existence and uniqueness of solutions are established under suitable over-specified and regularity assumptions. Special cases corresponding to particular choices of the memory kernel are analyzed in detail.

In \textbf{Chapter \ref{chapter:IFL}}, we investigate an IP for a space-time FDE with variable coefficients. The ill-posedness of the problem is addressed through a quasi-boundary value regularization. A finite difference discretization yields large two-by-two block linear system. Its spectral behavior is analyzed within the framework of GLT, which provides guidance for the construction of an effective Krylov subspace preconditioner. Numerical experiments demonstrate the accuracy of the reconstruction and the efficiency of the proposed preconditioning strategy.

\textbf{Chapter \ref{chapter:RFDE}} investigates an IP for a two-dimensional space-time FDE involving a tempered time-fractional derivative and Caputo fractional derivatives in space. The ill-posed problem is regularized by a quasi-boundary value method. A finite difference discretization yields large two-by-two block structured linear systems with multilevel Toeplitz-like structure, for which we develop an efficient preconditioning strategy combined with GMRES. A comprehensive GLT analysis of the coefficient matrix sequences guides the construction of the preconditioner. Numerical experiments confirm the stability, accuracy, and computational efficiency of the proposed approach in reconstructing the source term from noisy data.

\textbf{Chapter \ref{chapter:space-FDE}} studies the spectral properties and preconditioning of linear systems arising from the finite difference discretization of a time-dependent space-FDE with variable coefficient $a(x)$ on $\Omega \subset \mathbb{R}^{d}$, $d=1,2$. The equation is discretized by the shifted Grünwald formula in space and the $\theta$-method in time, leading to $(d+1)$-level Toeplitz-like all-at-once systems. Using GLT theory, the spectral distribution of the corresponding matrix sequences and of suitable preconditioned systems is analyzed, fully covering the case of nonconstant $a$. Numerical results confirm the theoretical findings.

\textbf{Chapter \ref{GM GLT2}} is devoted to the analysis of the spectral distribution of geometric means of Hermitian positive definite GLT matrix sequences to the maximal setting, allowing for degenerate symbols and non-commuting cases. It is shown that, when the associated GLT symbols commute, the geometric mean is again a GLT sequence, with symbol given by the geometric mean of the original symbols, without requiring invertibility almost everywhere. This result completes the proof of a long-standing conjecture in both the scalar and block settings, for arbitrary dimensions. The chapter also presents detailed numerical experiments illustrating that, in the presence of non-commuting symbols and simultaneous degeneracy, the classical relation may fail. These experiments reveal new spectral phenomena and suggest a natural candidate symbol for the geometric mean in such generalized situations.

\subsection*{Main Contributions of the Thesis}
The main contributions of this thesis are theoretical, analytical, and computational, and can be summarized as follows.
\begin{enumerate}
	\item \textbf{Inverse source problems for generalized diffusion equations.} 
	In Chapter~\ref{chapter:TFDE}, we formulate and analyze two ISPs for a two-dimensional generalized diffusion equation involving a Caputo-type time-fractional derivative with a general memory kernel and nonlocal boundary conditions. Using the spectral decomposition of the associated elliptic operator, we derive explicit series representations of the solutions in terms of bi-orthogonal systems of eigenfunctions. We establish rigorous existence and uniqueness results for regular solutions and demonstrate that both ISPs are ill-posed in the Hadamard sense. Several special cases, including tempered, multi-term, mono-term, and classical diffusion models, are obtained as particular choices of the memory kernel, illustrating the unifying nature of the generalized formulation.
	\item \textbf{Krylov preconditioning for inverse problems in space-time fractional diffusion equations.} 
	In Chapter~\ref{chapter:IFL}, we address an ISP for a space-time FDE involving a tempered fractional derivative in the time variable and a Caputo fractional derivative in the space variable. The ill-posed problem is regularized using the quasi-boundary value method. The considered equation is of great interest in applications, and a finite difference approximation of the resulting regularized problem leads to a two-by-two block linear system of large dimensions. A detailed spectral analysis of the associated matrix sequences is conducted using the GLT framework. Based on this analysis, we develop an efficient block triangular preconditioner for the GMRES method and prove that the preconditioned matrix sequence is weakly clustered at one in the GLT, singular value, and eigenvalue sense.
	\item \textbf{Efficient Krylov solvers for inverse source problems in two-dimensional space-time fractional diffusion equations.} 
	In Chapter~\ref{chapter:RFDE}, we extend the previous work to the two-dimensional setting, developing a fast algorithm to reconstruct the source term in a space-time FDE with variable coefficients. The finite difference discretization of the regularized problem leads to a large-scale all-at-once system with multilevel Toeplitz-like structure. We conduct a comprehensive GLT analysis of the coefficient matrix sequences and construct an effective preconditioner. Numerical experiments demonstrate the accuracy of the reconstruction and the efficiency of the proposed preconditioning strategy across various fractional orders and problem sizes. The convergence behavior depends mainly on the spatial discretization, with limited sensitivity to temporal refinement, confirming that the number of outliers depends only on the spatial dimension.
	\item \textbf{Structure-preserving preconditioning for variable-coefficient space-fractional diffusion equations.} 
	In Chapter~\ref{chapter:space-FDE}, we study the spectral properties and preconditioning of linear systems arising from the finite difference discretization of time-dependent space-FDEs with variable coefficients. The model involves a one-sided Riemann--Liouville fractional derivative multiplied by a variable coefficient, discretized using the shifted Grünwald formula in space and the $\theta$-method in time. The resulting all-at-once systems exhibit a $(d+1)$-level Toeplitz-like matrix structure. Using GLT theory, we provide a detailed analysis of the spectral distribution of both the original and preconditioned matrix sequences and propose structure-preserving preconditioners based on circulant approximations. The main novelty is that the analysis fully covers the nonconstant variable coefficient case, both in one and two spatial dimensions, extending several existing results in the literature which were often restricted to more simplified settings.
	\item \textbf{Maximal results for geometric means of GLT sequences.} 
	In Chapter~\ref{GM GLT2}, we address a fundamental open question concerning the spectral distribution of geometric means of Hermitian positive definite (HPD) GLT matrix sequences. We prove that when the GLT symbols commute, the invertibility condition on the symbols can be completely removed, resolving a long-standing conjecture. This result is established both in the scalar setting for any dimension and in the general block setting under the assumption that the symbols commute. We further demonstrate that the result is maximal: when the symbols do not commute and are simultaneously degenerate, the classical relation may fail and the geometric mean sequence may not coincide with the classical geometric mean of the symbols. Numerical experiments reveal new spectral phenomena explained by the theory of momentary GLT symbols, including secondary eigenvalue distributions and extremal spectral behavior connected to the order of vanishing of the associated GLT symbol.
\end{enumerate}
Collectively, these contributions provide a comprehensive framework for addressing ISPs in FDEs, combining analytical regularization techniques, spectral analysis of structured matrix sequences, and advanced preconditioning strategies. The thesis advances both the theoretical understanding of ISPs for fractional models and the development of efficient numerical solvers for large-scale structured linear systems arising from their discretization.

	\addcontentsline{toc}{chapter}{Introduction}
	\chapter{Fractional operators and Mittag-Leffler type functions}\label{chapter:Fractional Calculus}
This chapter presents the fundamental concepts of FC,
encompassing both fractional operators and Mittag--Leffler type functions,
which constitute the analytical foundation of this thesis. We review several widely used definitions of fractional integrals and derivatives, which play a central role in the
modeling, analysis, and numerical approximation of FDEs. 

FC generalizes classical differential calculus by allowing
differentiation and integration of non-integer order. Formally, this extension
is expressed as
\[
\frac{\partial^{n}u}{\partial z^{n}}
\;\longrightarrow\;
\frac{\partial^{\alpha}u}{\partial z^{\alpha}},
\qquad n\in\mathbb{N},\ \alpha\in\mathbb{C},
\]
where the resulting operators exhibit intrinsic nonlocality and memory effects.
These properties distinguish fractional derivatives from their integer order
counterparts and make them particularly suitable for modeling anomalous
transport and diffusion processes in heterogeneous and complex media
\cite{Podlubny1998,Samko1993,Metzler2000}.

The origins of FC can be traced back to the seventeenth
century, when L’Hôpital raised the question of interpreting derivatives of
non-integer order. Subsequent contributions by Euler, Lacroix, Fourier, and
Abel gradually clarified the analytical structure of such operators, while the
work of Liouville and Riemann in the nineteenth century provided a rigorous
integral-based framework. These developments led to the classical
Riemann--Liouville formulation and, later, to alternative definitions such as
the Caputo, Gr\"unwald--Letnikov, and Riesz derivatives, each tailored to
different analytical and modeling requirements
\cite{Miller1993,Podlubny1998,Kilbas2006}.

In addition to fractional operators, special functions play a fundamental role
in the theory of FDEs. In particular, Mittag--Leffler type functions naturally arise in the explicit representation of solutions of FDEs and multi-term FDEs. They generalize the exponential function and govern the temporal behavior, stability, and long-time dynamics of fractional models \cite{Gorenflo2014}. For this reason, this chapter also introduces Mittag--Leffler function (MLF) and generalized MLF, which provide essential analytical tools for the theoretical and
numerical developments presented in the subsequent chapters.

\section{Fractional operators}
\begin{definition}\label{def: RL_Integration}
	\cite{Podlubny1998} Let $\xi >0$, and $u(z)\in L^{1}([a,b])$. Then, the left- and right-sided Riemann–Liouville fractional integrals are defined as
	\begin{equation*}
		{\mathcal{D}_{a_+,z}^{-\xi}}u(z)=\frac{1}{\Gamma(\xi)}\int_{a}^{z}\frac{u(\sigma)}{(\,z-\sigma)\,^{1-\xi}}d{\sigma},
	\end{equation*}
	\begin{equation*}
		{\mathcal{D}_{b_{-},z}^{-\xi}}u(z)
		=\frac{1}{\Gamma(\xi)}\int_{z}^{b}\frac{u(\sigma)}{(\,\sigma-z)\,^{1-\xi}}d{\sigma}.
	\end{equation*}
	where $\Gamma(\cdot)$ is the gamma function.
	\smallskip 
	
	\noindent
	Since the development of FC, several definitions of
	non-integer order derivatives have been proposed, including the
	Riemann--Liouville, Caputo, and Gr\"unwald--Letnikov formulations
	\cite{Podlubny1998,Samko1993}. These definitions are generally not equivalent
	on bounded domains, unless additional regularity or boundary conditions are
	imposed on the function $u(z)$. In the present work, we first focus on the
	Riemann--Liouville fractional derivatives (R-LFD), which are introduced below.
\end{definition}
\begin{definition}\label{def: RL_derivative}
	\cite{Podlubny1998} Let $\xi \in (n-1, n)$, $n \in \mathbb{N^{+}}$, and let $u(z)\in AC([a,b])$ be absolutely continuous. Then, the left- and right-sided R-LFD are defined as
	\begin{equation*}
		\mathcal{D}_{a_+,z}^{\xi} u(z)
		=
		\frac{1}{\Gamma(n-\xi)}
		\frac{d^{\,n}}{dz^{\,n}}
		\int_{a}^{z}
		\frac{u(\sigma)}{(z-\sigma)^{\xi-n+1}}
		\, d\sigma,
	\end{equation*}
	
	\begin{equation*}
		\mathcal{D}_{b_-,z}^{\xi} u(z)
		=
		\frac{(-1)^{n}}{\Gamma(n-\xi)}
		\frac{d^{\,n}}{dz^{\,n}}
		\int_{z}^{b}
		\frac{u(\sigma)}{(\sigma-z)^{\xi-n+1}}
		\, d\sigma.
	\end{equation*}
	
\end{definition}
If $\alpha \in \mathbb{N}$, the R-LFDs reduce to the classical integer order derivatives. We now introduce the Riesz fractional derivative, which is closely related to the Riemann--Liouville formulation.

\begin{definition}\label{def:Riesz_fractional derivative}
	\cite{Samko1993}    
	For $\xi\in (n-1, n)$, $n \in \mathbb{N^{+}}$, the Riesz fractional derivative is defined as
	\begin{align}\label{eq:RFDE_LR}
		\frac{\partial^{\xi} u(z)\,}{\partial\vert z\vert^{\xi}}=c(\xi)\left(\mathcal{D}_{a_+,z}^{\xi} u(z)+\mathcal{D}_{b_-,z}^{\xi} u(z)\right),	
	\end{align}
	where $c(\,\xi)\,=-\frac{1}{2\cos{\left(\frac{\xi\pi}{2}\right)}}>0$ with $\xi\neq1$.
\end{definition}
In the 1D case, the Riesz fractional derivative coincides with the fractional Laplacian, which can also be defined as follows.
\begin{definition}\label{def: Hypersingular_integral}
	Let $u(z)$ be in the Schwartz space of rapidly decaying functions and $\alpha\in (0,1)$. The fractional Laplacian operator $(-\Delta)^{\xi}$ is defined as
	\begin{equation}\label{eq:PV}
		(-\Delta)^{\xi}{u}(z)=C(\xi) P.V.\int_{\mathbb{R}}{\frac{{u}(z_{1})-{u}(z_{2})}{|z_{1}-z_{2}|^{1+2\xi}}}d(z_{2}),
	\end{equation}
	where $C(\xi)=\frac{\xi2^{2\xi}\Gamma(\xi+\frac{1}{2})}{\sqrt{\pi}\Gamma(1-\xi)}$ is a normalization constant.
\end{definition}
A generalized class of fractional derivatives is obtained by introducing
exponential tempering, leading to the so-called tempered fractional
derivatives \cite{Meerschaert1997}. In the Riemann--Liouville framework, these operators are defined as follows.
\begin{definition}\label{def:RLTF}
	\cite{BaeumerMeerschaert2010} Let $\xi \in (n-1, n)$, $n \in \mathbb{N^{+}}$, $\vartheta\geq 0$ and let $u(z)\in AC([a,b])$. Then, the left- and right-sided Riemann–Liouville tempered fractional derivatives are defined as
	
	\begin{equation}\label{eq:LRRLD}
		{\mathcal{D}_{a_{+},z}^{\xi,\vartheta}}u(z)=\frac{{\rm{e}}^{-\vartheta z}}{\Gamma(n-\xi)}\frac{d^{n}}{dz^{n}}\int_{a}^{z}\frac{{\rm{e}}^{\vartheta \sigma}u(\sigma)}{(z-\sigma)\,^{\xi-n+1}}d{\sigma},
	\end{equation}
	\begin{equation}\label{eq:RRLD}
		{\mathcal{D}_{b_{-},z}^{\xi,\vartheta}}u(z)=\frac{(\,-1)\,^{n}{\rm{e}}^{\vartheta z}}{\Gamma(n-\xi)}\frac{d^{n}}{dz^{n}}\int_{z}^{b}\frac{{\rm{e}}^{-\vartheta \sigma}u(\sigma)}{(\sigma-z)\,^{\xi-n+1}}d{\sigma}.
	\end{equation}
\end{definition} 
It is evident that, for $\vartheta=0$, the exponential tempering disappears and the tempered fractional derivatives reduce to the classical
R-LFDs.


Another widely used definition is the Caputo fractional derivative (CFD), introduced by the Italian mathematician Michele Caputo \cite{Podlubny1998}, which is defined as follows.

\begin{definition}\label{def:Caputo_FD}
	\cite{Podlubny1998} Let $\xi \in (n-1, n)$, $n \in \mathbb{N^{+}}$, and assume that the $n^\text{th}$ derivative of $u(z)$ is absolutely integrable, then the left- and right-sided CFDs of order $\xi$ are defined as
	\begin{equation}\label{eq:CLFD}
		{}^{C}{\mathcal{D}_{a_{+},z}^{\xi}}u(z)=\frac{1}{\Gamma(n-\xi)}\int_{a}^{z}\frac{1}{(z-\sigma)\,^{\xi-n+1}}\frac{d^{n}u(\sigma)}{d\sigma^{n}}d{\sigma}, 
	\end{equation}
	\begin{equation}\label{eq:CRFD}
		{}^{C}{\mathcal{D}_{b_{-},z}^{\xi}}u(z)=\frac{(-1)^{n}}{\Gamma(n-\xi)}\int_{z}^{b}\frac{1}{(\sigma-z)\,^{\xi-n+1}}\frac{d^{n}u(\sigma)}{d\sigma^{n}}d{\sigma}.
	\end{equation}
\end{definition}
{
	\begin{remark}\label{remark:Caputo_RL}

		\rm{\cite{Kilbas2006}} If $\xi\notin\mathbb{N}$ and $u(x)$ are such that the Caputo and R-LFDs exist then they are linked through the following relation 
		\begin{align*}
			{\mathcal{D}_{a_{+},z}^{\xi}}u(z)&={{}^{C}{\mathcal{D}_{a_{+},z}^{\xi}}u(z)}+\sum_{i=0}^{n-1}{\frac{(z-a)^{i-\xi}}{\Gamma{(i-\xi+1)}}}\frac{d^{i}u(a)}{dz^{i}},\\\nonumber
			{\mathcal{D}_{b_{-},z}^{\xi}}u(z)&={{}^{C}{\mathcal{D}_{b_{-},z}^{\xi}}u(z)}+\sum_{i=0}^{n-1}{\frac{(-1)^{i}(b-z)^{i-\xi}}{\Gamma{(i-\xi+1)}}}\frac{d^{i}u(b)}{dz^{i}}. 
		\end{align*}
		In the case where $\frac{d^{i}u(a)}{dz^{i}}=\frac{d^{i}u(b)}{dz^{i}}=0$, then the R-LFDs \eqref{eq:LRRLD} and \eqref{eq:RRLD} coincide with the CFDs \eqref{eq:CLFD} and \eqref{eq:CRFD}.
	\end{remark}
}

A tempered extension of the CFD has been introduced to account for finite memory effects through exponential damping. This formulation, known as the Caputo tempered fractional derivative, combines the regularity advantages of the CFD with the flexibility of tempered fractional operators and is defined as follows
\cite{BaeumerMeerschaert2010,Samko1993}.

\begin{definition}
	\cite{BaeumerMeerschaert2010,Samko1993} Let $\xi \in (n-1, n)$, $n \in \mathbb{N^{+}}$, $\vartheta\geq 0$ and let $u(z)\in AC([a,b])$. Then, the left- and right-sided Caputo tempered fractional derivatives are defined as
	\begin{equation*}
		{}^{C}{\mathcal{D}_{a_{+},z}^{\xi}}u(z)=\frac{\rm{e}^{-\vartheta z}}{\Gamma(n-\xi)}\int_{a}^{z}\frac{1}{(z-\sigma)\,^{\xi-n+1}}\frac{d^{n}\left(\rm{e}^{\vartheta \sigma}u(\sigma)\right)}{d\sigma^{n}}d{\sigma}, 
	\end{equation*}
	\begin{equation*}
		{}^{C}{\mathcal{D}_{b_{-},z}^{\xi}}u(z)=\frac{(-1)^{n}{\rm{e}}^{\vartheta z}}{\Gamma(n-\xi)}\int_{z}^{b}\frac{1}{(\sigma-z)\,^{\xi-n+1}}\frac{d^{n}\left(\rm{e}^{-\vartheta \sigma}u(\sigma)\right)}{d\sigma^{n}}d{\sigma}.
	\end{equation*}
\end{definition} 
\section{Mittag-Leffler type functions}
This section is devoted to the introduction of the multinomial
MLF and to the presentation of its key analytical
properties, which will be used throughout the subsequent analysis.
\begin{definition}\cite{Luchko1999}
	For $\xi_i, \gamma>0,\ z_i\in \mathbb{C}, i=1,2,...,n,\ n\in \mathbb{N}$, the multinomial MLF is defined as
	$$E_{(\xi_1,\xi_2,...,\xi_n),\gamma}(z_1,z_2,...,z_n):=\sum_{k=0}^{\infty}\quad \sum_{\substack {l_1+l_2+...+l_n=k \\ l_1\geq 0,...,l_n\geq 0}}(k;l_1,...,l_n)
	\frac{\displaystyle\prod_{i=1}^{n}z_{i}^{l_i}}{\Gamma\biggl(\gamma+ \displaystyle\sum_{i=1}^{n}\xi_{i}l_{i}\biggr)},$$
	where $(k;l_1,...,l_n)=\displaystyle\frac{k!}{l_1! \times...\times l_n!}.$
\end{definition}
\begin{remark}\label{Multinomialeq}
	Notice that the multinomial MLF can be represented as
	\begin{align*}
		E_{(\xi_1,\xi_2,...,\xi_n),\gamma}&(z_1,z_2,...,z_n)\\=& \sum_{k=0}^{\infty}\sum_{l_1=0}^{k}\sum_{l_2=0}^{k-l_1}...
		\sum_{l_{n-1}=0}^{k-l_1-\cdots-l_{n-2}}\frac{k!}{l_1!l_2!...(k-l_1-l_2-...-l_{n-1})!}\\
		 \times & \frac{z_{1}^{k-l_{1}-l_{2}-...-l_{n-1}}z_{2}^{l_1}...z_{n}^{l_{n-1}}}{\Gamma(\gamma+\xi_{1}(k-l_{1}-l_{2}-...-l_{n-1})+\xi_{2}l_{1}+...+\xi_{n}l_{n-1})}.
	\end{align*}
\end{remark}
\begin{remark}
	The parameter of the multinomial MLF commutes, i.e.,
	\begin{equation}\label{symetriceq}
		E_{(\xi_1,\xi_2,...,\xi_n),\gamma}{(z_1,z_2,...,z_n)}=E_{(\xi_n,,...,\xi_2,\xi_1),\gamma}{(z_n,...,z_2,z_1)}.
	\end{equation}
\end{remark}
\begin{remark}
	For $z_{1}\neq0$, $z_{2}\neq0$, $z_{3}=...=z_{n}=0,\ n \in \mathbb{N}$  the multinomial MLF takes the following form
	\begin{eqnarray*}
		E_{(\xi_{1},\xi_{2}),\gamma}(z_{1},z_{2})&= &\sum_{k=0}^{\infty}\sum_{\substack {l_1+l_2=k \\ l_1\geq 0,\,l_2\geq 0}}\frac{k!}{l_{1}!l_{2}!}\frac{{z_{1}^{l_1}}z_{2}^{l_2}}{\Gamma(\gamma+\xi_{1}l_1+\xi_{2}l_2)},\\
		&=&\sum_{k=0}^{\infty}\sum_{i=0}^{k}\frac{k!}{i!(k-i)!}\frac{{z_{1}^{i}}z_{2}^{k-i}}{\Gamma(\gamma+\xi_{1}i+\xi_{2}(k-i))}.
	\end{eqnarray*}
\end{remark}
\begin{remark}
	For $z_{1}\neq0$ and $z_{2}=z_{3}=...=z_{n}=0,\ n \in \mathbb{N}$  the multinomial MLF takes the following form
	$$E_{\xi_{1},\gamma}(z_{1})= \sum_{k=0}^{\infty}\frac{{z_{1}^{k}}}{\Gamma(\gamma+\xi_{1}k)}.$$
\end{remark}
\begin{lemma}\label{Podlubny}\cite{Yamamoto2015}
	Let $0<\gamma<2$ and $0<\xi_n<\cdots<\xi_1<1$. Assume that the angle $K$ satisfies $\xi_1 \pi/2<
	K<\xi_1 \pi$, and that $K\leq|arg(z_1)|\leq\pi$ and $-K\leq z_{i}<0$, where $K>0, \;i=2,...,n$. Then there exists a positive constant $C_0$, depending only on $K$, $\xi_i$ $(i=1,\ldots,n)$, and $\gamma$, such that the multinomial MLF satisfies
	$$
	|E_{(\xi_1-\xi_m,...,\xi_1-\xi_2,\xi_1),\gamma}(z_n,...,z_2,z_1)|\leq \frac{C_0}{1+|z_1 |}.
	$$
\end{lemma}
\begin{lemma}\label{LaplaceT}\cite{Malik2021}
	For $\xi_i, \gamma, \tau ,\upsilon_i>0$,\; $i=1,2,...,n,\ n\in \mathbb{N}$ the Laplace transform of the multinomial MLF is given by
	\begin{equation*}
		\pounds{\{\tau^{\gamma-1}{E}_{(\xi_{1},\xi_{2},...,\xi_{n}),\gamma}(-\upsilon_{1}\tau^{\xi_{1}},...,-\upsilon_{n}\tau^{\xi_{n}})\}}\\
		=\frac{s^{-\gamma}}{1+\displaystyle\sum_{i=1}^{n}\upsilon_{i}s^{-\xi_{i}}}.
	\end{equation*}
	If \; $\Big|\sum_{i=1}^{n}\upsilon_{i}s^{-\xi_{i}}\Big|<1.$
\end{lemma}
\begin{lemma}\label{composite}\cite{Malik2021}
	For $\xi_{i}, \varsigma, \gamma, \tau, \upsilon_{i}>0,\;i=1,2,...,n, \ n\in \mathbb{N}$, we have the following relation
	\begin{align*}
		\tau^{\varsigma} \ast \tau^{\gamma-1}&{E}_{(\xi_{1},\xi_{2},...,\xi_{n}), \gamma}(-\upsilon_{1}\tau^{\xi_{1}},...,-\upsilon_{n}\tau^{\xi_{n}})\\=&\,\Gamma(\varsigma+1)\tau^{\gamma +\varsigma}{E}_{(\xi_{1},\xi_{2},...,\xi_{m}),\gamma +\varsigma+1}(-\upsilon_{1}\tau^{\xi_{1}},...,-\upsilon_{n}\tau^{\xi_{n}}).
	\end{align*}
\end{lemma}
\begin{remark}
	For  $\upsilon_{1}\neq 0,\upsilon_{2}=\upsilon_{3}=...=\upsilon_{n}=0, \ n \in \mathbb{N}$
	the above Lemma reduces the well known result
	$$\tau^{\varsigma} \ast \tau^{\gamma-1}{E}_{\xi_{1}, \gamma}(-\upsilon\tau^{\xi_1})=\Gamma(\varsigma+1)\tau^{\gamma+\varsigma}{E}_{\xi_{1}, \gamma+\varsigma+1}(-\upsilon\tau^{\xi_1}).$$
\end{remark}
\begin{proposition}\label{properties}\cite{Malik2021} The following identities hold for Mittag-Leffler type functions:
	\begin{itemize}
			\item $\tau^{2}{E}_{\xi_{1},3}(-\upsilon_{1}\tau^{\xi_1})=\frac{\tau^{2}}{\Gamma(3)}-\upsilon_{1} \tau^{2+\xi_1}{E}_{\xi_{1},3+\xi_{1}}(-\upsilon_{1}\tau^{\xi_1}),$
			\item $\tau^{2-\xi}{E}_{(\xi_{1},\xi_{1}-\xi_{2}),3-\xi_{2}}(-\upsilon_{1}\tau^{\xi_{1}},-\upsilon_{2}\tau^{\xi_{1}-\xi_{2}})\\+\upsilon_{2}\tau^{2+\xi_{1}-2\xi_{2}}{E}_{(\xi_{1},\xi_{1}-\xi_{2}),3+\xi_{1}-2\xi_{2}}(-\upsilon_{1}\tau^{\xi_{1}},-\upsilon_{2}\tau^{\xi_{1}-\xi_{2}})\\=\frac{\tau^{2-\xi_{2}}}{\Gamma(3-\xi_{2})}-\upsilon_{1} \tau^{2+\xi_{1}-\xi_{2}}{E}_{(\xi_{1},\xi_{1}-\xi_{2}),3+\xi_{1}-\xi_{2}}(-\upsilon_{1}\tau^{\xi_{1}},-\upsilon_{2}\tau^{\xi_{1}-\xi_{2}}).$
		\end{itemize}
	\end{proposition}




	\chapter{Matrix sequences: spectral distributions and asymptotic structures
}\label{chapter:Structured Matrix Sequences}
Matrix sequences constitute a fundamental tool in the analysis of numerical methods for differential equations. A matrix sequence is defined as a family of square matrices $\{A_n\}_{n}$, where each matrix $A_n$ has size $d_n \times d_n$ and the dimension $d_n$ increases with $n$, satisfying $d_{n_1} < d_{n_2}$ whenever $n_1 < n_2$. This increase reflects the refinement of a discretization process and naturally arises in the numerical approximation of continuous problems. Matrix sequences typically originate from the discretization of ordinary differential equation (ODEs), partial differential equations (PDEs), and FDEs. Discretizing a differential equation on a computational grid leads to a linear system whose size depends on the number of degrees of freedom. As the mesh is refined to improve accuracy, the number of unknowns increases and the associated coefficient matrices grow in size, producing a sequence indexed by the discretization parameter. These matrices are not arbitrary, but often inherit structural properties from the underlying differential operator, the equation coefficients, and the geometry of the domain.

The structure of a matrix sequence may be explicit, as in the case of Toeplitz, circulant, $\tau$, or Hankel matrices, or asymptotic and hidden, as described by the theory of GLT matrix sequences. The study of matrix sequences is primarily motivated by their asymptotic spectral behavior. The limiting distribution of eigenvalues and singular values provides essential insight into stability, conditioning, and convergence of numerical schemes, and plays a key role in the design of efficient algorithms, including preconditioned iterative solvers and multi-level methods.

This chapter introduces the core concepts underlying the study of
matrix sequences, which provide the spectral and algebraic foundation for the
numerical methods developed in this thesis. We first review eigenvalue and
singular value distributions as the main tools for describing the asymptotic
spectral behavior of large matrices arising from discretized PDEs and FDEs. We then present the notion of convergence for
matrix sequences through approximating classes of sequences, which allows one
to capture both explicit and hidden asymptotic structures. Finally, we discuss
multi-level matrix structures resulting from multi-dimensional discretizations.
These topics form the basis of the GLT framework used throughout the thesis;
see {\rm \cite[Chapter 2]{garoni2018}} for a comprehensive treatment.

\section{Spectral distributions of matrix sequences}\label{sig-eig}
In many relevant situations, important asymptotic information about a matrix sequence is captured by its eigenvalue or singular value distribution. In this section, we formally
introduce this concept, provide its intuitive interpretation, and outline its
role in characterizing the asymptotic spectral behavior of large matrices.

\begin{definition}\label{def:distrib}
	Let $\{A_n\}_n$ be a matrix sequence, with $A_n$ of size $d_n$, and let $\psi : D \subset \mathbb{R}^t \to \mathbb{C}^{r \times r}$ be a measurable function defined on a set $D$ with $0 < \mu_t(D) < \infty$.
	\begin{itemize}
		\item We say that $\{A_n\}_n$ has an (asymptotic) singular value distribution described by $\psi$, and we write $\{A_n\}_n \sim_\sigma \psi$, if
		\[
		\lim_{n \to \infty} \frac{1}{d_n} \sum_{i=1}^{d_n} F(\sigma_i(A_n)) = \frac{1}{\mu_t(D)} \int_D \frac{1}{r} \sum_{i=1}^r F(\sigma_i(\psi(\bm{x}))) \, d\bm{x}, \quad \forall F \in C_c(\mathbb{R}).
		\]
		
		\item We say that $\{A_n\}_n$ has an (asymptotic) spectral (or eigenvalue) distribution described by $\psi$, and we write $\{A_n\}_n \sim_\lambda \psi$, if
		\[
		\lim_{n \to \infty} \frac{1}{d_n} \sum_{i=1}^{d_n} F(\lambda_i(A_n)) = \frac{1}{\mu_t(D)} \int_D \frac{1}{r} \sum_{i=1}^r F(\lambda_i(\psi(\bm{x}))) \, d\bm{x}, \quad \forall F \in C_c(\mathbb{C}).
		\]
		
		\item If $\psi$ describes both the singular value and eigenvalue distribution of $\{A_n\}_n$, we write $\{A_n\}_n \sim_{\sigma, \lambda} \psi$.
	\end{itemize}
	In this case, the function $\psi$ is referred to as the \textit{eigenvalue (or spectral) symbol} of $\{A_n\}_n$.\\
	The same definition holds when the considered matrix sequence shows a multi-level structure. In that case $n$ is replaced by $\bm{n}$ uniformly in $A_n$ and $d_n$.
\end{definition}

From an intuitive viewpoint, when the symbol $\psi$ is continuous a.e. and the eigenvalues of $A_n$ can be associated with evaluations of
$\psi$ on an equispaced grid over the domain $D$, the notion of spectral
distribution reduces to the classical observation that a Riemann sum of a
Riemann-integrable function converges to the corresponding integral. More
generally, if $\psi$ is continuous a.e. and the discretization
parameter $n$ is sufficiently large, then the eigenvalues of $A_n$, after a
suitable ordering and possibly up to a negligible number of outliers of order
$o(d_n)$, behave like uniform samples of $\psi$ over $D$. An analogous
interpretation holds for singular value distributions. Within this framework,
the following result is therefore natural.

\begin{proposition}
	$\{A_n\}_n \sim_\lambda \psi$ and $\Lambda(A_n) \subseteq S$ for all $n\in\mathbb{N},$ then $\operatorname{ess\,Ran}\,\psi\subseteq\bar{S}$ a.e., where ${S}\subseteq \mathbb{C}$ and $\bar{S}$ denotes its closure.
\end{proposition}
In particular $\{A_n\}_n \sim_\lambda \psi$ and $A_n$ positive definite for all $n$ imply that $\psi$ is nonnegative definite almost everywhere, simply nonnegative almost everywhere if $r=1$. 
For singular values, the analogue of positivity is uniform invertibility. Suppose that
	\(\{A_n\}_n\sim_\sigma \psi\) and that there exists \(\delta>0\) such that
	\(\sigma_{\min}(A_n)\ge \delta\) for all \(n\). Then
	\(\sigma_{\min}(\psi)>0\) almost everywhere, i.e.
	\[
	\mu_t\bigl(\{x\in D:\sigma_{\min}(\psi(x))=0\}\bigr)=0.
	\]
	Conversely, if \(\sigma_{\min}(\psi)=0\) on a set of positive measure, then the sequence
	\(\{A_n\}_n\) cannot be uniformly invertible.
\smallskip

\noindent As will become apparent, determining the spectral distribution is generally more straightforward for Hermitian matrix sequences. Nevertheless, several effective techniques are available for treating sequences in which the non-Hermitian component is sufficiently small. In particular, the following proposition shows that perturbations of Hermitian sequences with small norm do not affect the asymptotic spectral distribution.

\begin{proposition}[{\cite[Corollary 2]{Barbarino2020}}]
	Let $\{A_{n}\}_n$ be a matrix sequence such that  $A_{n} = B_{n} + C_{n}$, where $B_{n}$ is Hermitian $\forall n \in \mathbb{N}$. Suppose that
	\begin{itemize}
		\item $\{B_{n}\}_n \sim_\lambda \psi$,
		\item $\|C_{n}\|=o(1)$.
	\end{itemize}
	Then $\{A_{n}\}_n \sim_\lambda \psi$.
\end{proposition}
An important aspect of spectral distributions for matrix sequences is their
intrinsic non-uniqueness. Once a function describing the asymptotic spectral
distribution has been identified, it is by no means unique. Indeed, there
exist infinitely many alternative representations that encode exactly the same
spectral information. In particular, any measurable function that satisfies
the same integral identities with respect to all compactly supported
continuous test functions yields an equivalent distribution. More precisely,
let $\psi$ be a distribution function defined on a measurable set $D$; then
any measurable mapping $\varphi : E \subset \mathbb{R}^s \to \mathbb{C}$ such
that
\[
\frac{1}{\mu_t(D)} \int_D F(\psi(\bm{x}))\, d\bm{x}
=
\frac{1}{\mu_s(E)} \int_E F(\varphi(\bm{x}))\, d\bm{x},
\quad \forall F \in C_c(\mathbb{C}),
\]
represents the same asymptotic distribution. These alternative
representations may differ substantially in form, including the dimension of
their domain or the number of variables involved. Moreover, by suitably extending the notion of distribution, the associated symbol may take matrix-valued values, provided that the corresponding integral characterization remains valid; see, for example, \cite{barbarino2020block1d,barbarino2020blockmulti}.

The notion of clustering can be regarded as a special case of spectral
distribution.

\begin{definition} \label{clustering}
	Let $\{A_n\}_n$ be a matrix sequence with $A_n \in \mathbb{C}^{d_n \times d_n}$,
	and let $S \subseteq \mathbb{C}$ be a nonempty closed set.
	\begin{itemize}
		\item The sequence $\{A_n\}_n$ is said to be \emph{strongly clustered at $S$ in
			the sense of the eigenvalues} if, for any $\varepsilon > 0$ and as
		$n \to +\infty$, the number of eigenvalues of $A_n$ lying outside the
		$\varepsilon$-neighborhood of $S$ is bounded independently of $n$, that is,
		\[
		\#\bigl\{ j \in \{1,\dots,d_n\} : \lambda_j(A_n) \notin B(S,\varepsilon) \bigr\}
		= O(1), \qquad n \to +\infty.
		\]
		
		\item The sequence $\{A_n\}_n$ is said to be \emph{weakly clustered at $S$ in the
			sense of the eigenvalues} if, for any $\varepsilon > 0$ and as
		$n \to +\infty$, the number of eigenvalues of $A_n$ lying outside
		$B(S,\varepsilon)$ is negligible with respect to the matrix size, namely,
		\[
		\#\bigl\{ j \in \{1,\dots,d_n\} : \lambda_j(A_n) \notin B(S,\varepsilon) \bigr\}
		= o(d_n), \qquad n \to +\infty.
		\]
	\end{itemize}
\end{definition}
\begin{remark} \label{remark:distrib-cluster}
	Let us clarify the relationship between the concepts of distribution and clustering. We recall that the essential range of a measurable function $g:D\subseteq\R^t\to\C$ is the set
	\begin{equation*}
		\{ z\in\C : \mu_t \big(\{ g \in B(z, \varepsilon) \} \big)> 0 \quad\forall\varepsilon > 0 \}.
	\end{equation*}
	As reported in \cite[Theorem 4.2]{GoSe}, it holds
	\begin{equation*}
		\{A_n\}_n \sim_\lambda \psi
		\qquad\Longrightarrow\qquad
		\{A_n\}_n \text{ is weakly clustered at the essential range of } \psi.
	\end{equation*}
	Furthermore, if the essential range of $\psi$ consists of a single number $s\in\C$, then
	\begin{equation*}
		\{A_n\}_n \sim_\lambda \psi
		\quad\iff\quad
		\{A_n\}_n \text{ is weakly clustered at } s \text{ in the sense of the eigenvalues}.
	\end{equation*}
	The latter case is typically of interest in the context of preconditioning. Corresponding statements can be given for the singular values, with the obvious suitable changes.
\end{remark}
A corresponding definition can be given for singular values by considering
sets $S \subseteq \mathbb{R}^+$. When the term clustering is used without
further specification, it is implicitly understood in the weak sense. The
special case of \emph{single-point clustering}, where $S=\{s\}$ consists of a
single element, plays a particularly important role in the theory of
preconditioning.

If a matrix sequence admits a continuous spectral distribution, clustering
naturally follows. More precisely, if $\{A_n\}_n$ has eigenvalue distribution
described by a function $\psi$, the eigenvalues can be interpreted as a
uniform sampling of $\psi$, and are therefore clustered around its range.
This intuition is formalized in the following result.

\begin{proposition}
	If $\{A_n\}_n \sim_{\lambda} \psi$, then $\{A_n\}_n$ is weakly clustered at the
	essential range of $\psi$ in the sense of the eigenvalues. Moreover, if
	$\operatorname{ess\,Ran}\psi = \{s\}$ with $s \in \mathbb{C}$ fixed, then
	$\{A_n\}_n \sim_{\lambda} \psi$ if and only if $\{A_n\}_n$ is weakly clustered
	at $s$ in the sense of the eigenvalues.
\end{proposition}

An analogous statement holds for singular values.

\begin{proposition}
	If $\{A_n\}_n \sim_{\sigma} \psi$, then $\{A_n\}_n$ is weakly clustered at the
	essential range of $\psi$ in the sense of the singular values. Furthermore, if
	$\operatorname{ess\,Ran}\psi = \{s\}$ with $s \in \mathbb{R}$ fixed, then
	$\{A_n\}_n \sim_{\sigma} \psi$ if and only if $\{A_n\}_n$ is weakly clustered
	at $s$ in the sense of the singular values.
\end{proposition}

\section{Approximating classes of sequences}\label{acs}
We introduce the notion of an approximating class of sequences (a.c.s.), which
plays a fundamental role in the modern analysis of matrix sequences. In
essence, an a.c.s.\ associated with a given matrix sequence $\{A_n\}_n$
consists of a family of matrix sequences that approximate $\{A_n\}_n$
asymptotically, up to perturbations that are negligible in a precise sense.
Such perturbations typically include a component of low rank relative to the
matrix dimension and another component whose spectral norm vanishes in the
limit. This framework provides a rigorous notion of convergence for
matrix sequences and induces a natural topology on their space. The formal
definition of a.c.s.\ and the corresponding convergence results originate from
earlier works and are presented in a comprehensive and modern form in
\cite{SerraCapizzano2001,Tilli1998,garoni2017,garoni2018}.

\begin{definition}\label{def:acs}
	Let $\{A_n\}_n$ be a matrix sequence and let $\{\{B_{n,j}\}_n\}_j$ be a class of matrix sequences, with $A_n$ and $B_{n,j}$ of size $d_n$. We say that $\{\{B_{n,j}\}_n\}_j$ is an a.c.s. for $\{A_n\}_n$ if the following condition is met: for every $j$, there exists $n_j$ such that, for every $n \geq n_j$,
	\[
	A_n = B_{n,j} + R_{n,j} + N_{n,j},
	\]
	\[
	\text{rank}(R_{n,j}) \leq c(j) d_n \quad \text{and} \quad \|N_{n,j}\| \leq \omega(j),
	\]
	where $n_j$, $c(j)$, and $\omega(j)$ depend only on $j$, and
	\[
	\lim_{j \to \infty} c(j) = \lim_{j \to \infty} \omega(j) = 0.
	\]
	$
	\{\{B_{n,j}\}_n\}_j \xrightarrow{\text{a.c.s. }} \{A_n\}_n
	$ denotes that $\{\{B_{n,j}\}_n\}_j$ is an a.c.s. for $\{A_n\}_n$.
\end{definition}
In intuitive terms, we say that the family $\{\{B_{n,j}\}_n\}_j$ forms an
a.c.s. for $\{A_n\}_n$ if, for all
sufficiently large values of $j$, the sequence $\{B_{n,j}\}_n$ provides an
accurate asymptotic approximation of $\{A_n\}_n$. More precisely, each matrix
$A_n$ can be decomposed into the sum of the corresponding matrix $B_{n,j}$, a
perturbation of rank negligible with respect to the dimension of $A_n$, and a
remainder whose norm is arbitrarily small. As $j$ increases, both the rank and
the norm of these perturbation terms become increasingly insignificant for
large values of $n$.

The following theorem represents the expression of a related convergence theory and it is a powerful tool
used, for example, in the construction of the GLT $\ast$-Algebra.
\begin{theorem}\label{th: acs}
	Let $\{A_n\}_n$, $\{B_{n,j}\}_n$, with $j, n \in \mathbb{N}$, be matrix sequences and let $\psi, \psi_j : D \subset \mathbb{R}^d \to \mathbb{C}^{r \times r}$ be measurable functions defined on a set $D$ with positive and finite Lebesgue measure. Suppose that:
	\begin{enumerate}
		\item $\{B_{n,j}\}_n \sim_{\sigma} \psi_j$ for every $j$;
		\item $\{\{B_{n,j}\}_n\}_j \xrightarrow{\text{a.c.s. }} \{A_n\}_n$;
		\item $\psi_j \to \psi$ in measure.
	\end{enumerate}
	Then
	\[
	\{A_n\}_n \sim_{\sigma} \psi.
	\]
\end{theorem}

\begin{theorem}\label{th: acs eig}
	Given a Hermitian matrix sequence $\{A_n\}_n$, let $\{B_{n,j}\}_n$, be a class of
	Hermitian matrix sequences and let $\psi, \psi_j : D \subset \mathbb{R}^d \to \mathbb{C}^{r \times r}$ be measurable functions defined on a set $D$ with positive and finite Lebesgue measure. Suppose that:
	\begin{enumerate}
		\item $\{B_{n,j}\}_n \sim_{\lambda} \psi_j$ for every $j$;
		\item $\{\{B_{n,j}\}_n\}_j \xrightarrow{\text{a.c.s. } } \{A_n\}_n$;
		\item $\psi_j \to \psi$ in measure.
	\end{enumerate}
	Then
	\[
	\{A_n\}_n \sim_{\lambda} \psi.
	\]
	
\end{theorem}
For clarity, we note that we will occasionally consider so-called
\emph{constant} families of matrix sequences of the form
$\{\{B_n\}_n\}_j$, where the underlying sequence $\{B_n\}_n$ is independent of the index $j$. In such cases, the parameter $j$ plays no effective role and serves only as a formal label. As an example, the family
$\{\{I_n\}_n\}_j$ represents the constant identity matrix sequence, since the element associated with each value of $j$ coincides with the same sequence $\{I_n\}_n$.

\section{Multi-level discretization matrices}\label{multi level} The discretization of PDEs and FDEs on multi-dimensional domains, such as rectangles, cubes, or higher-dimensional hyperrectangles, inherently results in the formation of multi-level matrices. These matrices are characterized by their hierarchical structure, reflecting the underlying geometry and grid of the computing domain. Multi-level matrices are not simply collections of numbers; they are recursively structured into blocks: each matrix is segmented into larger blocks, which can be further subdivided into smaller sub-blocks, and so on. When a PDE or FDE is discretized on a $d$-dimensional grid comprising $n_1, n_2, \ldots, n_d$ points in each spatial direction, the resultant matrix is generally of size $N \times N$, where $N = n_1 n_2 \cdots n_d$.

At the highest level, the matrix consists of $n_1^2$ blocks of size $N/n_1 $, corresponding to subdivisions along the first coordinate; each block at this level is further broken down into $n_2^2$ blocks of size $N/(n_1 n_2) $, according to the second coordinate, and this recursive process continues until, at the final level, we arrive at individual scalar entries. In general, at the $j$-th level, each block is partitioned into $n_j^2$ blocks of order $N/(n_1 \cdots n_j)$. This multi-level structure is not just a mathematical curiosity: it is vital for the efficient storage, computation, and theoretical analysis of the matrices that arise from multi-dimensional problems.

Using the conventional notation of two scalar indices $i, j = 1, \ldots, N$ to refer to the entries of a multi-level matrix quickly becomes impractical. To address this, mathematicians use multi-index notation, where each position in the matrix or vector is labeled by a $d$-tuple, or multi-index, $\bm{i} = (i_1, \ldots, i_d)$ in $\mathbb{Z}^{d}$, that records its location along every dimension of the grid. This approach not only makes the notation more concise and natural but also aligns perfectly with the recursive construction of multi-level matrices. Using multi-indices, we can systematically and clearly describe the positions of entries in both vectors and matrices: for example, a vector $\bm{x} \in \mathbb{C}^{N(\bm{n})}$ with $ N(\bm{n}) =n_1, \ldots, n_d$ can be indexed as $\bm{x} = [ x_{\bm{i}}]_{\bm{i=1}}^{\bm{n}}$, The first component is $x_{(1,\ldots,1)}$, the second is $x_{(1,\ldots,1,2)}$, and so on, until $x_{(n_1,\ldots,n_d)}$. For instance, the component in position $n_d + 3$ will be $x_{(1,\ldots,1,2,3)}$. Likewise, a $d$-level matrix $A \in \mathbb{C}^{{N(\bm{n})}\times{N(\bm{n})}}$ can be written as $A = [a_{\bm{i}\bm{j}}]_{\bm{i},\bm{j}=1}^{\bm{n}}$, where $\bm{n} = (n_1, \ldots, n_d)$.

In this setting, a multi-level matrix sequence is understood as a family of
matrices $\{A_{\boldsymbol{n}}\}_{\boldsymbol{n}}$, indexed by a multi-index $\boldsymbol{n}$ whose components tend to infinity and reflect the refinement of the discretization along each spatial direction, with
$A_{\boldsymbol{n}} \in \mathbb{C}^{{N(\bm{n})}\times{N(\bm{n})}}$. This notion captures the intrinsic multi-dimensional structure of matrices arising from
high-dimensional problems. With these definitions in place, the next chapter is devoted to the structured matrix classes that will be central to the subsequent analysis.

	\chapter{Structured matrices}\label{chap:Structured matrices}
Matrices arising from the discretization of differential equations frequently exhibit characteristic structural patterns that reflect the analytical properties of the underlying continuous operators. Identifying and exploiting such structures is a central theme in numerical linear algebra, as it often leads to significant reductions in computational cost and provides insight into the behavior of numerical schemes. Among the most prominent examples are Toeplitz matrices, which are defined by constant entries along each diagonal. This structure naturally emerges from the discretization of constant coefficient PDEs and FDEs, where translation invariance of the operator is preserved at the discrete level. In a one-dimensional setting, the use of uniform grids typically produces standard Toeplitz matrices, while higher-dimensional problems posed on rectangular domains give rise to multi-level Toeplitz structures. In this case, the discretization leads to block Toeplitz matrices whose blocks are themselves Toeplitz, with the number of recursive levels directly linked to the spatial dimension of the problem. Such structured matrices play a fundamental role in both theoretical analysis and algorithmic design, forming the basis for fast matrix--vector products, efficient solvers, and the spectral tools developed in the subsequent sections of this chapter.

Throughout this chapter, all the structures are introduced directly in the $d$-level, $r$-block setting, that is, for $\bm{n} \in \mathbb{N}^d$ and $r \in \mathbb{N}$ the matrices under consideration belong to $\mathbb{C}^{{N(\bm{n})r}\times{N(\bm{n})r}}$ and their entries are $r \times r$ matrices rather than scalars. This level of generality is required in the rest of the thesis, where blocks of size $r>1$ arise from systems of FDEs, from vector-valued unknowns and from all-at-once discretizations. The scalar case treated in the classical literature is recovered by setting $r=1$, in which case all the statements below reduce verbatim to their scalar counterparts.

This chapter is devoted to the systematic study of Toeplitz and related
structured matrices. In Section \ref{TM}, we introduce the formal definition of multi-level Toeplitz matrices and review their main algebraic and spectral properties, drawing on the extensive literature developed over the past decades. Unless stated otherwise, these results are collected from \cite[Chapter~3]{garoni2018}. In addition, we discuss several matrix classes that either arise as special cases of Toeplitz matrices or exhibit strong structural and spectral connections with them, which makes them particularly relevant in numerical applications. Circulant and $\omega$-circulant matrices are presented in Sections~\ref{cir mat} and \ref{omega cir mat}, respectively, while Sections~\ref{tau-mat} and \ref{hank matr} are devoted to $\tau$ and Hankel matrices.

\section{Toeplitz matrices}\label{TM}
Given $\bm{n} \in \mathbb{N}^d$ and $r \in \mathbb{N}$, a matrix of the form
\[
T_{\bm{n}} := [A_{\bm{i}-\bm{j}}]_{\bm{i},\bm{j}=\bm{1}}^{\bm{n}} \in \mathbb{C}^{{N(\bm{n})r}\times{N(\bm{n})r}},
\]
with blocks $A_{\bm{k}} \in \mathbb{C}^{r\times r}$, $\bm{k} \in \{-(\bm{n}-\bm{1}), \ldots, \bm{n}-\bm{1}\}$, where the block in position $(\bm{i}, \bm{j})$ depends only on the difference $\bm{i}-\bm{j}$, is called a multi-level block Toeplitz matrix, or, more precisely, a $d$-level $r$-block Toeplitz matrix. For $r=1$ the blocks are scalars and one recovers the classical $d$-level (scalar) Toeplitz matrices.

In the simplest case where $d = 1$ and $\bm{n} = n$, a 1-level $r$-block Toeplitz matrix is just a classical block Toeplitz matrix, whose blocks remain constant along the block diagonals:
\[
T_n = [A_{i-j}]_{i,j=1}^n =
\begin{bmatrix}
	A_0 & A_{-1} & A_{-2} & \cdots & \cdots & A_{1-n} \\
	A_1 & A_0 & A_{-1} & A_{-2}&  & \vdots \\
	A_2 & A_1 & \ddots &\ddots & \ddots & \vdots \\
	\vdots & A_2&\ddots&\ddots&\ddots&A_{-2}\\
	\vdots &  & \ddots & \ddots & \ddots &A_{-1} \\
	A_{n-1} & \cdots & \cdots & A_{2} & A_1&A_{0}
\end{bmatrix}
\in \mathbb{C}^{nr\times nr}.
\]
Note that $T_n$ is completely specified by its $2n - 1$ blocks $A_k \in \mathbb{C}^{r\times r}$, $k = -(n-1), \ldots, n-1$, or equivalently by its first block row and first block column.

For $d = 2$ and $\bm{n} = (n_1, n_2)$, we have a 2-level $r$-block Toeplitz matrix:

\[
T_{\bm{n}} = \bigg[\left[ A_{(i_1 - j_1,\, i_2 - j_2)} \right]_{i_2, j_2 = 1}^{n_2}\bigg]_{i_1, j_1 = 1}^{n_1}
= \left[ \mathcal{A}_{i_1 - j_1} \right]_{i_1, j_1 = 1}^{n_1}\]
\[
= 
\begin{bmatrix}
	\mathcal{A}_0 & \mathcal{A}_{-1} & \mathcal{A}_{-2} & \cdots & \cdots & \mathcal{A}_{1-n_{1}} \\
	\mathcal{A}_1 & \mathcal{A}_0 & \mathcal{A}_{-1} & \mathcal{A}_{-2}&  & \vdots \\
	\mathcal{A}_2 & \mathcal{A}_1 & \ddots &\ddots & \ddots & \vdots \\
	\vdots & \mathcal{A}_2&\ddots&\ddots&\ddots&\mathcal{A}_{-2}\\
	\vdots &  & \ddots & \ddots & \ddots &\mathcal{A}_{-1} \\
	\mathcal{A}_{n_{1}-1} & \cdots & \cdots & \mathcal{A}_{2} & \mathcal{A}_1&\mathcal{A}_{0}
\end{bmatrix}
\in \mathbb{C}^{n_1 n_2 r \times n_1 n_2 r},
\]
where each block $\mathcal{A}_k$ is itself a 1-level $r$-block Toeplitz matrix:\\
$
\mathcal{A}_k = \left[ A_{(k,\, i_2 - j_2)} \right]_{i_2, j_2 = 1}^{n_2}$

\[=
\begin{bmatrix}
	A_{(k, 0)} & A_{(k, -1)} & A_{(k, -2)} & \cdots & \cdots & A_{(k, 1-n_{2})} \\
	A_{(k, 1)} & A_{(k, 0)} & A_{(k, -1)} & A_{(k, -2)}&  & \vdots \\
	A_{(k, 2)} & A_{(k, 1)} & \ddots &\ddots & \ddots & \vdots \\
	\vdots & A_{(k, 2)}&\ddots&\ddots&\ddots&A_{(k, -2)}\\
	\vdots &  & \ddots & \ddots & \ddots &A_{(k, -1)} \\
	A_{(k, n_{2}-1)} & \cdots & \cdots & A_{(k, 2)} & A_{(k, 1)}&A_{(k, 0)}
\end{bmatrix}\in \mathbb{C}^{n_2 r\times n_2 r}.
\]
For $r=1$, this type of matrix is also known as a \emph{block Toeplitz matrix with Toeplitz blocks} or BTTB matrix; for a general $r$ it is a BTTB matrix whose innermost entries are $r \times r$ matrices.
\smallskip 

In the case where $d = 3$, a 3-level $r$-block Toeplitz matrix is made up of blocks, and each block has its own 2-level $r$-block Toeplitz structure.  A more general way to describe a $d$-level $r$-block Toeplitz matrix is as a block Toeplitz matrix whose blocks are $(d-1)$-level $r$-block Toeplitz matrices.  In this recursive design, the blocks in the first block row and first block column of the matrix tell you everything you need to know about it, just like in the case of a single-level matrix.
From the computational perspective, this makes the storage requirement for Toeplitz matrices drop to only $O(N(\bm{n})r^2)$, compared to $O\big((N(\bm{n})r)^2\big)$ for a generic dense matrix.\\
A compact expression for a multi-level block Toeplitz matrix can be formulated by defining, for \(n\in\mathbb{N}\) and \(k\in\mathbb{Z}\), define the matrix \(J_n^{(k)}\in\mathbb{R}^{n\times n}\) by
\[
[J_n^{(k)}]_{i,j}
:= \delta_{i-j,\,k}
=
\begin{cases}
	1, & i-j=k,\\[2mm]
	0, & \text{otherwise},
\end{cases}
\qquad i,j=1,\dots,n,
\]
where \(\delta_{p,q}\) is the Kronecker delta, i.e.
\[
\delta_{p,q}=
\begin{cases}
	1, & p=q,\\
	0, & p\neq q.
\end{cases}
\]
Equivalently, $J_n^{(k)}$ is a matrix with all ones on the $k$-th diagonal and zeros elsewhere. The main diagonal is numbered as 0, diagonals above as $-1, -2, ..., -n+1$, and those below as $1, 2, ..., n-1$. For the multi-level case, where $\bm{n} \in \mathbb{N}^{d}$ and $\bm{k} \in \mathbb{Z}^{d}$, are both $d$-dimensional vectors (for example, $\bm{n} = (n_1, ..., n_d)$ and $\bm{k} = (k_1, ..., k_d)$), we define its multi-level version:
\[
J_{\bm{n}}^{(\bm{k})} = J_{n_1}^{(k_1)} \otimes J_{n_2}^{(k_2)} \otimes \cdots \otimes J_{n_d}^{(k_d)}.\]

Then, $T_{\bm{n}}$ admits the following representation:
\begin{equation}\label{T comp}
	T_{\bm{n}} = [A_{\bm{i}-\bm{j}}]_{\bm{i},\bm{j}=\bm{1}}^{\bm{n}} = \sum_{\bm{k}=\bm{1-n}}^{\bm{n-1}} J_{\bm{n}}^{(\bm{k})} \otimes A_{\bm{k}}.
\end{equation}
It can be proven componentwise by explicitly computing the $(\bm{i}, \bm{j})$-th block
on the right-hand side of the equation and applying the properties of the
Kronecker product.
\smallskip 

Given a matrix-valued function $f : [-\pi, \pi]^d \to \mathbb{C}^{r \times r}$ belonging to the space of integrable functions $L^1([-\pi, \pi]^d)$, that is, such that all its entries are in $L^1([-\pi,\pi]^d)$, it is possible to construct a class of multi-level block Toeplitz matrices associated with $f$. Let $\hat{f}_{\bm{k}}$ denote the $\bm{k}$-th Fourier coefficient of $f$, defined for any $\bm{k} \in \mathbb{Z}^d$ as
\[
\hat{f}_{\bm{k}} := \frac{1}{(2\pi)^d} \int_{[-\pi,\pi]^d} f(\bm{\theta})\, e^{-i \bm{k} \cdot \bm{\theta}} \, d\bm{\theta} \in \mathbb{C}^{r \times r},
\]
where $i^2 = -1$ and $\bm{k} \cdot \bm{\theta} = k_1\theta_1 + \cdots + k_d\theta_d$, the integral being computed entrywise.

For any $\bm{n} \in \mathbb{N}^d$, the $\bm{n}$-th Toeplitz matrix associated with $f$ is the $d$-level $r$-block Toeplitz matrix, denoted by $T_{\bm{n}}(f)$, whose blocks are the Fourier coefficients of $f$, arranged along the block diagonals according to the following definition:
\[
T_{\bm{n}}(f) := [\hat{f}_{\bm{i}-\bm{j}}]_{\bm{i},\bm{j}=\bm{1}}^{\bm{n}} \in \mathbb{C}^{{N(\bm{n})r}\times{N(\bm{n})r}}.
\]
The family $\{T_{\bm{n}}(f)\}_{\bm{n} \in \mathbb{N}^d}$ is called the family of (multi-level block) Toeplitz matrices associated with $f$, which is referred to as the \emph{generating function}.
Using Equation (\ref{T comp}), $T_{\bm{n}}(f)$ can also be represented as
\[
T_{\bm{n}}(f) = \sum_{\bm{k} = \bm{1-n}}^{\bm{n-1}} J_{\bm{n}}^{(\bm{k})} \otimes \hat{f}_{\bm{k}} = \sum_{k_1 = 1 - n_1}^{n_1 - 1} \cdots \sum_{k_d = 1 - n_d}^{n_d - 1} \left( J_{n_1}^{(k_1)} \otimes \cdots \otimes J_{n_d}^{(k_d)} \right) \otimes \hat{f}_{(k_1, \ldots, k_d)}.
\]

\subsection{Relationship between the Toeplitz matrix and its generating function}
By constructing a Toeplitz matrix $T_{\bm{n}}(f)$ from a given matrix-valued function $f$ in $L^1([-\pi, \pi]^d)$, we are essentially creating a map that takes functions and produces matrices:
\[
T_{\bm{n}}(\cdot) : L^1([-\pi, \pi]^d, \mathbb{C}^{r \times r}) \to \mathbb{C}^{{N(\bm{n})r}\times{N(\bm{n})r}}
\]
for any fixed $\bm{n} \in \mathbb{N}^d$ and $r \in \mathbb{N}$. The properties of the function $f$ are closely related to those of the Toeplitz matrix.
To start, the operator $T_{\bm{n}}(\cdot)$ is linear. This follows directly from the linearity of the Fourier coefficient operation.
\begin{proposition}
	For any $\alpha, \beta \in \mathbb{C}$ and $f, g \in L^1([-\pi, \pi]^d, \mathbb{C}^{r \times r})$, it holds that
	\[
	T_{\bm{n}}(\alpha f + \beta g) = \alpha\, T_{\bm{n}}(f) + \beta\, T_{\bm{n}}(g).
	\]
	Moreover, further results follow from the identities
	\[
	\widehat{(f(-\bm{\theta}))}_{\bm{k}} = \hat{f}_{-\bm{k}}
	\quad \text{and} \quad
	\widehat{(f^{*})}_{\bm{k}} = \big(\hat{f}_{-\bm{k}}\big)^{*},
	\quad \text{for all } \bm{k} \in \mathbb{Z}^d,
	\]
	where $f^{*}(\bm{\theta}) := \overline{f(\bm{\theta})}^{\top}$ denotes the conjugate transpose of $f$.
\end{proposition}
\begin{proposition}
	For any $f \in L^1([-\pi, \pi]^d, \mathbb{C}^{r \times r})$, the following properties hold:
	\[
	T_{\bm{n}}(f(\bm{\theta}))^{\top} = T_{\bm{n}}\big(f(-\bm{\theta})^{\top}\big),
	\]
	\[
	T_{\bm{n}}(f)^* = T_{\bm{n}}(f^{*}),
	\]
	where $T_{\bm{n}}(f)^*$ is the Hermitian (conjugate transpose) of $T_{\bm{n}}(f)$.
\end{proposition}
In particular, if $f$ is Hermitian a.e., that is, $f(\bm{\theta}) = f(\bm{\theta})^{*}$ for almost every $\bm{\theta}$, then $T_{\bm{n}}(f)$ is Hermitian, so that Hermitian matrix-valued generating functions produce block Toeplitz matrices with real eigenvalues. More precisely, the spectrum of $T_{\bm{n}}(f)$ is contained within the essential infimum of the smallest eigenvalue of $f$ and the essential supremum of its largest eigenvalue.
\begin{theorem}\label{the 2.3}
	Given $f \in L^1([-\pi, \pi]^d, \mathbb{C}^{r \times r})$, Hermitian a.e., let
	\[
	m_f := \operatorname*{ess\,inf}_{\bm{\theta}} \lambda_{\min}\big(f(\bm{\theta})\big),
	\qquad
	M_f := \operatorname*{ess\,sup}_{\bm{\theta}} \lambda_{\max}\big(f(\bm{\theta})\big).
	\]
	Then, for any $\bm{n} \in \mathbb{N}^d$, it holds that
	\[
	\lambda_j \big(T_{\bm{n}}(f)\big) \in [m_f, M_f], \quad \forall\, j = 1, \ldots, N(\bm{n})r.
	\]
	If, in addition, $m_f < M_f$, or equivalently $f$ is not constant a.e., then it holds
	\[
	\lambda_j \big(T_{\bm{n}}(f)\big) \in (m_f, M_f), \quad \forall\, j = 1, \ldots, N(\bm{n})r.
	\]
\end{theorem}
A valuable addition to this theorem is the following result concerning the
extreme eigenvalues of $T_{\bm{n}}(f)$: for any fixed $j \ge 1$ and as $\bm{n} \to \infty$,
the $j$ smallest eigenvalues of $T_{\bm{n}}(f)$ converge to $m_f$,
while the $j$ largest eigenvalues converge to $M_f$.

\begin{theorem}
	Given $f \in L^1([-\pi, \pi]^d, \mathbb{C}^{r \times r})$, Hermitian a.e., let $m_f$ and $M_f$ be defined as in Theorem~\ref{the 2.3}, and let
	\[
	\lambda_1 \big(T_{\bm{n}}(f)\big) \leq \lambda_2 \big(T_{\bm{n}}(f)\big) \leq \cdots \leq \lambda_{N(\bm{n})r} \big(T_{\bm{n}}(f)\big)
	\]
	be the eigenvalues of $T_{\bm{n}}(f)$ arranged in increasing order. Then, for any fixed $j \geq 1$, it holds that
	\[
	\lim_{\bm{n} \to \infty} \lambda_j\big(T_{\bm{n}}(f)\big) = m_f, \qquad
	\lim_{\bm{n} \to \infty} \lambda_{N(\bm{n})r-j+1}\big(T_{\bm{n}}(f)\big) = M_f.
	\]
\end{theorem}
From Theorem~\ref{the 2.3}, we immediately deduce the positivity and monotonicity of the block Toeplitz operator.
\begin{corollary}
	For any $f, g \in L^1([-\pi, \pi]^d, \mathbb{C}^{r \times r})$, Hermitian a.e., the following holds:
	\[
	f \geq O_{r} \text{ a.e.} \iff T_{\bm{n}}(f) \geq O_{N(\bm{n})r}, \quad \forall \bm{n} \in \mathbb{N}^d,
	\]
	and, consequently,
	\[
	f \geq g \text{ a.e.} \iff T_{\bm{n}}(f) \geq T_{\bm{n}}(g), \quad \forall \bm{n} \in \mathbb{N}^d,
	\]
	where the inequalities between matrices are understood in the sense of the partial ordering induced by positive semidefiniteness.
\end{corollary}
These results remain valid even if all of the inequalities are replaced by strict inequalities. In other words, if you consider situations where $f \geq O_r$ a.e., or $f \geq g$ a.e., the corresponding conclusions about the block Toeplitz matrices still hold when using the strict inequality signs.

The spectral norm of a block Toeplitz matrix is connected with the $\infty$-norm of its generating function through the following inequality, proved in \cite{serra97}.

\begin{proposition}
	For any $f \in L^{\infty}([-\pi, \pi]^d, \mathbb{C}^{r \times r})$ and any $\bm{n} \in \mathbb{N}^d$, it holds that
	\[
	\|T_{\bm{n}}(f)\| \leq \|f\|_{\infty},
	\qquad
	\|f\|_{\infty} := \operatorname*{ess\,sup}_{\bm{\theta}} \|f(\bm{\theta})\|.
	\]
\end{proposition}
Finally, we present a useful relation between the separability of the function $f$ and the associated multi-level block Toeplitz matrices.
\begin{proposition}
	Let $f \in L^1([-\pi, \pi]^d, \mathbb{C}^{r \times r})$. If $f$ is separable, i.e.,
	\[
	f = f_1 \otimes \cdots \otimes f_d,
	\]
	with $f_1, \ldots, f_{d-1}$ scalar-valued and $f_d$ matrix-valued, then, for any $\bm{n} \in \mathbb{N}^d$, it holds that
	\[
	T_{\bm{n}}(f) = T_{n_1}(f_1) \otimes \cdots \otimes T_{n_d}(f_d).
	\]
\end{proposition}
Now, by allowing $\bm{n} \in \mathbb{N}^d$ to vary rather than keeping it fixed, we obtain a multi-level block Toeplitz matrix sequence $\{T_{\bm{n}}(f)\}_{\bm{n}}$.
The deep connection between the singular value and eigenvalue distributions of a Toeplitz sequence and its generating function has been established in
various forms over the years. We report here the version presented in \cite[Theorem~3.5]{garoni2018}, where an elegant proof based on the notion of
a.c.s. is provided. For the original formulation in the scalar case, see \cite{tyrtyshnikov1998}, while the case of matrix-valued generating functions considered here can be found in \cite{Tilli1998}.
\begin{theorem}
	With reference to Definition \ref{def:distrib}, if $f \in L^1([-\pi, \pi]^d, \mathbb{C}^{r \times r})$, then
	\[
	\{T_{\bm{n}}(f)\}_{\bm{n}} \sim_\sigma f.
	\]
	If, moreover, $f$ is Hermitian a.e., then
	\[
	\{T_{\bm{n}}(f)\}_{\bm{n}} \sim_\lambda f.
	\]
	Here the singular values and the eigenvalues of the matrix-valued symbol $f(\bm{\theta})$ play the role of the values of the scalar symbol.
\end{theorem}

\section{Circulant matrices}\label{cir mat}
Given $\bm{n} \in \mathbb{N}^d$ and $r \in \mathbb{N}$, a matrix of the form
\[
C_{\bm{n}} := [A_{(\bm{i}-\bm{j}) \bmod \bm{n}}]_{\bm{i},\bm{j}=\bm{1}}^{\bm{n}} \in \mathbb{C}^{{N(\bm{n})r}\times{N(\bm{n})r}},
\]
with blocks $A_{\bm{k}} \in \mathbb{C}^{r \times r}$, $\bm{k} \in \{\bm{0}, \ldots, \bm{n}-\bm{1}\}$, where the block in position $(\bm{i}, \bm{j})$ depends only on the difference
$(\bm{i} - \bm{j}) \bmod \bm{n}$, is called a \emph{multi-level block} or \emph{$d$-level $r$-block circulant matrix}.
In particular, since the $(\bm{i}, \bm{j})$-th block depends only on $\bm{i} - \bm{j}$, circulant matrices can be regarded as a special subclass of Toeplitz matrices.
\smallskip

When $d = 1$ and $\bm{n} = n$, this reduces to a classic $r$-block circulant matrix, where each block column is a one-position downward shift of the previous one:
\[
C_n = \bigr[A_{(i-j) \bmod n}\bigr]_{i,j=1}^n =
\begin{bmatrix}
	A_0 & A_{n-1} & A_{n-2} & \cdots & \cdots & A_{1} \\
	A_1 & A_0 & A_{n-1} & A_{n-2}&  & \vdots \\
	A_2 & A_1 & \ddots &\ddots & \ddots & \vdots \\
	\vdots & A_2&\ddots&\ddots&\ddots&A_{n-2}\\
	\vdots &  & \ddots & \ddots & \ddots &A_{n-1} \\
	A_{n-1} & \cdots & \cdots & A_{2} & A_1&A_{0}
\end{bmatrix}
\in \mathbb{C}^{nr \times nr}.
\]
$C_n$ is entirely determined by the $n$ blocks in its first block column, $\{A_k\}_{k=0}^{n-1}$. Similarly, each block row is a one-element rightward shift of the previous block row, so that the matrix is fully specified by its first block row.

For $d = 2$ and $\bm{n} = (n_1, n_2)$, we obtain a 2-level $r$-block circulant matrix:
\[
C_{\bm{n}} = \Bigr[\bigr[ A_{((i_1-j_1) \bmod n_1,\, (i_2-j_2) \bmod n_2)} \bigr]_{i_2,j_2=1}^{n_2}\Bigr]_{i_1,j_1=1}^{n_1}
= \bigr[ \mathcal{A}_{(i_1-j_1) \bmod n_1} \bigr]_{i_1,j_1=1}^{n_1}\]
\[
=
\begin{bmatrix}
	\mathcal{A}_0 & \mathcal{A}_{n_1-1} & \mathcal{A}_{n_1-2} & \cdots & \cdots & \mathcal{A}_{1} \\
	\mathcal{A}_1 & \mathcal{A}_0 & \mathcal{A}_{n_1-1} & \mathcal{A}_{n_1-2}&  & \vdots \\
	\mathcal{A}_2 & \mathcal{A}_1 & \ddots &\ddots & \ddots & \vdots \\
	\vdots & \mathcal{A}_2&\ddots&\ddots&\ddots&\mathcal{A}_{n_1-2}\\
	\vdots &  & \ddots & \ddots & \ddots &\mathcal{A}_{n_1-1} \\
	\mathcal{A}_{n_{1}-1} & \cdots & \cdots & \mathcal{A}_{2} & \mathcal{A}_1&\mathcal{A}_{0}
\end{bmatrix}
\in \mathbb{C}^{n_1 n_2 r\times n_1 n_2 r},
\]
in which each block $\mathcal{A}_k$ is itself a 1-level $r$-block circulant matrix:
\[
\mathcal{A}_k =
\left[ A_{(k,\ (i_2-j_2) \bmod n_2 )} \right]_{i_2, j_2=1}^{n_2}\]
\[
=
\begin{bmatrix}
	A_{(k, 0)} & A_{(k, n_2-1)} & A_{(k, n_2-2)} & \cdots & \cdots & A_{(k, 1)} \\
	A_{(k, 1)} & A_{(k, 0)} & A_{(k, n_2-1)} & A_{(k, n_2-2)}&  & \vdots \\
	A_{(k, 2)} & A_{(k, 1)} & \ddots &\ddots & \ddots & \vdots \\
	\vdots & A_{(k, 2)}&\ddots&\ddots&\ddots&A_{(k, n_2-2)}\\
	\vdots &  & \ddots & \ddots & \ddots &A_{(k, n_2-1)} \\
	A_{(k, n_{2}-1)} & \cdots & \cdots & A_{(k, 2)} & A_{(k, 1)}&A_{(k, 0)}
\end{bmatrix}\in \mathbb{C}^{n_2 r\times n_2 r}.
\]
For $r=1$ this structure is often called a \emph{block circulant matrix with circulant blocks (BCCB)}.

For $d = 3$, a 3-level $r$-block circulant matrix is constructed as a block circulant matrix, where each block is itself a 2-level $r$-block circulant matrix. More generally, a $d$-level $r$-block circulant matrix can be defined recursively as a block circulant matrix whose blocks have $(d-1)$-level $r$-block circulant structure. Just as in the classic (uni-level) case, the entire matrix is fully specified by the blocks in its first block column, or, equivalently, by those in its first block row. 

A compact expression for a uni-level $r$-block circulant matrix $C_{n,r}$ can be derived by using
\[
[Z_n]_{i,j} := \delta_{(i-j-1) \bmod n} =
\begin{cases}
	1 & \text{if } i - j = 1 \mod{n}, \\
	0 & \text{otherwise},
\end{cases}
\qquad i, j = 1, \ldots, n,
\]
where $n \in \mathbb{N}$ and $\delta$ is the Dirac delta function.
More explicitly, $Z_n$ is the matrix that has all ones on the first and $(1-n)$-th diagonals, with zeros elsewhere:
\[
Z_n =
\begin{bmatrix}
	0      &        &        &        &       &   1    \\
	1      & \ddots &        &        &       &        \\
	& \ddots & \ddots &        &       &        \\
	&        & \ddots & \ddots &       &        \\
	&        &        & \ddots & \ddots&        \\
	&        &        &        &1      &   0    \\
\end{bmatrix}
\in \mathbb{C}^{n\times n}.
\]
Raising $Z_n$ to the power $k$, for $k = 1, \ldots, n-1$, cyclically shifts the diagonals by $k$ positions, with entries given by
\[
[Z_n^k]_{i,j} := \delta_{(i-j-k) \bmod n} =
\begin{cases}
	1 & \text{if } i-j = k \mod{n}, \\
	0 & \text{otherwise},
\end{cases}
\qquad i, j = 1, \ldots, n.
\]

More explicitly,
\[
Z_n^2 =
\begin{bmatrix}
	0      &        &        &        &  1    &   0   \\
	0      & \ddots &        &        &       &   1   \\
	1      & \ddots & \ddots &        &       &       \\
	&  \ddots& \ddots & \ddots &       &      \\
	&        & \ddots & \ddots & \ddots&       \\
	&        &        &    1   &   0   &   0   \\
\end{bmatrix}, \qquad
Z_n^3 =
\begin{bmatrix}
	0      &        &        &     1  &    0  &   0   \\
	0      & \ddots &        &        &    1  &   0   \\
	0      & \ddots & \ddots &        &       &   1   \\
	1      & \ddots & \ddots & \ddots &       &      \\
	& \ddots & \ddots & \ddots & \ddots&       \\
	&        &  1      &    0   &   0   &   0   \\
\end{bmatrix}, \quad \cdots
\]
Hence, any uni-level $r$-block circulant matrix can be expressed as a linear combination of Kronecker products of these powers with the blocks $A_k$:
\begin{equation}\label{2.2} 
	C_{n,r} = \sum_{k=0}^{n-1} Z_n^k \otimes A_k. 
\end{equation}
Conversely, since $Z_n^n = I_n$ and, more generally, $Z_n^k = Z_n^{k \bmod n}$ for any $k \in \mathbb{Z}$, any finite sum of the form
\[
\sum_{k \in K} Z_n^k \otimes A_k, \quad A_k \in \mathbb{C}^{r\times r},\, K \subset \mathbb{Z},
\]
is an $r$-block circulant matrix. For this reason, $Z_n$ is commonly referred to as the generator of circulant matrices of order $n$.

The multi-level generalization of circulant matrices is straightforward. For $\bm{n} \in \mathbb{N}^d$ and $\bm{k} \in \mathbb{Z}^d$, we define
\[
Z_{\bm{n}}^{\bm{k}} := Z_{n_1}^{k_1} \otimes Z_{n_2}^{k_2} \otimes \cdots \otimes Z_{n_d}^{k_d}.
\]
Using this definition, any multi-level block circulant matrix can be written as
\[
C_{\bm{n},r} = \sum_{\bm{k}=\bm{0}}^{\bm{n-1}} Z_{\bm{n}}^{\bm{k}} \otimes A_{\bm{k}}.
\]
This formula is directly analogous to the uni-level case and demonstrates that, just as for classic circulant matrices, the structure of multi-level block circulant matrices is fully captured by these shift matrices and the set of blocks $A_{\bm{k}}$. 

\subsection{Spectral structure and fast solvers for circulant matrices}
Circulant matrices are widely used in numerical analysis due to their deep connection with the Discrete Fourier Transform (DFT) and, consequently, with the Fast Fourier Transform (FFT). Recall that the DFT is the linear mapping that associates to a vector $\bm{x} \in \mathbb{C}^n$ its product with the Fourier
matrix $F_n$, namely $F_{n}\bm{x}$. The FFT is an efficient algorithm for computing the DFT in $\mathcal{O}(n \log n)$ operations, in contrast to the $\mathcal{O}(n^2)$ operations required by a standard matrix--vector multiplication \cite{vanloan1992}. The fact that block circulant matrices are uniformly block-diagonalizable by the Fourier matrix makes them particularly attractive for practical applications.

We start by introducing the Fourier matrix of order $n \in \mathbb{N}$:
\begin{equation}\label{fourier-matrix}
F_n := \frac{1}{\sqrt{n}}\left[ e^{-2\pi i \frac{rs}{n}} \right]_{r,s=0}^{n-1} \in \mathbb{C}^{n\times n},\qquad i=\sqrt{-1}.
\end{equation}
For convenience, this expression involves a slight abuse of notation, as the row and column indices start from $0$ rather than $1$. It is straightforward
to verify that $F_n$ is both symmetric and unitary. Extending this definition, the multi-level Fourier matrix for $\bm{n} \in \mathbb{N}^d$ is given by
\[
F_{\bm{n}} := F_{n_1} \otimes \cdots \otimes F_{n_d},
\]
and its $r$-block counterpart is defined as
\[
F_{\bm{n},r} := F_{\bm{n}} \otimes I_r.
\]
The matrix $F_{\bm{n}}$ preserves the symmetry and unitarity of its constituent matrices $F_{n_i}$, $i = 1, \ldots, d$, as a direct consequence of the
structural properties of the Kronecker product; the same holds for $F_{\bm{n},r}$, which is unitary as well.

$F_{\bm{n}}$ inherits its symmetry and unitarity from its components $F_{n_i}$, $i = 1, \ldots, d$, by virtue of the properties of the Kronecker product.

For any multi-level block circulant matrix, we have the following block-diagonalization.
\begin{theorem}\label{thm:circulant}
	For all $\bm{n} \in \mathbb{N}^d$ and $r \in \mathbb{N}$,
	\begin{equation}\label{2.3}
		C_{\bm{n},r} = \left[ A_{(\bm{i-j}) \bmod \bm{n}} \right]_{\bm{i,j=1}}^{\bm{n}} = F_{\bm{n},r}\, \Lambda_{\bm{n},r}\, F_{\bm{n},r}^*, 
	\end{equation}
	where $\Lambda_{\bm{n},r}$ is the block diagonal matrix
	\[
	\Lambda_{\bm{n},r} := \underset{\bm{j=0,\ldots,n-1}}{\operatorname*{diag}}\left( p\left( \frac{2\pi \bm{j}}{\bm{n}} \right) \right),
	\]
	whose $r \times r$ diagonal blocks are the values of
	\[
	p(\bm{\theta}) := \sum_{\bm{k=0}}^{\bm{n-1}} A_{\bm{k}}\, e^{i \bm{k} \cdot \bm{\theta}} \in \mathbb{C}^{r\times r}, \qquad \bm{k} \cdot \bm{\theta} = k_1 \theta_1 + \cdots + k_d \theta_d.
	\]
\end{theorem}

The proof of this theorem relies on finding a diagonalization of the uni-level generator $Z_n$, which can be obtained explicitly by computing its eigenvalues and eigenvectors. The block-diagonalization for $C_{n,r}$ then follows directly from equation~(\ref{2.2}), and the extension to the multi-level case is achieved by exploiting the properties of the Kronecker product.

Theorem \ref{thm:circulant} offers numerous significant discoveries. 
Every multi-level $r$-block circulant matrix is unitarily similar to the block diagonal matrix $\Lambda_{\bm{n},r}$; in particular, it is normal if and only if each $r \times r$ block $p(2\pi \bm{j}/\bm{n})$ is normal, which for $r=1$ is always the case. Secondly, the spectrum of \( C_{\bm{n},r} \) is the union of the spectra of these blocks,
\[ \Lambda\big( C_{\bm{n},r}\big) = \bigcup_{\bm{j}=\bm{0}}^{\bm{n}-\bm{1}} \Lambda\big( p(\bm{\theta}_{\bm{j}})\big), \qquad p(\bm{\theta}_{\bm{j}}) = \sum_{\bm{k=0}}^{\bm{n-1}} A_{\bm{k}} e^{i \bm{k} \cdot \bm{\theta}_{\bm{j}}}, \quad \bm{\theta_j} := \frac{2\pi \bm{j}}{\bm{n}}, \quad \bm{j} = \bm{0}, \ldots, \bm{n} - \bm{1}, \]
which means that the blocks $p(\bm{\theta_j})$ can be computed by performing the DFT of the first block column of $C_{\bm{n},r}$. The corresponding invariant subspaces are spanned by the columns of $F_{\bm{n},r}$. Moreover, for $r=1$ circulant matrices form an
algebra that is closed under inversion, and every matrix in this algebra is uniformly diagonalized by $F_{\bm{n}}$; for $r>1$ the block circulant matrices are simultaneously block-diagonalized by $F_{\bm{n},r}$, and the algebraic operations are inherited from those of the blocks.

Last but not least, many computations involving block circulant matrices can be significantly {accelerated by using the FFT.}  The computational cost of the matrix-vector product $C_{\bm{n},r} \bm{x}$, where $\bm{x} \in \mathbb{C}^{N(\bm{n})r}$, is reduced from the conventional $O\big((N(\bm{n})r)^2\big)$ operations necessary for dense matrices to only $O\big(r N(\bm{n}) \log N(\bm{n})+r^2 N(\bm{n})\big)$ operations, by using the following approach:
\begin{enumerate}
	
	\item  Compute \( \bm{y} = F_{\bm{n},r}^* \bm{x} \) via the FFT applied to each of the $r$ components, in \( O(r N(\bm{n}) \log N(\bm{n})) \) flops; 
	\item  Compute \( \bm{z} = \Lambda_{\bm{n},r} \bm{y} \), that is, $N(\bm{n})$ products of size $r$, in \( O(r^2 N(\bm{n})) \) flops; 
	\item  Compute \( F_{\bm{n},r} \bm{z} \) via the FFT, in \( O(r N(\bm{n}) \log N(\bm{n})) \) flops.
\end{enumerate}
In the same way, $C_{\bm{n},r}^{-1} \bm{x}$ can be computed as $F_{\bm{n},r} \Lambda_{\bm{n},r}^{-1} F_{\bm{n},r}^* \bm{x}$, where the inversion of $\Lambda_{\bm{n},r}$ amounts to solving $N(\bm{n})$ independent linear systems of size $r$. Similar techniques apply to the most relevant matrix operations, such as
matrix--matrix products, matrix inversion, the solution of linear systems, and the computation of the spectrum, which makes block circulant matrices
exceptionally efficient and powerful in practical applications. For a fixed $r$, all the above costs reduce to $O(N(\bm{n}) \log N(\bm{n}))$, with a hidden constant depending on $r$.

\begin{remark}
	This algorithm can be adapted to accommodate general multi-level block Toeplitz matrices through the use of an immersion (or embedding) technique. The central idea is to embed the original Toeplitz matrix into a larger block circulant matrix, and then perform the matrix-vector multiplication using the circulant structure as described above. The required result is then obtained by appropriately extracting the relevant portion from the output. By adopting this strategy, the computational complexity of the block Toeplitz matrix-vector product is also reduced to $O(r N(\bm{n}) \log N(\bm{n})+r^2 N(\bm{n}))$.
\end{remark}
\begin{remark}
	Finally, it is important to note that each step of the process is fully parallelizable. Therefore, in a parallel model of computation, the theoretical computational cost of all the aforementioned algorithms can be further reduced to $O(\log N(\bm{n}))$, for a fixed block size $r$.
\end{remark}

\section{\texorpdfstring{$\omega$-circulant matrices}{omega-circulant matrices}}\label{omega cir mat}
By introducing a parameter $\omega \in \mathbb{C} \setminus \{0\}$, the class of uni-level circulant matrices can be extended to the so-called
\emph{$\omega$-circulant matrices}. In this thesis, we restrict our attention to the uni-level case and do not consider multi-level $\omega$-circulant matrices.
These matrices have the following structure:
\[
\left[C^{\omega}_n\right]_{i,j} :=
\begin{cases}
	a_{(i-j) \bmod n}      & \text{if } i > j, \\
	\omega\, a_{(i-j) \bmod n} & \text{if } i \leq j,
\end{cases}
\qquad i, j = 1, \ldots, n.
\]
where $a_k \in \mathbb{C}$ for $k = 0, \ldots, n-1$. Here, the parameter $\omega$ acts as a scaling factor for certain matrix entries, setting these matrices apart from the standard circulant case.

Essentially, all entries above and on the main diagonal are multiplied by the factor $\omega$, resulting in the following structure:
\[
C_n^\omega =
\begin{bmatrix}
	a_0 & \omega a_{n-1} & \omega a_{n-2} & \cdots & \cdots & \omega a_{1} \\
	a_1 & a_0 &\omega a_{n-1} & \omega a_{n-2}&  & \vdots \\
	a_2 & a_1 & \ddots &\ddots & \ddots & \vdots \\
	\vdots & a_2&\ddots&\ddots&\ddots& \omega a_{n-2}\\
	\vdots &  & \ddots & \ddots & \ddots & \omega a_{n-1} \\
	a_{n-1} & \cdots & \cdots & a_{2} & a_1&a_{0}
\end{bmatrix}
\in \mathbb{C}^{n\times n}.
\]
Note that $C_n^{\omega}$ still retains a Toeplitz structure, and that for $\omega = 1$ it reduces to the standard circulant matrix. To generalize equation (\ref{2.2}), we define the $\omega$-circulant generator as
\[
\left[ Z_n^{\omega} \right]_{i,j} :=
\begin{cases}
	\delta_{(i-j-1) \bmod n}      & \text{if } i > j, \\
	\omega \, \delta_{(i-j-1) \bmod n} & \text{if } i \leq j,
\end{cases}
\qquad i, j = 1, \ldots, n.
\]
Equivalently,
\[
\left[ Z_n^{\omega} \right]_{i,j} =
\begin{cases}
	1     & \text{if } i - j = 1, \\
	\omega & \text{if } i = 1 \wedge j = n, \\
	0     & \text{otherwise},
\end{cases}
\qquad i, j = 1, \ldots, n.
\]
In other words, $Z_n^{\omega}$ has all ones on the first superdiagonal and $\omega$ in the top-right corner:
\[
Z_n^{\omega} =
\begin{bmatrix}
	0      &        &        &        &       &   \omega   \\
	1      & \ddots &        &        &       &        \\
	& \ddots & \ddots &        &       &        \\
	&        & \ddots & \ddots &       &        \\
	&        &        & \ddots & \ddots&        \\
	&        &        &        &1      &   0    \\
\end{bmatrix}
\in \mathbb{C}^{n\times n}.
\]
In this framework, $C_n^{\omega}$ can be expressed as a linear combination of powers of $Z_n^{\omega}$:
\[
C_n^{\omega} = \sum_{k=0}^{n-1} a_k (Z_n^{\omega})^k.
\]
A diagonalization formula, paralleling equation (\ref{2.3}), naturally extends to the case of $\omega$-circulant matrices:
\[
C_n^{\omega} = F_n^{\omega} \Lambda_n^{\omega} (F_n^{\omega})^{-1},
\]
where the diagonal matrix containing the eigenvalues of $C_n^{\omega}$ is given by
\[
\Lambda_n^{\omega} := \underset{j=0,\ldots,n-1}{\operatorname*{diag}}\left( p\left( \frac{2\pi j}{n} \right) \right)
\]
with
\[
p(\theta) := \sum_{k=0}^{n-1} \omega^{\frac{k}{n}} a_k e^{ik\theta},
\]
and the adapted Fourier matrix containing  the eigenvectors is defined as
\[
F_n^{\omega} := D_n^{\omega} F_n, \qquad D_n^{\omega} := \underset{k=0,\ldots,n-1}{\operatorname*{diag}}\left( \omega^{-\frac{k}{n}} \right).
\]
Observe that the eigenvalues of $C^{\omega}_n$ are obtained by computing the
matrix--vector product $F^{\omega}_n \mathbf{c}$, where $\mathbf{c}$ denotes the first column of
$C^{\omega}_n$.

Extending the computation of the matrix--vector product to $\omega$-circulant matrices, as well as adapting all the algorithms discussed in the previous
subsection, is immediate and requires only minor modifications:
\begin{enumerate}
	\item Compute $\bm{y} = (D_n^{\omega} F_n)^{-1} \bm{x} = F_n^* (D_n^{\omega})^{-1} \bm{x}$ using the FFT, which requires $O(n \log n)$ flops;
	\item Compute $\bm{z} = \Lambda_n^{\omega} \bm{y}$, in $O(n)$ flops;
	\item Compute $D_n^{\omega} F_n \bm{z}$ via the FFT, again in $O(n \log n)$ flops.
\end{enumerate}
Overall, this sequence of operations leads to a computational cost of $O(n \log n)$ in a sequential implementation, while in a parallel environment it requires $O(\log n)$ steps. Nevertheless, the presence of the parameter $\omega$ may affect numerical stability. Since $F_n$ is unitary, the conditioning of the transformation $D_n^{\omega} F_n$ depends exclusively on the diagonal scaling $D_n^{\omega}$. More precisely, as $n \to +\infty$, the condition number of $D_n^{\omega}$ approaches $\max\{|\omega|,|\omega|^{-1}\}$, and therefore becomes
unbounded whenever $\omega \to 0$ or $|\omega| \to \infty$.
\smallskip

\noindent Recent developments have extended the use of \(\omega\)-circulant preconditioners to more complex applications. For instance, Fung and Hon \cite{FungHon2024} introduced block \(\omega\)-circulant preconditioners for parabolic optimal control problems, demonstrating their effectiveness in accelerating the convergence of iterative solvers for large-scale all-at-once systems. More recently, Hon, Fung, and Lin \cite{HonFungLin2026} developed an optimal diagonalization-based preconditioner for parabolic optimal control problems, further showcasing the versatility of circulant-based approaches. These works highlight the continued relevance of \(\omega\)-circulant structures in modern numerical linear algebra.
\begin{remark}
		In the standard circulant case of Theorem~\ref{thm:circulant}, the Fourier matrix
		\(F_n\) is unitary, so \(F_n^{*}=F_n^{-1}\) and the two notations
		\(F_n\Lambda_n F_n^{*}\) and \(F_n\Lambda_n F_n^{-1}\) coincide.
		For the \(\omega\)-circulant case, the adapted Fourier matrix
		\(F_n^{\omega}=D_n^{\omega}F_n\) is not unitary unless \(|\omega|=1\),
		since \(D_n^{\omega}\) is a non-trivial diagonal scaling. Hence
		\((F_n^{\omega})^{*}\neq(F_n^{\omega})^{-1}\) in general, and the
		diagonalization must be written as
		\[
		C_n^{\omega}=F_n^{\omega}\Lambda_n^{\omega}(F_n^{\omega})^{-1}.
		\]
\end{remark}
\subsection{\texorpdfstring{ Block $\omega$-circulant matrices}{Block omega-circulant matrices}}
To set the stage for certain practical applications to the solution of linear systems, we briefly extend the definition and diagonalization of $\omega$-circulant matrices to the block case. For $n,r \in \mathbb{N}$, we define an $r$-block $\omega$-circulant matrix $C_{n,r}^{\omega} \in \mathbb{C}^{^{nr \times nr}}$ as follows:
\[
\left[ C_{n,r}^{\omega} \right]_{i,j} :=
\begin{cases}
	A_{(i-j) \bmod n}         & \text{if } i > j, \\
	\omega\, A_{(i-j) \bmod n} & \text{if } i \leq j,
\end{cases}
\qquad i, j = 1, \ldots, n,
\]
where $A_k \in \mathbb{C}^{r\times r}$ for $k = 0, \ldots, n-1$ are $n$ general (not necessarily structured) matrices. Explicitly:

\[
C_{n,r}^\omega =
\begin{bmatrix}
	A_0 & \omega A_{n-1} & \omega A_{n-2} & \cdots & \cdots & \omega A_{1} \\
	A_1 & A_0 &\omega A_{n-1} & \omega A_{n-2}&  & \vdots \\
	A_2 & A_1 & \ddots &\ddots & \ddots & \vdots \\
	\vdots & A_2&\ddots&\ddots&\ddots& \omega A_{n-2}\\
	\vdots &  & \ddots & \ddots & \ddots & \omega A_{n-1} \\
	A_{n-1} & \cdots & \cdots & A_{2} & A_1& A_{0}
\end{bmatrix}.
\]
A compact representation for $C_{n,r}^{\omega}$ is given by
\[
C_{n,r}^{\omega} = \sum_{j=0}^{n-1} (Z_n^{\omega})^j \otimes A_j.
\]
Moreover, $C_{n,r}^{\omega}$ can be block-diagonalized as
\[
C_{n,r}^{\omega} = F_{n,r}^{\omega} \Lambda_{n,r}^{\omega} (F_{n,r}^{\omega})^{-1},
\]
where
\[
F_{n,r}^{\omega} := F_n^{\omega} \otimes I_r,
\]
and
\[
\Lambda_{n,r}^{\omega} :=  \underset{j=0,\ldots,n-1}{\operatorname*{diag}}\left( p\left( \frac{2\pi j}{n} \right) \right)
\]
with
\[
p(\theta) := \sum_{k=0}^{n-1} \omega^{\frac{k}{n}} A_k e^{ik\theta} \in \mathbb{C}^{r\times r}.
\]
Note that $C_{n,r}^{\omega}$ is not a two-level $\omega$-circulant matrix, since the matrices $A_k$ do not necessarily possess a specific structure and the function $p(\theta)$ is an $r \times r$ matrix rather than a scalar. Nevertheless, when $r$ is a fixed integer, the matrix-vector product can still be computed in $\mathcal{O}(n \log n)$ floating-point operations in a sequential setting, or in $\mathcal{O}(\log n)$ steps in a parallel computational model, where the hidden constant in the $\mathcal{O}$-notation depends on $r$. If $r$ varies and tends to infinity and each inner matrix $A_k$ admits a diagonalization with cost $\mathcal{O}(r \log r)$, then the overall computational cost becomes $\mathcal{O}(nr \log(nr))$ flops in the sequential case, or $\mathcal{O}(\log(nr))$ parallel steps. For a review on recent applications concerning \(r\)-block \(\omega\)-circulant structures, refer to \cite{Diaz} and references therein.
\section{\texorpdfstring{$\tau$}{tau} matrices}\label{tau-mat}
A uni-level $r$-block $\tau$ matrix is defined as a square matrix
$S_{n,r} \in \mathbb{R}^{nr \times nr}$, partitioned into $r \times r$ blocks $S_{i,j}$, whose blocks satisfy the so-called \emph{cross-sum condition}
\[
S_{i-1,j} + S_{i+1,j} = S_{i,j-1} + S_{i,j+1},
\qquad i,j = 1,\ldots,n,
\]
where we assume the boundary conditions
\[
S_{n+1,j} = S_{i,n+1} = S_{0,j} = S_{i,0} = O_r.
\]
For $r=1$ the blocks reduce to the scalar entries $s_{i,j}$ and the classical definition is recovered. This relation implies the block-centrosymmetry of $S_{n,r}$, that is, symmetry with respect to the center of the matrix at the block level. As a consequence, $S_{n,r}$ is both symmetric and persymmetric, namely symmetric with respect to the
main anti-diagonal, whenever the generating blocks are symmetric. Furthermore, all the blocks of $S_{n,r}$ are uniquely determined by the blocks in its first block row or, equivalently, in its first block column. Additional details on this definition and the resulting
properties can be found in~\cite{BiniCapovani1983}.

Much like circulant matrices, there is a generator for the set of $\tau$ matrices, defined by
\[
[W_n]_{i,j} := \delta_{|i-j|-1} =
\begin{cases}
	1 & \text{if } i - j = \pm 1, \\
	0 & \text{otherwise},
\end{cases}
\qquad i, j = 1, \ldots, n,
\]
or, in explicit matrix form
\[
W_n =
\begin{bmatrix}
	0      &     1   &        &        &       &       \\
	1      & 0 &  1      &        &       &        \\
	& 1 & \ddots &   \ddots     &       &        \\
	&        & \ddots & \ddots & \ddots      &        \\
	&        &        & \ddots & \ddots& 1       \\
	&        &        &        &1      &   0    \\
\end{bmatrix}
\in\mathbb{R}^{n\times n}.
\]
The $r$-block $\tau$ algebra is defined as the set of matrices generated by $W_n$ in the block sense, meaning that every $r$-block $\tau$ matrix can be represented as
\[
S_{n,r} = \sum_{k=0}^{n-1} W_n^{k} \otimes A_k, \qquad A_k \in \mathbb{R}^{r\times r},
\]
that is, as a linear combination, with matrix coefficients, of the powers $W_n^0 = I_n$, $W_n$, $W_n^2$, $\ldots$, $W_n^{n-1}$. This representation is equivalent to the cross-sum condition stated above, since the latter acts independently on each of the $r^2$ entry positions inside the blocks. As is characteristic of matrix algebras, for $r=1$ this set is closed under inversion, a property that follows from the Hamilton--Cayley theorem~\cite{Bhatia1997}. Furthermore, for any $n, k \in \mathbb{N}$, the matrix $W_n^k$ is symmetric and is invertible if and only if $n$ is even.

This construction can be extended to the multi-level setting by defining the multi-level generators for $\bm{n} \in \mathbb{N}^d$ as
\begin{equation}
	\begin{aligned}
		W_{\bm{n}}^{(1)} &= W_{n_{1}} \otimes \left( \bigotimes_{l=2}^{d} I_{n_{l}} \right), \\[6pt]
		W_{\bm{n}}^{(2)} &= I_{n_{1}} \otimes W_{n_{2}} \otimes \left( \bigotimes_{l=3}^{d} I_{n_{l}} \right), \\[6pt]
		&\;\;\vdots \\[6pt]
		W_{\bm{n}}^{(d)} &= \left( \bigotimes_{l=1}^{d-1} I_{n_{l}} \right) \otimes W_{n_{d}}.
	\end{aligned}
\end{equation}
Accordingly, a multi-level $r$-block $\tau$ matrix admits the representation
\[
S_{\bm{n},r} = \sum_{\bm{k}=\bm{0}}^{\bm{n-1}} \left( W_{n_1}^{k_1} \otimes \cdots \otimes W_{n_d}^{k_d} \right) \otimes A_{\bm{k}}, \qquad A_{\bm{k}} \in \mathbb{R}^{r\times r}.
\]

Matrices $S_{\bm{n},r}$ that belong to the multi-level $r$-block $\tau$ algebra are referred to as multi-level block $\tau$ matrices.

\subsection{Spectral properties and the discrete sine transform}
In close analogy with the diagonalization theory for circulant matrices, $\tau$ matrices can also be diagonalized by means of a related fast transform.
\smallskip
More precisely, the Discrete Sine Transform (DST) is the linear mapping that associates with each vector $\bm{x} \in \mathbb{C}^n$ the vector $Q_n \bm{x}$, where the matrix $Q_n$, known as the sine transform matrix,
is defined as

\[
Q_n := \sqrt{\frac{2}{n + 1}} \left[ \sin \left( \frac{r s \pi}{n + 1} \right) \right]_{r,s = 1}^{n}.
\]
It is closely connected to the Fast Sine Transform, which enables the DST to be computed in $\mathcal{O}(n \log n)$ real arithmetic operations, requiring roughly half the computational cost of the
well-known FFT \cite{vanloan1992}. The matrix $Q_n$ is real, symmetric, and orthogonal. Its multilevel extension, for $\bm{n} \in \mathbb{N}^d$,
is defined as
\[
Q_{\bm{n}} := Q_{n_1} \otimes \cdots \otimes Q_{n_d},
\]
and the $d$-level DST is given by the matrix--vector product $Q_{\bm{n}}\bm{x}$, with $\bm{x} \in \mathbb{C}^{N(\bm{n})}$. By exploiting the properties of the Kronecker product, it follows that $Q_{\bm{n}}$ is a real, symmetric, and orthogonal matrix, and therefore
$
Q_{\bm{n}} = Q_{\bm{n}}^{T} = Q_{\bm{n}}^{-1}.
$
Its $r$-block counterpart is defined as
\[
Q_{\bm{n},r} := Q_{\bm{n}} \otimes I_r,
\]
which inherits the same properties, so that $Q_{\bm{n},r} = Q_{\bm{n},r}^{T} = Q_{\bm{n},r}^{-1}$.
As a consequence, we obtain the following spectral decomposition for multi-level block $\tau$ matrices.

\begin{theorem}
	For all $\bm{n} \in \mathbb{N}^d$ and $r \in \mathbb{N}$, it holds that
	\[
	S_{\bm{n},r} = Q_{\bm{n},r}\, \Lambda_{\bm{n},r}\, Q_{\bm{n},r},
	\]
	where $\Lambda_{\bm{n},r}$ is block diagonal and its $r \times r$ diagonal blocks are given by the $d$-level DST of the first block column (or, equivalently, the first block row) of $S_{\bm{n},r}$.
\end{theorem}
All the properties established after Theorem~\ref{thm:circulant} for block circulant matrices extend to the multi-level block $\tau$ setting, once the necessary adaptations are
made. In particular, any multi-level block $\tau$ matrix is real and, when the generating blocks $A_{\bm{k}}$ are symmetric, its spectrum is the union of the spectra of the $r \times r$ diagonal blocks of $\Lambda_{\bm{n},r}$, which are obtained in explicit form by applying the DST to its first block column, while the corresponding invariant subspaces are
spanned by the columns of $Q_{\bm{n},r}$. Moreover, the fundamental matrix operations, including matrix--vector multiplication, inversion, the solution
of linear systems, and spectral computation, can be performed with a sequential computational cost of $O(N(\bm{n}) \log N(\bm{n}))$, or in
$O(\log N(\bm{n}))$ steps within a parallel computational framework, with hidden constants depending on the block size $r$.
\vspace{0.3cm}

\noindent Recently, \(\tau\) preconditioners have gained renewed attention for solving FDEs. Li and Hon \cite{LiHon2024} developed multilevel tau preconditioners for symmetrized multilevel Toeplitz systems arising from space-FDEs, providing a comprehensive GLT-based analysis of their spectral properties. Building on this work, Hon, Li, Sormani, Krause, and Serra-Capizzano \cite{HonLiSormani2025} introduced symbol-based multilevel block \(\tau\) preconditioners for general multilevel block Toeplitz systems, establishing a rigorous GLT framework for analyzing their asymptotic spectral behavior. These developments demonstrate that \(\tau\)-preconditioners remain a powerful tool for structured linear systems, particularly in the fractional diffusion setting.
\section{Hankel matrices}\label{hank matr}
Given $\bm{n} \in \mathbb{N}^d$ and $r \in \mathbb{N}$, a matrix of the form
\[
H_{\bm{n}} := \left[ A_{\bm{i}+\bm{j}} \right]_{\bm{i},\bm{j}=\bm{1}}^{\bm{n}}
\in \mathbb{C}^{{N(\bm{n})r}\times{N(\bm{n})r}},
\]
with blocks $A_{\bm{k}} \in \mathbb{C}^{r \times r}$, where the block in position $(\bm{i}, \bm{j})$ depends only on the sum $\bm{i} + \bm{j}$, is known as a multi-level block (or $d$-level $r$-block) Hankel matrix.

In practice, it can be interpreted as a mirrored Toeplitz structure. In the one-dimensional case, i.e., for $d=1$ and $\bm{n}=n$, a one-level (uni-level) $r$-block Hankel matrix is defined by the property that its blocks are constant along the block anti-diagonals:
\[
H_n =
\big[ A_{i+j} \big]_{i,j=1}^{n}
=
\begin{bmatrix}
	A_2 & A_3 & A_4 & \cdots & A_{n+1} \\
	A_3 & A_4 &  & \reflectbox{$\ddots$} & \vdots \\
	A_4 &  & \reflectbox{$\ddots$} &  & A_{2n-2} \\
	\vdots & \reflectbox{$\ddots$} &  & A_{2n-2} & A_{2n-1} \\
	A_{n+1} & \cdots & A_{2n-2} & A_{2n-1} & A_{2n}
\end{bmatrix}
\in \mathbb{C}^{nr \times nr}.
\]
Just like Toeplitz matrices, $H_n$ is completely determined by its $2n-1$ blocks $A_k$ for $k=2,\ldots,2n$, or, equivalently, by its first block row and last block column.

For $d=2$ and $\bm{n} = (n_1, n_2)$, a 2-level $r$-block Hankel matrix takes the form
\[
H_{\bm{n}} =
\big[ A_{\bm{i}+\bm{j}} \big]_{\bm{i,j=1}}^{\bm{n}}
=
\bigg[ \left[A_{(i_1+j_1,\, i_2+j_2)} \right]_{i_2,j_2=1}^{n_2}\bigg]_{i_1,j_1=1}^{n_1}
=
\big[ \mathcal{A}_{i_1+j_1} \big]_{i_1,j_1=1}^{n_1}
\]
\[
=
\begin{bmatrix}
	\mathcal{A}_2 & \mathcal{A}_3 & \mathcal{A}_4 & \cdots & \mathcal{A}_{n_1+1} \\
	\mathcal{A}_3 & \mathcal{A}_4 &  & \reflectbox{$\ddots$} & \vdots \\
	\mathcal{A}_4 &  & \reflectbox{$\ddots$} &  & \mathcal{A}_{2n_1-2} \\
	\vdots & \reflectbox{$\ddots$} &  & \mathcal{A}_{2n_1-2} & \mathcal{A}_{2n_1-1} \\
	\mathcal{A}_{n_1+1} & \cdots & \mathcal{A}_{2n_1-2} & \mathcal{A}_{2n_1-1} & \mathcal{A}_{2n_1}
\end{bmatrix}
\in \mathbb{C}^{n_{1}n_{2}r \times n_{1}n_{2}r},
\]
where each block $\mathcal{A}_k$ is itself a 1-level $r$-block Hankel matrix,\\
\[\mathcal{A}_k =
\big[ A_{(k,\, i_2+j_2)} \big]_{i_2,j_2=1}^{n_2}\]

\[
=
\begin{bmatrix}
	A_{(k,2)} & A_{(k,3)} & A_{(k,4)} & \cdots & A_{(k,n_2+1)} \\
	A_{(k,3)} & A_{(k,4)} &  & \reflectbox{$\ddots$} & \vdots \\
	A_{(k,4)} &  & \reflectbox{$\ddots$} &  & A_{(k,2n_2-2)} \\
	\vdots & \reflectbox{$\ddots$} &  & A_{(k,2n_2-2)} & A_{(k,2n_2-1)} \\
	A_{(k,n_2+1)} & \cdots & A_{(k,2n_2-2)} & A_{(k,2n_2-1)} & A_{(k,2n_2)}
\end{bmatrix}
\in \mathbb{C}^{n_2 r \times n_2 r}.
\]
In the multi-level setting, a $d$-level $r$-block Hankel matrix admits a recursive block structure, in the sense that it can be interpreted as a block Hankel matrix whose blocks are $(d-1)$-level $r$-block Hankel matrices. As in the
one-level case, the entire matrix is uniquely specified by the blocks in its first block row and last block column.

To obtain a compact representation of $H_{\bm{n}}$, we introduce, for $n \in \mathbb{N}$ and $k \in \mathbb{Z}$, the matrix
\[
\big[ K_n^{(k)} \big]_{i,j} := \delta_{i+j - k} =
\begin{cases}
	1 & \text{if } i + j = k, \\
	0 & \text{otherwise},
\end{cases}
\qquad i, j = 1, \ldots, n,
\]
which places ones on the $(k-1)$-th anti-diagonal and zeros elsewhere.

The multi-level generalization for $\bm{n} \in \mathbb{N}^d$ and $\bm{k} \in \mathbb{Z}^d$ is given by
\[
K_{\bm{n}}^{(\bm{k})} := K_{n_1}^{(k_1)} \otimes K_{n_2}^{(k_2)} \otimes \cdots \otimes K_{n_d}^{(k_d)}.
\]
This leads to the following compact expression for $H_{\bm{n}}$:

\begin{equation}\label{2.4}
	H_{\bm{n}} = \left[ A_{\bm{i}+\bm{j}} \right]_{\bm{i},\bm{j}=\bm{1}}^{\bm{n}}
	= \sum_{\bm{k}=\bm{2}}^{2\bm{n}} K_{\bm{n}}^{(\bm{k})} \otimes A_{\bm{k}}. 
\end{equation}
In analogy with the Toeplitz and circulant cases, the above representation follows from a direct verification of the individual matrix blocks.
\smallskip

\noindent Block Hankel matrices are closely related to block Toeplitz matrices through the anti-identity matrix; see \cite{FasinoTilli2000} for spectral properties and applications.

	\chapter{GLT matrix sequences and their algebraic structure}\label{GLT}
We are now in a position to introduce the theory of GLT matrix sequences, which provides the central spectral framework for the analysis carried out in this thesis. The GLT machinery
offers a unified language for describing the asymptotic behavior of large structured matrices generated by the discretization of differential and nonlocal operators and, at the same time, supplies the theoretical tools required for
the design of efficient numerical solvers.

In the present setting we consider the multi-level block case, corresponding to $d$-level GLT matrix sequences with $r \times r$ matrix-valued symbols, where $d \ge 1$ denotes the number of structural levels and $r \ge 1$ is the block size. For fixed $d$ and $r$, the associated GLT class is characterized by
symbols given by measurable functions
\[
\kappa : [0,1]^d \times [-\pi,\pi]^d \longrightarrow \mathbb{C}^{r\times r},
\]
which encode the asymptotic spectral distribution of the matrix sequence. While the scalar case is thoroughly treated in the monographs \cite{garoni2017,garoni2018}, the matrix-valued framework adopted here follows the more general and flexible approach developed in
\cite{barbarino2020block1d,barbarino2020blockmulti}. This setting is essential for the study of block discretizations and coupled systems arising in PDEs and FDEs. Moreover, the theory continues to evolve, as witnessed by recent developments on reduced and rectangular GLT
sequences \cite{Barbarino2022,Barbarino2022RectangularGLT}.

Rather than presenting the full constructive definition of the multi-level GLT class, we focus on the fundamental operative properties that uniquely characterize it and that are most relevant for numerical applications. These properties, commonly referred to as the \emph{GLT axioms}, provide a
complete algebraic and spectral description of the class and will be systematically used in the subsequent analysis. They are introduced and discussed in Section~\ref{ast algebra}. In addition, we highlight how the GLT framework naturally guides the construction of effective preconditioning strategies, thereby establishing the link between asymptotic spectral theory and fast iterative solvers.

\section{The GLT class}\label{GLT Clas}
Let us first review some basic aspects to better understand GLT sequences. The theory of LT sequences provides a powerful framework for analyzing the asymptotic singular value and eigenvalue distributions of structured matrix sequences that arise in numerical analysis. In the context of PDE discretization, the matrices resulting from numerical approximations often fail to retain the classical Toeplitz structure due to the presence of variable coefficients. Traditionally, Toeplitz matrices have played a central role in problems involving translation-invariant differential operators, where constant coefficients ensure a globally repeating structure. Instead, the resulting matrices exhibit a Toeplitz-like structure that varies smoothly along the diagonals as the matrix size increases, leading to the broader class of LT sequences \cite{Tilli1998, Tyrtyshnikov1996}.\\
A classical Toeplitz matrix generated by a function $f: [-\pi, \pi]^d \to \mathbb{C}^{r \times r}$ belonging to $L^1([-\pi, \pi]^d)$, where the constant diagonal entries correspond to the Fourier coefficients of $f$. However, for differential operators with non-constant coefficients, the resulting matrices do not maintain this global Toeplitz structure but instead belong to the LT class. Specifically, an LT matrix $A_{\bm{n}}$ is obtained by modulating the entries of a Toeplitz matrix using an almost everywhere (a.e.) continuous function $a(x): [0,1]^d \to \mathbb{C}^{r \times r}$, scaling the entries along the same diagonal
$ a\left(\frac{\bm{i}}{\bm{n}}\right)$ for $\bm{i}= \bm{1},\dots,\bm{n}.$ As $\bm{n} \to \infty$, the differences between consecutive diagonal elements vanish due to continuity, making the structure asymptotically Toeplitz-like. More precisely, this structure can be represented as the product of a Toeplitz matrix and a diagonal matrix, where the diagonal entries correspond to evaluations of the function $a(x)$, such that $A_{\bm{n}} = D_{\bm{n}}(a)T_{\bm{n}},$
where $D_{\bm{n}}(a)$ is a diagonal matrix whose entries are given by $D_{\bm{n}}(a) = \mathrm{diag}\left(a\left(\frac{\bm{i}}{\bm{n}}\right)\right)_{\bm{i}=\bm{1}}^{\bm{n}}.$
This formulation allows LT sequences to accurately capture local variations in variable-coefficient differential operators while preserving the structural properties of Toeplitz matrices. Formally, a LT sequence \(\{A_{\mathbf{n}}\}_{\mathbf{n}}\) can
be expressed through the product of the weight function \(a(x)\) and the
generating function \(f(\theta)\), leading to the LT symbol
\[
a \otimes f : [0,1]^d \times [-\pi, \pi]^d \to \mathbb{C}^{r \times r},
\]
where
\[
(a\otimes f)(x,\theta) := a(x)\,f(\theta),
\qquad (x,\theta)\in[0,1]^d\times[-\pi,\pi]^d.
\]
Here the symbol ``\(\otimes\)'' in \(a\otimes f\) denotes the pointwise matrix
product of the two matrix-valued functions, not the Kronecker product of
the two matrices. Since \(a(x),f(\theta)\in\mathbb{C}^{r\times r}\), the
product is well defined and belongs to \(\mathbb{C}^{r\times r}\), in
accordance with the target space of the LT symbol. This is
consistent with the use of ``\(\otimes\)'' in GLT~\ref{glt:kron}, where at the matrix level
``\(\otimes\)'' denotes the Kronecker product, while at the symbol level it
denotes the pointwise product of the symbols.

Building on this concept, a GLT matrix sequence is defined as the limit, in the sense of the a.c.s., as described in Definition~\ref{def:acs}, of a finite sum of LT sequences. The construction of GLT sequences was formalized through a.c.s. limits in \cite{SerraCapizzano2003} and refined to a \textbf{d-level $*$-algebra of structured matrix sequences} satisfying specific axioms, listed in \cite[Chapter~6, pp.~118--120]{garoni2018}. The theory of GLT sequences provides a fundamental framework for analyzing structured matrix sequences arising from discretizations of differential equations.   
\vspace{0.3cm}

\noindent A particularly relevant development in the spectral analysis of structured matrix sequences concerns symmetrized Toeplitz structures. Hon, Mursaleen, and Serra-Capizzano \cite{Hon2019} provided a note on the spectral distribution of symmetrized Toeplitz sequences, establishing precise distribution results for this important class of matrices. Building on this work, Ferrari, Furci, Hon, Mursaleen, and Serra-Capizzano \cite{FerrariFurciHon2019} studied the eigenvalue distribution of special \(2\times 2\) block matrix sequences, with applications to symmetrized Toeplitz structures. These results are particularly relevant for understanding the spectral behavior of preconditioned matrix sequences and for analyzing the convergence of iterative methods applied to non-Hermitian Toeplitz systems. The analysis presented in later chapters of this thesis, particularly concerning nonsymmetric Toeplitz systems, draws upon this theoretical framework.

\section{Matrix sequences with explicit or hidden (asymptotic) structure}
In this section, we introduce three types of matrix structures that serve as the fundamental building blocks of the GLT $*$-algebras. Specifically, for any positive integers $d$ and $r$, we consider the set of $d$-level $r$-block GLT matrix sequences. This set forms a $*$-algebra of matrix sequences, which is both maximal and isometrically equivalent to the maximal $*$-algebra of $2d$-variate $r \times r$ matrix-valued measurable functions (with respect to the Lebesgue measure) that are naturally defined over $[0,1]^d \times [-\pi,\pi]^d$; see~\cite{barbarino2020block1d, barbarino2020blockmulti, garoni2017, garoni2018} and references therein.

The reduced version of this structure plays a crucial role in approximating integro-differential operators, including their fractional versions, particularly when defined over general (non-Cartesian) domains. This concept was initially introduced in~\cite{SerraCapizzano2003,SerraCapizzano2006} and later extensively developed in~\cite{Barbarino2022}, where GLT symbols are again defined as measurable functions over $\Omega \times [-\pi, \pi]^d$, with $\Omega$ being Peano-Jordan measurable and contained within $[0,1]^d$. Additionally, the reduced versions also form maximal $*$-algebras that are isometrically equivalent to their corresponding maximal $*$-algebras of measurable functions.

These GLT $*$-algebras provide a rich framework of hidden (asymptotic) structures, built upon three fundamental classes of explicit algebraic structures: $d$-level $r$-block Toeplitz matrix sequences and sampling diagonal matrix sequences (discussed in Sections~\ref{blcktop} and~\ref{blckdiag}), along with asymptotic structures described by zero-distributed matrix-sequences (see Section~\ref{zero dist}). Notably, the latter class serves an analogous role to compact operators in relation to bounded linear operators, forming a two-sided ideal of matrix sequences within any of the GLT $*$-algebras.

\section{Zero-distributed sequences}\label{zero dist}

Zero-distributed sequences are defined as matrix sequences $\{A_n\}_n$ such that $\{A_n\}_n \sim_\sigma 0$.
Note that, for any $r \geq 1$, $\{A_n\}_n \sim_\sigma 0$ is equivalent to $\{A_n\}_n \sim_\sigma O_r$, where $O_r$ is the $r \times r$ zero matrix.
The following theorem, taken from~\cite{SerraCapizzano2001,garoni2017}, provides a useful characterization for detecting this type of sequence.

\begin{theorem}\label{lemma:zero-distr}
	Let $\{A_n\}_n$ be a matrix sequence, with $A_n$ of size $d_n$, and let $p \in [1,\infty]$, with $\|X\|_p$ being the Schatten $p$-norm of $X$, that is, the $\ell^p$ norm of the vector of its singular values. Let $\|\cdot\| = \|\cdot\|_\infty$ be the spectral norm. With the natural convention $1/\infty = 0$, we have:
	\begin{itemize}
		\item $\{A_n\}_n \sim_\sigma 0$ if and only if $A_n = R_n + N_n$ with $\mathrm{rank}(R_n)/d_n \to 0$ and $\|N_n\| \to 0$ as $n \to \infty$;
		\item $\{A_n\}_n \sim_\sigma 0$ if there exists $p \in [1,\infty]$ such that
		\[
		\frac{\|A_n\|_p}{(d_n)^{1/p}} \to 0 \quad \text{as } n \to \infty.
		\]
	\end{itemize}
\end{theorem}

As in Section \ref{acs}, the same definition can be given and the corresponding result (with obvious changes) holds when the involved matrix sequences show a multi-level structure. In that case, $n$ is replaced by $\bm{n}$ uniformly in $A_n$, $N_n$, $R_n$, $d_n$.

Further results connecting Schatten $p$-norms, unitarily invariant norms (see~\cite{bhatia2007}), variational characterization, and Toeplitz/Hankel structures can be found in~\cite{SerraCapizzanoTilli-LPO}.\\
\section{Multi-level block Toeplitz matrices}\label{blcktop}

Given $\bm{n} \in \mathbb{N}^d$, a matrix of the form
\[
[A_{\bm{i}-\bm{j}}]_{\bm{i},\bm{j}=\bm{1}}^{\bm{n}} \in \mathbb{C}^{{N(\bm{n})r}\times{N(\bm{n})r}},
\]
with blocks $A_{\bm{k}} \in \mathbb{C}^{r \times r}$, $\bm{k} \in \{- (\bm{n} - \bm{1}), \ldots, \bm{n} - \bm{1}\}$, is called a \textit{multi-level block Toeplitz matrix}, or, more precisely, a $d$-level $r$-block Toeplitz matrix.

Given a matrix-valued function $f : [ -\pi, \pi ]^d \to \mathbb{C}^{r \times r}$ belonging to $L^1( [ -\pi, \pi ]^d )$, the $\bm{n}$-th Toeplitz matrix associated with $f$ is defined as
\[
T_{\bm{n}}(f) := [ \hat{f}_{\bm{i}-\bm{j}} ]_{\bm{i},\bm{j}=\bm{1}}^{\bm{n}} \in \mathbb{C}^{{N(\bm{n})r}\times{N(\bm{n})r}} ,
\]
where
\[
\hat{f}_{\bm{k}} = \frac{1}{(2\pi)^d} \int_{ [ -\pi, \pi ]^d } f(\theta) e^{- i (\bm{k} \cdot \bm{\theta})} d\bm{\theta} \in \mathbb{C}^{r \times r}, \quad \bm{k} \in \mathbb{Z}^d,
\]
are the Fourier coefficients of $f$.

The family $\{ T_{\bm{n}}(f) \}_{ \bm{n} \in \mathbb{N}^d }$ is called the family of (multi-level block) Toeplitz matrices associated with $f$, which is referred to as the \emph{generating function}.
\section{Block diagonal sampling matrices}\label{blckdiag}
Given $d \geq 1$, $\bm{n} \in \mathbb{N}^d$ and a function $a : [0,1]^d \to \mathbb{C}^{r \times r}$, we define the multi-level block diagonal sampling matrix $D_{\bm{n}}(a)$ as the block diagonal matrix
\begin{equation}\label{eq:da}
	D_{\bm{n}}(a) = \underset{\bm{i}=\bm{1},\ldots,\bm{n}}{\operatorname*{diag}}\left( a\left(\frac{\bm{i}}{\bm{n}}\right) \right) \in \mathbb{C}^{{N(\bm{n})r}\times{N(\bm{n})r}}.    
\end{equation}
Essentially, $D_{\bm{n}}(a)$ is a diagonal matrix whose entries are the evaluations of $a(x)$ at the points of a uniform mesh over the domain $[0, 1]^d$, with the points arranged in lexicographical order. 
\section{The \texorpdfstring{$\ast$}{*}-algebra of \texorpdfstring{$d$}{d}-level \texorpdfstring{$r$}{r}-block GLT matrix sequences}\label{ast algebra}
Let $r \geq 1$ be a fixed integer. A multi-level $r$-block GLT sequence, or simply a GLT sequence if $r$ does not need to be specified, is a special multi-level $r$-block matrix sequence equipped with a measurable function $\kappa : [0, 1]^d \times [-\pi, \pi]^d \to \mathbb{C}^{r \times r}$, $d \geq 1$, called the \emph{symbol}. We write
$
\{A_n\}_n \sim_{\mathrm{GLT}} \kappa
$
to denote that \(\{A_n\}_n\) is a GLT sequence with symbol \(\kappa\). A fundamental property of the GLT symbol is that it is unique: if two symbols are associated with the same matrix sequence, then they must coincide almost everywhere.
\begin{proposition}
	If $\{A_{\bm{n}}\}_{\bm{n}} \sim_{\mathrm{GLT}} \kappa$ and $\{B_{\bm{n}}\}_{\bm{n}} \sim_{\mathrm{GLT}} \xi$, then $
	\kappa=\xi$ \text{a.e.}
\end{proposition}
The symbol plays a central role in determining the singular value distribution and, in the Hermitian setting, the eigenvalue distribution of GLT matrix sequences.

It can be proven that the set of multilevel block GLT sequences is the $\ast$-algebra generated by the three classes of sequences defined in Sections \ref{zero dist}, \ref{blcktop}, \ref{blckdiag}: zero-distributed sequences, multi-level block Toeplitz sequences, and block diagonal sampling matrix sequences. The GLT class satisfies several algebraic and topological properties, which are discussed in detail in \cite{barbarino2020block1d, barbarino2020blockmulti, garoni2017, garoni2018}. Here, we focus on the main operative properties listed below, which provide a complete characterization of GLT sequences, equivalent to the full constructive definition.
\begin{gltax} \label{glt:distrib}
	If $\{A_n\}_n \sim_\glt \kappa$, then $\{A_n\}_n \sim_\sigma \kappa$ in the sense of Definition \ref{def:distrib}, with $t=2d$ and $D=[0,1]^{d}\times[-\pi,\pi]^{d}$. If moreover each $A_n$ is Hermitian, then $\{A_n\}_n \sim_\lambda \kappa$, again in the sense of Definition \ref{def:distrib}, with $t=2d$.
\end{gltax}
\begin{gltax} \label{glt:symbols}
	It holds
	\begin{itemize}
		\item $\{T_\mathbf{n}(f)\}_\mathbf{n} \sim_\glt \kappa(\mathbf{x},\bm{\theta}) = f(\bm{\theta})$ for any $f: [-\pi,\pi]^d \to \mathbb{C}^{r \times r}$ with 
		$f \in L^1([-\pi,\pi]^d)$;
		\item $\{D_\mathbf{n}(a)\}_\mathbf{n} \sim_\glt \kappa(\mathbf{x},\bm{\theta}) = a(\mathbf{x})$ for any Riemann-integrable $a:[0,1]^d \to \mathbb{C}^{r \times r}$;
		\item $\{Z_n\}_n \sim_\glt \kappa(\mathbf{x},\bm{\theta}) = O_r$ if and only if $\{Z_n\}_n \sim_\sigma 0$, i.e., $\{Z_n\}_n$ is zero-distributed.
	\end{itemize}
\end{gltax}

\begin{gltax} \label{glt:algebra}
	If $\{A_n\}_n \sim_\glt \kappa$ and $\{B_n\}_n \sim_\glt \varsigma$, then
	\begin{itemize}
		\item $\{A_n^*\}_n \sim_\glt \bar{\kappa}$, where $A_n^*$ denotes the conjugate transpose of $A_n$ and $\bar{\kappa}$ the complex conjugate of $\kappa$;
		\item $\{\alpha A_n + \beta B_n\}_n \sim_\glt \alpha\kappa +\beta\varsigma$ for any $\alpha,\beta\in\C$;
		\item $\{A_n B_n\}_n \sim_\glt \kappa\varsigma$;
		\item $\{A_n^\dagger\}_n \sim_\glt\kappa^{-1}$ for any $\kappa$ invertible a.e., where $A_n^\dagger$ denotes the Moore-Penrose pseudoinverse of $A_n$.
	\end{itemize}
\end{gltax}
\begin{gltax}\label{glt:acs}
	\(\{A_{\bm{n}}\}_{\bm{n}} \sim_{\mathrm{GLT}} \kappa\) if and only if there exist \(\{B_{{\bm{n}},j}\}_{\bm{n}} \sim_{\text{GLT}} \kappa_j\) such that \(\{\{B_{{\bm{n}},j}\}_{\bm{n}}\}_j \xrightarrow{\text{a.c.s.wrt }j} \{A_{\bm{n}}\}_{\bm{n}}\) and \(\kappa_j \rightarrow \kappa\) in measure.
\end{gltax}
\begin{gltax}\label{glt:hermitian}
If $\{A_{\bm{n}}\}_n \sim_{\mathrm{GLT}} \kappa$ and $A_{\bm{n}} = X_{\bm{n}} + Y_{\bm{n}}$, where
\begin{itemize}
	\item every $X_{\bm{n}}$ is Hermitian,
	\item $ ||X_{\bm{n}}||, \, ||Y_{\bm{{n}}}|| \leq C$ for some constant $C$ independent of $n$,
	\item $N(\bm{n})^{-1} \|Y_{\bm{n}}\|_1 \to 0$,
\end{itemize}
then $\{A_{\bm{n}}\}_n \sim_\lambda \kappa$.
\end{gltax}
\begin{gltax}\label{glt:function}
	If $\{A_{\bm{n}}\}_n \sim_{\mathrm{GLT}} \kappa$ and each $A_{\bm{n}}$ is Hermitian, then $\{f(A_{\bm{n}})\}_n \sim_{\mathrm{GLT}} f(\kappa)$ for every continuous function $f : \mathbb{C} \to \mathbb{C}$.
\end{gltax}
\begin{gltax} \label{glt:kron}
	Let $\{A _{\bm{n}}\}_{\bm{n}} \sim_\glt \kappa$ and $\{B _{\bm{m}}\}_{\bm{m}} \sim_\glt \xi$ with
	\begin{align*}
		&\kappa(x_1,\theta_1) : [0,1]^{d_1}\times[-\pi,\pi]^{d_1} \to\C,
		\\
		&\xi(x_2,\theta_2) : [0,1]^{d_2}\times[-\pi,\pi]^{d_2}\to\C.
	\end{align*}
	
	Then, it holds 
	$$\{A _{\bm{n}}\}_{\bm{n}}\otimes\{B _{\bm{m}}\}_{\bm{m}}=\{A _{\bm{n}} \otimes B_{\bm{m}}\}_{\bm{n}\bm{m}} \sim_\glt \kappa\otimes\xi,$$ 
	where $(\kappa\otimes\xi): [0,1]^{d_1 + d_2}\times[-\pi,\pi]^{d_1 + d_2}\to\C$ is given by
	\begin{equation*}
		(\kappa\otimes\xi)(x_1,x_2,\theta_1,\theta_2) := \kappa(x_1,\theta_1)\xi(x_2,\theta_2),
	\end{equation*}
	i.e., the tensor product is a \((d_1+d_2)\)-level GLT matrix sequence with
	scalar symbol. In particular, when \(d_1=d_2=1\), the tensor product is a
	two-level GLT matrix sequence.
\end{gltax}
\vspace{0.3cm}

Note that, by \textbf{GLT \ref{glt:distrib}}, it is always possible to obtain the singular value distribution from the GLT symbol, while the eigenvalue distribution can only be deduced either if the involved
matrices are Hermitian or the related matrix-sequence is quasi-Hermitian in the sense of \textbf{GLT \ref{glt:hermitian}}.

\subsection{GLT theory and preconditioning: a practical perspective}
Before proceeding to the numerical applications in Part~II, it is useful to clarify how the abstract GLT framework directly informs preconditioner design. The central idea is as follows: if a preconditioned matrix sequence has GLT symbol equal to the constant function \(1\), then, by GLT~\ref{glt:distrib}, its singular values are asymptotically distributed as \(1\). In other words, the singular values cluster at \(1\) in the weak sense. If, in addition, the matrices are Hermitian, the eigenvalues also cluster at \(1\). For non-Hermitian matrices, eigenvalue clustering at \(1\) may still be obtained under suitable additional assumptions, such as normality or the existence of a low-rank perturbation structure that allows the application of perturbation arguments \cite{Barbarino2020, GoSe}. Such clustering is of significant practical importance: when the eigenvalues or singular values of the preconditioned matrix are clustered around \(1\), Krylov subspace methods such as GMRES or MINRES typically exhibit fast convergence. Thus, the GLT symbol serves as a powerful predictive tool: it not only describes the asymptotic spectral behavior of the coefficient matrices but also guides the construction of preconditioners whose symbols match those of the original matrices, thereby ensuring favorable convergence properties. Throughout the following chapters, this principle will be applied systematically to design and analyze preconditioners for ISPs of the FDEs.
	\chapter{Inverse source problems for a 2D generalized diffusion equation}\label{chapter:TFDE}

Having introduced in the previous chapters the analytical foundations of fractional operators, the spectral framework of matrix sequences, the main classes of structured matrices, and the GLT theory, we now turn to the study of inverse source problems (ISPs) for a two-dimensional generalized diffusion equation. This chapter is based on the work \cite{IlyasSerraCapizzano2025}, with additional details, examples, and reformulations included for completeness of the thesis presentation. This chapter is devoted to the identification of unknown source coefficients in diffusion equation governed by generalized time-fractional derivative of Caputo-type and nonlocal boundary conditions. Such models arise naturally in the description of anomalous transport processes in heterogeneous and complex media, where memory and spatial nonlocality play a dominant role \cite{Metzler2000,Sandev1,Kochubei}.

We investigate two ISPs based on the given final time and integral observational data and derive explicit series representations of the solutions using a bi-orthogonal system of functions. \cite{Kapustin,Kapustin1,Bazhlekova,IlyasSaif2023}. Rigorous existence and uniqueness results are established, and both ISPs are shown to be ill-posed in the sense of Hadamard, since small perturbations in the measured data may lead to unbounded reconstruction errors \cite{Prilepko}. Special cases of the diffusion equation are obtained by appropriate choices of the memory kernel in the time-fractional derivative, allowing the proposed framework to encompass tempered, multi-term, mono-term, and classical diffusion models within a unified setting \cite{Sandev2,Podlubny1998}. The results presented in this chapter provide a rigorous mathematical foundation for ISPs in generalized FDEs and contribute to a deeper understanding of anomalous transport processes governed by memory dependent and nonlocal dynamics.

The remainder of this chapter is organized as follows. Section~\ref{Bi-orth} presents the spectral problem associated with the system \eqref{problem}--\eqref{initialcdn}, which forms the analytical foundation for the subsequent inverse analysis. Section~\ref{Main IP-I} is devoted to the first inverse source problem, where the existence and uniqueness of the regular solution are established and its ill-posed nature is rigorously demonstrated. Section~\ref{INVERSEPROBLEM-II} investigates the second inverse source problem, proving the existence and uniqueness of a regular solution and analyzing its inherent instability. Section~\ref{exmp} provides several numerical examples that illustrate and validate the theoretical results.

\section{Mathematical formulation and preliminaries}
In this section, we introduce the generalized time-fractional derivative of Caputo-type and recall its main analytical properties. We then formulate the two dimensional generalized diffusion equation under consideration, discuss several important special cases of this equation, and define the mathematical formulation of the ISPs.
\subsection{Generalized Caputo-type fractional derivative}
In this study, we utilize the generalized time-fractional derivative of Caputo-type as formulated by Kochubei \cite{Kochubei}, presented in the following form
\begin{equation}\label{Kach}_{G}\mathcal{D}_{t}^{(\varrho)}u(x,y,t):=\int_{0}^{t}\varrho(t-\tau) \frac{\partial}{\partial \tau} u(x,y,\tau)d\tau,\end{equation}
where $\varrho(t)\in L^{1}_{loc} (\mathbb{R_{+}})$ represents a locally integrable, non-negative function commonly referred to as the memory kernel. The  Laplace transform of $\varrho(t)$, denoted $\widehat{\varrho}(s)$, exists for all $s>0$ and satisfies
\begin{equation}\label{sf}
	\widehat{\varrho}(s)\in \mathbb{SF},\quad\qquad\mbox{and}\qquad\quad \lim\limits_{s\to\infty} s \widehat{\varrho}(s)=+\infty,
\end{equation}
where $\mathbb{SF}$ denotes the class of Stieltjes functions, as discussed in \cite{Bazhlekova}. The class of $\mathbb{SF}$ includes all functions on $\mathbb{R_+}$ that have the following representation
\begin{equation}\label{repres}
	\widehat{\phi}(s)=\frac{a}{s}+b+\int_{0}^{\infty}\mathrm{e}^{-s\tau}\varphi(\tau)d\tau,\qquad s>0,
\end{equation}
where $a,b\geq0$, $\varphi$ is the completely monotone function and the Laplace transform of $\varphi$ exists for any $s>0$. A function $\phi$ defined on $\mathbb{R}_+$ is said to be a complete Bernstein function $(\mathbb{CBF})$ if and only if
\begin{equation}\label{bern}
	\frac{\widehat{\phi}(s)}{s} \in \mathbb{SF}, \quad s > 0.
\end{equation}
Assuming $\widehat{\varrho}(s) \in \mathbb{SF}$ implies that $\widehat{\varrho}(s)$ can be represented in the specific form given by Eq. \eqref{repres}.  Therefore, the kernel $\varrho(t)$ is completely monotone for $t > 0$ and can have an additional term $b\delta(t)$, where $b \geq 0$, and $\delta(\cdot)$ denotes the Dirac delta function. Moreover, the condition $\hat{\varrho}(s) \in \mathbb{SF}$ is equivalent to the condition $s\hat{\varrho}(s) \in \mathbb{CBF}$, see Eq. \eqref{bern}.

\vspace{0.4cm}

\subsection{Statement of the problem}
In this work, we investigate a two-dimensional generalized diffusion equation with memory kernel of the form
\begin{equation}\label{problem}
	{}_{G}\mathcal{D}_{t}^{(\varrho)}u(x,y,t)=u_{xx}(x,y,t)+u_{yy}(x,y,t)+F(x,y,t), \quad (x,y,t) \in \Omega,
\end{equation}
within the domain
$$\Omega:=(0,1)\times(0,1)\times(0,T],$$
alongside the following boundary conditions
\begin{align}
	u(0,y,t)&=0,\;\; u_{x}(1,y,t)-d\big(u_{xx}(1,y,t)+u_{yy}(1,y,t)\big)=0,\label{boundarycdn}\\
	& \;\;\;y\in[0,1], \;\;  t\in[0,T],\;\;\; {d=const}>0,\nonumber\\
	u(x,0,t)&=0,\qquad u_{y}(x,1,t)=0,\quad \quad x\in[0,1], \quad \quad t\in[0,T],\label{boundarycdn2}
\end{align}
and the initial condition
\begin{equation}\label{initialcdn}
	u(x,y,0)=\psi(x,y), \quad\;\; \quad y\in[0,1],\qquad\;\;\quad x\in[0,1],
\end{equation}
where ${}_{G}\mathcal{D}_{t}^{(\varrho)}$ is given by Eq. (\ref{Kach}).
\subsection{Special cases}
We examine specific cases of Eq. (\ref{problem}) based on the properties of $\varrho(t)$.
\smallskip 

\noindent\textbf{$n$-Tempered FDE:}
Assume the following power-law memory kernel $$\varrho(t)=\sum_{\xi=1}^{n}\frac{P_{\xi}}{\Gamma(1-\beta_{\xi})}e^{-pt}t^{-\beta_{\xi}},$$ where  $0<\beta_n<...<\beta_2<\beta_{1}<1$,\; $D_{\xi}>0,\; \xi=1,2,...n,$ and $p>0$ represents the truncation parameter. Applying this kernel in (\ref{problem}), we derive the corresponding $n$-tempered FDE
\begin{equation}\label{temperedeq}
	\displaystyle\sum_{\xi=1}^{n}P_{\xi} \;\mathcal{D}_{t}^{\beta_\xi, p}\;u(x,y,t)=u_{xx}(x,y,t)+u_{yy}(x,y,t)+F(x,y,t),
\end{equation}
where $\mathcal{D}_{t}^{\beta_\xi, p}$ denotes the tempered fractional derivative of order $\beta_\xi\in(0,1)$ and $p>0$.
\vspace{0.3cm}

Now, we present few special cases of Eq. (\ref{temperedeq}) and Eq. (\ref{problem}):
\smallskip

\noindent \textbf{Case-I:} Setting $p=0$ in Eq. (\ref{temperedeq}), we obtain the following multi-term FDE
\begin{equation}\label{bieq}
	\displaystyle\sum_{\xi=1}^{n}P_{\xi} \;{}_C\mathcal{D}_{t}^{\beta_\xi}u(x,y,t)=u_{xx}(x,y,t)+u_{yy}(x,y,t)+F(x,y,t),
\end{equation} 
where $0<\beta_n<...<\beta_2<\beta_{1}<1$,\; $P_{\xi}>0,\; \xi=1,2,...n$. ${}_C\mathcal{D}_{t}^{\beta_\xi}$ represents the CFD of order $\beta_\xi\in(0,1)$.
\vspace{0.3cm}

\smallskip

\noindent \textbf{Case-II:} Take $P_1=1,\; P_{2}=P_3=...=P_{n}=0,\;$ and $p=0$ in Eq. (\ref{temperedeq}), the following mono-FDE is obtained
\begin{equation}\label{monoeq}
	_C\mathcal{D}_{t}^{\beta_{1}}u(x,y,t)=u_{xx}(x,y,t)+u_{yy}(x,y,t)+F(x,y,t),
\end{equation}
where $0<\beta_{1}<1$.
\vspace{0.2cm}
\smallskip

\noindent \textbf{Case-III:}
By choosing $\beta_{1}=1,\;P_1=1,\; P_{2}=P_3=...=P_{n}=0,\;$ and $p=0$ in Eq. (\ref{temperedeq}), we derive the SDE as follows 
\begin{equation}\label{standardeq}
	u_{t}(x,y,t)=u_{xx}(x,y,t)+u_{yy}(x,y,t)+F(x,y,t).
\end{equation}
\vspace*{0.3cm}
\subsection{Mathematical formulation of the ISPs}
In connection with the system (\ref{problem})-(\ref{initialcdn}), we explore two ISPs, defining each and addressing their solutions.

\smallskip

\noindent \textbf{In the first-ISP}, the source term in Eq. (\ref{problem}) depends only on spatial variables, denoted as $F(x,y,t) = f(x,y)$. The first-ISP is to determine both $f(x,y)$ and $u(x,y,t)$ simultaneously for the system (\ref{problem})-(\ref{initialcdn}), given an additional condition 
\begin{eqnarray}\label{overdeterminationcdn} 
	{u(x,y,T)}=g(x,y),\qquad x\in[0,1]\quad y\in[0,1].
\end{eqnarray} 
The additional condition in the first-ISP denotes the concentration of diffusive particles at time $T$ in a diffusion process and the final temperature distribution at time $T$ in heat conduction.
The solution of first-ISP, expressed as $\{u(x,y,t), f(x,y)\}$, is termed regular if it satisfies the system (\ref{problem})-(\ref{initialcdn}) and the additional condition (\ref{overdeterminationcdn}), with the requirements:  $$f(x,y)\in C([0,1]\times[0,1]); \quad u(x,y,t)\in C(\overline{\Omega}); \;\quad\overline{\Omega}=[0,1]\times[0,1]\times[0,T];  $$
$$ u_{xx}(\cdot,y,t) \in C^2([0,1]);\quad u_{yy}(x,\cdot,t) \in C^2([0,1]);\quad {}_{G}\mathcal{D}_{t}^{(\varrho)}u(x,y,\cdot)\in C([0,T]).$$
\smallskip

\noindent \textbf{In the second-ISP}, we consider the source term in the form $F(x,y,t)=w(t)f(x,y,t)$, {where $f(x,y,t)$ is given}. The objective of the second-ISP is to identify both $w(t)$ and $u(x,y,t)$ for the system defined by Eqs. (\ref{problem})-(\ref{initialcdn}). To ascertain $w(t)$, we impose the following integral condition
\begin{equation}\label{integraloverdeterminationcdn}
	\int_{0}^{1} \int_{0}^{1}u(x,y,t) dxdy =\mathcal{E}(t), \quad \quad t\in [0,T],
\end{equation}
where $\mathcal{E}(t)\in AC([0,T])$ and AC denotes the space of absolutely continuous functions. The integral-type condition serves as essential information for recovering the source term. For instance, in diffusion processes, Eq. \eqref{problem} describes the anomalous diffusion of particles in turbulent plasma, where $\mathcal{E}(t)$ represents the total concentration of particles. Similarly, in the case of heat conduction, $\mathcal{E}(t)$ corresponds to the total heat energy in a heated plate. The regular solution of second-ISP, represented as $\{u(x,y,t), w(t)\}$, meets the conditions: $$u(x,y,t)\in C(\overline{\Omega});\quad w(t)\in C([0,T]); \quad  u_{xx}(\cdot,y,t) \in C^2([0,1]);$$  $$\quad u_{yy}(x,\cdot,t) \in C^2([0,1]);\qquad {}_{G}\mathcal{D}_{t}^{(\varrho)}u(x,y,\cdot)\in C([0,T]).$$ 

\section{Spectral problem}\label{Bi-orth}
This section recalls the spectral analysis of the problem outlined in (\ref{problem})-(\ref{initialcdn}), as explored in \cite{Kapustin}. Let us consider the spectral problem corresponding to (\ref{problem})-(\ref{initialcdn})
\begin{equation}\label{sproblem2D}
	\left . \textstyle\begin{array}{l}
		X_{xx}(x,y)+X_{yy}(x,y)+\lambda X(x,y)=0,\quad  (x,y)\in(0,1)\times(0,1),\\
		X(0,y)=X(1,y),\;\;\qquad\qquad\qquad\qquad\qquad y\in[0,1] \\ X_{x}(1,y)-d\left(X_{xx}(1,y)+X_{yy}(1,y)\right)=0,\;\;\;\;\;\;\;\;y\in[0,1],\;\; d>0, \\
		X(x,0)=0,\qquad X_{y}(x,1)=0,\qquad\qquad \qquad x\in[0,1]. \end{array}\displaystyle \right \}\end{equation}
The eigenfunctions associated with this spectral problem are given by
\begin{equation}\label{x0}
	X_{j,k}(x,y) = 2 \sin\left(\sqrt{\nu^{k}_{j}} x\right) \sin(\mu_k y) ,\;  j = 0, 1, 2, \dots, \;; k = 1, 2,\dots, 
\end{equation}
with corresponding eigenvalues given by 
$$\lambda_{j,k} = \nu^{k}_{j}+\mu_{k}^{2},$$
where $\mu_k = \frac{\pi}{2} (2k - 1)$, and  $\nu^{k}_{j}$ denotes the roots of the characteristic equation
\begin{equation}\label{root}
	\tan\big(\sqrt{s}\big) = \frac{\sqrt{s}}{d(s + \mu_k^2)},
\end{equation}
listed in ascending order. The square root branch is taken on the right half-plane, ensuring $\nu^{k}_{j} > 0$ for $j \geq 1$, with further behavior outlined as $\lim_{j \to \infty}\left(\sqrt{\nu^{k}_{j}} - \pi j\right) = 0$, $\nu^{k}_{j} > 0$ if $d \mu_k^2 < 1$, $\nu^{k}_{0} < 0$ if $d \mu_k^2 > 1$, and $\nu^{k}_{0} = 0$ if $d \mu_k^2 = 1$. In the latter case, the eigenfunction takes the form $$X_{0,k}(x,y) = 2x \sin(\mu_k y).$$ It is verified that the system $Y_{j,k}(x,y), j = 1, 2, \dots; k = 1, 2, \dots$, defined as
$$
Y_{j,k}(x,y) = \frac{2}{	\chi_{j,k}} \left[\sin\left(\sqrt{\nu^{k}_{j}} x\right) - \frac{\sin\left(\sqrt{\nu^{k}_{j}}\right)}{\sin\left(\sqrt{\nu^{k}_{0}}\right)} \sin\left(\sqrt{\nu^{k}_{0}} x\right)\right] \sin(\mu_k y),
$$
is bi-orthogonal in $L^2((0,1)\times(0,1))$ to the system $X_{j,k}(x,y),\; j, k = 1, 2, \dots$. The denominator $	\chi_{j,k}$ for $Y_{j,k}(x,y)$ is given by
\begin{eqnarray}\label{ajk}
	\chi_{j,k} &=& 1 +\frac{ d \sin^2\left(\sqrt{\nu^{k}_{j}}\right)}{\sqrt{\nu^{k}_{j}}} \left(\sqrt{\nu^{k}_{j}} - \frac{\mu^{2}_{k}}{\sqrt{\nu^{k}_{j}}}\right) \nonumber\\&=& d \sin^2\left(\sqrt{\nu^{k}_{j}}\right) + \int_0^1 \sin^2\left(\sqrt{\nu^{k}_{j}} x\right) dx.
\end{eqnarray}
After performing the integral and using Eq. (\ref{root}), we get
\begin{equation}\label{sss}
	\sqrt{\nu^{k}_{j}} = \pi j + \frac{\pi k}{d \gamma_{j,k}} + \sigma^{k}_{j},
\end{equation}
where $\gamma_{j,k} = (\pi j)^2 + \mu_k^2$, and $\left|\sigma^{j}_{k}\right| \leq A \gamma_{j,k}^{-\frac{3}{2}}$, for a constant $A>0$.

The following theorem applies to the basis property of the eigenfunctions.
\begin{theorem} \cite{Kapustin}
	The system of functions $\{X_{j,k}(x,y)\}$, with $j = 1, 2, \dots;$\\$ k = 1, 2, \dots$, form a Riesz basis in $L^2((0,1)\times(0,1))$.
\end{theorem}

\begin{lemma}\label{regularties}
	Consider $\phi(x,y)\in C^{4}([0,1]\times[0,1])$ that meets the following conditions
	\begin{equation}\label{conditions}
		\left . \textstyle\begin{array}{l} \phi(x,0)=0,\quad \phi(0,y)=0,\quad \phi(1,y)=0 \quad \phi_{y}(x,1)=0,\quad \phi_{x}(0,y)=0,\\\phi_{x}(1,y)=0,\quad\phi_{xx}(x,0)=0,\quad {\phi_{xx}(0,y)}=0,\quad {\phi_{xx}(1,y)}=0,\\ {\phi_{xyx}(1,y)}=0,\quad
			\displaystyle\int_{0}^{1}\phi(x,y)\sin\left(\sqrt{\nu^{j}_{0}} x\right)dx=0,\quad y\in [0,1], \quad j=1,2,\cdots \end{array}\displaystyle \right \}.\end{equation}
	then the following inequalities 
	\begin{eqnarray}
		\sum_{j=1}^{\infty}\sum_{k=1}^{\infty}\lambda_{j,k}|\phi_{j,k}|&\leq& C\|\phi^{(iv)}\|_{C^{4}([0,1]\times[0,1])},\label{equalti1}\\
		|\phi_{j,k}|&\leq&\frac{1}{\nu^{j}_{k}\mu^{2}_{j}}\|\phi_{xyxy}\|_{C^{4}([0,1]\times[0,1])} \label{equality2}
	\end{eqnarray}
	hold, where $$\phi_{j,k}:=\frac{2}{{\chi_{j,k}}}\int_{0}^{1}\int_{0}^{1}\phi(x,y)Y_{j,k}(x,y)dy dx.$$
\end{lemma} 
\begin{proof} 
	The inequality (\ref{equalti1}) is proved in \cite{Ismail}. To prove the inequality (\ref{equality2}), consider  
	\begin{equation*}
		\phi_{j,k}=\langle\phi(x,y),Y_{j,k}(x,y)\rangle=\frac{2}{{\chi_{j,k}}}\int_{0}^{1}\int_{0}^{1}\phi(x,y)Y_{j,k}(x,y)dxdy.
	\end{equation*}
	Using integration by parts w.r.t $x$ followed by	$y$, we obtain
	\begin{equation*}
		|\phi_{j,k} |\leq\frac{2}{{\chi_{j,k}}\sqrt{\nu^{k}_{j}}\; \mu_{k}}\int_{0}^{1}\int_{0}^{1}\phi_{xy}(x,y)Y_{j,k}(x,y)dy dx.
	\end{equation*}
	Similarly, by integrating by parts twice w.r.t $x$  and once w.r.t $y$, we find
	\begin{equation*}
		|\phi_{j,k} |\leq\frac{2}{{\chi_{j,k}}\; {\nu^{k}_{j}}\; \mu^{2}_{k}}\int_{0}^{1}\int_{0}^{1}\phi_{xyxy}(x,y)Y_{j,k}(x,y)dy dx.
	\end{equation*}
	Using Cauchy–Schwarz inequality, we deduce
	\begin{equation*}
		|\phi_{j,k}|\leq\frac{1}{\nu^{k}_{j}\;\mu^{2}_{k}}\|\phi_{xyxy}\|_{C^{4}([0,1]\times[0,1])},
	\end{equation*}
	where $\|Y_{j,k}(x,y)\|_{C^{4}([0,1]\times[0,1])}\leq 1$.
\end{proof}
\section{First-ISP}\label{Main IP-I}
The aim of this section is to recover a space-dependent source term, namely $F(x,y,t)=f(x,y)$ and $u(x,y,t)$, from the systems (\ref{problem})-(\ref{initialcdn}) and (\ref{overdeterminationcdn}).
\subsection{Series solution construction for first-ISP}
The series solution for first-ISP in relation to the system defined by Eqs. (\ref{problem})-(\ref{initialcdn}) and (\ref{overdeterminationcdn}) can be written as follows
\begin{eqnarray}
	u(x,y,t) &= &\sum_{j=1}^{\infty}\sum_{k=1}^{\infty} U_{j,k}(t)X_{j,k}(x,y) \label{vsolution},\end{eqnarray}and\begin{eqnarray}
	f(x,y)&=&\sum_{j=1}^{\infty}\sum_{k=1}^{\infty}f_{j,k}X_{j,k}(x,y),\label{fsolution}
\end{eqnarray}
where $U_{j,k}(t)$ and $f_{j,k}$, for $j,k=1,2,\cdots$, are the unknowns, governed by the differential equation system
\begin{equation}\label{LapT1T2T3}
	_{G}\mathcal{D}_{t}^{(\varrho)}U_{j,k}(t)=-\lambda_{j,k} U_{j,k}(t)+f_{j,k}.
\end{equation} 
\smallskip

\noindent Taking into account the Laplace transform of (\ref{LapT1T2T3}) and using the initial condition in (\ref{initialcdn}), we obtain
\begin{equation}\label{T1T2T3}
	U_{j,k}(t)=\psi_{j,k}\mathbf{W}_{j,k}(t)+f_{j,k}\mathbf{Z}_{j,k}(t),
\end{equation}
where $$\psi_{j,k}=\big\langle \psi(x,y), Y_{j,k}(x,y)\big\rangle:=\int_{0}^{1}\int_{0}^{1}\psi(x,y)Y_{j,k}(x,y)dxdy.$$
The functions $\mathbf{W}_{j,k}(t)$ and $\mathbf{Z}_{j,k}(t)$  are defined as follows
\begin{align}\label{AS-val}
	\mathbf{W}_{j,k}(t)=&1-\lambda_{j,k}\int_{0}^{t}\mathbf{V}_{j,k}(\tau)d\tau  \quad\mbox{and}\quad
	\mathbf{Z}_{j,k}(t)=\int_{0}^{t}\mathbf{V}_{j,k}(\tau)d\tau,
\end{align}
where 
\begin{equation}\label{v-values}
	\mathbf{V}_{j,k}(t)=\mathcal{L}^{-1}\left( \frac{1}{s\widehat{\varrho}(s)+\lambda_{j,k}}\right),
\end{equation}
and $\mathbf{V}_{j,k}(t)$ is the impulse-response function. From Lemma 4.1 in \cite{Bazhlekova}, we establish that ${\mathbf{V}_{j,k}}(t)$ is a completely monotone function, meaning that its Laplace transform $\widehat{\mathbf{V}_{j,k}}(s)$ is a Stieltjes function. 
\smallskip

\noindent In the Laplace domain, the representations of $\mathbf{W}_{j,k}(t)$ and $\mathbf{Z}_{j,k}(t)$ are given as follows
\begin{eqnarray}\label{ghi values}
	\left.
	\textstyle\begin{array}{l}
		\widehat{\mathbf{W}_{j,k}}(s)= \frac{1}{s}- \frac{\lambda_{j,k}}{s(s\widehat{\varrho}(s)+\lambda_{j,k})},   \\
		\widehat{\mathbf{Z}_{j,k}}(s)=\frac{1}{s(s\widehat{\varrho}(s)+\lambda_{j,k})}.
	\end{array}\displaystyle\right\}
\end{eqnarray}
The additional condition (\ref{overdeterminationcdn}) provides the following expression for $f_{j,k}$
\begin{equation}
	f_{j,k}=\frac{1}{\mathbf{Z}_{j,k}(T)}\Big(g_{j,k}-\psi_{j,k} \mathbf{W}_{j,k}(T)\Big),\label{f0}
\end{equation}
where $g_{j,k}=\big\langle g(x,y), Y_{j,k}(x,y)\big\rangle.$
\smallskip

\noindent By plugging the value of $f_{j,k}$ in Eq. (\ref{T1T2T3}), we obtain
\begin{equation}
	U_{j,k}(t)=\psi_{j,k} \mathbf{W}_{j,k}(t)+\frac{ \mathbf{Z}_{j,k}(t)}{\mathbf{Z}_{j,k}(T)}\Big(g_{j,k}-\psi_{j,k}\mathbf{W}_{j,k}(T)\Big). \label{T0}
\end{equation}
Now, we examine specific cases of Eqs. (\ref{T1T2T3}) and (\ref{f0}) based on the behavior of $\varrho(t).$
\smallskip

\noindent To analyze the Laplace transform of the power-law memory kernel, given by $\widehat{\varrho}(s)=\sum_{\xi=1}^{n}\frac{P_{\xi}}{(s+p)^{1-\beta_{\xi}}}$ and apply it in Eq. (\ref{ghi values}), we obtain the following estimates
\begin{equation}\label{AS values}
	\left.
	\textstyle\begin{aligned}
		\mathbf{W}_{j,k}(t)=&\frac{1}{P_{1}}-\frac{\lambda_{j,k}}{D_1}\displaystyle\int_{0}^{t}e^{-p\tau}\tau^{\beta_1-1}{E}_{\bar{\beta},{\beta_1}}(\cdots\mathbf{\lambda \tau}\cdots)d\tau\\=&\frac{1}{P_{1}}e^{-pt}{E}_{\bar{\beta},1}(\cdots\mathbf{\lambda t}\cdots),\\
		\mathbf{Z}_{j,k}(t)=&\frac{1}{P_{1}}\displaystyle\int_{0}^{t}e^{-p\tau}\tau^{\beta_1-1}{E}_{\bar{\beta},{\beta_1}}(\cdots\mathbf{\lambda \tau}\cdots)d\tau\\
		=&\frac{1}{P_{1}}e^{-pt}t^{\beta_{1}}{E}_{\bar{\beta},\beta_{1}+1}(\cdots\mathbf{\lambda t}\cdots),
	\end{aligned}\displaystyle\right\}, 
\end{equation}
where $E_{\bar{\beta},\beta_{1}}(\cdot)$ denotes the multinomial (generalized) Mittag--Leffler function, with parameters   $$\bar{\beta}=(\beta_1-\beta_n,...,\beta_1-\beta_2,\beta_1),$$ and argument vector
$$(\cdots\mathbf{\lambda t}\cdots)=\biggl(-\frac{P_n}{P_1}t^{\beta_1-\beta_n},...,-\frac{P_2}{P_1}t^{\beta_1-\beta_2},-\frac{\lambda _{j,k} }{P_1}t^{\beta_1}\biggr).$$
Plugging the values of Eq. (\ref{AS values}) into Eq. (\ref{T1T2T3}), we infer the expression of $U_{j,k}(t)$ for the $n$-tempered FDE (\ref{temperedeq}). Its precise formulation is given by
\begin{align}\label{tempT value}
	U_{j,k}(t)=& \frac{1}{P_{1}}\Bigg(\psi_{j,k}e^{-pt}{E}_{\bar{\beta},1}(\cdots\mathbf{\lambda t}\cdots)+ f_{j,k}e^{-pt}t^{\beta_{1}}{E}_{\bar{\beta},\beta_{1}+1}(\cdots\mathbf{\lambda t}\cdots)\Bigg),	
\end{align}
where
\begin{equation}\label{tempf value}
	f_{j,k}=\frac{P_{1}g_{j,k}-\psi_{j,k} e^{-pT}{E}_{\bar{\beta},1}(\cdots\mathbf{\lambda T}\cdots)}{e^{-pT}T^{\beta_{1}}{E}_{\bar{\beta},\beta_{1}+1}(\cdots\mathbf{\lambda T}\cdots)}.
\end{equation}
\smallskip

\noindent $\bullet$ Substituting $p=0$ in Eqs. (\ref{tempT value}) and (\ref{tempf value}), the following expression of $U_{j,k}(t)$ for the solution of multi-term FDE (\ref{bieq}) is obtained
\begin{align*}
	U_{j,k}(t)=& \frac{1}{P_{1}}\Big(\psi_{j,k} {E}_{\bar{\beta},1}(\cdots\mathbf{\lambda t}\cdots)+ f_{j,k}t^{\beta_{1}}{E}_{\bar{\beta},\beta_{1}+1}(\cdots\mathbf{\lambda t}\cdots)\Big),	
\end{align*}
where
\begin{equation*}
	f_{j,k}=\frac{P_{1}g_{j,k}-\psi_{j,k} {E}_{\bar{\beta},1}(\cdots\mathbf{\lambda T}\cdots)}{T^{\beta_{1}}{E}_{\bar{\beta},\beta_{1}+1}(\cdots\mathbf{\lambda T}\cdots)}.
\end{equation*}  
\smallskip

\noindent $\bullet$ By setting $p=0$,\; $P_{1}=1,$ and $P_{2}=P_3=...=P_{n}=0$ in Eqs. (\ref{tempT value}) and (\ref{tempf value}), we get the following expression for $U_{j,k}(t)$, which provides a solution for a mono-FDE (\ref{monoeq})
\begin{equation*}
	U_{j,k}(t)=\psi_{j,k} {E}_{{\beta_1},1}(-\lambda_{j,k}t^{\beta_{1}})+ f_{j,k}t^{\beta_{1}}{E}_{{\beta_1},\beta_{1}+1}(-\lambda_{j,k}t^{\beta_{1}}),	
\end{equation*}
where
\begin{equation*}
	f_{j,k}=\frac{1}{T^{\beta_{1}}{E}_{{\beta_1},\beta_{1}+1}(-\lambda_{j,k}T^{\beta_{1}})}\left(g_{j,k}-\psi_{j,k} {E}_{{\beta_1},1}(-\lambda_{j,k}T^{\beta_{1}})\right).
\end{equation*}  
\smallskip

\noindent $\bullet$ Taking $\beta_{1}=1,\;p=0$,\; $P_{1}=1,$ and $P_{2}=P_3=...=P_{n}=0$ in Eqs. (\ref{tempT value}) and (\ref{tempf value}), the expression for $U_{j,k}(t)$ in the solution of the SDE (\ref{standardeq}) is given by
$$	U_{j,k}(t)= \psi_{j,k} \Big(1-\lambda_{j,k}\int_{0}^{t} e^{- \lambda_{j,k} \zeta}d\zeta\Big)+ f_{j,k}\int_{0}^{t} e^{- \lambda_{j,k}\zeta}d\zeta,$$
where
$$f_{j,k}=\frac{g_{j,k}-\psi_{j,k} \Big(1-\lambda_{j,k}\displaystyle\int_{0}^{t} e^{- \lambda_{j,k} \zeta}d\zeta\Big)}{\displaystyle\int_{0}^{t} e^{- \lambda_{j,k}\zeta}d\zeta}.$$
\subsection{Ill-posedness of first-ISP}
In this subsection, we see the ill-posedness of first-ISP.	Reviewing some relevant results will be the first step before moving forward.
\begin{lemma}\label{AS}\cite{IlyasMalik2024CAM} Given that $\lambda_{j,k} > 0$, the following holds: \begin{equation}\label{relation} \frac{d}{dt}\mathbf{W}_{j,k}(t) = -\lambda_{j,k}\mathbf{V}_{j,k}(t). \end{equation} Moreover, for $T > 0$, we have the following bound: \begin{equation} \label{estimateS} \mathbf{Z}_{j,k}(t) = \int_{0}^{T}\mathbf{V}_{j,k}(\zeta) d\zeta \leq \frac{1}{\lambda_{j,k}}. \end{equation} \end{lemma}
\vspace{0.4cm}
\smallskip

\noindent
Consider the following example as a way to demonstrate the ill-posedness of first-ISP. The final and initial data are set as follows $$\tilde{g}(x,y)=\frac{1}{\sqrt{\lambda_{j,k}}}Y_{j,k}(x,y),\qquad \mbox{and} \qquad\tilde{\psi}(x,y)=0,$$ where $ j,k \in \mathbb{N}$.  With these choices, the source term becomes
\begin{equation}
	\tilde{f}(x,y)=\frac{1}{\sqrt{\lambda_{j,k}} \mathbf{Z}_{j,k}(T)}X_{j,k}(x,y) . 
\end{equation}
Now, by taking an alternative final data and the initial data
$${g}(x,y)=0,\qquad \mbox{and} \qquad{\psi}(x,y)=0,$$
then, the source term becomes $f(x,y)=0$. We then calculate the $L^{2}$-norm between the two final data $\tilde{g}$ and $g$
\begin{eqnarray*}
	\|\tilde{g}-g\|_{L^{2}((0,1)\times(0,1))}&=&\bigg\|\frac{1}{\sqrt{\lambda_{j,  k}}}X_{j,k}(x,y)\bigg\|_{L^{2}((0,1)\times(0,1))}.
\end{eqnarray*}
Due to (\ref{sss}), the equality 
\begin{equation*}
	\|\tilde{g}-g\|_{L^{2}((0,1)\times(0,1))}=\frac{1}{ \sqrt{\lambda_{j,k}}}
\end{equation*}
is obtained, leading to the limit relation
\begin{equation}\label{errorfinal}
	\lim_{j, k \to +\infty}\|\tilde{g}-g\|_{L^{2}((0,1)\times(0,1))}=\lim_{j, k \to +\infty}\frac{1}{\sqrt{\lambda_{j,k}}}=0.
\end{equation}
To compute the error in the $L^{2}$-norm between the two source terms $\tilde{f}$ and $f$, we have 
\begin{eqnarray*}
	\|\tilde{f}-f\|_{L^{2}((0,1)\times(0,1))}&=&\bigg\|\frac{X_{j,k}(x,y)}{\sqrt{\lambda_{j, k}} \mathbf{Z}_{j,k}(T)}\bigg\|_{L^{2}((0,1)\times(0,1))}
	=\frac{1}{\sqrt{\lambda_{j, k}}\mathbf{Z}_{j,k}(T)}.
\end{eqnarray*}
Using the estimate from Lemma \ref{AS}, we arrive at
\begin{equation*}
	\|\tilde{f}-f\|_{L^{2}((0,1)\times(0,1))} \geq  {\sqrt{\lambda_{j, k}} }.
\end{equation*}
This leads to
\begin{equation}
	\lim_{j, k \to +\infty}\|\tilde{f}-f\|_{L^{2}((0,1)\times(0,1))} > \lim_{j, k \to +\infty}  { \sqrt{\lambda_{j, k}} }=+\infty \label{errorf}.
\end{equation}
From Eqs. (\ref{errorfinal}) and (\ref{errorf}), we conclude that first-ISP is ill-posed.

\subsection{The main results of first-ISP}
In this subsection, we outline the existence and the uniqueness of the solution of the first-ISP. To demonstrate the regularity of solutions for both first-ISP and second-ISP, we present the following lemma:

\begin{lemma}\cite{Bazhlekova}\label{*} Consider $\mathbf{Z}_{j,k}(t)$ as defined by Eq. (\ref{AS-val}), which is continuous on $[0, \infty)$, vanish at $t=0$, positive, and nondecreasing on $\mathbb{R_{+}}$. For any $t, T>0$, the following inequalities 
	\begin{itemize}
		\item [(a)] $\mathbf{Z}_{j,k}(t)\leq\frac{1}{\lambda_{j,k}}$,
		\item  [(b)] $\frac{1}{\mathbf{Z}_{j,k}(t)}\leq C_{0}\lambda_{j,k}$,
	\end{itemize}
	hold. where $C_{0}=C_{0}(T)>0$. \end{lemma}
\noindent $\bullet$ From the expression for $ \mathbf{W}_{j,k}(t) $ in Eq. (\ref{AS-val}), we can written as
\begin{equation*}
	\mathbf{W}_{j,k}(t) = 1 - \lambda_{j,k} \int_{0}^{t} \mathbf{V}_{j,k}(\zeta) \, d\zeta.
\end{equation*}
The form of $ \mathbf{V}_{j,k}(t)$ here aligns with the expression provided in \cite{Bazhlekova} (Lemma 4.1, Eq. (42)). Given this, the estimate for $ \mathbf{V}_{j,k}(t)$ can be written as
\begin{equation}\label{es2}
	0 < c_{0} \leq \lambda_{j,k} \int_{0}^{t} \mathbf{V}_{j,k}(\zeta) \, d\zeta < 1,
\end{equation}
based on the inequalities in (\ref{es2}), we deduce that
\begin{equation}\label{hi}
	\mathbf{W}_{j,k}(t) \leq C_{1},
\end{equation}
where $ c_{0} $ and $C_{1} $ are positive constants.
\begin{theorem}\label{existenceisp1Theorem} Suppose that \(\psi(y, x)\) and \(g(y, x)\) are given functions in \(C^4([0,1] \times [0,1])\) and meet the conditions given in \eqref{conditions}. Then, a regular solution exists for first-ISP associated with the system \eqref{problem}–\eqref{initialcdn} and \eqref{overdeterminationcdn}. \end{theorem}
\begin{proof} This theorem is proved in two parts.
	\smallskip
	
	\noindent\textbf{Part-I:}	To establish the existence of a solution for first-ISP, we need to demonstrate the continuity of the associated series $$f(x,y),\quad u(x,y,t), \quad u_{xx}(x,y,t),\quad u_{yy}(x,y,t),\quad \mbox{and}\quad {}_{G}\mathcal{D}_{t}^{(\varrho)}u(x,y,t).$$
	\smallskip
	
	\noindent First, we focus on verifying the continuity of $f(x,y)$. Using Lemma \ref{*}, along with estimate (\ref{hi}) and Eq. (\ref{f0}), we obtain
	\begin{equation*}
		|f_{j,k}|\leq C_{0} \lambda_{j,k}(|g_{j,k}|-C_{1}|\psi_{j,k}|).
	\end{equation*}
	By invoking Lemma \ref{regularties}, we deduce
	\begin{equation}\label{convf1n}
		|f_{j,k}|\leq C_{0}C\Big(\|g^{(iv)}\|_{C^4([0,1]\times[0,1])}-C_{1}\|\psi^{(iv)}\|_{C^4([0,1]\times[0,1])}\Big)
	\end{equation}
	From inequality (\ref{convf1n}), we deduce that the series of $f(x,y)$ as given by Eq. (\ref{fsolution}) is bounded. Thus, the continuity of $f(x,y)$ is established by the WMT.
	
	\noindent  Now, our goal is to establish the continuity of $u(x,y,t)$. By applying Lemma \ref{*}, Lemma \ref{regularties}, and Eq. (\ref{T0}), we find that
	\begin{eqnarray}
		|U_{j,k}(t)|&\leq& \frac{C}{\lambda_{j,k} }\Big(K\|\psi^{(iv)}\|_{C^4([0,1]\times[0,1])}+C_{0}\|g^{(iv)}\|_{C^4([0,1]\times[0,1])}\Big),\label{Bispaceconvert1n}
	\end{eqnarray}
	where $K=C_{1}(1-C_{0})$ is a positive constant. Based on Eq. (\ref{Bispaceconvert1n}), we deduce that the series of $u(x,y,t)$ in Eq. (\ref{vsolution}) is bounded above. Hence, the continuity of $u(x,y,t)$ follows by applying the WMT.
	\smallskip
	
	\noindent Our next task is to demonstrate the continuity of $u_{xx}(x,y,t)$. To achieve this, we differentiate each term of the series solution provided in Eq. (\ref{vsolution}) w.r.t $x$
	\begin{eqnarray}\label{vxxsolutionx}
		u_{xx}(x,y,t) &=& \sum_{j=1}^{\infty} \sum_{k=1}^{\infty}
		U_{j,k}(t)X_{xx, j, k}(x,y),
	\end{eqnarray}
	where
	$$
	X_{xx,j, k}(x,y)=-\nu^{k}_{j}X_{j, k}(x,y).
	$$	Indeed, we can have $$|X_{j, k}(x,y)|<1,$$ which implies
	\begin{eqnarray} \label{estmatimates xx2}
		\Big|X_{xx, j, k}(x,y)\Big|\leq \nu^{k}_{j}.
	\end{eqnarray}
	Based on Lemma \ref{regularties}, Eq. (\ref{Bispaceconvert1n}), and inequality (\ref{estmatimates xx2}), we derive the following expression
	\begin{equation}\label{partxxx}
		u_{xx}(x,y,t) = \sum_{j=1}^{\infty} \sum_{k=1}^{\infty}
		\frac{C}{\mu^{2}_{k}\lambda_{j,k} }\Big(K\|\psi^{(iv)}\|_{C^4([0,1]\times[0,1])}+C_{0}\|g^{(iv)}\|_{C^4([0,1]\times[0,1])}\Big).
	\end{equation}
	Eq. (\ref{partxxx}) demonstrates that the series form of $u_{xx}(x,y,t)$ is bounded above. Consequently, by WMT, $u_{xx}(x,y,t)$ is continuous.
	\smallskip
	
	\noindent Next, we establish the continuity of $u_{yy}(x,y,t)$. By applying term-by-term differentiation w.r.t $y$ on the series solution from Eq. (\ref{vsolution}), we obtain
	\begin{equation}\label{vyysolutionx}
		u_{yy}(x,y,t) = -\sum_{j=1}^{\infty} \sum_{k=1}^{\infty} \mu^{2}_{k}U_{j, k}(t)X_{j,k}(x,y).
	\end{equation}
	Using Lemma \ref{regularties} and  Eq. (\ref{Bispaceconvert1n}), we obtain
	\begin{equation}\label{partyyy}
		u_{yy}(x,y,t) = \sum_{j=1}^{\infty} \sum_{k=1}^{\infty}
		\frac{C}{\nu^{k}_{j}\lambda_{j,k} }\Big(K\|\psi^{(iv)}\|_{C^4([0,1]\times[0,1])}+C_{0}\|g^{(iv)}\|_{C^4([0,1]\times[0,1])}\Big).
	\end{equation}
	Thus, $u_{yy}(x,y,t)$, as given in Eq. (\ref{partyyy}), is also bounded above, ensuring its continuity by WMT.
	\smallskip
	
	\noindent In the similar way, we establish the continuity of ${}_{G}\mathcal{D}_{ t}^{(\varrho)}u(x,y,t)$. For the continuity of ${}_{G}\mathcal{D}_{t}^{(\varrho)}u(x,y,t)$, using Eq. (\ref{vsolution}), we have
	\begin{equation*}
		{}_{G}\mathcal{D}_{t}^{(\varrho)}u(x,y,t) = \sum_{j=1}^{\infty}\sum_{k=1}^{\infty}{}_{G}\mathcal{D}_{0_+, t}^{(\varrho)} U_{j, k}(t)X_{j,k}(x,y).
	\end{equation*}
	By virtue of Eq. (\ref{LapT1T2T3}), we obtain
	\begin{align*}
		\sum_{j=1}^{\infty}\sum_{k=1}^{\infty}{}_{G}\mathcal{D}_{t}^{(\varrho)}U_{j,k}(t)X_{j,k}(x,y)=&\sum_{j=1}^{\infty}\sum_{k=1}^{\infty}\Big(-\lambda_{j,k} U_{j,k}(t)+f_{j,k}\Big)X_{j,k}(x,y)\\=& -\sum_{j=1}^{\infty}\sum_{k=1}^{\infty}\lambda_{j,k} U_{j,k}(t)X_{j,k}(x,y)+f(x,y).
	\end{align*}	
	Using Lemma \ref{*}, Lemma \ref{regularties} and Eq. (\ref{Bispaceconvert1n}), we can prove the continuity of the series, i.e., $$\sum_{j=1}^{\infty}\sum_{k=1}^{\infty}\lambda_{j,k} U_{j,k}(t)X_{j,k}(x,y).$$ Furthermore, we have discussed already the continuity of $f(x,y).$ Hence, we can say that ${}_{G}\mathcal{D}_{t}^{(\varrho)}u(x,y,t)$ is a continuous function due to WMT.
	\vspace{0.4cm}
	\smallskip
	
	\noindent\textbf{Part-II:} 
	Assume $$\{{}_{(1)}u(x,y,t), {}_{(1)}f(x,y)\}\qquad\qquad \mbox{and} \qquad\qquad\{{}_{(2)}u(x,y,t), {}_{(2)}f(x,y)\}$$ are two regular solution sets of first-ISP. Suppose that for some $x_{0},y_{0}\in(0,1)\times(0,1)$, we have ${}_{(1)}u(x_{0},y_{0},t)={}_{(2)}u(x_{0},y_{0},t)$. Then, it follows that $${}_{(1)}u(x,y,t)={}_{(2)}u(x,y,t),\;\quad\qquad (x,y,t)\in\Omega, $$and similarly$$ {}_{(1)}f(x,y)={}_{(2)}f(x,y),\;\quad\qquad\forall\;x,y\in(0,1)\times(0,1).$$ 
	To prove the uniqueness for the solution of first-ISP, let's consider the following expressions for ${}_{(1)}U_{j,k}(t)$ and ${}_{(2)}U_{j,k}(t)$
	\begin{eqnarray}
		{}_{(1)}U_{j,k}(t)&=&\int_{0}^{1}\int_{0}^{1}
		{}_{(1)}u(x,y,t)X_{j,k}(x,y)dxdy,\label{integreleq}\\ {}_{(2)}{U}_{j,k}(t)&=&\int_{0}^{1}\int_{0}^{1}{}_{(2)}u(x,y,t)X_{j,k}(x,y)dxdy.\label{Vtildetintegral}
	\end{eqnarray}
	Applying $_{G}\mathcal{D}_{t}^{(\varrho)}$ to both sides of Eq. (\ref{Vtildetintegral}) and using Eq. (\ref{problem}), along with the initial condition (\ref{initialcdn}), yields
	\begin{equation*}
		{}_{(2)}U_{j,k}(t)=\psi_{j,k} \mathbf{W}_{j,k}(t)+{}_{(2)}f_{j,k}\mathbf{Z}_{j,k}(t).
	\end{equation*}
	Similarly, an expression for ${}_{(2)}{U}_{j,k}(t)$ is obtained.
	\begin{equation*}
		{}_{(1)}U_{j,k}(t)=\psi_{j,k} \mathbf{W}_{j,k}(t)+{}_{(1)}f_{j,k}\mathbf{Z}_{j,k}(t).
	\end{equation*}
	Since ${}_{(1)}u(x_{0},y_{0},t)={}_{(2)}u(x_{0},y_{0},t)$ for some $x_{0},y_{0}\in (0,1)\times(0,1)$, it follows that 
	$${}_{(1)}U_{j,k}(t)={}_{(2)}{U}_{j,k}(t),$$ and therefore
	\begin{equation*}
		{}_{(1)}f_{j,k}\mathbf{Z}_{j,k}(t)={}_{(2)}f_{j,k}\mathbf{Z}_{j,k}(t), \quad \Rightarrow\quad ({}_{(1)}f_{j,k}-{}_{(2)}f_{j,k})\mathbf{Z}_{j,k}(t)=0.
	\end{equation*}
	Applying the Laplace transform, we get
	\begin{equation}\label{A1}
		\frac{{}_{(1)}f_{j,k}-{}_{(2)}f_{j,k}}{s(s\widehat{\varrho}(s)+\lambda_{j,k})}=0, \quad Re\; s>0,\quad \Rightarrow\quad \frac{{}_{(1)}f_{j,k}-{}_{(2)}f_{j,k}}{\Delta+\lambda_{j,k}}=0,
	\end{equation} 
	where $\Delta=s\widehat{\varrho}(s)$.
	\smallskip
	
	\noindent By selecting the disk for $j=1,\; k=1  $ such that $\lambda_{1,1}$
	lies within it, and applying the Cauchy Integral Theorem to Eq. (\ref{A1}), we obtain
	\begin{equation*}
		{}_{(1)}f_{1,1}={}_{(2)}f_{1,1}.
	\end{equation*}
	Similarly, by considering different disks for $j,k=2,3,\cdots,$  we find that
	\begin{equation*}
		{}_{(1)}{f}_{j,k}={}_{(2)}{f}_{j,k}.
	\end{equation*}
	Therefore, we conclude that $${}_{(1)}f(x,y)={}_{(2)}{f}(x,y).$$
\end{proof}
\section{Second-ISP}\label{INVERSEPROBLEM-II}
This section examines the second-ISP for the system described by Eqs. (\ref{problem})-(\ref{initialcdn}) and (\ref{integraloverdeterminationcdn}), incorporating the source term $F(x,y,t) = w(t)f(x,y,t)$. We will discuss how to recover $w(t)$ along with $u(x,y,t)$, given the additional condition specified by (\ref{integraloverdeterminationcdn}).
\subsection{Series solution construction for second-ISP}
The solution to second-ISP for the system given by Eqs. (\ref{problem})-(\ref{initialcdn}) and (\ref{integraloverdeterminationcdn}) can be expressed as
\begin{equation}
	u(x,y,t) = \sum_{j=1}^{\infty}\sum_{k=1}^{\infty} U_{j,k}(t)X_{j,k}(x,y)\label{usolutionISP2},
\end{equation}
Here, $U_{j,k}(t)$ for $j,k\in\mathbf{N}$ are unknown that satisfy the following differential equation
\begin{equation}\label{LapT}
	_{G}\mathcal{D}_{t}^{(\varrho)}U_{j,k}(t)=-\lambda_{j,k} U_{j,k}(t)+w(t)f_{j,k}(t), 
\end{equation}
Taking the Laplace transform of Eq. (\ref{LapT}) and Applying the initial condition provided in (\ref{initialcdn}), we derive
\begin{equation}\label{TISP2}
	U_{j,k}(t)=\psi_{j,k} \mathbf{W}_{j,k}(t)+w(t)f_{j,k}(t)\ast \mathbf{V}_{j,k}(t).
\end{equation}
Hence, solution of the second-ISP can be written as
\begin{equation}
	u(x,y,t) = \sum_{j=1}^{\infty}\sum_{k=1}^{\infty} \left(\psi_{j,k} \mathbf{W}_{j,k}(t)+w(t)f_{j,k}(t)\ast \mathbf{V}_{j,k}(t)\right)X_{j,k}(x,y)\label{solution-second},
\end{equation}
where $w(t)$ remains to be determined. Next, we examine some specific cases of Eq. (\ref{TISP2}) based on the memory kernel $\varrho(t)$. 
\smallskip

\noindent
\textbf{Solution of second-ISP for $n$-Tempered FDE:}
From Eq. (\ref{AS values}), we derive the following expression for $U_{j,k}(t)$ in the second-ISP for the $n$-tempered FDE (\ref{temperedeq})
\begin{align}\label{tem T}
	U_{j,k}(t)=&\frac{1}{P_{1}}\psi_{j,k} e^{-pt}{E}_{\bar{\beta},1}(\cdots\mathbf{\lambda t}\cdots)+ \frac{w(t)}{P_{1}}f_{j,k}(t)\ast e^{-pt}t^{\beta_{1}-1}{E}_{\bar{\beta},\beta_{1}}(\cdots\mathbf{\lambda t}\cdots).
\end{align}
\begin{itemize}
	\item For \( p=0 \), the equation above matches the expression of \( U_{j,k}(t) \) in second-ISP for solving the multi-term FDE (\ref{bieq}).
	\item When \( p=0 \), \( P_{1}=1 \), and \( P_{2}=P_3=\dots=P_n=0 \), Eq. (\ref{tem T}) simplifies to the expression of \( U_{j,k}(t) \) in second-ISP for the solution of a mono-FDE.
	\item The expression of \( U_{j,k}(t) \) in second-ISP for the SDE (\ref{standardeq}) can be derived by setting \( \beta_{1}=1 \), \( p=0 \), \( P_{1}=1 \), and \( P_{2}=P_3=\dots=P_n=0 \) in Eq. (\ref{tem T}).
\end{itemize}
\subsection{The main results of second-ISP}
The following section presents a theorem that establishes the conditions necessary for problem second-ISP to have a regular solution.
\smallskip

\noindent
\textbf{$\bullet$ Existence of second-ISP}	
\begin{theorem}\label{main rslt}
	Assume the following:
	\begin{itemize}
		\item[(a)] $\psi(x,y) \in C^4([0,1] \times [0,1])$ satisfies conditions (\ref{conditions}).\\
		\item[(b)]$f(x,y,t) \in C^{0,2,2}(\bar{\Omega})$ meets the conditions $$f(x,0,t)=0,\;\frac{\partial}{\partial y}f(x,1,t)=0,\;\frac{\partial}{\partial x}f(0,y,t)=0,$$$$\frac{\partial}{\partial x}f(1,y,t)-d\left(\frac{\partial^{2}}{\partial x^{2}}f(1,y,t)+\frac{\partial^{2}}{\partial y^{2}}f(1,y,t)\right)=0.$$
		Additionally, if
		$$
		\bigg(\displaystyle \int_{0}^{1}\int_{0}^{1}f(x,y,t)dxdy\bigg)^{-1}\leq Q_{0},
		$$
		for some constant $Q_{0}>0$. \\
		\item[(c)] $\mathcal{E}(t)\in AC([0,T])$, and satisfies the consistency condition, i.e.,		$$
		\int_{0}^{1}\int_{0}^{1}{\psi(x,y)}dxdy={\mathcal{E}(0)},$$\\
		then a regular solution for second-ISP exists.
	\end{itemize}
\end{theorem}
\begin{proof}
	To construct the $w(t)$, we employ an additional condition given by (\ref{integraloverdeterminationcdn}), which yields
	\begin{equation*}
		\int_{0}^{1}\int_{0}^{1} {}_{G}\mathcal{D}_{t}^{(\varrho)}u(x,y,t) dxdy={}_{G}\mathcal{D}_{t}^{(\varrho)} \mathcal{E}(t).
	\end{equation*}
	According to Eq. (\ref{problem}), this results in
	\begin{equation*}
		\int_{0}^{1}\int_{0}^{1}\Big(u_{xx}(x,y,t)+u_{yy}(x,y,t)+w(t)f(x,y,t)\Big)dxdy={}_{G}\mathcal{D}_{t}^{(\varrho)} \mathcal{E}(t),
	\end{equation*}
	From this, the following expression for $w(t)$ is obtained
	\begin{equation}\label{time}
		w(t)=\frac{{}_{G}\mathcal{D}_{t}^{(\varrho)} \mathcal{E}(t)-\mathcal{T}(t)-\displaystyle\int_{0}^{t}\mathcal{G}(t,\tau)w(\tau)d\tau}{\displaystyle\int_{0}^{1}\int_{0}^{1}f(x,y,t)dxdy}.
	\end{equation}
	Setting
	\begin{eqnarray}
		\mathcal{T}(t)&=&-2\sum_{j=1}^{\infty}\sum_{k=1}^{\infty}\psi_{j,k} \mathbf{W}_{j,k}(t)\left(1-\cos\left(\sqrt{\nu^{k}_{j}}\right)\right)\left(\frac{\sqrt{\nu^{k}_{j}}}{\mu_{k}}+\frac{2\mu_{k}}{\sqrt{\nu^{k}_{j}}}\right)\label{timefeka},\end{eqnarray}	\begin{equation}
		\int_{0}^{1}\int_{0}^{1}f(x,y,t)dxdy= -2\sum_{j=1}^{\infty}\sum_{k=1}^{\infty}
		f_{j,k}(t)\left(1-\cos\left(\sqrt{\nu^{k}_{j}}\right)\right)\left(\frac{\sqrt{\nu^{k}_{j}}}{\mu_{k}}+\frac{2\mu_{k}}{\sqrt{\nu^{k}_{j}}}\right)\label{timefxt},\end{equation}	\begin{equation}
		\mathcal{G}(t,\tau)=-2\sum_{j=1}^{\infty}\sum_{k=1}^{\infty}f_{j,k}(\tau) \mathbf{V}_{j,k}(t-\tau)\left(1-\cos\left(\sqrt{\nu^{k}_{j}}\right)\right)\left(\frac{\sqrt{\nu^{k}_{j}}}{\mu_{k}}+\frac{2\mu_{k}}{\sqrt{\nu^{k}_{j}}}\right).
	\end{equation}
	Based on Lemma 4.1 (refer to \cite{Bazhlekova}, Eq. (42) on page 10) and Lemma \ref{*}, there are positive constants $Q_{1}>0$ and $Q_{2}>0$ such that 
	\begin{equation}\label{tg-extimates}
		\|\mathcal{T}(t)\|_{C([0,T])}\leq Q_{1},\qquad \mbox{and} \qquad
		\|\mathcal{G}(t, \tau)\|_{C([0,T]\times[0,T])}\leq Q_{2}.
	\end{equation}
	In order to define the mapping $\mathcal{I}:C([0,T])\rightarrow C([0,T])$, let's write it as follows
	\begin{equation}\label{B(a(t))}
		\mathcal{I}(w(t)):=\frac{{}_{G}\mathcal{D}_{t}^{(\varrho)} \mathcal{E}(t)-\mathcal{T}(t)-\displaystyle\int_{0}^{t}\mathcal{G}(t,\tau)w(\tau)d\tau}{\displaystyle\int_{0}^{1}\int_{0}^{1}f(x,y,t)dxdy}.
	\end{equation}
	As a result of the uniform convergence of $$\|\mathcal{T}(t)\|_{C([0,T])}\qquad \mbox{and} \qquad\|\mathcal{G}(t, \zeta)\|_{C([0,T]\times[0,T])},$$ the mapping $\mathcal{I}(w(t))$ is well-defined.
	\smallskip
	
	\noindent Now, we demonstrate that the mapping $\mathcal{I}(w(t))$ is a contraction. Consider
	\begin{equation*}
		|\mathcal{I}(w_{1}(t))-\mathcal{I}(w_{2}(t))| =\left(\displaystyle\int_{0}^{1}\int_{0}^{1}f(x,y,t)dxdy\right)^{-1} {\displaystyle\int_{0}^{t}|\mathcal{G}(t, \tau)||w_{1}(\tau)-w_{2}(\tau)|d\tau}.
	\end{equation*}
	Using the estimate of $\mathcal{G}(t, \zeta)$ given in \eqref{tg-extimates} and assumptions of Theorem \ref{main rslt}, we find that
	\begin{eqnarray*}
		\max_{0\leq t\leq T}|\mathcal{I}(w_{1}(t))-\mathcal{I}(w_{2}(t))| &\leq& TQ_{0}Q_{2}\max_{0\leq t\leq T}|w_{1}(\tau)-w_{2}(\tau)|,
	\end{eqnarray*}
	which implies
	\begin{eqnarray*}
		\|\mathcal{I}(w_{1})-\mathcal{I}(w_{2})\|_{C([0,T])} &\leq& TQ_{0}Q_{2}\|w_{1}-w_{2}\|_{C([0,T])}.
	\end{eqnarray*}
	This demonstrates that the mapping $\mathcal{I}(w(t))$ is a contraction for  $\displaystyle T<\frac{1}{Q_{0}Q_{2}}$, by applying the Banach fixed-point theorem, which guarantees the unique existence of $w(t)$.
	Since $w\in C([0,T])$, we have
	$$\|w\|_{C([0,T])}\leq Q, $$
	where $Q$ is a constant. 
	\vspace{0.3cm}
	
	\noindent 
	To show the existence of a solution for the second-ISP, we will verify the continuity of $$u(x,y,t),\quad u_{xx}(x,y,t),\quad u_{yy}(x,y,t),\quad\mbox{and}\quad{}_{G}\mathcal{D}_{t}^{(\varrho)}u(x,y,t).$$ First, note that the norm estimate for $f(x,y,t)$ in the function space $C^{0,2,2}(\bar{\Omega}))$ is given as follows
	\begin{eqnarray*}
		\|f\|^2_{C^{0,2,2}(\bar{\Omega}))}&=&\int_{0}^{T}\|f\|^2_{C^{2}([0,1]\times[0,1])}dt\\&=&\int_{0}^{T}\sum_{j=1}^{\infty}\sum_{k=1}^{\infty}|f_{j,k}(t)|^{2}dt\\&=&\sum_{j=1}^{\infty}\sum_{k=1}^{\infty}\|f_{j,k}(t)\|^{2}_{C^2([0,T])}.
	\end{eqnarray*}
	Due to Young's inequality for integral convolution and using Eq. (42) from \cite{Bazhlekova}, we obtain
	\begin{eqnarray}\label{festimate}
		\Big\|\int_{0}^{t}\mathbf{V}_{j,k}(t-\tau)f_{j,k}(\tau)d\tau\Big\|^{2}_{C^2([0,T])}&\leq& \Big(\int_{0}^{T}\mathbf{V}_{j,k}(t)dt\Big)^{2}\int_{0}^{T}|f_{j,k}(t)|^{2}dt\nonumber\\&\leq&\frac{1}{\lambda^2_{j,k}}\int_{0}^{T}|f_{j,k}(t)|^{2}dt.
	\end{eqnarray}
	By virtue of Lemmas \ref{regularties}, \ref{*}, and inequality (\ref{festimate}), we have
	\begin{equation}
		|U_{j,k}(t)|\leq \frac{C_{0}}{|\lambda_{j,k}|^{2 }}\left(\|\psi^{(iv)}\|_{C^4([0,1]\times [0,1])}+Q\|f_{j,k}(t)\|_{C^2([0,T])}\right)\label{Biconvert1n}.
	\end{equation}
	Based on inequality (\ref{Biconvert1n}), the uniform convergence of $u(x,y,t)$,  defined by Eq. (\ref{vsolution}), is bounded. Therefore, applying the WMT established the continuity of $u(x,y,t)$.
	\smallskip
	
	\noindent
	We now proceed to demonstrate the continuity of the second partial derivatives, $u_{xx}(x,y,t)$ and $u_{yy}(x,y,t)$. Using the relationships established in Eqs. (\ref{Biconvert1n}) and (\ref{estmatimates xx2}), it follows that $u_{xx}(x,y,t)$ and $u_{yy}(x,y,t)$ are bounded above. Therefore, by WMT, $u_{xx}(x,y,t)$ and $u_{yy}(x,y,t)$ are indeed continuous functions.
	\smallskip
	
	\noindent In a similar manner, we can show the continuity of ${}_{G}\mathcal{D}_{t}^{(\varrho)}u(x,y,t)$.
\end{proof}
\vspace{0.4cm}
\smallskip

\noindent
\textbf{$\bullet$ Uniqueness of second-ISP}	
\begin{theorem}
	A regular solution to second-ISP is uniquely determined under the assumptions of Theorem \ref{main rslt}.
\end{theorem}
\begin{proof}
	The uniqueness of $w(t)$ has been previously established by means of the Banach fixed-point theorem. Now, we prove the
	uniqueness of $u(x,y,t)$. Define 
	\begin{equation*}
		\overline{u}(x,y,t)=u(x,y,t)-v(x,y,t),\qquad\quad \overline{w}(t)=w_{1}(t)-w_{2}(t).
	\end{equation*}
	Then, the function $\overline{u}(x,y,t)$ satisfies the following equation
	\begin{equation*}
		{}_{G}\mathcal{D}_{t}^{(\varrho)}\overline{u}(x,y,t)=\overline{u}_{xx}(x,y,t)+\overline{u}_{yy}(x,y,t), \qquad\; (x,y,t) \in \Omega,
	\end{equation*}
	alongside the following boundary conditions
	\begin{align*}
		\overline{u}(0,y,t)&=0,\;\; \overline{u}_{x}(1,y,t)-d\big(\overline{u}_{xx}(1,y,t)+\overline{u}_{yy}(1,y,t)\big)=0,\\ &\qquad\;\;\;y\in[0,1], \;\;  t\in[0,T],\;\;\; {d=const}>0,\\
		&\overline{u}(x,0,t)=0,\qquad \overline{u}_{y}(x,1,t)=0,\quad \quad x\in[0,1], \quad \quad t\in[0,T],
	\end{align*}
	and the initial condition 
	\begin{equation}\label{initialcdn-2ndISP}
		\overline{u}(x,y,0)=0, \quad\;\; \quad y\in[0,1],\qquad\;\;\quad x\in[0,1].
	\end{equation}
	Consider the further function
	\begin{eqnarray}
		\overline{U}_{j,k}(t)&=&\int_{0}^{1}\int_{0}^{1}
		{}_{(1)}\overline{u}(x,y,t)X_{j,k}(x,y)dxdy.  \label{integrel-2ndISP}
	\end{eqnarray}
	Applying $_{G}\mathcal{D}_{t}^{(\varrho)}$ to both sides of Eq. (\ref{Vtildetintegral}), we have
	\begin{eqnarray*}
		_{G}\mathcal{D}_{t}^{(\varrho)}\overline{U}_{j,k}(t)&=&_{G}\mathcal{D}_{t}^{(\varrho)}\int_{0}^{1}\int_{0}^{1}
		{}_{(1)}\overline{u}(x,y,t)X_{j,k}(x,y)dxdy.   
	\end{eqnarray*}
	Due to Eq. (\ref{problem}), we obtain
	\begin{equation*}
		\overline{U}_{j,k}(t)= \overline{U}_{j,k}(0)\mathbf{W}_{j,k}(t)+\overline{w}(t)f_{j,k}(t)\ast \mathbf{V}_{j,k}(t).
	\end{equation*}
	By applying the initial conditions \eqref{initialcdn-2ndISP}, taking into account the uniqueness of $\overline{w}(t)$ and the completeness of $Y_{j,k}(x,y)$, we infer
	\begin{equation*}
		\overline{U}_{j,k}(t)= 0.
	\end{equation*}
	As a consequence, we deduce 
	$$u(x,y,t)=v(x,y,t).$$
\end{proof}

\subsection{Ill-posedness of second-ISP}
In this subsection, we discuss the ill-posedness of second-ISP. 
\vspace{0.6cm}

The following example illustrates the ill-posed nature of second-ISP. Let us consider 
\begin{eqnarray*}
	\hat{f}(x,y,t)&=&{\lambda_{j,k}}\varrho(t)\Big(C_{1}(t)+\lambda_{j,k}  C_{2}(t)\Big)X_{j,k}(x,y),
\end{eqnarray*}
where $$C_{1}(t)=\sum_{r=1}^{n}\frac{\Gamma(1-\beta_{r})\Gamma\Big(R-\displaystyle\sum_{\xi=2}^{n}\beta_{\xi} \Big)t^{\beta_1}}{\Gamma\Big(R-\beta_{r}-\displaystyle\sum_{\xi=2}^{n}\beta_{\xi}\Big)},\qquad R\geq2,$$and$$C_{2}(t)=\sum_{r=1}^{n}\frac{\Gamma(1-\beta_{r})}{P_{r}}t^{\beta_{1}+\beta_{r}}.$$ The initial condition given in Eq. (\ref{initialcdn}) is taken to be zero, and let the additional condition be represented as follows
\begin{equation*}
	\int_{0}^{1}\int_{0}^{1}\hat{u}(x,y,t)dxdy=4\pi^{2} j(2k-1)\;\cdot\;{t}^{\Big(R-1-\displaystyle\sum_{\xi=2}^{n}\beta_{\xi}\Big)}.
\end{equation*}
We know $\lambda_{j,k}=\nu^{k}_{j}+\mu^{2}_{k}$, where $\mu^{2}_{j}=\pi/2(2j-1)$ and $\nu^{j}_{k}$ is given in (\ref{sss}).
\smallskip

\noindent For the given data, the solution \eqref{solution-second} can be written as 
\begin{eqnarray*}
	\hat{u}(x,y,t) &= &\sum_{j=1}^{\infty}\sum_{k=1}^{\infty} w(t)\hat{f}_{j,k}(t)\ast \mathbf{V}_{j,k}(t)X_{j,k}(x,y),
\end{eqnarray*}
where $\hat{f}_{j,k}(t)=\langle \hat{f}(x,y,t),Y_{j,k}(x,y)\rangle.$
\smallskip

\noindent
By applying {Lemma \ref{composite} and Proposition \ref{properties}}, we derive the solution $\hat{u}(x,y,t)$ in the form
\begin{eqnarray*}
	\hat{u}(x,y,t)&=&{\lambda_{j,k}}\;{t}^{\Big(R-1-\displaystyle\sum_{\xi=2}^{n}\beta_{\xi}\Big)}X_{j,k}(x,y)).
\end{eqnarray*}
Additionally, $w(t)$ given in Eq. (\ref{time}) can be written as
\begin{equation*}
	w(t)=\frac{{}_{G}\mathcal{D}_{t}^{(\varrho)} E(t)-\mathcal{T}(t)-\displaystyle\int_{0}^{t}\mathcal{G}(t,\tau)w(\tau)d\tau}{\displaystyle\int_{0}^{1}\int_{0}^{1}\hat{f}(x,y,t)dxdy},
\end{equation*}
where
\begin{eqnarray*}
	\int_{0}^{1}\int_{0}^{1}\hat{f}(x,y,t)dxdy&=&4\pi^{2} j(2k-1)\Big(C_{1}(t)+\lambda_{j,k}  C_{2}(t)\Big),\end{eqnarray*}$$
{}_{G}\mathcal{D}_{t}^{(\varrho)}E(t)=4\pi^{2} j(2k-1) \;{}_{G}\mathcal{D}_{t}^{(\varrho)}\;{t}^{\Big(R-1-\displaystyle\sum_{\xi=2}^{n}\beta_{\xi}\Big)},\qquad \mathcal{T}(t)=0,$$
\begin{eqnarray*}
	\mathcal{G}(t,\tau)&=&4\pi^{2} j(2k-1)\Big(C_{1}(t)+\lambda_{j,k}  C_{2}(t)\Big)\mathbf{V}_{j,k}(t-\tau).
\end{eqnarray*}
Thus, the expression for $w(t)$ becomes
$$w(t)={t}^{\Big(R-1-\displaystyle\sum_{\xi=1}^{n}\beta_{\xi}\Big)}.$$
Consider another source term like $f(x,y,t)=0$, then the solution corresponding to $f(x,y,t)$ is $u(x,y,t)=0.$
\\
\\
The error in $L^{2}(0,T;L^{2}((0,1)\times(0,1)))$ between the two solutions is
\begin{eqnarray*}
	\|\hat{u}-u\|_{L^{2}(0,T;L^{2}((0,1)\times(0,1)))}&=&\int_{0}^{T}\|\hat{u}-u\|_{L^{2}((0,1)\times(0,1))}dt \\&\leq&  {\lambda_{j,k}}\int_{0}^{T}{t}^{\Big(R-1-\displaystyle\sum_{\xi=2}^{n}\beta_{\xi}\Big)} \;dt\\&=& \frac{{\lambda_{j,k}}}{R-\displaystyle\sum_{\xi=2}^{n}\beta_{\xi}} {T}^{\Big(R-\displaystyle\sum_{\xi=2}^{n}\beta_{\xi}\Big)}.
\end{eqnarray*}
Taking the limit as ${(j,k) \to \infty}$, we find
\begin{eqnarray} \label{uillposed}
	\lim_{j,k \to \infty}\|\hat{u}-u\|_{L^{2}(0,T;L^{2}((0,1)\times(0,1)))}=\frac{{T}^{\Big(R-\displaystyle\sum_{\xi=2}^{n}\beta_{\xi}\Big)}}{R-\displaystyle\sum_{\xi=2}^{n}\beta_{\xi}}  \lim_{j,k \to \infty} \;{\lambda_{j,k}}=+\infty.
\end{eqnarray}
Therefore, based on Eq. (\ref{uillposed}), we guarantee that second-ISP is ill-posed.
\section{Examples}\label{exmp}
This section provides examples illustrating first-ISP and second-ISP.\\
\noindent\textbf{Example 01:}
For the specific case of first-ISP, consider the initial data 
$$\psi(x,y)={2	\chi_{j,k} x \sin\left(\sqrt{\nu^{k}_{j}} x\right)}\sin(\mu_k y),$$ 
where $	\chi_{j,k}$ is given by equation (\ref{ajk}). The final data is defined as
$$g(x,y)= \frac{2T^{\beta_{1}} \sin\left(\sqrt{\nu^{k}_{j}} x\right)\sin(\mu_k y)}{\Gamma(\beta_{1}+2)}.$$  
Applying the spectral expansion for $\psi(x,y)$ and  , the following expression are obtained
\begin{align*}
	\phi_{1,1}=&\langle\psi(x,y),Y_{1,1}(x,y)\rangle,\qquad j, k=1,\\
	\psi_{1,1}=&
	\frac{1}{\mu_k } \left( \frac{\sin\left(\sqrt{\nu^{1}_{1}}\right) - \sqrt{\nu^{1}_{1}} \cos\left(\sqrt{\nu^{1}_{1}}\right)}{\nu^{1}_{1}}\right) \\
	-&
	\frac{1}{\mu_k } \left(  \frac{\sin\left(\sqrt{\nu^{1}_{1}}\right)}{\sin\left(\sqrt{\nu^{1}_{0}}\right)} \cdot \frac{\sin\left(\sqrt{\nu^{1}_{0}}\right) - \sqrt{\nu^{1}_{0}} \cos\left(\sqrt{\nu^{1}_{0}}\right)}{\nu^{1}_{0}} \right),
\end{align*}
Using the spectral expansion for $g(x,y)$, we have
\begin{align*} 
	g_{1,1}=&\langle g(x,y),Y_{1,1}(x,y)\rangle,\qquad j, k=1,\\
	g_{1, 1}=&\frac{T^{\beta_{1}}}{\Gamma(\beta_{1}+2)}.
\end{align*}
Substituting the values of $\psi_{1,1}$ and $g_{1, 1}$ into Eq. (\ref{f0}), the solution $f(x,y)$ in Eq. (\ref{fsolution}) is given by
\begin{equation*}
	f(x,y)= 
	\frac{g_{1,1}}{\mathbf{Z}_{1,1}(T)}X_{1,1}(x,y) + \frac{\psi_{1,1}\mathbf{W}_{1,1}(T)}{\mathbf{Z}_{1,1}(T)}X_{1,1}(x,y).	\end{equation*}
Similarly, applying Eqs. (\ref{T0}) and (\ref{vsolution}), we obtain
\begin{align*}
	u(x,y,t)=& \frac{g_{1,1}\mathbf{Z}_{1,1}(t)}{\mathbf{Z}_{1,1}(T)}X_{1,1}(x,y)+ {\psi_{1,1}}\Big(\mathbf{W}_{1,1}(t)-\frac{ \mathbf{Z}_{1,1}(t)}{\mathbf{Z}_{1,1}(T)}\mathbf{W}_{1,1}(T)\Big)X_{1,1}(x,y).
\end{align*}
\smallskip
\noindent To find the solution for the mono-FDE with $\varrho(t)=\frac{t^{-\beta_{1}}}{\Gamma(1-\beta_{1})}$, we have
\begin{eqnarray*}
	f(x,y)&=& 
	\frac{g_{1,1}}{\displaystyle\int_{0}^{T} \tau^{\beta_{1}-1}{E}_{\beta_{1},\beta_{1}}(-\lambda_{1,1} \tau^{\beta_{1}})d\tau }X_{1,1}(x,y)\\&-&\frac{\psi_{1,1}{E}_{{\beta_1},1}(-\lambda_{1,1}T^{\beta_{1}})}{\displaystyle\int_{0}^{T} \tau^{\beta_{1}-1}{E}_{\beta_{1},\beta_{1}}(-\lambda_{1,1} \tau^{\beta_{1}})d\tau}X_{1,1}(x,y).	\end{eqnarray*}

\begin{eqnarray*}
	u(x,y,t)&=& \frac{g_{1,1}\displaystyle\int_{0}^{t} \tau^{\beta_{1}-1}{E}_{\beta_{1},\beta_{1}}(-\lambda_{1,1} \tau^{\beta_{1}})d\tau} {\displaystyle\int_{0}^{T} \tau^{\beta_{1}-1}{E}_{\beta_{1},\beta_{1}}(-\lambda_{1,1} \tau^{\beta_{1}})d\tau}X_{1,1}(x,y)+ \psi_{1,1}\bigg({E}_{{\beta_1},1}(-\lambda_{1,1}t^{\beta_{1}})\\&-&\frac{ \displaystyle\int_{0}^{t} \tau^{\beta_{1}-1}{E}_{\beta_{1},\beta_{1}}(-\lambda_{1,1} \tau^{\beta_{1}})d\tau}{\displaystyle\int_{0}^{T} \tau^{\beta_{1}-1}{E}_{\beta_{1},\beta_{1}}(-\lambda_{1,1} \tau^{\beta_{1}})d\tau}{E}_{{\beta_1},1}(-\lambda_{1,1}T^{\beta_{1}})\bigg)X_{1,1}(x,y).
\end{eqnarray*}
\smallskip
\noindent
\textbf{Example 02:}
In this example, we define the source term as
\begin{equation*}
	{f}(x,y,t)=2\varrho(t)\Big(C_{1}(t)+\lambda_{1,1}  C_{2}(t)\Big)  \sin\left(\sqrt{\nu^{1}_{1}} x\right)\sin(\mu_1 y),
\end{equation*}
for $j,k=1.$	Here, we have $$C_{1}(t)=\sum_{r=1}^{n}\frac{\Gamma(1-\beta_{r})\Gamma\Big(R-\displaystyle\sum_{\xi=2}^{n}\beta_{\xi} \Big)t^{\beta_1}}{\Gamma\Big(R-\beta_{r}-\displaystyle\sum_{\xi=2}^{n}\beta_{\xi}\Big)},\qquad R\geq2,$$and$$C_{2}(t)=\sum_{r=1}^{n}\frac{\Gamma(1-\beta_{r})}{P_{r}}t^{\beta_{1}+\beta_{r}}.$$ The initial condition in Eq. (\ref{initialcdn}) is set to zero, while an additional condition is
\begin{equation*}
	\int_{0}^{1}\int_{0}^{1}{u}(x,y,t)dxdy=\frac{8{t}^{\Big(R-1-\displaystyle\sum_{\xi=2}^{n}\beta_{\xi}\Big)}}{\pi^{2} j(2k-1)}, 
\end{equation*}
{Indeed, by employing equation \eqref{solution-second}, the solution to the system is given by
	\begin{equation*}
		U_{j,k}(t)=w(t)f_{1,1}(t)\ast \mathbf{V}_{1,1}(t),
	\end{equation*}
	where $f_{1,1}(t)=\langle f(x,y,t),Y_{1,1}(x,y)\rangle.$}
\smallskip

\noindent Due to {Lemma \ref{composite} and Proposition \ref{properties}}, we obtain the solution ${u}(x,y,t)$ as
\begin{equation*}
	{u}(x,y,t)=2{t}^{\Big(R-1-\displaystyle\sum_{\xi=2}^{n}\beta_{\xi}\Big)} \sin\left(\sqrt{\nu^{1}_{1}} x\right)\sin(\mu_1 y).
\end{equation*}
Moreover, $w(t)$ in (\ref{time}) is defined as follows
\begin{equation*}
	w(t)=\frac{{}_{G}\mathcal{D}_{t}^{(\varrho)} E(t)-\mathcal{T}(t)-\displaystyle\int_{0}^{t}\mathcal{G}(t,\tau)w(\tau)d\tau}{\displaystyle\int_{0}^{1}\displaystyle\int_{0}^{1}f(x,y,t)dxdy},
\end{equation*}
where
\begin{eqnarray*}
	\int_{0}^{1}\int_{0}^{1}f(x,y,t)dxdy&=&\frac{8\varrho(t)}{\pi^{2} j(2k-1)}\Big(C_{1}(t)+\lambda_{1,1}  C_{2}(t)\Big),\\
	{}_{G}\mathcal{D}_{t}^{(\varrho)}E(t)&=&\frac{8\varrho(t)}{\pi^{2} j(2k-1)}{}_{G}\mathcal{D}_{t}^{(\varrho)}\;{t}^{\Big(R-1-\displaystyle\sum_{\xi=2}^{n}\beta_{\xi}\Big)}, \quad \mathcal{T}(t)=0,\\
	\mathcal{G}(t,\tau)&=&\frac{8\varrho(t)}{\pi^{2} j(2k-1)}\Big(C_{1}(t)+\lambda_{1,1}  C_{2}(t)\Big)\mathbf{V}_{1,1}(t-\tau).
\end{eqnarray*}
Hence, 
$$w(t)={t}^{\Big(R-1-\displaystyle\sum_{\xi=1}^{n}\beta_{\xi}\Big)}.$$
\smallskip
\noindent For the mono-FDE solution, we use $\varrho(t)=\frac{t^{-\beta_{1}}}{\Gamma(1-\beta_{1})}$  with $C_{1}(t)$ and $C_{2}(t)$ given by
$$C_{1}(t)=\frac{\Gamma(1-\beta_{1})\Gamma(2)t^{\beta_1}}{\Gamma(2-\beta_{1})},\qquad C_{2}(t)={\Gamma(1-\beta_{1})}t^{2 \beta_{1}}.$$ 
Then, the solution takes the form
\begin{eqnarray*}
	{f}(x,y,t)=\Big(\frac{\Gamma(2)}{\Gamma(2-\beta_{1})}+\lambda_{1,1} t^{\beta_{1}}\Big) 2 \sin\left(\sqrt{\nu^{1}_{1}} x\right)\sin(\mu_1 y),
\end{eqnarray*}
\begin{equation*}
	{u}(x,y,t)=2t  \sin\left(\sqrt{\nu^{1}_{1}} x\right)\sin(\mu_1 y),
\end{equation*}
with
$$w(t)={t}^{(1-\beta_{1})}.$$
\vspace{0.3cm}

\noindent The analytical results presented in this chapter establish a rigorous foundation for ISPs in a two-dimensional generalized diffusion equation. We have shown that these problems are severely ill-posed in the Hadamard sense, requiring regularization to restore stability. The series representations derived here, based on bi-orthogonal systems of eigenfunctions, provide explicit formulas for the solutions and highlight the challenges that arise when attempting to reconstruct unknown source coefficients from noisy observational data. In the following chapters, we shift our focus from the analytical treatment of ISPs to their numerical solution. The regularized formulations of these ISPs will be discretized using finite difference schemes, leading to large structured linear systems of the all-at-once type. The GLT framework, developed in Part~I, will then become essential for understanding the spectral behavior of these systems and for designing efficient preconditioners that accelerate the convergence of Krylov subspace methods. In particular, Chapters~\ref{chapter:IFL} and~\ref{chapter:RFDE} address the numerical solution of ISPs for space-time FDEs, building upon the regularization and analytical techniques introduced in this chapter.

	\chapter{Krylov preconditioning for inverse problems in space-time fractional diffusion equations}\label{chapter:IFL}

Having introduced in the previous chapters the analytical foundations of fractional operators, the spectral framework of matrix sequences, and the theory of GLT sequences, we now turn to the numerical solution of an IP for a space-time FDE. This chapter is based on the work \cite{Ilyas2026a}, with additional details, numerical experiments, and reformulations included for completeness of the thesis presentation. This chapter investigates the reconstruction of a space-dependent source term in a model that incorporates a tempered fractional derivative in time and a Caputo fractional derivative in space, both of order in $(0,1)$. Such equations are particularly valuable for modeling complex transport phenomena where diffusion transitions from anomalous to normal behavior over time, a feature that cannot be captured by standard fractional or classical models \cite{Samko1993, Sandev2017}.

The inverse problem is inherently ill-posed in the Hadamard sense, as small perturbations in the final-time data can lead to large errors in the recovered source term. To address this instability, we apply a quasi-boundary value regularization method, which reformulates the problem into a well-posed form \cite{Denche2005, Hao2008}. The finite difference discretization of the regularized system leads to a large two-by-two block linear system. We conduct a detailed spectral analysis of the associated matrix sequences using the GLT framework, which provides a rigorous description of their asymptotic spectral behavior despite the presence of the variable coefficient that spoils the pure Toeplitz structure \cite{garoni2017, garoni2018}. Based on this analysis, we design an efficient preconditioner for the GMRES method. Numerical experiments confirm the effectiveness of the proposed preconditioning strategy in reducing iteration counts and computational cost.

The aim of this work is to numerically solve the IP of identifying the source term $f(x)$ in the space-time FDE
\begin{equation} \label{eq:problem}
	\begin{dcases}
		{}^T\!\mathcal{D}_t^{\xi,\rho} u(x,t) + a(x) {}^C\!\mathcal{D}_x^{\eta} u(x,t) = q(t)f(x),
		& x\in (0,1), \; t\in (0,T),
		\\
		u(0,t) = u(1,t) = 0,
		& t\in (0,T),
		\\
		u(x,0) = \varphi(x),
		& x\in (0,1),
		\\
		u(x,T) = \phi(x),
		& x\in (0,1),
	\end{dcases}
\end{equation}
where $T>0$ is the final time, $a(x)$ is a continuous a.e.\! variable coefficient such that $a(x)\neq 0$, $q(t)$ is a given function, $\varphi(x)$ is the initial condition, and $\phi(x)$ is the final time data. Moreover, ${}^{T}\!\mathcal{D}_{t}^{\xi,\rho}$ represents the tempered fractional derivative in the time variable of order $\xi\in(0,1)$, defined for $\rho>0$ as follows
\begin{equation*}
	{}^T\!\mathcal{D}_t^{\xi,\rho} u(x,t) := \frac{1}{\Gamma(1-\xi)} \int_0^t e^{-\rho(t-\tau)} (t-\tau)^{-\xi} \frac{\partial}{\partial\tau} u(x,\tau) \,\mathrm{d} \tau,
\end{equation*}
while ${}^C\!\mathcal{D}_x^{\eta}$ denotes the Caputo fractional derivative in the space variable of order $\eta\in(0,1)$, defined as
\begin{equation*}
	{}^C\!\mathcal{D}^\eta_x u(x,t) := \frac{1}{\Gamma(1-\eta)} \int_0^x {(x-\zeta)^{-\eta}} \frac{\partial}{\partial\zeta} u(\zeta,t) \,\mathrm{d}\zeta.
\end{equation*}
To address the ill-posedness in the sense of Hadamard of the inverse problem \eqref{eq:problem}, we apply the so-called nonlocal or quasi-boundary value method, transforming it into a well-posed formulation. First, we observe that the measurement of the final time data is inevitably contaminated by noise, therefore we denote the noisy data as $\phi_\varepsilon(x)$, where $\varepsilon>0$ is the noise level, and assume that
\begin{equation*}
	\norm{\phi - \phi_\varepsilon}_{L^2} \leq \varepsilon,
\end{equation*}
where $\norm{\cdot}_{L^2}$ is the $L^2$ norm on $(0,1)$. Applying the quasi-boundary value method leads to the following well-posed regularized problem
\begin{equation} \label{eq:quasi-problem}
	\begin{dcases}
		{}^T\!\mathcal{D}_t^{\xi,\rho} u(x,t) + a(x) {}^C\!\mathcal{D}_x^{\eta} u(x,t) = q(t)f_{\varepsilon,\lambda}(x),
		& x\in (0,1), \; t\in (0,T),
		\\
		u(0,t) = u(1,t) = 0,
		& t\in (0,T),
		\\
		u(x,0) = \varphi(x),
		& x\in (0,1),
		\\
		u(x,T) = \phi_\varepsilon(x) - \lambda\, a(x) {}^C\!\mathcal{D}_x^{\eta} f_{\varepsilon,\lambda}(x),
		& x\in (0,1),
	\end{dcases}
\end{equation}
where $\lambda>0$ is a regularization parameter, $f_{\varepsilon,\lambda}$ approximates $f(x)$ and the solution $u(x,t)$ represents an approximation of the solution of the original problem \eqref{eq:problem}.
\vspace{0.2cm}

The remainder of this chapter is organized as follows. Section \ref{sec:discre} provides the discretization of the regularized problem by using a finite difference scheme in space and time and presents the resulting linear system, emphasizing the structural features of the resulting matrices. In Section \ref{sec:glt-analysis} the coefficient matrix sequence is studied via GLT tools, while in Section \ref{sec:precond} the same type of analysis is performed for the preconditioned matrix sequences. In Section \ref{sec:num7}, numerical experiments are conducted to assess the effectiveness of the proposed preconditioner.

\section{Discretization and matrix structures} \label{sec:discre}

We now present the discretization of the regularized problem \eqref{eq:quasi-problem}, using a finite difference scheme \cite{Sun2020}. Then we describe with care the global structure of the arising coefficient matrix.


Let $n$ and $m$ be positive integers and denote the corresponding spatial and temporal meshes as
\begin{alignat}{3}
	x_i &= i\Delta x, && \qquad\Delta x = \frac{1}{n+1}, && \qquad i=0,1,\ldots,n, \label{eq:space-mesh} \\
	t_j &= j\Delta t, && \qquad\Delta t=\frac{T}{m}, && \qquad j=0,1,\ldots,m. \label{eq:time-mesh1}
\end{alignat}

To deal with the tempered fractional derivative, we use the $L1$ formula:
\begin{align*}
	{}^T\!\mathcal{D}_{L1}^{\xi,\rho} u(x,t_j) =& \frac{(\Delta t)^{-\xi}}{\Gamma(2-\xi)} \biggl[ b_0 u(x,t_j) - \sum_{l=1}^{j-1} \left( b_{j-l-1} - b_{j-l} \right) e^{\rho(t_l-t_j)} u(x,t_l) \\-& b_{j-1} e^{\rho(t_0-t_j)} u(x,t_0) \biggr],
\end{align*}
in which
\begin{equation*}
	b_l = (l+1)^{1-\xi}-l^{1-\xi},
	\qquad l = 0,1,\ldots
\end{equation*}
Denoting
\begin{equation*}
	\alpha_m := (\Delta t)^\xi\Gamma(2-\xi)
	= \bigg(\frac{T}{m}\bigg)^\xi\Gamma(2-\xi)
\end{equation*}
and
\begin{equation}
	\begin{split} \label{eq:time-coeff}
		\gamma_0 &:= b_0, \\
		\gamma_l &:= b_l - b_{l-1} = (l+1)^{1-\xi} - 2l^{1-\xi} + (l-1)^{1-\xi}, \qquad l=1,\ldots,m-1,
	\end{split}
\end{equation}
and using $t_j = j\Delta t$, the $L1$ approximation of the tempered fractional derivative is expressed as
\begin{equation} \label{eq:tem-discre}
	{}^T\!\mathcal{D}_{L1}^{\xi,\rho} u(x,t_j) = \frac{1}{\alpha_m} \left[ \gamma_0 u(x,t_j) + \sum_{l=1}^{j-1} \gamma_{j-l} e^{-(j-l)\rho\Delta t} u(x,t_l) - b_{j-1} e^{-j \rho \Delta t} u(x,t_0) \right].
\end{equation}
The following error estimate holds, {see \cite[page 10]{Zhao2021}}. Assuming that $u(\cdot,t)\in C^2(0,T)$. For $t_j\in (0,T)$, we have
\begin{align*}
	\abs{{}^T\!\mathcal{D}_t^{\xi,\rho} u(x,t)|_{t=t_j}-{}^T\!\mathcal{D}_{L1}^{\xi,\rho}u(x,t_j)} \leq& \chi_\xi
	\biggr[ \rho^2 \max_{t_0\leq t\leq t_j} \abs{u(x,t)} + 2\rho \max_{t_0\leq t\leq t_j} \abs{u_t(x,t)}\\ +& \Delta t^{2-\xi} \max_{t_0\leq t\leq t_j} \abs{u_{tt}(x,t)} \biggl],
\end{align*}
where $u_t$ and $u_{tt}$ denote the first and second partial derivative of $u$ with respect to $t$, and
\begin{equation*}
	\chi_\xi :=  \frac{1}{2\Gamma(1-\xi)} \left[ \frac{1}{4} + \frac{\xi}{(1-\xi)(2-\xi)} \right].
\end{equation*}

To discretize the CFD, we exploit the spatial counterpart of the $L1$ formula:
\begin{equation} \label{eq:Caputo-discre}
	{}^C\!\mathcal{D}_{L1}^\eta u(x_i,t) = \frac{1}{\beta_n} \left[ \delta_0 u(x_i,t) + \sum_{k=1}^{i-1} \delta_{i-k} u(x_k,t) - d_{i-1} u(x_0,t) \right],
\end{equation}
where
\begin{equation*}
	\beta_n := (\Delta x)^\eta \Gamma(2-\eta) = \bigg(\frac{1}{n+1}\bigg)^\eta \Gamma(2-\eta)
\end{equation*} and
\begin{equation} \label{eq:space-coeff}
	\begin{split}
		d_k &:= (k+1)^{1-\eta}-k^{1-\eta},
		\qquad k = 0,1,\ldots \\
		\delta_0 &:= d_0, \\
		\delta_k &:= d_k - d_{k-1} = (k+1)^{1-\eta} - 2k^{1-\eta} + (k-1)^{1-\eta}, \qquad k=1,\ldots,n-1.
	\end{split}
\end{equation}
The corresponding error estimate takes the form, {see \cite[page 32]{Sun2020}}. Assume that $u(x,\cdot)\in C^2(0,1)$. For $x_i\in (0,1)$, we have
\begin{equation*}
	\abs{{}^C\!\mathcal{D}_x^\eta u(x,t)|_{x=x_i}-{}^C\!\mathcal{D}_{L1}^\eta u(x_i,t)} \leq \chi_\eta
	\Delta x^{2-\eta} \max_{0\leq x \leq x_i} \abs{u_{xx}(x,t)},
\end{equation*}
where $u_{xx}$ denote the second partial derivative of $u$ with respect to $x$ and 
\begin{equation*}
	\chi_\eta :=  \frac{1}{2\Gamma(1-\eta)} \left[ \frac{1}{4} + \frac{\eta}{(1-\eta)(2-\eta)} \right].
\end{equation*}
\vspace{0.3cm}

Substituting the fractional derivatives in problem \eqref{eq:quasi-problem} with their respective $L1$ approximations \eqref{eq:tem-discre} and \eqref{eq:Caputo-discre}, and omitting the truncation errors, the finite difference scheme is constructed as follows: in $(x_i,t_j)$ for $i=1,\ldots,n$ and $j=1,\ldots,m$ we have
\begin{equation*}
	\begin{dcases}
		&\frac{1}{\alpha_m} \left[ \gamma_0 u_i^{(j)} + \sum_{l=1}^{j-1} \gamma_{j-l} e^{-(j-l)\rho\Delta t} u_i^{(l)} - b_{j-1} e^{-j\rho\Delta t} u_i^{(0)} \right]\\
		&+ \frac{a_i}{\beta_n}
		\left[ \delta_0 u_i^{(j)} + \sum_{k=1}^{i-1} \delta_{i-k} u_k^{(j)} - d_{i-1} u_0^{(j)} \right]
		= q^{(j)} f_i,
		\\
		&u_0^{(j)} = u_{n+1}^{(j)} = 0,
		\\
		&u_i^{(0)} = \varphi_i,
		\\
		&u_i^{(m)} = \phi_{\varepsilon,i} - \lambda \frac{a_i}{\beta_n}
		\left[ \delta_0 f_i + \sum_{k=1}^{i-1} \delta_{i-k} f_k - d_{i-1} f_0 \right],
	\end{dcases}
\end{equation*}
where
\begin{alignat*}{3}
	u_i^{(j)} &:= u(x_i,t_j), & \qquad
	a_i &:= a(x_i), & \qquad
	f_i &:= f_{\varepsilon,\lambda}(x_i),
	\\
	\varphi_i &:= \varphi(x_i), & \qquad
	\phi_{\varepsilon,i} &:= \phi_\varepsilon(x_i), & \qquad
	q^{(j)} &:= q(t_j).
\end{alignat*}
By defining the corresponding vectors
\begin{alignat*}{3}
	\mathbf{u}^{(j)} &:= \begin{bmatrix} u_i^{(j)} \end{bmatrix}_{i=1}^n, & \qquad
	\mathbf{a} &:= \begin{bmatrix} a_i \end{bmatrix}_{i=1}^n, & \qquad
	\mathbf{f} &:= \begin{bmatrix} f_i \end{bmatrix}_{i=1}^n,
	\\
	\boldsymbol{\varphi} &:= \begin{bmatrix} \varphi_i \end{bmatrix}_{i=1}^n, & \qquad
	\boldsymbol{\phi}_\varepsilon &:= \begin{bmatrix} \phi_{\varepsilon,i} \end{bmatrix}_{i=1}^n,
\end{alignat*}
we can write the numerical scheme in matrix form as
\begin{equation*}
	\begin{dcases}
		&\gamma_0 \mathbf{u}^{(j)} + \sum_{l=1}^{j-1} \gamma_{j-l} e^{-(j-l)\rho\Delta t} \mathbf{u}^{(l)}
		+ \frac{\alpha_m}{\beta_n} D_n(a) B_n \mathbf{u}^{(j)}
		\\&- \alpha_m q^{(j)} \mathbf{f}
		= b_{j-1} e^{-j\rho\Delta t} \boldsymbol{\varphi},
		\qquad j=1,\ldots,m, \\
		&\mathbf{u}^{(m)} + \frac{\lambda}{\beta_n} D_n(a) B_n \mathbf{f} = \boldsymbol{\phi}_\varepsilon,
	\end{dcases}
\end{equation*}
where $D_n(a)$ is the diagonal sampling matrix of size $n$ associated to $a(x)$, while $B_n$ is the following lower triangular Toeplitz matrix
\begin{equation*}
	B_n := \begin{bmatrix}
		\delta_0 & 0 & \cdots & \cdots & 0 \\
		\delta_1 & \delta_0 & 0 & & \vdots \\
		\delta_2 & \delta_1 & \delta_0 & \ddots & \vdots \\
		\vdots & \ddots & \ddots & \ddots& 0  \\
		\delta_{n-1} & \cdots & \delta_2 & \delta_1 & \delta_0
	\end{bmatrix}.
\end{equation*}
By gathering the equations for all the time levels and the regularization equation, we obtain the all-at-once linear system
\begin{equation} \label{eq:lin-syst}
	A_N\mathbf{v}_\lambda = \mathbf{z}_\lambda,
\end{equation}
where $N = N(m,n) := (m+1)n$ and
\begin{gather*}
	\renewcommand{\arraystretch}{1.5}
	A_N := \left[ \begin{array}{cccc|c}
		H_{n} & \mathrm{O}_n & \cdots & \mathrm{O}_n & -\alpha_m q^{(1)} I_n \\
		\gamma_1 e^{-\rho \Delta t} I_n & H_{n} & \ddots & \vdots & -\alpha_m q^{(2)} I_n \\
		\vdots & \ddots & \ddots & \mathrm{O}_n & \vdots \\
		\gamma_{m-1} e^{-(m-1)\rho\Delta t} I_n & \cdots & \gamma_1 e^{-\rho \Delta t} I_n & H_{n} & -\alpha_m q^{(m)} I_n \\
		\hline
		\mathrm{O}_n & \cdots & \mathrm{O}_n & I_n & \frac{\lambda}{\beta_n} G_n
	\end{array} \right]
	\in \C^{N\times N},
	\\
	\renewcommand{\arraystretch}{1.5}
	\mathbf{v}_{\lambda} :=
	\left[ \begin{array}{c}
		\mathbf{u}^{(1)} \\ \mathbf{u}^{(2)} \\ \vdots \\ \mathbf{u}^{(m)} \\ \hline \mathbf{f}
	\end{array} \right]
	\in\C^N,
	\qquad
	\mathbf{z}_{\lambda} :=
	\left[ \begin{array}{c}
		b_0 e^{-\rho\Delta t} \boldsymbol{\varphi} \\
		b_1 e^{-2\rho\Delta t} \boldsymbol{\varphi} \\
		\vdots \\
		b_{m-1} e^{-m\rho\Delta t} \boldsymbol{\varphi} \\
		\hline \boldsymbol{\phi}_\varepsilon
	\end{array} \right]
	\in\C^N,
\end{gather*}
in which $I_n$ is the $n\times n$ identity matrix, $\mathrm{O}_n$ is the $n\times n$ zero matrix, and
\begin{equation} \label{eq:Gm-def}
	H_{n}:=\gamma_0 I_n + \frac{\alpha_m}{\beta_n} G_n, \quad \textit{with}\quad G_n := D_n(a) B_n.
\end{equation}
Moreover, setting
\begin{equation*}
	U_m := \begin{bmatrix}
		\gamma_0 & 0 & \cdots & 0 \\
		\gamma_1 e^{-\rho\Delta t} & \gamma_0 & \ddots & \vdots \\
		\vdots & \ddots & \ddots & 0 \\
		\gamma_{m-1} e^{-\rho(m-1)\Delta t} & \cdots & \gamma_1 e^{-\rho\Delta t} & \gamma_0
	\end{bmatrix},
\end{equation*}
the coefficient matrix can be expressed into the compact form
\begin{equation} \label{eq:coeff-mat}
	\renewcommand{\arraystretch}{1.5}
	A_N = \left[ \begin{array}{c|c}
		U_m \otimes I_n + I_m\otimes \frac{\alpha_m}{\beta_n} G_n & -\alpha_m \mathbf{q}\otimes I_n \\
		\hline
		\mathbf{e}_m^T \otimes I_n & \frac{\lambda}{\beta_n} G_n
	\end{array} \right],
\end{equation}
in which
\begin{equation*}
	\mathbf{e}_m := \begin{bmatrix}
		0 \\ \vdots \\ 0 \\ 1
	\end{bmatrix} \in\R^m,
	\qquad
	\mathbf{q} := \begin{bmatrix}
		q^{(1)} \\ q^{(2)} \\ \vdots \\ q^{(m)}
	\end{bmatrix} \in\R^m.
\end{equation*}
\vspace{0.3cm}

{
	For concluding the present section, two remarks are useful, the first for clarifying the approximation order of the used scheme and the second concerning the non singular character of the coefficient matrix $A_N$.}

\begin{remark}
	{We remark that the tempered fractional derivative is discretized by the $L1$ finite difference scheme. Under the regularity assumption $u(\cdot,t)\in C^{2}(0,T)$, the $L1$ approximation satisfies an error estimate of order $O(\Delta t^{2-\xi})$, while the spatial $L1$-type discretization of the CFD yields an error of order $O(\Delta x^{2-\eta})$. These estimates are standard and well documented for FDEs. However, due to the weak initial-time singularity typically exhibited by solutions of time-FDEs, the effective temporal convergence rate on uniform grids is often reduced to first order. Since the IP is regularized through the quasi-boundary value method, the resulting formulation is well posed. Consequently, the consistency of the space-time discretization, together with stability of the discrete scheme, ensures convergence of the numerical solution of the linear system~\eqref{eq:lin-syst} to 	the solution of the regularized problem~\eqref{eq:quasi-problem} as $\Delta x,\Delta t \to 0$, despite the first-order accuracy in time.
	}
\end{remark}
\begin{remark} \label{remark:structure vs invertibility}
	{The structure of the matrix $A_N$ is very interesting because, up to a rank $m$ correction, it is a block lower triangular matrix.
		This suggests two things: the first is either an easy to invert preconditioner or a plain Sherman-Morrison-Woodbury (SMW) direct solver plus possibly iterative refinement, the second is a technique for checking the invertibility of the matrix $A_N$. As reported in the conclusions, the first item represents an additional research line which we leave for a future work, while the second is treated in detail here.
	}
	
	{
		We start from the expression of $A_N$ in (\ref{eq:coeff-mat}). Let us consider the symmetric permutation matrix
		\begin{equation*}
			\Pi_N \coloneqq
			\left[
			\begin{array}{ccc|c}
				& & I_{ n } & \\
				& \iddots & & \\
				I_{ n } & & & \\
				\hline
				& & & I_{ n }
			\end{array}
			\right]
		\end{equation*}
		Hence $A_N$ is invertible if and only if
		\begin{equation*}
			\renewcommand{\arraystretch}{1.3}
			\Pi_N A_N \Pi_N =
			\left[
			\begin{array}{c|c}
				Q_N & -\alpha_m\,\tilde{\mathbf{q}}\otimes I_n \\
				\hline
				\mathbf{e}_1^T \otimes I_n & \frac{\lambda}{\beta_n} G_n
			\end{array}
			\right]
		\end{equation*}
		is invertible with $\mathbf{\tilde{q}} \coloneqq \left[  q^{(m)}, \dots, q^{(1)} \right]^{T}$ and
		\begin{gather*}
			Q_N \coloneqq
			\begin{bmatrix}
				Q_{m,m} & Q_{m,m-1} & \cdots & \cdots & Q_{m,1} \\
				\mathrm{O}_n & Q_{m-1,m-1} & Q_{m-1,m-2} & \cdots & Q_{m-1,1} \\
				\mathrm{O}_n & \mathrm{O}_n & Q_{m-2,m-2} & \cdots & Q_{m-2,1} \\
				\vdots & & \ddots  & \ddots & \vdots \\
				\mathrm{O}_n & \cdots & \cdots & \mathrm{O}_n  & Q_{1,1}
			\end{bmatrix},
		\end{gather*}
		\begin{equation*}
			Q_{j,k} \coloneqq
			\begin{cases}
				u_{j,k} I_n
				& \text{for } j<k, \\
				u_{j,j} I_n + \frac{\alpha_m}{\beta_n} G_n
				& \text{for } j=k, \\
				0
				& \text{elsewhere},
			\end{cases}
		\end{equation*}
		in which $u_{j,k} = [U_n]_{j,k}$ are the entries of the matrix $U_n$ appearing in \eqref{eq:coeff-mat}.
		For the first step of the block Gaussian elimination, let us define
		\begin{equation*}
			P_N^{(1)} \coloneqq
			\begin{bmatrix}
				I_n & & \\
				& \ddots & \\
				- Q_{m,m}^{-1} & & I_n
			\end{bmatrix}.
		\end{equation*}
		As a consequence,
		\begin{equation*}
			\renewcommand{\arraystretch}{1.3}
			A_N^{(1)} \coloneqq P_N^{(1)} A_N =
			\left[
			\begin{array}{c|c}
				Q_N & -\alpha_m\,\tilde{\mathbf{q}}\otimes I_n \\
				\hline
				\boldsymbol{\alpha}^{(1)} & \boldsymbol{\beta}^{(1)}
			\end{array}
			\right]
		\end{equation*}
		with
		\begin{align*}
			\boldsymbol{\alpha}^{(1)} & \coloneqq
			\left[
			\begin{array}{c|c|c|c|c}
				\mathrm{O}_n & -Q_{m,m}^{-1} Q_{m,m-1} & -Q_{m,m}^{-1} Q_{m,m-2} & \cdots & -Q_{m,m}^{-1} Q_{m,1}
			\end{array}
			\right],
			\\
			\boldsymbol{\beta}^{(1)} & \coloneqq \frac{\lambda}{\beta_n} G_n +\alpha_m q^{(m)} Q_{m,m}^{-1}.
		\end{align*}
		Iterating the block Gaussian steps, at each step $k=2,\ldots,m$ we define
		\begin{equation*}
			P_N^{(k)} \coloneqq
			\left[
			\begin{array}{ccc|c}
				I_n & & & \\
				& \ddots & & \\
				& & I_n \\
				\hline
				\mathbf{p}_1^{(k)} & \cdots & \mathbf{p}_m^{(k)}
				& I_n
			\end{array}
			\right],
			\;\;\;\;
			\mathbf{p}^{(k)}_j \coloneqq
			\begin{cases}
				\mathrm{O}_n
				& \text{if } j\neq k \\
				-\boldsymbol{\alpha}_k^{(k-1)} Q_{m+1-k,m+1-k}^{-1}
				& \text{if } j = k,
			\end{cases}
		\end{equation*}
		so that
		\begin{equation*}
			A_N^{(k)} \coloneqq P_N^{(k)} A_N^{(k-1)} =
			\left[
			\begin{array}{c|c}
				Q_N & \tilde{\mathbf{{q}}} \\
				\hline
				\begin{matrix}
					\boldsymbol{\alpha}_1^{(k)} & \cdots & \boldsymbol{\alpha}_m^{(k)}
				\end{matrix}
				& \boldsymbol{\beta}^{(k)}
			\end{array}
			\right],
		\end{equation*}
		where
		\begin{align*}
			\boldsymbol{\alpha}_j^{(k)} & \coloneqq
			\begin{cases}
				\mathrm{O}_n
				& \text{for } j=1,\dots,k \\
				\mathbf{p}_k^{(k)} Q_{m+1-k,m+1-j} + \boldsymbol{\alpha}_j^{(k-1)}
				& \text{for } j=k+1,\dots m,
			\end{cases}
			\\
			\boldsymbol{\beta}^{(k)} & \coloneqq
			-\alpha_m\,q^{(m+1-k)} \mathbf{p}_k^{(k)} + \boldsymbol{\beta}^{(k-1)}.
		\end{align*}
		In conclusion, for $k=m$ we find
		\begin{equation*}
			A_N^{(m)} =
			\left[
			\begin{array}{c|c}
				Q_N & -\alpha_m\,\tilde{\mathbf{q}}\otimes I_n \\
				\hline
				\mathrm{O}_{n\times mn } & \boldsymbol{\beta}^{(m)}
			\end{array}
			\right].
		\end{equation*}
		We are left to analyze the spectrum of $\boldsymbol{\beta}^{(m)}$, which by expanding the recursive definition takes the form
		\begin{equation*}
			\boldsymbol{\beta}^{(m)} = -\sum_{k=1}^m \alpha_m\, q^{(m+1-k)} \mathbf{p}_k^{(k)} + \frac{\lambda}{\beta_n} G_n.
		\end{equation*}
		It is evident that $\boldsymbol{\beta}^{(m)} $ is a sum of products and inverses of lower triangular matrices, so that it is lower triangular as well. Its non singular character holds if and only if the matrix $A_N$ is invertible, and the check can be done on the diagonal entries of $\boldsymbol{\beta}^{(m)}$. This criterion is easy to check and computationally inexpensive, since it only involves element-wise operations, and in fact it has been checked to hold true in our numerical experiments in Section \ref{sec:num7}.}
\end{remark}

\section{GLT analysis} \label{sec:glt-analysis}

This section is devoted to the spectral analysis of the coefficient matrix sequence via the GLT theory, in order to gather information on the asymptotic behavior of the matrices. We start with the sequence associated with the space block $G_n$, then proceed to the time block $U_m$, and conclude with the overall coefficient matrix sequence.

\begin{theorem} \label{thm:GLT-Bm}
	There exists a continuous $2\pi$-periodic function $g_\eta$ such that
	\begin{equation*}
		\{B_n\}_n \sim_{\glt,\sigma} g_\eta(\theta),
		\qquad \theta\in [-\pi,\pi].
	\end{equation*}
	Moreover,
	\begin{equation*}
		\{B_n\}_n \sim_\lambda \delta_0.
	\end{equation*}
\end{theorem}
\begin{proof}
	To demonstrate that $\{B_n\}_n$ is a GLT sequence, we claim that $B_n$ is a Toeplitz matrix generated by a measurable function, i.e., there exists $g_\eta$ such that $B_n = T_n(g_\eta)$ for all $n\in\N$. The  idea is simple: we prove that the series
	\begin{equation*}
		g_\eta(\theta) := \sum_{k=0}^{\infty} \delta_k e^{\iu k\theta},
	\end{equation*}
	where the coefficients $\delta_k$ are given by \eqref{eq:space-coeff}, belongs to the Wiener class, so that $g_\eta$ is well defined, continuous and $2\pi$-periodic, and the corresponding Toeplitz matrix sequence is precisely $\{B_n\}_n$ by definition.
	
	Let us consider the auxiliary function $\psi_\eta(x) := x^{1-\eta}$. From \eqref{eq:space-mesh}, for any index $k$ we have $k = (n+1)x_k$. Hence,
	\begin{align*}
		\delta_k &= (k+1)^{1-\eta} - 2k^{1-\eta} + (k-1)^{1-\eta}
		\\
		&= (n+1)^{1-\eta} \big[ x_{k+1}^{1-\eta} - 2 x_k^{1-\eta} + x_{k-1}^{1-\eta} \big]
		\\
		&= (n+1)^{1-\eta} \big[ \psi_\eta(x_{k+1}) - 2\psi_\eta(x_k) + \psi_\eta(x_{k-1}) \big]
		\\
		&= (n+1)^{1-\eta} \left[ \frac{ \psi_\eta(x_{k+1}) - 2\psi_\eta(x_k) + \psi_\eta(x_{k-1}) }{\Delta x^2} \right] \Delta x^2
		\\
		&= \frac{1}{(n+1)^{1+\eta}} \big[ \psi_\eta''(x_k) + O(\Delta x^2) \big] \qquad \text{as } \Delta x\to 0,
	\end{align*}
	where we used a standard second order central finite difference scheme and the fact that $\Delta x = \frac{1}{n+1}$. From
	\begin{equation*}
		\psi_\eta''(x) = \eta(\eta-1) \frac{1}{x^{1+\eta}},
	\end{equation*}
	we have
	\begin{align*}
		\delta_k
		&= \frac{\eta(\eta-1)}{(n+1)^{1+\eta}} \cdot  \frac{1}{x_k^{1+\eta}} + \Delta x^{1+\eta} O(\Delta x^2) \\
		&= \frac{\eta(\eta-1)}{(n+1)^{1+\eta}} \cdot \frac{(n+1)^{1+\eta}}{k^{1+\eta}} + O(\Delta x^{3+\eta}) \\
		&= \eta(\eta-1) \frac{1}{k^{1+\eta}} + O(\Delta x^{3+\eta}).
	\end{align*}
	Since $\eta\in (0,1)$, we deduce that the series
	\begin{equation*}
		\sum_{k=0}^{\infty} \abs{\delta_k}
	\end{equation*}
	is convergent, and consequently $g_\eta$ belongs to the Wiener class. Therefore, {since $g_\eta$ belongs to the Wiener class, it is a continuous and hence Lebesgue integrable function, so that $B_n = T_n(g_\eta)$ and Axiom GLT \ref{glt:symbols} and Axiom GLT \ref{glt:distrib}} imply
	\begin{equation*}
		\{B_n\}_n \sim_{\glt,\sigma} g_\eta(\theta).
	\end{equation*}
	Finally, all the eigenvalues of $B_n$ are trivially equal to $\delta_0$ because of the triangular structure, so that
	\begin{equation*}
		\{B_n\}_n \sim_\lambda \delta_0
	\end{equation*}
	follows immediately from the definition of distribution.
\end{proof}

\begin{corollary}
	If $a:[0,1]^2\to\R$ is continuous a.e, then
	\begin{equation*}
		\{G_n\}_n  \sim_{\glt,\sigma} a(x) g_\eta(\theta),
		\qquad x\in [0,1], \;\theta\in [-\pi,\pi].
	\end{equation*}
	with $g_\eta$ as in Theorem \ref{thm:GLT-Bm}. Moreover,
	\begin{equation*}
		\{G_n\}_n  \sim_\lambda \delta_0 a(x),
		\qquad x\in [0,1].
	\end{equation*}
\end{corollary}
\begin{proof}
	Recalling \eqref{eq:Gm-def}, by GLT \ref{glt:symbols} we have
	\begin{equation*}
		\{D_n(a)\}_n \sim_\glt a(x), \qquad x\in [0,1].
	\end{equation*}
	Then, by Theorem \ref{thm:GLT-Bm} and the $\ast$-algebra properties detailed in GLT \ref{glt:algebra}, we deduce
	\begin{equation*}
		\{G_n\}_n \sim_\glt a(x)g_\eta(\theta),
		\qquad x\in [0,1], \;\theta\in [-\pi,\pi].
	\end{equation*}
	By GLT \ref{glt:distrib}, we immediately have
	\begin{equation*}
		\{G_n\}_n \sim_\sigma a(x)g_\eta(\theta),
		\qquad x\in [0,1], \;\theta\in [-\pi,\pi].
	\end{equation*}
	Finally, we notice that $G_n$ is lower triangular and the main diagonal is given by $\begin{bmatrix} \delta_0 a(x_i) \end{bmatrix}_{i=1}^n$. In conclusion
	\begin{equation*}
		\{G_n\}_n \sim_\lambda \delta_0 a(x),
		\qquad x\in [0,1].
	\end{equation*}
\end{proof}

We proceed to the analysis of $\{U_m\}_m$, relying on the same technique employed in the proof of Theorem \ref{thm:GLT-Bm}.

\begin{theorem} \label{thm:GLT-Un}
	There exists a continuous $2\pi$-periodic function $h_\xi$ such that
	\begin{equation*}
		\{U_m\}_m \sim_{\glt,\sigma} h_\xi(\theta),
		\qquad \theta\in [-\pi,\pi].
	\end{equation*}
	Moreover,
	\begin{equation*}
		\{U_m\}_m \sim_\lambda \gamma_0.
	\end{equation*}
\end{theorem}

\begin{proof}
	We follow the proof of Theorem \ref{thm:GLT-Bm}, showing that there exists $h_\xi$ such that $U_m = T_m(h_\xi)$ for all $m\in\N$. We prove that the series    \begin{equation*}
		h_\xi(\theta) := \sum_{k=0}^{\infty} \gamma_k e^{-\rho k \Delta t} e^{\iu k\theta},
	\end{equation*}
	with $\gamma_k$ defined by \eqref{eq:time-coeff}, belongs to the Wiener class. Considering $\psi_\xi(t) := t^{1-\xi}$ and using \eqref{eq:time-mesh1}, we have
	\begin{align*}
		\gamma_k &= (k+1)^{1-\xi} - 2k^{1-\xi} + (k-1)^{1-\xi}
		\\
		&= \left( \frac{m}{T} \right)^{1-\xi} \left[ t_{k+1}^{1-\xi} - 2 t_k^{1-\xi} + t_{k-1}^{1-\xi} \right]
		\\
		&= \left( \frac{m}{T} \right)^{1-\xi} \big[ \psi_\xi(t_{k+1}) - 2\psi_\xi(t_k) + \psi_\xi(t_{k-1}) \big]
		\\
		&= \left( \frac{m}{T} \right)^{1-\xi} \left[ \frac{ \psi_\xi(t_{k+1}) - 2\psi_\xi(t_k) + \psi_\xi(t_{k-1}) }{\Delta t^2} \right] \Delta t^2
		\\
		&= \left(\frac{T}{m}\right)^{1+\xi} \big[ \psi_\xi''(t_k) + O(\Delta t^2) \big] \qquad \text{as } k\to\infty,
	\end{align*}
	using $\Delta t = \frac{T}{m}$ and a standard central finite difference scheme. From
	\begin{equation*}
		\psi_\xi''(x) = \xi(\xi-1) \frac{1}{t^{1+\xi}},
	\end{equation*}
	we have
	\begin{align*}
		\gamma_k
		&= \left(\frac{T}{m}\right)^{1+\xi} \xi(\xi-1) \frac{1}{t_k^{1+\xi}} + \Delta t^{1+\xi} O(\Delta t^2) \\
		&= \left(\frac{T}{m}\right)^{1+\xi} \frac{\xi(\xi-1)}{k^{1+\xi}} \left(\frac{m}{T}\right)^{1+\xi} + O(\Delta t^{3+\xi}) \\
		&= \frac{\xi(\xi-1)}{k^{1+\xi}} + O(\Delta t^{3+\xi}).
	\end{align*}
	Since $\xi\in (0,1)$ and $\rho>0$, the series
	\begin{equation*}
		\sum_{k=0}^{\infty} \abs{\gamma_k e^{-\rho k \Delta t}}
	\end{equation*}
	is convergent, and $h_\xi$ belongs to the Wiener class by definition. {Since $h_\xi$ belongs to the Wiener class, it is a continuous and hence Lebesgue integrable function, so that $U_m = T_m(h_\xi)$. As a consequence of Axiom GLT \ref{glt:symbols} and of Axiom GLT \ref{glt:distrib}} we conclude
	\begin{equation*}
		\{U_m\}_m \sim_{\glt,\sigma} h_\xi(\theta),
	\end{equation*}
	while the eigenvalue distribution is obvious from the triangular structure.
\end{proof}
\begin{remark}
	{As it is evident from the final part of the proofs of Theorem \ref{thm:GLT-Bm} and Theorem \ref{thm:GLT-Un}, the fact that the involved generating functions belong to the Wiener class is not a goal, but it is instrumental in applying the GLT axioms for deducing in a fast way the distribution in the sense of the singular values.}
\end{remark}

We are now ready to compute the GLT symbol of the sequence associated with the full coefficient matrix. To provide meaningful results, we introduce a mild assumption regarding the rate at which the discretization parameters tend to zero with respect to one another.

\begin{theorem} \label{thm:coeff-mat-sym}
	With $A_N$ defined in \eqref{eq:coeff-mat}, assume that
	\begin{equation*}
		\frac{\alpha_m}{\beta_n} \to \nu, \qquad \text{as } m,n\to\infty,
	\end{equation*}
	with $\nu\neq 0$. Then, it holds
	\begin{equation*}
		\{A_N\}_N \sim_{\glt,\sigma} h_\xi(\theta_1) + \nu a(x)g_\eta(\theta_2), \qquad (\theta_1,\theta_2)\in [-\pi,\pi]^2, \; x\in [0,1],
	\end{equation*}
	where $h_\xi$ and $g_\eta$ are given in Theorems \ref{thm:GLT-Bm}-\ref{thm:GLT-Un}.
\end{theorem}

\begin{proof}
	We begin by decomposing \eqref{eq:coeff-mat} as
	\begin{align*}
		\renewcommand{\arraystretch}{1.5}
		A_N = &\underbrace{ \left[ \begin{array}{cc}
				U_m \otimes I_n & \mathrm{O}_{mn\times n} \\
				\mathrm{O}_{n\times mn} & \mathrm{O}_n
			\end{array} \right] }_{=:\, C_N^{(1)}}
		+ \underbrace{ \left[ \begin{array}{cc}
				I_m\otimes \frac{\alpha_m}{\beta_n} G_n & \mathrm{O}_{mn\times n} \\
				\mathrm{O}_{n\times mn} & \frac{\lambda}{\beta_n} G_n
			\end{array} \right] }_{=:\, C_N^{(2)}}
		+ \underbrace{ \left[ \begin{array}{cc}
				\mathrm{O}_{mn} & -\alpha_m\mathbf{q}\otimes I_n \\
				\mathbf{e}_m^T \otimes I_n & \mathrm{O}_n
			\end{array} \right] }_{=:\, C_N^{(3)}},
	\end{align*}
	so we can focus on each component individually. For the first one, we have
	\begin{equation*}
		\renewcommand{\arraystretch}{1.3}
		C_N^{(1)}
		= U_{m+1} \otimes I_n
		- \underbrace{\left[ \begin{array}{c|c}
				\mathrm{O}_{mn} & \mathrm{O}_{mn\times n} \\
				\hline
				\begin{matrix}
					\gamma_m e^{-\rho m\Delta t} I_n & \cdots & \gamma_1 e^{-\rho\Delta t} I_n
				\end{matrix} & \gamma_0 I_n
			\end{array} \right]}_{=: R^{(1)}_N},
	\end{equation*}
	where $\rk R^{(1)}_N \leq n = o(N)$ as $m,n\to\infty$. By Lemma \ref{lemma:zero-distr}, $\{R^{(1)}_N\}_N$ is zero-distributed, hence by GLT \ref{glt:symbols}
	\begin{equation*}
		\{R^{(1)}_N\}_N\sim_\glt 0.
	\end{equation*}
	Moreover, we have $U_{m+1} \otimes I_n = T_{m+1}(h_\xi) \otimes T_n(1)$, therefore applying GLT \ref{glt:symbols}--\ref{glt:kron} leads to
	\begin{equation} \label{eq:S1}
		\big\{C_N^{(1)}\big\}_N \sim_\glt h_\xi(\theta_1).
	\end{equation}
	
	Concerning $C_N^{(2)}$, we notice that
	\begin{align*}
		C_N^{(2)} =& \left[ \begin{array}{cc}
			I_m\otimes \frac{\alpha_m}{\beta_n} G_n & \mathrm{O}_{mn\times n} \\
			\mathrm{O}_{n\times mn} & \frac{\alpha_m}{\beta_n} G_n + \frac{\lambda-\alpha_m}{\beta_n} G_n
		\end{array} \right]= I_{m+1} \otimes \frac{\alpha_m}{\beta_n} G_n +
		\underbrace{\left[ \begin{array}{cc}
				\mathrm{O}_{mn} & \mathrm{O}_{mn\times n} \\
				\mathrm{O}_{n\times mn} & \frac{\lambda-\alpha_m}{\beta_n} G_n
			\end{array} \right]}_{=: R_N^{(2)}},
	\end{align*}
	in which once again $\rk R_N^{(2)} = n = o(N)$ as $m,n\to\infty$, so that
	\begin{equation*}
		\{R_N^{(2)}\}_N\sim_\glt 0.
	\end{equation*}
	On the other hand, with the assumption that $\frac{\alpha_m}{\beta_n}\to \nu$, by GLT \ref{glt:symbols} and GLT \ref{glt:kron} we have
	\begin{equation*}
		\left\{I_{m+1} \otimes \frac{\alpha_m}{\beta_n} G_n\right\}_N = \left\{ T_{m+1}(1) \otimes \frac{\alpha_m}{\beta_n} D_n(a)T_n(g_\eta) \right\}_N \sim_\glt \nu a(x) g_\eta(\theta_2).
	\end{equation*}
	We conclude by GLT \ref{glt:algebra} that
	\begin{equation} \label{eq:S2}
		\big\{C_N^{(2)}\big\}_N \sim_\glt \nu a(x) g_\eta(\theta_2).
	\end{equation}
	
	Finally, we have
	\begin{equation*}
		C_N^{(3)} = \left[ \begin{array}{cc}
			\mathrm{O}_{mn} & -\alpha_m\mathbf{q}\otimes I_n \\
			\mathrm{O}_{n\times mn} & \mathrm{O}_n
		\end{array} \right]
		+ \left[ \begin{array}{cc}
			\mathrm{O}_{mn} & \mathrm{O}_{mn\times n} \\
			\mathbf{e}_m^T \otimes I_n & \mathrm{O}_n
		\end{array} \right],
	\end{equation*}
	where each of these components has rank at most $m=o(N)$. Therefore,
	\begin{equation} \label{eq:S3}
		\big\{C_N^{(3)}\big\}_N \sim_\glt 0.
	\end{equation}
	
	By summing the contributions in \eqref{eq:S1}, \eqref{eq:S2} and \eqref{eq:S3}, we find the GLT symbol of the overall coefficient matrix sequence $\{A_N\}_N$ via GLT \ref{glt:algebra}, and then the singular value distribution via GLT \ref{glt:distrib}.
	
\end{proof}
\vspace{0.3cm}

\noindent Since the matrices \(A_N\) are non-Hermitian, the GLT symbol obtained in Theorem~\ref{thm:coeff-mat-sym} directly implies the singular value distribution via GLT~\ref{glt:distrib}. However, the GLT theory alone does not guarantee the eigenvalue distribution. In the present setting, the eigenvalue distribution can be established because the preconditioned matrix sequence is weakly clustered at one, as shown in Theorem~\ref{glt_prec1}. This clustering result relies on the fact that the perturbation between the preconditioner and the coefficient matrix is of low rank relative to the matrix size, which is sufficient to characterize the eigenvalue behavior in the sense of weak clustering. In general, for non-Hermitian matrix sequences, eigenvalue distributions require additional assumptions, such as Hermitian structure, normality, or suitable perturbation arguments; see \cite{Barbarino2020, GoSe} and references therein.

\section{GLT preconditioning} \label{sec:precond}

In the present section, we propose preconditioners $P_N$ for the linear system \eqref{eq:lin-syst}, with the aim of clustering at 1 the eigenvalues of $\{ P_N^{-1} A_N \}_N $. Let us define
\begin{equation} \label{eq:precond}
	\renewcommand{\arraystretch}{1.5}
	P_N := \left[ \begin{array}{c|c}
		U_m \otimes I_n + I_m\otimes \frac{\alpha_m}{\beta_n} G_n & \mathrm{O}_{nm \times n} \\
		\hline
		\mathbf{e}_m^T \otimes I_n & \frac{\lambda}{\beta_n} G_n
	\end{array} \right].
\end{equation}
Note that both $U_m$ and $G_n$ are invertible, in consequence of their lower triangular structure with non-zero diagonal elements, and the same is true for $U_m \otimes I_n + I_m\otimes \frac{\alpha_m}{\beta_n} G_n$, which is also lower triangular. Because of the block triangular structure, the eigenvalues of $P_N$ are the eigenvalues of the diagonal blocks, therefore $P_N$ is invertible.

The following theorem proves that $\{P_N\}_N$ has the same GLT symbol, and hence the same singular value distribution, as the coefficient matrix sequence.

\begin{theorem} \label{thm:precond-sym}
	With $P_N$ defined in \eqref{eq:precond}, assume that
	\begin{equation*}
		\frac{\alpha_m}{\beta_n} \to \nu, \qquad \text{as } m,n\to\infty,
	\end{equation*}
	with $\nu\neq 0$. Then, it holds
	\begin{equation*}
		\{P_N\}_N \sim_{\glt,\sigma} h_\xi(\theta_1) + \nu a(x)g_\eta(\theta_2), \qquad (\theta_1,\theta_2)\in [-\pi,\pi]^2, \; x\in [0,1],
	\end{equation*}
	where $h_\xi$ and $g_\eta$ are given in Theorems \ref{thm:GLT-Bm}-\ref{thm:GLT-Un}.
\end{theorem}

\begin{proof}
	Immediate by following the same reasoning as in the proof of Theorem \ref{thm:coeff-mat-sym}.
\end{proof}

The GLT symbol and singular value distribution of the preconditioned matrix sequence can be derived as a corollary of Theorems \ref{thm:coeff-mat-sym} and \ref{thm:precond-sym}, simply by exploiting the algebra properties in GLT \ref{glt:algebra}. However, it can also be computed more in general avoiding the assumption on the discretization parameters, as we show in the following theorem. Moreover, we can derive the distribution of the eigenvalues.

\begin{theorem} \label{glt_prec1}
	With $A_N$ as in \eqref{eq:coeff-mat} and $P_N$ as in \eqref{eq:precond}, it holds
	\begin{equation*}
		\{P_N^{-1}A_N\}_N \sim_{\glt,\sigma} 1.
	\end{equation*}
	Moreover,
	\begin{equation*}
		\{P_N^{-1}A_N\}_N \sim_\lambda 1.
	\end{equation*}
\end{theorem}

\begin{proof}
	Setting $R_N := A_N - P_N$, we have
	\begin{equation*}
		P_N^{-1} A_N = I_N + P_N^{-1}R_N.
	\end{equation*}
	From the fact that $\rk R_N\leq n$, it follows $\rk(P_N^{-1}R_N) \leq n = o(N)$, therefore
	\begin{equation*}
		\{P_N^{-1}R_N\}_N \sim_\glt 0
	\end{equation*}
	by Lemma \ref{lemma:zero-distr}. Hence, by GLT \ref{glt:algebra} and the fact that $\{I_N\}_N \sim_\glt 1$,
	\begin{equation*}
		\{P_N^{-1} A_N\}_N \sim_\glt 1.
	\end{equation*}
	The singular value distribution follows immediately from GLT \ref{glt:distrib}. Concerning the eigenvalue distribution, it is obvious that
	\begin{equation*}
		\lambda_j (P_N^{ -1 } A_N ) = 1 + \lambda_j (P_N^{-1} R_N),
		\qquad j=1,\ldots,N,
	\end{equation*}
	therefore
	\begin{equation*}
		\# \left\{ j=1,\ldots,N  \;\colon\; \lambda_j (P_N^{-1} A_N) \neq 1 \right\} \leq n = o(N).
	\end{equation*}
	In other words, the eigenvalues of $\{ P_N^{-1} A_N \}_N$ are weakly clustered at 1, in the sense of Definition \ref{clustering}. Recalling Remark \ref{remark:distrib-cluster}, this is equivalent to the fact that the sequence is spectrally distributed as 1.
\end{proof}
\vspace{0.2cm}

\noindent Recent developments in preconditioning for all-at-once systems and inverse problems for FDEs are closely related to the present work. Li, Lin, Hon, and Wu \cite{LiLinHon2024} developed a preconditioned MINRES method for block lower triangular Toeplitz systems, which are ubiquitous in evolutionary PDE discretizations. More recently, Lin, Ke, Hon, and Ng \cite{LinKe2026} proposed a matching Schur complement preconditioning technique specifically designed for discrete time fractional diffusion ISPs, demonstrating its effectiveness in accelerating the convergence of Krylov methods. These works highlight the growing interest in structure-preserving preconditioners for all-at-once formulations, and the block triangular preconditioner proposed in this chapter builds upon this foundation.

\begin{remark}\label{rem:comput1}
	In the special case where $a(x)=a$ is constant, so that $D_n(a)=aI_n$ and $G_n = aB_n$, a more efficient preconditioner can be constructed by approximating each Toeplitz matrix appearing in $P_N$ with the corresponding Strang preconditioner, as follows:
	\begin{equation*}
		\renewcommand{\arraystretch}{1.5}
		S_N := \left[ \begin{array}{c|c}
			S(U_m) \otimes I_n + I_m\otimes \frac{\alpha_m}{\beta_n} a S(B_n) & \mathrm{O}_{nm \times n} \\
			\hline
			\mathbf{e}_m^T \otimes I_n & \frac{\lambda}{\beta_n} a S(B_n)
		\end{array} \right],
	\end{equation*}
	where $S(U_m)$, $S(B_n)$ are the standard Strang preconditioners of $U_m$ and $B_n$, respectively. Since circulant matrices are uniformly diagonalizable through the Fourier matrix $F_m$, we can write
	\begin{equation*}
		S(U_m) = F_m \Lambda_m^{(1)} F_m^{*},
		\qquad
		S(B_n) = F_n \Lambda_n^{(2)} F_n^{*},
	\end{equation*}
	where $\Lambda_m^{(i)}$, $i=1,2$, are the diagonal matrices containing the eigenvalues and can be computed starting from the coefficients of $U_m$, $B_n$ via the FFT. Then, $S_N$ becomes
	\begin{equation*}
		\renewcommand{\arraystretch}{1.5}
		S_N := \left[ \begin{array}{c|c}
			(F_m\otimes F_n)
			\big[\Lambda_m^{(1)}\otimes I_n + I_m\otimes \frac{\alpha_m}{\beta_n} a \Lambda_n^{(2)} \big](F_m\otimes F_n)^{*} & \mathrm{O}_{nm \times n} \\
			\hline
			\mathbf{e}_m^T \otimes I_n & \frac{\lambda}{\beta_n} a F_n \Lambda_n^{(2)} F_n^{*}
		\end{array} \right].
	\end{equation*}
	With this approach, the computational cost of inverting $S_N$ is of $O\big(nm\log(nm)\big)$ arithmetic operations, since all the necessary computations can be performed via the FFT. The same reasoning can be carried out using any suitable preconditioner for Toeplitz matrices, drawn from a matrix algebra which is uniformly diagonalizable through a fast transform.
\end{remark}

\begin{remark}\label{rem:comput2}
	From an implementation viewpoint, we can define a custom function that inverts $P_N$ through a forward substitution method, thanks to the block triangular structure of the matrix and the triangular structure of the blocks themselves. The computational cost amounts to $O(n^2m^2)$ arithmetic operations.
	{However, a linear system with coefficient matrix $\frac{\lambda}{\beta_n} D_n(a)B_n=\frac{\lambda}{\beta_n} G_n$ and variable coefficient $a$ can be still solved within $O\big(n\log n\big)$ arithmetic operations as in Remark \ref{rem:comput1}, while the more difficult additive structure in the other diagonal blocks $\gamma_0 I_n +\frac{\alpha_m}{\beta_n}G_n$ can be treated e.g. via preconditioned GMRES and Strang preconditioner $\frac{\alpha_m}{\beta_n} \bar a S(B_n)$, with $\bar a$ being an average of $a(x)$, again within $O\big(n\log n\big)$ arithmetic operations.	
	}
\end{remark}

\section{Numerical results}\label{sec:num7}

We conclude by assessing the efficiency of our preconditioning strategy in terms of iterations required by the GMRES.

Let us consider problem \eqref{eq:quasi-problem} in which we set
\begin{gather*}
	T = 1, \qquad a(x) = x, \qquad
	q(t) = t^2, \qquad \varphi(x) = 0.
\end{gather*}

To construct the noisy final data, we select the source function
\begin{equation*}
	f(x) = x\sin(\pi x)
\end{equation*}
and solve the linear system arising from the corresponding direct problem
\begin{equation*}
	\begin{dcases}
		{}^T\!\mathcal{D}_t^{\xi,\rho} u(x,t) + a(x) {}^C\!\mathcal{D}_x^{\eta} u(x,t) = q(t)f(x),
		& x\in (0,1), \; t\in (0,T),
		\\
		u(0,t) = u(1,t) = 0,
		& t\in (0,T),
		\\
		u(x,0) = \varphi(x),
		& x\in (0,1).
	\end{dcases}
\end{equation*}
With the same discretization technique adopted in Section \ref{sec:discre}, we obtain the linear system
\begin{equation*}
	\tilde{A}_{m,n}\mathbf{u} = \mathbf{b},
\end{equation*}
where
\begin{gather*}
	\renewcommand{\arraystretch}{1.5}
	\tilde{A}_{m,n} := \left[ \begin{array}{cccc}
		H_{n} & \mathrm{O}_n & \cdots & \mathrm{O}_n \\
		\gamma_1 e^{-\rho \Delta t} I_n & H_{n} & \ddots & \vdots \\
		\vdots & \ddots & \ddots & \mathrm{O}_n \\
		\gamma_{m-1} e^{-(m-1)\rho\Delta t} I_n & \cdots & \gamma_1 e^{-\rho \Delta t} I_n & H_{n}
	\end{array} \right]
	\in \C^{mn\times mn},
	\\
	\renewcommand{\arraystretch}{1.5}
	\mathbf{u} :=
	\left[ \begin{array}{c}
		\mathbf{u}^{(1)} \\
		\mathbf{u}^{(2)} \\
		\vdots \\
		\mathbf{u}^{(m)}
	\end{array} \right]
	\in\C^{mn},
	\qquad
	\mathbf{b} :=
	\left[ \begin{array}{c}
		b_0 e^{-\rho\Delta t} \boldsymbol{\varphi} + \alpha_m q^{(1)} \mathbf{f} \\
		b_1 e^{-2\rho\Delta t} \boldsymbol{\varphi} + \alpha_m q^{(2)} \mathbf{f} \\
		\vdots \\
		b_{m-1} e^{-m\rho\Delta t} \boldsymbol{\varphi} + \alpha_m q^{(n)} \mathbf{f}
	\end{array} \right]
	\in\C^{nm}.
\end{gather*}
By solving this system, we construct final time data $\mathbf{u}^{(m)}$, which we perturb by adding a random quantity $\varepsilon\delta (x)$, where $\delta (x)$ is randomly selected in the interval $(-1,1)$. In this way, we simulate the contamination of $\phi(x)$ by noise, obtaining $\boldsymbol{\phi}_\varepsilon$.

Then, we apply the GMRES method to \eqref{eq:lin-syst} and report in Table \ref{table:gmres} the number of iterations {and the CPU timings (see Remark \ref{remark:cpu-timing})} required by both the nonpreconditioned and preconditioned linear systems. We set the regularization parameter to $\lambda=5\cdot10^{-3}$, the noise level to $\varepsilon=0.01$, and vary both fractional orders. The initial guess is the zero vector, the tolerance is $\texttt{tol}=10^{-8}$, and the maximum number of iterations is equal to the size of the matrix, $N=(m+1)n$. The experiments are performed on Intel(R), Xeon(R), Silver 4114, CPU@2.20GHz.

Comparing the columns corresponding to the preconditioned and nonpreconditioned linear system, we observe that, for each choice of the fractional orders, $P_N$ is indeed effective in reducing the number of iterations, confirming the theoretical predictions. It also appears to be stable, maintaining good results as the problems sizes increase and across all choices of fractional orders. Indeed, $P_N$ demonstrates the best performance when the iteration counts for the nonpreconditioned linear system are particularly high, i.e, for $\xi=0.2$, $\eta=0.8$.

Furthermore, for each choice of $\xi$ and $\eta$, the convergence rate in the preconditioned framework seems to depend almost exclusively on the spatial discretization finesse: in fact, for each fixed $n$, the number of iterations does not vary significantly as we increase $m$. This is reasonable, considering that the number of outliers in the cluster at 1 of the preconditioned matrix depends only on $n$.

Looking more in detail, we observe that the growth in the parameter $n$ is less than linear, i.e., less than a constant times the square root of the global size, when $n$ and $m$ grow proportionally. The described behaviour agrees with the negative results in~\cite{MultiNo2,MultiNo}, where theoretical barriers to the optimality of preconditioned Krylov solvers are proven.

\begin{table}[ht] 
	\centering
	\caption{Iterations {and CPU times} required for solving the linear system with the GMRES method, without and with the preconditioner, for various choices of the fractional orders.}
	\label{table:gmres}
	\begin{adjustbox}{width=\textwidth}
		\begin{tabular}{cccccccccccccc}
			\toprule
			& & \multicolumn{4}{c}{$\xi=0.2$} & \multicolumn{4}{c}{$\xi=0.5$} & \multicolumn{4}{c}{$\xi=0.8$}
			\\
			& & \multicolumn{4}{c}{$\eta=0.8$} & \multicolumn{4}{c}{$\eta=0.5$} & \multicolumn{4}{c}{$\eta=0.2$}
			\\
			\cmidrule(lr){3-6}\cmidrule(lr){7-10}\cmidrule(lr){11-14}
			& & \multicolumn{2}{c}{Iter} & \multicolumn{2}{c}{CPU} & \multicolumn{2}{c}{Iter} & \multicolumn{2}{c}{CPU} & \multicolumn{2}{c}{Iter} & \multicolumn{2}{c}{CPU}
			\\
			\cmidrule(lr){3-4}\cmidrule(lr){5-6}\cmidrule(lr){7-8}\cmidrule(lr){9-10}\cmidrule(lr){11-12}\cmidrule(lr){13-14}
			$n$ & $m$ & - & $P_N$ & - & $P_N$ & - & $P_N$ & - & $P_N$ & - & $P_N$ & - & $P_N$
			\\
			\midrule
			$2^{4}$ & $2^{4}$ & 66 & 15 & 0.11 & 0.03 & 53 & 17 & 0.02 & 0.01 & 45 & 16 & 0.03 & 0.01 \\
			$2^{4}$ & $2^{5}$ & 62 & 16 & 0.04 & 0.02 & 60 & 17 & 0.04 & 0.02 & 63 & 16 & 0.03 & 0.01 \\
			$2^{4}$ & $2^{6}$ & 64 & 17 & 0.06 & 0.03 & 67 & 17 & 0.07 & 0.03 & 86 & 17 & 0.08 & 0.03 \\
			$2^{4}$ & $2^{7}$ & 61 & 17 & 0.11 & 0.05 & 67 & 17 & 0.13 & 0.06 & 141 & 17 & 0.30 & 0.05 \\
			$2^{4}$ & $2^{8}$ & 61 & 17 & 0.37 & 0.15 & 79 & 17 & 0.49 & 0.15 & 230 & 16 & 2.02 & 0.14 \\
			\midrule
			$2^{5}$ & $2^{4}$ & 107 & 26 & 0.04 & 0.01 & 75 & 30 & 0.02 & 0.02 & 49 & 28 & 0.01 & 0.01 \\
			$2^{5}$ & $2^{5}$ & 110 & 27 & 0.08 & 0.04 & 79 & 29 & 0.05 & 0.03 & 68 & 27 & 0.04 & 0.03 \\
			$2^{5}$ & $2^{6}$ & 110 & 27 & 0.16 & 0.06 & 86 & 29 & 0.12 & 0.06 & 90 & 27 & 0.13 & 0.06 \\
			$2^{5}$ & $2^{7}$ & 103 & 27 & 0.50 & 0.15 & 89 & 29 & 0.41 & 0.16 & 142 & 27 & 0.77 & 0.15 \\
			$2^{5}$ & $2^{8}$ & 100 & 27 & 1.20 & 0.42 & 96 & 30 & 1.13 & 0.47 & 230 & 28 & 3.76 & 0.43 \\
			\midrule
			$2^{6}$ & $2^{4}$ & 213 & 36 & 0.17 & 0.03 & 101 & 49 & 0.06 & 0.04 & 51 & 40 & 0.02 & 0.03 \\
			$2^{6}$ & $2^{5}$ & 211 & 35 & 0.32 & 0.05 & 103 & 48 & 0.11 & 0.07 & 71 & 39 & 0.06 & 0.06 \\
			$2^{6}$ & $2^{6}$ & 207 & 36 & 1.34 & 0.14 & 110 & 47 & 0.43 & 0.20 & 94 & 39 & 0.31 & 0.15 \\
			$2^{6}$ & $2^{7}$ & 193 & 37 & 2.05 & 0.37 & 109 & 48 & 0.84 & 0.47 & 145 & 41 & 1.30 & 0.39 \\
			$2^{6}$ & $2^{8}$ & 194 & 36 & 4.48 & 1.13 & 117 & 49 & 2.27 & 1.54 & 231 & 40 & 5.65 & 1.26 \\
			\midrule
			$2^{7}$ & $2^{4}$ & 503 & 46 & 1.46 & 0.08 & 136 & 71 & 0.16 & 0.12 & 56 & 57 & 0.05 & 0.13 \\
			$2^{7}$ & $2^{5}$ & 493 & 46 & 6.68 & 0.18 & 136 & 71 & 0.66 & 0.32 & 76 & 56 & 0.25 & 0.24 \\
			$2^{7}$ & $2^{6}$ & 458 & 46 & 8.78 & 0.36 & 137 & 70 & 1.07 & 0.59 & 98 & 56 & 0.60 & 0.45 \\
			$2^{7}$ & $2^{7}$ & 440 & 46 & 11.96 & 0.88 & 136 & 72 & 1.80 & 1.64 & 146 & 56 & 2.42 & 1.44 \\
			$2^{7}$ & $2^{8}$ & 436 & 48 & 28.55 & 3.20 & 138 & 71 & 5.55 & 5.67 & 232 & 56 & 10.76 & 3.90 \\
			\midrule
			$2^{8}$ & $2^{4}$ & 1237 & 60 & 37.41 & 0.31 & 181 & 101 & 1.10 & 0.62 & 59 & 77 & 0.18 & 0.43 \\
			$2^{8}$ & $2^{5}$ & 1165 & 59 & 50.73 & 0.53 & 174 & 100 & 1.53 & 0.99 & 78 & 76 & 0.39 & 0.73 \\
			$2^{8}$ & $2^{6}$ & 1149 & 58 & 68.13 & 1.10 & 171 & 99 & 2.30 & 1.99 & 100 & 76 & 1.01 & 1.43 \\
			$2^{8}$ & $2^{7}$ & 1044 & 59 & 92.23 & 3.03 & 168 & 100 & 4.67 & 5.41 & 150 & 76 & 3.92 & 3.85 \\
			$2^{8}$ & $2^{8}$ & 985 & 60 & 151.27 & 8.66 & 166 & 101 & 11.27 & 15.08 & 234 & 76 & 17.51 & 10.48 \\
			\bottomrule
		\end{tabular}
	\end{adjustbox}
\end{table}


\begin{remark} \label{remark:cpu-timing}
	{
		Table \ref{table:gmres} also reports the measured CPU timings, even if, in this setting, they do not compare the preconditioned and non-preconditioned cases in a totally meaningful way. In fact, since the algorithms are implemented in MATLAB, we are limited in our ability to optimize the code compared with what is possible in more sophisticated programming languages. Therefore, the measured CPU time also reflects implementation-specific overhead rather than the intrinsic efficiency of the algorithms, and may therefore be misleading.}
\end{remark}

	\chapter{Efficient Krylov solvers for inverse source problem in 2D space-time fractional diffusion equation}\label{chapter:RFDE}

The motivation for this work stems from the increasing importance of space-time FDEs, which offer a versatile framework for modeling complex, multiscale phenomena in a rich variety of scientific disciplines. By incorporating fractional derivatives, FDEs can capture essential features such as memory effects, non-local interactions, and anomalous diffusion behaviors. These characteristics make FDEs ideal for modeling processes that involve complex transport phenomena, including anomalous diffusion, subdiffusion, and superdiffusion in heterogeneous media.  In particular, the Caputo fractional derivative is widely used for modeling spatial heterogeneities, as it incorporates physically meaningful initial conditions and avoids singularities at the boundaries. Meanwhile, tempered fractional derivatives are effective in describing transient anomalous diffusion, where the diffusion initially follows a fractional power law and later transitions to normal diffusion \cite{Sandev2017}. 
\vspace{0.15cm}

\noindent In this chapter, based on the results presented in \cite{IlyasCapizzano2026}, we develop a fast algorithm to reconstruct the source term \(f(x,y)\) in the context of a two-dimensional space-time FDE, generalizing the work in \cite{Ilyas2026a} devoted to the uni-dimensional space version. Additional details, numerical experiments, and reformulations are included for completeness of the thesis presentation. Since IPs are often ill-posed in the sense of Hadamard, small perturbations in the input data can lead to significant errors in the solution \cite{IlyasMalik2024CAM, Jin2015, Tuan2017, IlyasIqbal2024, Ilyas2026b}. To address this, we apply the quasi-boundary value method, a nonlocal regularization technique that has been successfully utilized to handle similar ill-posed problems in both parabolic and FDEs \cite{Denche2005, Hao2008}. The finite difference discretization of the regularized problem leads to a two-by-two block linear system. We design symbol-based approximations \cite{Ilyas2026a,HonLiSormani2025} to the resulting coefficient matrices and construct an effective preconditioner for the GMRES method \cite{Saad1986GMRES}. A detailed spectral analysis of the associated matrix sequences is carried out using the GLT theory, which provides a rigorous description of their asymptotic spectral behavior despite the presence of the variable coefficient \( a(x,y) \) that spoils the pure Toeplitz structure \cite{garoni2017, garoni2018, barbarino2020blockmulti, barbarino2020block1d}. Numerical experiments confirm the effectiveness of the proposed preconditioning strategy in reducing iteration counts and computational cost.

\subsection*{Problem Formulation}

We consider the inverse source problem of identifying the source term \( f(x,y) \) from the space-time fractional diffusion equation
\begin{equation}\label{eq:problem2}
	\begin{dcases} &{}^{T}\mathcal{D}_{t}^{\xi,\rho}u(x,y,t)+ a(x,y)\left({}^{C}\mathcal{D}_{x}^{\eta_{1}}+{}^{C}\mathcal{D}_{y}^{\eta_{2}}\right)u(x,y,t)=q(t)f(x,y), \; (x,y,t) \in \Omega_{T},\\ &	u(0,y,t)=0=u(1,y,t), \qquad y \in (0,1), \quad t\in(0,T),\\
		&u(x,0,t)=0=u(x,1,t), \qquad x \in (0,1), \quad t\in(0,T),\\
		&u(x,y,0) = \varphi(x,y), \qquad x \in (0,1), \quad y \in (0,1), \\
		&u(x,y,T)=\phi(x,y),\qquad x\in(0,1), \quad y\in(0,1),\end{dcases}\end{equation}
where \(\Omega_{T} := (0,1) \times (0,1) \times (0,T)\), \({}^{T}\mathcal{D}_{t}^{\xi,\rho}\) denotes the tempered fractional derivative of order \(\xi \in (0,1)\) and \(\rho > 0\) is the truncation parameter, \(a(x,y)\) is a continuous variable coefficient such that \(a(x,y) \neq 0\), \({}^{C}\mathcal{D}_{x}^{\eta_{1}}\) and \({}^{C}\mathcal{D}_{y}^{\eta_{2}}\) denote the Caputo fractional derivatives of orders \(\eta_{1}, \eta_{2} \in (0,1)\), \(q(t)\) is a given function, \(\varphi(x,y)\) is the initial condition, and \(\phi(x,y)\) is the final time datum.
\smallskip

\noindent The tempered fractional derivative in time variable of order $\xi\in(0,1)$ is defined as follows
$${}^{T}\mathcal{D}_{t}^{\xi,\rho}u(x,y,t):=\frac{1}{\Gamma(1-{\xi})}\int_{0}^{t}e^{-\rho(t-\tau)}(t-\tau)^{-{\xi}}\frac{\partial}{\partial \tau}u(x,y,\tau)d\tau.$$
The Caputo fractional derivative in the space variable of orders $\eta\in(0,1)$ is defined as follows
\[{}^{C}\mathcal{D}^{\eta}_{x}u(x,y,t)=\frac{1}{\Gamma(1-\eta)}\int\limits^{x}_{0}{(x-\zeta)^{-\eta}}\frac{\partial}{\partial \zeta}{u(x,y,\zeta)}\;d\zeta,\qquad x >0,\]

In practice, the measurement of the final time datum is inevitably affected by noise. Let us denote the noisy datum by \(\phi_{\varepsilon}(x,y)\), where \(\varepsilon > 0\) represents the noise level. Assuming the noise is sufficiently small, we impose the condition

\begin{equation*}
	\norm{\phi - \phi_{\varepsilon}}_{L^2} \leq \varepsilon,
\end{equation*}

where \(\norm{\cdot}_{L^2}\) is the \(L^2\) norm on \((0,1) \times (0,1)\). Applying the quasi-boundary value method yields a class of parameter-dependent well-posed regularized problems
\begin{equation}\label{eq:quasi-problem2}
	\begin{dcases} &{}^{T}\mathcal{D}_{t}^{\xi,\rho}u(x,y,t)+ a(x,y)\left({}^{C}\mathcal{D}_{x}^{\eta_{1}}+{}^{C}\mathcal{D}_{y}^{\eta_{2}}\right)u(x,y,t)=q(t)f_{\varepsilon,\lambda}(x,y), \; (x,y,t) \in \Omega_{T},\\	&u(0,y,t)=0=u(1,y,t), \qquad y \in (0,1), \quad t\in(0,T),\\
		&u(x,0,t)=0=u(x,1,t), \qquad  x \in (0,1), \qquad t\in(0,T),\\
		&u(x,y,0) = \varphi(x,y), \qquad x \in (0,1), \qquad y \in (0,1), \\
		&u(x,y,T)=\phi_{\varepsilon}(x,y) - \lambda\, a(x,y) \left({}^{C}\mathcal{D}_{x}^{\eta_{1}}+{}^{C}\mathcal{D}_{y}^{\eta_{2}}\right) f_{\varepsilon,\lambda}(x,y), x\in(0,1), \quad y\in(0,1),\end{dcases}
\end{equation}
which approximate the original FDE in \eqref{eq:problem2}. Here, \(\lambda > 0\) is a regularization parameter, \(f_{\varepsilon,\lambda}\) represents the regularized approximation of the source term \(f(x,y)\), and the solution \(u(x,y,t) \equiv u_{\varepsilon,\lambda}(x,y,t)\) provides an approximation to the original solution of problem \eqref{eq:problem2}. The regularization term \(\lambda \, a(x,y) \left({}^{C}\mathcal{D}_{x}^{\eta_{1}} + {}^{C}\mathcal{D}_{y}^{\eta_{2}}\right)\) ensures that the reconstructed source term is stable and well-behaved, even in the presence of noise.
\vspace{0.2cm}

The structure of the paper is as follows. Section \ref{sec:discre2} outlines the discretization of the regularized problem using a finite difference scheme for both space and time, and discusses the structural properties of the resulting linear system. In Section \ref{sec:glt-analysis2}, we analyze the coefficient matrix sequence using GLT theory. Section \ref{sec:precond2} extends this analysis to the preconditioned matrix sequences. In Section \ref{sec:num2}, we conduct numerical experiments to evaluate the performance of the proposed preconditioner and we report a detailed analysis of the computational costs.

\section{Discretization and matrix structures} \label{sec:discre2}

In the current section, we present the discretization of the regularized problem \eqref{eq:quasi-problem2} using a finite difference scheme \cite{Sun2020}. We then carefully discuss the global structure of the resulting coefficient matrix.


Let $n_{1}$, $n_{2}$ and $m$ be positive integers and denote the corresponding spatial and temporal meshes as
\begin{alignat}{3}
	x_i &= i\Delta x, && \qquad\Delta x = \frac{1}{n_{1}+1}, && \qquad i=0,1,\ldots,n_{1}, \label{eq:1space-mesh} \\
	y_j &= j\Delta y, && \qquad\Delta y = \frac{1}{n_{2}+1}, && \qquad j=0,1,\ldots,n_{2}, \label{eq:2space-mesh} \\
	t_s &= s\Delta t, && \qquad\Delta t=\frac{T}{m}, && \qquad s=0,1,\ldots,m. \label{eq:time-mesh}
\end{alignat}
The approximate values of the functions at the grid points are denoted as follows
\begin{alignat*}{3}
	u_{i,j}^{s} &:= u(x_i,y_{j},t_s), & \qquad
	a_{i,j} &:= a(x_i,y_{j}), & \qquad
	f_{i,j} &:= f_{\varepsilon,\lambda}(x_i,y_{j}),
	\\
	\varphi_{i,j} &:= \varphi(x_i,y_{j}), & \qquad
	\phi_{i,j} &:= \phi_\varepsilon(x_i,y_{j}), & \qquad
	q^{s} &:= q(t_s).
\end{alignat*}
To deal with the tempered fractional derivative, we use the $L1$ formula:
\begin{equation*}
	{}^T\!\mathcal{D}_{L1}^{\xi,\rho} u(x_{i},y_{j},t_{s}) = \frac{(\Delta t)^{-\xi}}{\Gamma(2-\xi)} \left[ b_0 u_{i,j}^{s} - \sum_{l=1}^{s-1} \left( b_{s-l-1} - b_{s-l} \right) e^{\rho(t_l-t_s)} u_{i,j}^{l} - b_{s-1} e^{\rho(t_0-t_s)} u_{i,j}^{0} \right],
\end{equation*}
in which
\begin{equation*}
	b_l = (l+1)^{1-\xi}-l^{1-\xi},
	\qquad l = 0,1,\ldots
\end{equation*}
Denoting
\begin{equation*}
	\alpha_{m} := (\Delta t)^\xi\Gamma(2-\xi)
	= \bigg(\frac{T}{m}\bigg)^\xi\Gamma(2-\xi)
\end{equation*}
and
\begin{equation}
	\begin{split} \label{eq:time-coeff2}
		\gamma_0 &:= b_0, \\
		\gamma_l &:= b_l - b_{l-1} = (l+1)^{1-\xi} - 2l^{1-\xi} + (l-1)^{1-\xi}, \qquad l=1,\ldots,m-1,
	\end{split}
\end{equation}
and using $t_s = s\Delta t$, the $L1$ approximation of the tempered fractional derivative is expressed as
\begin{equation} \label{eq:tem-discre2}
	{}^T\!\mathcal{D}_{L1}^{\xi,\rho} u(x_{i},y_{j},t_s) = \frac{1}{\alpha_m} \left[ \gamma_0 u_{i,j}^{s} + \sum_{l=1}^{s-1} \gamma_{s-l} e^{-(s-l)\rho\Delta t} u_{i,j}^{l} - b_{s-1} e^{-s \rho \Delta t} u_{i,j}^{0} \right].
\end{equation}

To discretize the Caputo fractional derivative with respect to $x$, we exploit the $L1$ formula:
\begin{equation} \label{eq:Caputo-xdiscre}
	{}^C\!\mathcal{D}_{L1}^{\eta_{1}} u(x_i,y_{j},t_{s}) = \frac{1}{\beta_{n_{1}}} \left[ \delta_0 u_{i,j}^{s} + \sum_{k_{1}=1}^{i-1} \delta_{i-k_{1}} u_{k_{1},j}^{s} - d_{i-1} u_{0,j}^{s} \right],
\end{equation}
where
\begin{equation}\label{eq:beta-n1}
	\beta_{n_1} := (\Delta x)^{\eta_{1}} \Gamma(2-\eta_{1}) = \bigg(\frac{1}{n_{1}+1}\bigg)^{\eta_{1}} \Gamma(2-\eta_{1})
\end{equation}
and
\begin{equation} \label{eq:space-xcoeff}
	\begin{split}
		d_{k_{1}} &:= (k_{1}+1)^{1-\eta_{1}}-k_{1}^{1-\eta_{1}},
		\qquad k_{1} = 0,1,\ldots \\
		\delta_0 &:= d_0, \\
		\delta_{k_{1}} &:= d_{k_{1}} - d_{k_{1}-1} = (k_{1}+1)^{1-\eta_{1}} - 2k_{1}^{1-\eta_{1}} + (k_{1}-1)^{1-\eta_{1}}, \qquad k_{1}=1,\ldots,n_{1}-1.
	\end{split}
\end{equation}
\vspace{0.3cm}

Now, we proceed to discretize the Caputo fractional derivative with respect to $y$ by employing the $L1$ approximation formula
\begin{equation} \label{eq:Caputo-ydiscre}
	{}^C\!\mathcal{D}_{L1}^{\eta_{2}} u(x_i,y_{j},t_{s}) = \frac{1}{\beta_{n_{2}}} \left[ \sigma_0 u_{i,j}^{s} + \sum_{k_{2}=1}^{j-1} \sigma_{j-k_{2}} u_{i,k_{2}}^{s} - c_{j-1} u_{i,0}^{s} \right],
\end{equation}
where
\begin{equation}\label{eq:beta-n2}
	\beta_{n_2} := (\Delta x)^{\eta_{2}} \Gamma(2-\eta_{2}) = \bigg(\frac{1}{n_{2}+1}\bigg)^{\eta_{2}} \Gamma(2-\eta_{2})
\end{equation}
and
\begin{equation} \label{eq:space-ycoeff}
	\begin{split}
		c_{k_{2}} &:= (k_{2}+1)^{1-\eta_{2}}-k_{2}^{1-\eta_{2}},
		\qquad k_{2} = 0,1,\ldots \\
		\sigma_0 &:= c_0, \\
		\sigma_{k_2} &:= c_{k_2} - c_{k_{2}-1} = (k_{2}+1)^{1-\eta_{2}} - 2k_{2}^{1-\eta_{2}} + (k_{2}-1)^{1-\eta_{2}}, \qquad k_{2}=1,\ldots,n_{2}-1.
	\end{split}
\end{equation}
\vspace{0.3cm}

Substituting the fractional derivatives in problem \eqref{eq:quasi-problem2} with their respective $L1$ approximations \eqref{eq:tem-discre2} and \eqref{eq:Caputo-xdiscre} and \eqref{eq:Caputo-ydiscre}, and omitting the truncation errors, the finite difference scheme is constructed as follows:
\begin{equation*}
	\begin{dcases}
		\frac{1}{\alpha_m} \Bigg[ \gamma_0 u_{i,j}^{s} + \sum_{l=1}^{s-1} \gamma_{s-l} e^{-(n-l)\rho\Delta t} u_{i,j}^{l} - b_{s-1} e^{-s \rho \Delta t} u_{i,j}^{0} \Bigg]
		+ \frac{a_{i,j}}{\beta_{n_{1}}} \Bigg[ \delta_0 u_{i,j}^{s} + \sum_{k_{1}=1}^{i-1} \delta_{i-k_{1}} u_{k_{1},j}^{s} \\- d_{i-1} u_{0,j}^{s} \Bigg]+\frac{a_{i,j}}{\beta_{n_{2}}} \left[ \sigma_0 u_{i,j}^{s} + \sum_{k_{2}=1}^{j-1} \sigma_{j-k_{2}} u_{i,k_{2}}^{s} - c_{j-1} u_{j,0}^{s} \right]
		= q^{s} f_{i,j},
		\\
		u_{0,j}^{s} = u_{n_{1}+1,j}^{s} = 0,\qquad u_{i,0}^{s} = u_{i,n_{2}+1}^{s} = 0,
		\\
		u_{i,j}^{0} = \varphi_{i,j},
		\\
		u_{i,j}^{m} = \phi_{i,j} - \lambda \frac{a_{i,j}}{\beta_{n_{1}}}
		\Bigg[ \delta_0 f_{i,j} + \sum_{k_{1}=1}^{i-1} \delta_{i-k_{1}} f_{k_{1},j} - d_{i-1} f_{0,j} \Bigg] - \lambda \frac{a_{i,j}}{\beta_{n_{2}}}
		\Bigg[ \sigma_0 f_{i,j} + \sum_{k_{2}=1}^{j-1} \sigma_{j-k_{2}} f_{i,k_{2}}\\ - c_{j-1} f_{i,0} \Bigg].
	\end{dcases}
\end{equation*}
By defining the corresponding vectors
\begin{alignat*}{3}
	\mathbf{u}^{s} &:= \begin{bmatrix} u_{i,j}^{s} \end{bmatrix}_{i=1,\ldots,n_{1},\,j=1,\ldots,n_{2}}, & \quad
	\mathbf{a} &:= \begin{bmatrix} a_{i,j} \end{bmatrix}_{i=1,\ldots,n_{1},\,j=1,\ldots,n_{2}}, & \quad
	\mathbf{f} &:= \begin{bmatrix} f_{i,j} \end{bmatrix}_{i=1,\ldots,n_{1},\,j=1,\ldots,n_{2}},
	\\
	\boldsymbol{\varphi} &:= \begin{bmatrix} \varphi_{i,j} \end{bmatrix}_{i=1,\ldots,n_{1},\,j=1,\ldots,n_{2}}, & \quad
	\boldsymbol{\phi}&:= \begin{bmatrix} \phi_{i,j} \end{bmatrix}_{i=1,\ldots,n_{1},\,j=1,\ldots,n_{2}},
\end{alignat*}
Now, let $\mathbf{n}=n_{1}n_{2}$, then we can write the numerical scheme in matrix form as
\begin{equation*}
	\begin{dcases}
		\gamma_0 \mathbf{u}^{s} + \sum_{l=1}^{s-1} \gamma_{s-l} e^{-(s-l)\rho\Delta t} \mathbf{u}^{l}
		+  D_{\mathbf{n}}(a) \left(B_{\mathbf{n}}^{x}+B_{\mathbf{n}}^{y}\right) \mathbf{u}^{s}
		- \alpha_m q^{s} \mathbf{f}
		= b_{s-1} e^{-s\rho\Delta t} \boldsymbol{\varphi},
		\\ \quad s=1,\ldots,m, \\
		\mathbf{u}^{m} + {\lambda} D_{\mathbf{n}}(a) \left(B_{\mathbf{n}}^{x}+B_{\mathbf{n}}^{y}\right) \mathbf{f} = \boldsymbol{\phi},
	\end{dcases}
\end{equation*}
where $D_{\mathbf{n}}(a)$ is the diagonal sampling matrix of size $M$ associated to $a(x,y)$, while $B_{M}^{x}$ are $B_{M}^{y}$ are defined as follows
\begin{equation}\label{eq:xB}
	B_{\mathbf{n}}^{x}=I_{n_{2}} \otimes\frac{\alpha_{m}}{\beta_{n_{1}}} B^{x}_{n_{1}},
\end{equation}
and
\begin{equation}\label{eq:yB}
	B_{\mathbf{n}}^{y}=\frac{\alpha_{m}}{\beta_{n_{2}}}B^{y}_{n_{2}} \otimes I_{n_{1}}.
\end{equation}
Here $B_{n_{1}}$ are $B_{n_{2}}$  are the following lower triangular Toeplitz matrices
\begin{equation*}
	B^{x}_{n_{1}} := \begin{bmatrix}
		\delta_0 & 0 & \cdots & \cdots & 0 \\
		\delta_1 & \delta_0 & 0 & & \vdots \\
		\delta_2 & \delta_1 & \delta_0 & \ddots & \vdots \\
		\vdots & \ddots & \ddots & \ddots& 0  \\
		\delta_{n_{1}-1} & \cdots & \delta_2 & \delta_1 & \delta_0
	\end{bmatrix},
\end{equation*}
and
\begin{equation*}
	B^{y}_{n_{2}} := \begin{bmatrix}
		\sigma_0 & 0 & \cdots & \cdots & 0 \\
		\sigma_1 & \delta_0 & 0 & & \vdots \\
		\sigma_2 & \delta_1 & \delta_0 & \ddots & \vdots \\
		\vdots & \ddots & \ddots & \ddots& 0  \\
		\sigma_{n_{2}-1} & \cdots & \delta_2 & \sigma_1 & \sigma_0
	\end{bmatrix}.
\end{equation*}
By gathering the equations for all the time levels and the regularization equation, we obtain the all-at-once linear system
\begin{equation} \label{eq:lin-syst2}
	A_{\widetilde{\mathbf{N}}}\mathbf{v}_\lambda = \mathbf{z}_\lambda,
\end{equation}
where $\widetilde{\mathbf{N}} = \widetilde{\mathbf{N}}(m,\mathbf{n}) := (m+1)\nu(\mathbf{n})$ and
\begin{gather*}
	\renewcommand{\arraystretch}{1.5}
	A_{\widetilde{\mathbf{N}}} := \left[ \begin{array}{cccc|c}
		\gamma_0 I_{\mathbf{n}} +  G_\mathbf{n} & \mathrm{O}_\mathbf{n} & \cdots & \mathrm{O}_\mathbf{n} & -\alpha_m q^{1} I_\mathbf{n} \\
		\gamma_1 e^{-\rho \Delta t} I_\mathbf{n} & \gamma_0 I_\mathbf{n} +  G_\mathbf{n} & \ddots & \vdots & -\alpha_m q^{2} I_\mathbf{n} \\
		\vdots & \ddots & \ddots & \mathrm{O}_\mathbf{n} & \vdots \\
		\gamma_{m-1} e^{-(m-1)\rho\Delta t} I_\mathbf{n} & \cdots & \gamma_1 e^{-\rho \Delta t} I_\mathbf{n} & \gamma_0 I_\mathbf{n} +  G_\mathbf{n} & -\alpha_m q^{m} I_\mathbf{n} \\
		\hline
		\mathrm{O}_\mathbf{n} & \cdots & \mathrm{O}_\mathbf{n} & I_\mathbf{n} & \frac{{\lambda}}{\alpha_{m}} G_\mathbf{n}
	\end{array} \right]
	\in \C^{\widetilde{\mathbf{N}}\times \widetilde{\mathbf{N}}},
	\\
	\renewcommand{\arraystretch}{1.5}
	\mathbf{v}_{\lambda} :=
	\left[ \begin{array}{c}
		\mathbf{u}^{1} \\ \mathbf{u}^{2} \\ \vdots \\ \mathbf{u}^{m} \\ \hline \mathbf{f}
	\end{array} \right]
	\in\C^{\widetilde{\mathbf{N}}},
	\qquad
	\mathbf{z}_{\lambda} :=
	\left[ \begin{array}{c}
		b_0 e^{-\rho\Delta t} \boldsymbol{\varphi} \\
		b_1 e^{-2\rho\Delta t} \boldsymbol{\varphi} \\
		\vdots \\
		b_{m-1} e^{-m\rho\Delta t} \boldsymbol{\varphi} \\
		\hline \boldsymbol{\phi}
	\end{array} \right]
	\in\C^{\widetilde{\mathbf{N}}},
\end{gather*}
in which $I_\mathbf{n}$ is the $\mathbf{n}\times\mathbf{n}$ identity matrix, $\mathrm{O}_\mathbf{n}$ is the $\mathbf{n}\times \mathbf{n}$ zero matrix, and
\begin{equation} \label{eq:Gm-def2}
	G_\mathbf{n} := D_\mathbf{n}(a) \left(B_{\mathbf{n}}^{x}+B_{\mathbf{n}}^{y}\right).
\end{equation}
Moreover, setting
\begin{equation*}
	U_m := \begin{bmatrix}
		\gamma_0 & 0 & \cdots & 0 \\
		\gamma_1 e^{-\rho\Delta t} & \gamma_0 & \ddots & \vdots \\
		\vdots & \ddots & \ddots & 0 \\
		\gamma_{m-1} e^{-\rho(m-1)\Delta t} & \cdots & \gamma_1 e^{-\rho\Delta t} & \gamma_0
	\end{bmatrix},
\end{equation*}
the coefficient matrix can be expressed into the compact form
\begin{equation} \label{eq:coeff-mat2}
	\renewcommand{\arraystretch}{1.5}
	A_{\widetilde{\mathbf{N}}} = \left[ \begin{array}{c|c}
		U_m \otimes I_\mathbf{n} + I_m\otimes G_\mathbf{n} & -\alpha_m \mathbf{q}\otimes I_\mathbf{n} \\
		\hline
		\mathbf{e}_{m}^T \otimes I_\mathbf{n} & \frac{{\lambda}}{\alpha_{m}} G_\mathbf{n}
	\end{array} \right],
\end{equation}
in which
\begin{equation*}
	\mathbf{e}_m := \begin{bmatrix}
		0 \\ \vdots \\ 0 \\ 1
	\end{bmatrix} \in\R^m,
	\qquad
	\mathbf{q} := \begin{bmatrix}
		q^{1} \\ q^{1} \\ \vdots \\ q^{m}
	\end{bmatrix} \in\R^m.
\end{equation*}
\vspace{0.3cm}

\section{GLT analysis} \label{sec:glt-analysis2}
This section presents the spectral analysis of the coefficient matrix sequence via GLT theory to describe its asymptotic behavior. We first report the symbols of $\{B^{x}_{n}\}_{n}$ and the time block $\{U_m\}_{m}$ given in Chapter \ref{chapter:IFL}, then analyze the space block $\{G_{\mathbf{n}}\}_{\mathbf{n}}$, and finally consider the full coefficient matrix sequence.

\begin{theorem} \cite{Ilyas2026a}\label{thm:GLT-Bm2}
	There exist continuous $2\pi$-periodic functions $g_{\eta_{1}}$ and $g_{\eta_{2}}$ such that
	\begin{align*}
		\{B^{x}_{n_{1}}\}_{n_{1}} \sim_{\glt,\sigma}& g_{\eta_{1}}(\theta_{1}),
		\qquad \theta_{1}\in [-\pi,\pi],\\
		\{B^{y}_{n_{2}}\}_{n_{2}} \sim_{\glt,\sigma}& g_{\eta_{2}}(\theta_{2}),
		\qquad \theta_{2}\in [-\pi,\pi].
	\end{align*}
	Moreover,
	\begin{align*}
		\{B^{x}_{n_{1}}\}_{n_{1}} \sim_\lambda& \delta_0,\\
		\{B^{y}_{n_{2}}\}_{n_{2}} \sim_\lambda& \sigma_0.
	\end{align*}
\end{theorem}
\begin{theorem} \cite{Ilyas2026a}\label{thm:GLT-Un2}
	There exists a continuous $2\pi$-periodic function $h_\xi$ such that
	\begin{equation*}
		\{U_m\}_m \sim_{\glt,\sigma} h_\xi(\theta),
		\qquad \theta\in [-\pi,\pi].
	\end{equation*}
	Moreover,
	\begin{equation*}
		\{U_m\}_m \sim_\lambda \gamma_0.
	\end{equation*}
\end{theorem}

\begin{corollary}\label{cor:distr}
	If $a:[0,1]^2\to\R$ is continuous a.e. and if
	\begin{equation*}
		\frac{\alpha_m}{\beta_{n_{1}}} \to \nu_{1},\quad \frac{\alpha_m}{\beta_{n_{2}}} \to \nu_{2}, \qquad \text{as } m,n_{1},n_{2}\to\infty,
	\end{equation*}
	with $\nu_{1}\neq 0$ or $\nu_{2}\neq 0$. Then, we have
	\begin{equation*}
		\{G_\mathbf{n}\}_\mathbf{n}  \sim_{\glt,\sigma} a(x,y) \left(\nu_{1}g_{\eta_{1}}(\theta_{1})+\nu_{2}g_{\eta_{2}}(\theta_{2})\right),
		\; (x,y)\in [0,1]^{2}, \;(\theta_{1},\, \theta_{2})\in [-\pi,\pi]^{2},
	\end{equation*}
	with $g_{\eta_{1}}$ and $g_{\eta_{2}}$ as in Theorem \ref{thm:GLT-Bm2}. Moreover,
	\begin{equation*}
		\{G_\mathbf{n}\}_\mathbf{n}  \sim_\lambda  a(x,y)\left(\delta_0+\sigma_{0}\right),
		\qquad (x,y)\in [0,1]^{2}.
	\end{equation*}
\end{corollary}
\begin{proof}
	From the definition in \eqref{eq:Gm-def2}, we have
	\begin{equation*}
		G_\mathbf{n} := D_\mathbf{n}(a) \left(B_{\mathbf{n}}^{x}+B_{\mathbf{n}}^{y}\right).
	\end{equation*}
	Thus, by GLT \ref{glt:symbols}, we deduce that
	\begin{equation*}
		\{D_\mathbf{n}(a)\}_{\mathbf{n}} \sim_\glt a(x,y), \qquad (x,y)\in [0,1]^{2}.
	\end{equation*}
	Then, using \eqref{eq:xB} and assuming $\frac{\alpha_m}{\beta_{n_{1}}}\to \nu_{1}$, by GLT \ref{glt:symbols}, and GLT \ref{glt:kron}, we obtain
	\begin{align*}
		\{B^{x}_{\mathbf{n}}\}_{\mathbf{n}}=\left\{I_{n_{2}} \otimes\frac{\alpha_{m}}{\beta_{n_{1}}} B^{x}_{n_{1}}\right\}_{\mathbf{n}}=\left\{T_{n_{2}}(1) \otimes\frac{\alpha_{m}}{\beta_{n_{1}}} T_{n_{1}(g_{\eta_{1}})}\right\}_{\mathbf{n}}  \sim_{\glt}& \,\nu_{1}g_{\eta_{1}}(\theta_{1}).
	\end{align*}
	Similarly, using \eqref{eq:yB} and assuming $\frac{\alpha_m}{\beta_{n_{2}}}\to \nu_{2}$, by GLT \ref{glt:symbols}, and GLT \ref{glt:kron}, we obtain
	\begin{align*}
		\{B^{y}_{\mathbf{n}}\}_{\mathbf{n}}=\left\{\frac{\alpha_{m}}{\beta_{n_{2}}} B^{y}_{n_{2}}  \otimes I_{n_{1}}\right\}_{\mathbf{n}}=\left\{\frac{\alpha_{m}}{\beta_{n_{2}}} T_{n_{2}(g_{\eta_{2}})} \otimes T_{n_{1}}(1)\right\}_{\mathbf{n}}  \sim_{\glt}& \,\nu_{2}g_{\eta_{2}}(\theta_{2}).
	\end{align*}
	Invoking the $\ast$-algebra properties detailed in GLT \ref{glt:algebra}, we obtain that
	\begin{equation*}
		\{B^{x}_{\mathbf{n}}+B^{y}_{\mathbf{n}}\}_{\mathbf{n}}   \sim_{\glt} \nu_{1}g_{\eta_{1}}(\theta_{1})+\nu_{2}g_{\eta_{2}}(\theta_{2}),
		\qquad  \;(\theta_{1},\, \theta_{2})\in [-\pi,\pi]^{2}.
	\end{equation*}
	Using the product property of GLT \ref{glt:algebra} applied to the sequences $\{D_\mathbf{n}(a)\}_\mathbf{n}$ and $\{B^{x}_\mathbf{n}+B^{y}_\mathbf{n}\}_\mathbf{n}$, we conclude that
	\begin{equation*}
		\{G_\mathbf{n}\}_\mathbf{n}  \sim_{\glt} a(x,y) \left(\nu_{1}g_{\eta_{1}}(\theta_{1})+\nu_{2}g_{\eta_{2}}(\theta_{2})\right),
		\qquad (x,y)\in [0,1]^{2}, \;(\theta_{1},\, \theta_{2})\in [-\pi,\pi]^{2},
	\end{equation*}
	By GLT \ref{glt:distrib}, we immediately obtain
	\begin{equation*}
		\{G_\mathbf{n}\}_\mathbf{n}  \sim_{\sigma} a(x,y) \left(\nu_{1}g_{\eta_{1}}(\theta_{1})+\nu_{2}g_{\eta_{2}}(\theta_{2})\right),
		\qquad (x,y)\in [0,1]^{2}, \;(\theta_{1},\, \theta_{2})\in [-\pi,\pi]^{2},
	\end{equation*}
	Finally, we observe that $B^{x}_\mathbf{n}$ and $B^{y}_\mathbf{n}$ are lower triangular matrices. Hence $\left(B^{x}_\mathbf{n} + B^{y}_\mathbf{n}\right)$ is lower triangular, and since $D_{\mathbf{n}}(a)$ is diagonal,
	\begin{equation*}
		G_\mathbf{n} := D_\mathbf{n}(a) \left(B_{\mathbf{n}}^{x}+B_{\mathbf{n}}^{y}\right),
	\end{equation*}
	is also lower diagonal. Therefore, the eigenvalues of $G_\mathbf{n}$ coincide with its diagonal entries.
	Since the main diagonal entries of $B^{x}_\mathbf{n}$ and $B^{y}_\mathbf{n}$ are $\delta_{0}$ and $\sigma_{0}$, the diagonal entries of $G_\mathbf{n}$ are $\{a(x_{i},y_{j})\left(\delta_0+\sigma_{0}\right)\}_{i,j=1}^{n_{1},n_{2}}$.
	Hence
	\begin{equation*}
		\{G_\mathbf{n}\}_\mathbf{n}  \sim_\lambda  a(x,y)\left(\delta_0+\sigma_{0}\right),
		\qquad (x,y)\in [0,1]^{2}.
	\end{equation*}
\end{proof}

We are now in a position to compute the GLT symbol of the sequence associated with the full coefficient matrix $A_{\widetilde{\mathbf{N}}}$. To ensure meaningful results, we assume a mild relation among the rates at which the discretization parameters tend to zero.

\begin{theorem} \label{thm:coeff-mat-sym2}
	With $A_{\widetilde{\mathbf{N}}}$ defined in \eqref{eq:coeff-mat2}, assume that
	\begin{equation*}
		\frac{\alpha_m}{\beta_{n_{1}}} \to \nu_{1},\quad \frac{\alpha_m}{\beta_{n_{2}}} \to \nu_{2}, \qquad \text{as } m,n_{1},n_{2}\to\infty,
	\end{equation*}
	with $\nu_{1}\neq 0$ or $\nu_{2}\neq 0$. Then, it holds
	\begin{equation*}
		\{A_{\widetilde{\mathbf{N}}}\}_{\widetilde{\mathbf{N}}} \sim_{\glt,\sigma} h_\xi(\theta_1) + a(x,y) \left(\nu_{1}g_{\eta_{1}}(\theta_{2})+\nu_{2}g_{\eta_{2}}(\theta_{3})\right),  
	\end{equation*}
	$$(\theta_1,\theta_2,\theta_{3})\in [-\pi,\pi]^3,\qquad \;\;(x,y)\in [0,1]^{2},$$
	where $h_\xi$ is given in Theorems \ref{thm:GLT-Bm2} and $g_{\eta_{1}}$ and $g_{\eta_{2}}$  are given in Theorems \ref{thm:GLT-Un2}.
\end{theorem}

\begin{proof}
	We proceed by decomposing \eqref{eq:coeff-mat2} as follows
	\begin{equation*}
		\renewcommand{\arraystretch}{1.5}
		A_{\widetilde{\mathbf{N}}} = \underbrace{ \left[ \begin{array}{cc}
				U_m \otimes I_\mathbf{n} & \mathrm{O}_{m\mathbf{n}\times \mathbf{n}} \\
				\mathrm{O}_{\mathbf{n}\times m\mathbf{n}} & \mathrm{O}_\mathbf{n}
			\end{array} \right] }_{=:\, C_{\widetilde{\mathbf{N}}}^{(1)}}
		+ \underbrace{ \left[ \begin{array}{cc}
				I_m\otimes  G_\mathbf{n} & \mathrm{O}_{m\mathbf{n}\times \mathbf{n}} \\
				\mathrm{O}_{\mathbf{n}\times m\mathbf{n}} & \frac{{\lambda}}{\alpha_m} G_\mathbf{n}
			\end{array} \right] }_{=:\, C_{\widetilde{\mathbf{N}}}^{(2)}}
		+ \underbrace{ \left[ \begin{array}{cc}
				\mathrm{O}_{m\mathbf{n}} & -\alpha_m\mathbf{q}\otimes I_\mathbf{n} \\
				\mathbf{e}_m^T \otimes I_\mathbf{n} & \mathrm{O}_\mathbf{n}
			\end{array} \right] }_{=:\, C_{\widetilde{\mathbf{N}}}^{(3)}}.
	\end{equation*}
	Thus, we analyze each component separately. For the first component, we have
	\begin{equation*}
		\renewcommand{\arraystretch}{1.3}
		C_{\widetilde{\mathbf{N}}}^{(1)}
		= U_{m+1} \otimes I_\mathbf{n}
		- \underbrace{\left[ \begin{array}{c|c}
				\mathrm{O}_{m\mathbf{n}} & \mathrm{O}_{m\mathbf{n}\times \mathbf{n}} \\
				\hline
				\begin{matrix}
					\gamma_m e^{-\rho m\Delta t} I_\mathbf{n} & \cdots & \gamma_1 e^{-\rho\Delta t} I_\mathbf{n}
				\end{matrix} & \gamma_0 I_m
			\end{array} \right]}_{=: R^{(1)}_{\widetilde{\mathbf{N}}}},
	\end{equation*}
	where $\rk R^{(1)}_{\widetilde{\mathbf{N}}} \leq \nu(\mathbf{n}) = o({\widetilde{\mathbf{N}}})$ as $m,\mathbf{n}\to\infty$. According to Theorem \ref{lemma:zero-distr}, the sequence $\{R^{(1)}_{\widetilde{\mathbf{N}}}\}_{\widetilde{\mathbf{N}}}$ is zero-distributed; hence, by GLT \ref{glt:symbols}, we deduce that
	\begin{equation*}
		\{R^{(1)}_{\widetilde{\mathbf{N}}}\}_{\widetilde{\mathbf{N}}}\sim_\glt 0.
	\end{equation*}
	Moreover, exploiting the identity $U_{m+1} \otimes I_\mathbf{n}= T_{m+1}(h_\xi) \otimes T_\mathbf{n}(1)$ and applying the axioms GLT \ref{glt:symbols}--\ref{glt:kron} yields
	\begin{equation} \label{eq:S12}
		\big\{C_{\widetilde{\mathbf{N}}}^{(1)}\big\}_{\widetilde{\mathbf{N}}} \sim_\glt h_\xi(\theta_1).
	\end{equation}
	Regarding $C_{\widetilde{\mathbf{N}}}^{(2)}$, we observe that
	\begin{equation*}
		C_{\widetilde{\mathbf{N}}}^{(2)} = \left[ \begin{array}{cc}
			I_m\otimes  G_\mathbf{n} & \mathrm{O}_{m\mathbf{n}\times\mathbf{n}} \\
			\mathrm{O}_{\mathbf{n}\times m\mathbf{n}} &  G_\mathbf{n} + \left(\frac{\lambda}{\alpha_m}-1\right) G_\mathbf{n}
		\end{array} \right]
		= I_{m+1} \otimes   G_\mathbf{n} +
		\underbrace{\left[ \begin{array}{cc}
				\mathrm{O}_{m\mathbf{n}} & \mathrm{O}_{m\mathbf{n}\times \mathbf{n}} \\
				\mathrm{O}_{\mathbf{n}\times m\mathbf{n}} & \left(\frac{\lambda}{\alpha_m}-1\right) G_\mathbf{n}
			\end{array} \right]}_{=: R_{\widetilde{\mathbf{N}}}^{(2)}},
	\end{equation*}
	where, once again, $\rk R_{\widetilde{\mathbf{N}}}^{(2)} = \nu(\mathbf{n}) = o(\widetilde{\mathbf{N}})$ as $m,\mathbf{n}\to\infty$. By Theorem \ref{lemma:zero-distr}, the sequence $\{R^{(2)}_{\widetilde{\mathbf{N}}}\}_{\widetilde{\mathbf{N}}}$ is zero-distributed; therefore, by GLT \ref{glt:symbols}, it follows that
	\begin{equation*}
		\{R_{\widetilde{\mathbf{N}}}^{(2)}\}_{\widetilde{\mathbf{N}}}\sim_\glt 0.
	\end{equation*}
	Invoking the assumptions $\frac{\alpha_m}{\beta_{n_1}}\to \nu_{1}$ and $\frac{\alpha_m}{\beta_{n_2}}\to \nu_{2}$ together with GLT \ref{glt:symbols} yields
	\begin{align*}
		\left\{I_{m+1} \otimes G_\mathbf{n}\right\}_{\widetilde{\mathbf{N}}} =&\left\{I_{m+1} \otimes  D_{\mathbf{n}}(a)\left(B_{\mathbf{n}}^{x}+B_{\mathbf{n}}^{y}\right)\right\}_{\widetilde{\mathbf{N}}}\\=& \left\{ T_{m+1}(1) \otimes  D_\mathbf{n}(a)\left(\nu_{1}T_\mathbf{n}(g_{\eta_{1}})+\nu_{2}T_\mathbf{n}(g_{\eta_{2}})\right) \right\}_{\widetilde{\mathbf{N}}} .
	\end{align*}
	Then, by GLT \ref{glt:kron}, we obtain
	\begin{equation*}
		\left\{I_{m+1} \otimes G_\mathbf{n}\right\}_{\widetilde{\mathbf{N}}}\sim_\glt a(x,y) \left(\nu_{1}g_{\eta_{1}}(\theta_{2})+\nu_{2}g_{\eta_{2}}(\theta_{3})\right).
	\end{equation*}
	By GLT \ref{glt:algebra}, we conclude that
	\begin{equation} \label{eq:S22}
		\big\{C_{\widetilde{\mathbf{N}}}^{(2)}\big\}_{\widetilde{\mathbf{N}}} \sim_\glt a(x,y) \left(\nu_{1}g_{\eta_{1}}(\theta_{2})+\nu_{2}g_{\eta_{2}}(\theta_{3})\right).
	\end{equation}
	We conclude by observing that
	\begin{equation*}
		C_{\widetilde{\mathbf{N}}}^{(3)} = \left[ \begin{array}{cc}
			\mathrm{O}_{m\mathbf{n}} & -\alpha_m\mathbf{q}\otimes I_\mathbf{n} \\
			\mathrm{O}_{\mathbf{n}\times m\mathbf{n}} & \mathrm{O}_\mathbf{n}
		\end{array} \right]
		+ \left[ \begin{array}{cc}
			\mathrm{O}_{m\mathbf{n}} & \mathrm{O}_{m\mathbf{n}\times \mathbf{n}} \\
			\mathbf{e}_{m}^T \otimes I_\mathbf{n} & \mathrm{O}_\mathbf{n}
		\end{array} \right],
	\end{equation*}
	where each component has rank at most $\mathbf{n}=o({\widetilde{\mathbf{N}}})$. Hence, we conclude that
	\begin{equation} \label{eq:S32}
		\big\{C_{\widetilde{\mathbf{N}}}^{(3)}\big\}_{\widetilde{\mathbf{N}}} \sim_\glt 0.
	\end{equation}
	By combining the contributions in \eqref{eq:S12}, \eqref{eq:S22}, and \eqref{eq:S32}, and invoking GLT \ref{glt:algebra}, we obtain the GLT symbol of the matrix sequence $\{A_{\widetilde{\mathbf{N}}}\}_{\widetilde{\mathbf{N}}}$
	\begin{equation*}
		\{A_{\widetilde{\mathbf{N}}}\}_{\widetilde{\mathbf{N}}} \sim_{\glt} h_\xi(\theta_1) + a(x,y) \left(\nu_{1}g_{\eta_{1}}(\theta_{2})+\nu_{2}g_{\eta_{2}}(\theta_{3})\right), 
	\end{equation*}
	By GLT \ref{glt:distrib}, the singular value distribution is
	\begin{equation*}
		\{A_{\widetilde{\mathbf{N}}}\}_{\widetilde{\mathbf{N}}} \sim_{\sigma} h_\xi(\theta_1) + a(x,y) \left(\nu_{1}g_{\eta_{1}}(\theta_{2})+\nu_{2}g_{\eta_{2}}(\theta_{3})\right), 
	\end{equation*}
\end{proof}
\vspace{0.3cm}

\noindent As in the one-dimensional case, the matrices \(A_{\widetilde{\mathbf{N}}}\) are non-Hermitian, and therefore the GLT symbol obtained in Theorem~\ref{thm:coeff-mat-sym2} directly yields the singular value distribution but does not automatically imply the eigenvalue distribution. The eigenvalue distribution in this setting follows from Theorem~\ref{glt_prec2}, where we prove that the preconditioned matrix sequence is weakly clustered at one. This clustering result is based on the low-rank perturbation structure of the preconditioned system, which is sufficient to characterize the eigenvalue behavior. For non-Hermitian matrix sequences in general, eigenvalue distributions require additional assumptions, such as Hermitian structure or suitable perturbation arguments; see \cite{Barbarino2020, GoSe}. In the present work, the weak clustering result ensures fast convergence of the preconditioned GMRES method despite the non-Hermitian nature of the coefficient matrices.

\begin{remark} \label{remark:dimensional parameters}
	The assumption $\nu_{1}\neq 0$ or $\nu_{2}\neq 0$ used in Corollary \ref{cor:distr} and in Theorem \ref{thm:coeff-mat-sym2} deserves a bit of attention and few comments since it is used also in the study of the preconditioning; see Theorem \ref{thm:precond-sym2}. In fact, if we assume $n_1\sim n_2$ i.e.
	$cn_1\le n_2 \le Cn_1$ with $c, C$ universal constants independent of $n_1,n_2$, the we find the following cases:
	\begin{itemize}
		\item when $\eta_{1}> \eta_2$, taking into account the expression in (\ref{eq:beta-n1}), (\ref{eq:beta-n2}), the given assumption implies
		$\nu_{1}> 0$ and $\nu_{2}= 0$;
		\item when $\eta_{2}> \eta_1$, taking into account the expression in (\ref{eq:beta-n1}), (\ref{eq:beta-n2}), the given assumption implies
		$\nu_{1}= 0$ and $\nu_{2}> 0$;
		\item when $\eta_{2}= \eta_1$, taking into account the expression in (\ref{eq:beta-n1}), (\ref{eq:beta-n2}), the given assumption implies
		$\nu_{1}> 0$ and $\nu_{2}> 0$ with $\nu_{1}=\nu_{2}> 0$ if $n_2=n_1+o(n_1)$;
	\end{itemize}
	Finally, it is worth observing that when $n_1\sim n_2$ is violated, the picture becomes more complicate, but it is also a degenerate case since it implies that one space direction is approximated with much more precision than the other.
\end{remark}

\section{GLT preconditioning} \label{sec:precond2}

The current section is divided into two parts. The first contains the proposal and the second the related spectral analysis.

\subsection{Block Triangular preconditioning}
In this subsection, we propose block triangular preconditioners $P_{\widetilde{\mathbf{N}}}$ for the linear system \eqref{eq:lin-syst2}, aiming to cluster the eigenvalues of $\{ P_{\widetilde{\mathbf{N}}}^{-1} A_{\widetilde{\mathbf{N}}} \}_{\widetilde{\mathbf{N}}} $ around $1$. We define
\begin{equation} \label{eq:precond2}
	\renewcommand{\arraystretch}{1.5}
	P_{\widetilde{\mathbf{N}}} := \left[ \begin{array}{c|c}
		U_m \otimes I_\mathbf{n} + I_m\otimes G_\mathbf{n} & \mathrm{O}_{\mathbf{n}m \times \mathbf{n}} \\
		\hline
		\mathbf{e}_m^T \otimes I_\mathbf{n} & \frac{{\lambda}}{\alpha_{m}} G_\mathbf{n}
	\end{array} \right].
\end{equation}
Note that both $U_m$ and $G_\mathbf{n}$ are invertible, since they are lower triangular matrices with nonzero diagonal entries. The same holds for $U_m \otimes I_\mathbf{n} + I_m \otimes G_\mathbf{n}$, which is also block lower triangular. Moreover, due to the block triangular structure, the eigenvalues of $P_{\widetilde{\mathbf{N}}}$ coincide with those of its diagonal blocks. Therefore, $P_{\widetilde{\mathbf{N}}}$ is invertible.

The linear systems with matrix $P_{\widetilde{\mathbf{N}}}$ that have to be solved when applying the preconditioned GMRES to the original linear systems can be solved by standard block forward elimination. However we have to specify how the linear systems having as coefficient matrices the diagonal blocks of $P_{\widetilde{\mathbf{N}}}$ are solved. Indeed the diagonal blocks have the expression $\gamma_0 I_{\mathbf{n}} + D_{\mathbf{n}}(a) B_{\mathbf{n}}$ and for them we can use several standard preconditioners such as $\gamma_0 I_{\mathbf{n}} + \hat a B_{\mathbf{n}}$ which is of Toeplitz lower triangular form, with $\hat a$ being the integral average of the function $a$, or the standard circulant preconditioner constructed as $\gamma_0 I_{\mathbf{n}} + \hat a S(B_{\mathbf{n}})$, $S(B_{\mathbf{n}})$ being a Strang-like circulant approximation of $B_{\mathbf{n}}$, constructed by preserving the first column of $B_{\mathbf{n}}$. The computational features of these choice are discussed at the end of Section \ref{sec:num2}.
\vspace{0.4cm}

\noindent The preconditioning strategies developed in this chapter are informed by recent advances in preconditioning for all-at-once systems arising from evolutionary PDEs. Huang, Fung, Hon, and Lin \cite{Huang2025} introduced an efficient preconditioner for evolutionary PDEs discretized with the $\theta$-method in time, providing a comprehensive analysis of its spectral properties. Their approach, which leverages the structure of the time-stepping matrix, is complementary to the block triangular preconditioner proposed here for the ISP in a two-dimensional FDE. Additionally, the matching Schur complement preconditioning technique of Lin, Ke, Hon, and Ng \cite{LinKe2026} for discrete time fractional diffusion ISPs demonstrates the effectiveness of exploiting the block structure of all-at-once systems in the fractional setting.

\subsection{GLT preconditioning analysis}
The following theorem establishes that $P_{\widetilde{\mathbf{N}}}$ has the same GLT symbol, and consequently the same singular value distribution, as the original coefficient matrix sequence.

\begin{theorem} \label{thm:precond-sym2}
	With $P_{\widetilde{\mathbf{N}}}$ defined in \eqref{eq:precond2}, assume that
	\begin{equation*}
		\frac{\alpha_m}{\beta_{n_{1}}} \to \nu_{1},\quad \frac{\alpha_m}{\beta_{n_{2}}} \to \nu_{2}, \qquad \text{as } m,n_{1},n_{2}\to\infty,
	\end{equation*}
	with $\nu_{1}\neq 0$ or $\nu_{2}\neq 0$. Then, it holds
	\begin{equation*}
		\{P_{\widetilde{\mathbf{N}}}\}_{\widetilde{\mathbf{N}}} \sim_{\glt,\sigma} h_\xi(\theta_1) + a(x,y) \left(\nu_{1}g_{\eta_{1}}(\theta_{2})+\nu_{2}g_{\eta_{2}}(\theta_{3})\right),  
	\end{equation*}
	$$ (\theta_1,\theta_2,\theta_{3})\in [-\pi,\pi]^3,\qquad (x,y)\in [0,1]^{2},$$
	where $h_\xi$ is given in Theorems \ref{thm:GLT-Bm2} and $g_{\eta_{1}}$ and $g_{\eta_{2}}$  are given in Theorems \ref{thm:GLT-Un2}.
\end{theorem}

\begin{proof}
	The result follows from the same GLT arguments used in the proof of Theorem \ref{thm:coeff-mat-sym2}.
\end{proof}

The GLT symbol and the singular value distribution of the preconditioned matrix sequence follow as a corollary of Theorems \ref{thm:coeff-mat-sym2} and \ref{thm:precond-sym2}, by exploiting the algebra properties in GLT \ref{glt:algebra}. However, they can also be derived in a more general setting, without assuming specific conditions on the discretization parameters, as shown in the following theorem. Moreover, the eigenvalue distribution can be characterized.

\begin{theorem}\label{glt_prec2}
	With $A_{\widetilde{\mathbf{N}}}$ as in \eqref{eq:coeff-mat2} and $P_{\widetilde{\mathbf{N}}}$ as in \eqref{eq:precond2}, it holds
	\begin{equation*}
		\{P_{\widetilde{\mathbf{N}}}^{-1}A_{\widetilde{\mathbf{N}}}\}_{\widetilde{\mathbf{N}}} \sim_{\glt,\sigma} 1.
	\end{equation*}
	Moreover,
	\begin{equation*}
		\{P_{\widetilde{\mathbf{N}}}^{-1}A_{\widetilde{\mathbf{N}}}\}_{\widetilde{\mathbf{N}}} \sim_\lambda 1.
	\end{equation*}
\end{theorem}

\begin{proof}
	Define $R_{\widetilde{\mathbf{N}}} := A_{\widetilde{\mathbf{N}}} - P_{\widetilde{\mathbf{N}}}$. Then,
	\begin{equation*}
		P_{\widetilde{\mathbf{N}}}^{-1} A_{\widetilde{\mathbf{N}}} = I_{\widetilde{\mathbf{N}}} + P_{\widetilde{\mathbf{N}}}^{-1}R_{\widetilde{\mathbf{N}}}.
	\end{equation*}
	Since $\rk R_{\widetilde{\mathbf{N}}}\leq \nu(\mathbf{n})$, it follows that $\rk(P_{\widetilde{\mathbf{N}}}^{-1}R_{\widetilde{\mathbf{N}}}) \leq \nu(\mathbf{n}) = o({\widetilde{\mathbf{N}}})$ as $m,\mathbf{n}\to\infty$. By Theorem \ref{lemma:zero-distr}, the sequence $\{P_{\widetilde{\mathbf{N}}}^{-1}R_{\widetilde{\mathbf{N}}}\}_{\widetilde{\mathbf{N}}}$ is zero-distributed; therefore, by GLT \ref{glt:symbols}, it follows that
	\begin{equation*}
		\{P_{\widetilde{\mathbf{N}}}^{-1}R_{\widetilde{\mathbf{N}}}\}_{\widetilde{\mathbf{N}}} \sim_\glt 0.
	\end{equation*}
	Hence, by exploiting part 2 of axiom GLT \ref{glt:algebra} and by the fact that $\{I_{\widetilde{\mathbf{N}}}\}_{\widetilde{\mathbf{N}}} \sim_\glt 1$, we infer
	\begin{equation*}
		\{P_{\widetilde{\mathbf{N}}}^{-1} A_{\widetilde{\mathbf{N}}}\}_{\widetilde{\mathbf{N}}} \sim_\glt 1.
	\end{equation*}
	Thus, by GLT \ref{glt:distrib}, we immediately obtain
	\begin{equation*}
		\{P_{\widetilde{\mathbf{N}}}^{-1} A_{\widetilde{\mathbf{N}}}\}_{\widetilde{\mathbf{N}}} \sim_\sigma 1.
	\end{equation*}
	Turning to the eigenvalue distribution, we observe that
	\begin{equation*}
		\lambda_n (P_{\widetilde{\mathbf{N}}}^{ -1 } A_{\widetilde{\mathbf{N}}} ) = 1 + \lambda_n (P_{\widetilde{\mathbf{N}}}^{-1} R_{\widetilde{\mathbf{N}}}),
		\qquad n=1,\ldots,{\widetilde{\mathbf{N}}},
	\end{equation*}
	therefore
	\begin{equation*}
		\# \left\{ n=1,\ldots,{\widetilde{\mathbf{N}}}  \;\colon\; \lambda_n (P_{\widetilde{\mathbf{N}}}^{-1} A_{\widetilde{\mathbf{N}}}) \neq 1 \right\} \leq \mathbf{n} = o({\widetilde{\mathbf{N}}}).
	\end{equation*}
	Thus, the eigenvalues of
		\(\{P_{\widetilde{\mathbf{N}}}^{-1}A_{\widetilde{\mathbf{N}}}\}_{\widetilde{\mathbf{N}}}\)
		are weakly clustered at \(1\) in the sense of Definition \ref{clustering}, which, by
		Remark  \ref{remark:distrib-cluster}, is equivalent to saying that the sequence is spectrally	distributed as \(1\).
\end{proof}

\begin{remark} \label{remark:future directions}
	We emphasize that the given sequence of theoretical steps for studying the distribution and clustering in the eigenvalue and singular value sense can be adapted with little effort to a great variety of different situations. Here we mention two of them, one concerning a generalization of problem \eqref{eq:problem2} and one regarding other classes of discretizations for the regularized FDEs in \eqref{eq:quasi-problem2}.
	\begin{enumerate}
		\item If the square $[0,1]\times [0,1]$ is replaced by a bounded $\Omega\subset \mathbb{R}^2$, with $\Omega$ simply Peano-Jordan measurable, then the domain $\Omega_T$ in \eqref{eq:problem2}--\eqref{eq:quasi-problem2} is replaced by $\Omega\times (0,T)$. Hence, simply adapting the steps in Section \ref{sec:glt-analysis2} and Section \ref{sec:precond2}, the same conclusions are obtained, using the reduced GLT theory \cite{Barbarino2022} in place of the classical one; for the reduced GLT theory refer to \cite{SerraCapizzano2003}[pp. 395–399] and \cite{SerraCapizzano2006}[Section 3.1.4] for an example, the initial proposal, the related terminology, and \cite{Barbarino2022} for a complete theoretical development.
		\item If the discretization formulae used in Section \ref{sec:discre2} are replaced by high order finite elements in space of order $p>1$ or by the isogeometric analysis with degree $p$ and regularity $k<p-1$, then we would end up with block multilevel matrix sequences, the dimensionality $r$ of the GLT symbol would by either $p^d$ or $(p-k)^d$ with $d=2$. Again the analysis given here can be performed almost verbatim with minimal technical changes, thanks to the uniform structure of the GLT axioms.
	\end{enumerate}
	Of course the directions indicated in 1) and 2) can be combined and even extended to higher space dimensionality with $\Omega\subset \mathbb{R}^d$, $d\ge 2$, so having the potential of providing a very general picture.
\end{remark}

\section{Numerical results and study of the computational cost}\label{sec:num2}
We evaluate the efficiency of the proposed preconditioner by measuring the number of GMRES iterations required for convergence.
Let us consider the two-dimensional problem on the domain
\(\Omega = (0,1)^2\) and \(t \in (0,T)\). We set
\[
T = 1, \qquad a(x,y) = x + y, \qquad q(t) = t^2, \qquad \varphi(x,y) = 0.
\]
To construct the noisy final data, we select the source function
\[
f(x,y) = x y \sin(\pi x)\sin(\pi y).
\]
and solve the linear system arising from the corresponding direct problem
\begin{equation*}
	\begin{dcases} {}^{T}\mathcal{D}_{t}^{\xi,\rho}u(x,y,t)+ a(x,y)\left({}^{C}\mathcal{D}_{x}^{\eta_{1}}+{}^{C}\mathcal{D}_{y}^{\eta_{2}}\right)u(x,y,t)=q(t)f(x,y), & (x,y,t) \in \Omega_{T},\\	u(0,y,t)=0=u(1,y,t), & y \in (0,1), \quad t\in(0,T),\\
		u(x,0,t)=0=u(x,1,t), & x \in (0,1), \quad t\in(0,T),\\
		u(x,y,0) = \varphi(x,y), & x \in (0,1), \quad y \in (0,1).\end{dcases}\end{equation*}
With the same discretization technique adopted in Section \ref{sec:discre2}, we obtain the linear system
\begin{equation*}
	\tilde{A}_{m,\mathbf{n}}\tilde{\mathbf{u}} = \tilde{\mathbf{b}},
\end{equation*}
where
\begin{gather*}
	\renewcommand{\arraystretch}{1.5}
	\tilde{A}_{m,\mathbf{n}} := \left[ \begin{array}{cccc}
		\gamma_0 I_\mathbf{n} +  G_\mathbf{n} & \mathrm{O}_\mathbf{n} & \cdots & \mathrm{O}_\mathbf{n} \\
		\gamma_1 e^{-\rho \Delta t} I_\mathbf{n} & \gamma_0 I_\mathbf{n} + G_\mathbf{n} & \ddots & \vdots \\
		\vdots & \ddots & \ddots & \mathrm{O}_m \\
		\gamma_{m-1} e^{-(m-1)\rho\Delta t} I_\mathbf{n} & \cdots & \gamma_1 e^{-\rho \Delta t} I_\mathbf{n} & \gamma_0 I_\mathbf{n} +  G_\mathbf{n}
	\end{array} \right]
	\in \C^{m\mathbf{n}\times m\mathbf{n}},
	\\
	\renewcommand{\arraystretch}{1.5}
	\tilde{\mathbf{u}} :=
	\left[ \begin{array}{c}
	\mathbf{u}^{(1)} \\
	\mathbf{u}^{(2)} \\
	\vdots \\
	\mathbf{u}^{(m)}
	\end{array} \right]
	\in\C^{m\mathbf{n}},
	\qquad
	\tilde{\mathbf{b}} :=
	\left[ \begin{array}{c}
		b_0 e^{-\rho\Delta t} \boldsymbol{\varphi} + \alpha_m q^{1} \mathbf{f} \\
		b_1 e^{-2\rho\Delta t} \boldsymbol{\varphi} + \alpha_m q^{2} \mathbf{f} \\
		\vdots \\
		b_{m-1} e^{-m\rho\Delta t} \boldsymbol{\varphi} + \alpha_m q^{m} \mathbf{f}
	\end{array} \right]
	\in\C^{m\mathbf{n}}.
\end{gather*}
The final time datum $\mathbf{u}^{m}$ is obtained by solving the forward problem. To model measurement noise, we introduce a perturbation of the form $\varepsilon \,\delta(x)$, where $\delta(x)$ is a random function uniformly distributed in $(-1,1)$. The resulting noisy datum is denoted by $\boldsymbol{\phi}$, representing a contaminated version of the exact final time observation.

Then, we apply the GMRES method to \eqref{eq:lin-syst2} and report in Table \ref{table:gmres11} the iteration counts and CPU times required by both the nonpreconditioned and preconditioned systems. We set the regularization parameter to $\lambda=5\cdot10^{-3}$, the noise level to $\varepsilon=0.01$, and consider different choices of the fractional orders $\xi$, $\eta_1$, and $\eta_2$. The initial guess is the zero vector, the tolerance is $\texttt{tol}=10^{-8}$, and the maximum number of iterations is equal to the size of the matrix.

Comparing the results obtained for the preconditioned and nonpreconditioned systems, we observe that the preconditioner $P_N$ consistently reduces the number of GMRES iterations for all tested values of the fractional orders. The reduction is particularly significant for $\xi=0.2$, $\eta_1=0.2$, and $\eta_2=0.8$, where the iteration counts of the nonpreconditioned system are the largest. For example, when $n_1=n_2=m=2^5$, the number of iterations decreases from $186$ to $39$, corresponding to a reduction of almost $80\%$. These results confirm the effectiveness of the proposed preconditioning strategy and are in agreement with the theoretical spectral analysis.

Furthermore, the preconditioner exhibits a remarkable robustness with respect to the problem dimensions. Indeed, for each fixed value of $n_1=n_2$, the number of iterations required by the preconditioned system remains nearly constant as $m$ increases. For instance, when $n_1=n_2=2^5$ and $(\xi,\eta_1,\eta_2)=(0.2,0.2,0.8)$, the preconditioned iteration counts vary only between $39$ and $40$, while the corresponding nonpreconditioned iteration counts remain close to $190$. Similar behaviour can be observed for all considered choices of the fractional orders.

In addition, for each fixed value of $m$, the number of iterations required by the preconditioned system grows moderately as the spatial discretization parameters $n_1$ and $n_2$ increase. In contrast, the nonpreconditioned iteration counts increase much more rapidly. This indicates that the preconditioner successfully mitigates the deterioration of the spectral properties associated with mesh refinement, leading to a substantially improved convergence behaviour.

The CPU timings generally reflect the reduction in iteration counts achieved by the preconditioner. Although the gain in computational time is less pronounced in some cases due to the additional cost associated with applying the preconditioner, the preconditioned framework remains competitive and considerably more scalable than the nonpreconditioned one as the problem dimensions increase.

Overall, the numerical experiments demonstrate that the proposed preconditioner provides a substantial acceleration of the GMRES method and maintains a robust performance across different discretization levels and fractional-order parameters.

\begin{table}[ht]
	\centering
	\caption{Iteration counts and CPU times required by the GMRES method for solving the linear system, with and without preconditioning, for different values of the fractional orders.}
	\label{table:gmres11}
		\begin{adjustbox}{width=\textwidth}
	\begin{tabular}{cccccccccccccc}
		\toprule
		& & \multicolumn{4}{c}{$\xi=0.2$} & \multicolumn{4}{c}{$\xi=0.5$} & \multicolumn{4}{c}{$\xi=0.8$}
		\\
		& & \multicolumn{4}{c}{$\eta_{1}=0.2,\;\;\eta_{2}=0.8$} & \multicolumn{4}{c}{$\eta_{1}=0.5,\;\;\eta_{2}=0.5$} & \multicolumn{4}{c}{$\eta_{1}=0.8,\;\;\eta_{2}=0.2$}
		\\
		\cmidrule(lr){3-6}\cmidrule(lr){7-10}\cmidrule(lr){11-14}
		& & \multicolumn{2}{c}{Iter} & \multicolumn{2}{c}{CPU} & \multicolumn{2}{c}{Iter} & \multicolumn{2}{c}{CPU} & \multicolumn{2}{c}{Iter} & \multicolumn{2}{c}{CPU}
		\\
		\cmidrule(lr){3-4}\cmidrule(lr){5-6}\cmidrule(lr){7-8}\cmidrule(lr){9-10}\cmidrule(lr){11-12}\cmidrule(lr){13-14}
		$n_{1}$ = $n_{2}$ & $m$ & - & $P_{\widetilde{\mathbf{N}}}$ & - & $P_{\widetilde{\mathbf{N}}}$ & - & $P_{\widetilde{\mathbf{N}}}$ & - & $P_{\widetilde{\mathbf{N}}}$ & - & $P_{\widetilde{\mathbf{N}}}$ & - & $P_{\widetilde{\mathbf{N}}}$
		\\
		\midrule
		$2^{2}$ & $2^{2}$ &32  & 8 & 0.02 &0.01  & 29 & 8 & 0.01 & 0.01 & 30 & 8 & 0.00 &0.00  \\
		$2^{2}$ & $2^{3}$ &34  & 8 &0.01  &0.00  & 33 & 8 & 0.01 & 0.01 &35  &8  & 0.00 &0.00  \\
		$2^{2}$ & $2^{4}$ & 34 & 8 & 0.01 & 0.00 &38  & 8 & 0.01 & 0.00 &47  & 8 &0.01  &0.01  \\
		$2^{2}$ & $2^{5}$ & 34 & 8 & 0.01 & 0.00 &42  & 8 & 0.01 & 0.01 &60 &8  & 0.01 &0.01  \\
		\midrule
		$2^{3}$ & $2^{2}$ & 59 & 17 & 0.01 &0.00  & 48 & 15 & 0.01 & 0.00 & 49 &16  &0.00  & 0.00 \\
		$2^{3}$ & $2^{3}$ & 57 & 16 & 0.01 & 0.01 & 50 & 16 & 0.01 & 0.01 &55  & 16 & 0.01 &0.00  \\
		$2^{3}$ & $2^{4}$ &  61& 16 & 0.01 & 0.01 & 54 & 16 & 0.01 & 0.01 &66 &15  &0.02  &0.01  \\
		$2^{3}$ & $2^{5}$ & 58 & 17 & 0.03 & 0.01 & 59 & 16 & 0.02 & 0.02 &  82& 16 &0.03  & 0.01 \\
		\midrule
		$2^{4}$ & $2^{2}$ & 98 & 27 & 0.04 & 0.02 &75  & 28 & 0.04 & 0.02 & 78 &27  &0.03  &0.01  \\
		$2^{4}$ & $2^{3}$ & 99 & 27 & 0.06 & 0.03 & 74 & 28 & 0.04 & 0.02 & 83 &27  &  0.04& 0.03 \\
		$2^{4}$ & $2^{4}$ & 97 & 28 & 0.12 & 0.04 & 79 & 28 & 0.09 & 0.04 &  92& 27 &0.11  &0.04  \\
		$2^{4}$ & $2^{5}$ & 98 & 27 & 0.22 & 0.08 & 84 & 28 & 0.16 & 0.07 &  113& 27 & 0.27 &0.08  \\
		\midrule
		$2^{5}$ & $2^{2}$ & 189 & 39 & 1.71 & 0.69 & 115 & 43 & 0.96 & 0.83 & 146 &38  & 1.18 &0.63  \\
		$2^{5}$ & $2^{3}$ & 189 & 40 & 2.00 & 0.82 & 114 & 42 & 1.02 & 0.83 &  138& 38 &1.58  & 1.63 \\
		$2^{5}$ & $2^{4}$ & 185 & 39 & 2.87 & 1.16 &115  & 41 & 1.40 & 1.13 & 144 &37  &2.04  &2.11  \\
		$2^{5}$ & $2^{5}$ &186  & 39 & 3.14 &1.10  & 116 & 40 & 2.13 & 1.15 & 149 &36  &3.10  &3.09  \\
		\bottomrule
	\end{tabular}
	\end{adjustbox}
\end{table}

\subsection*{Study of the computational costs and further proposals}

In this subsection, we study in detail the computational cost of our preconditioned GMRES method and we discuss further potential proposals regarding the preconditioning.

When applying one step of the preconditioned GMRES with preconditioner \( P_{\widetilde{\mathbf{N}}} \), we have to solve a linear system with coefficient matrix \( P_{\widetilde{\mathbf{N}}} \) and we have to perform matrix-vector products with the original matrix \( A_{\widetilde{\mathbf{N}}} \). The related arithmetic costs are the following
\begin{itemize}
	\item When performing one matrix-vector multiplication with matrix \( A_{\widetilde{\mathbf{N}}} \), we perform $\frac{m^2}{2} +O(m)$ products of a constant times a vector of size
	$\nu(\mathbf{n})$, $\frac{m^2}{2} +O(m)$ additions between vectors of size $\nu(\mathbf{n})$, and $m+1$ products of the matrix $D_{\mathbf{n}}(a)B_{\mathbf{n}}$ times a vector: hence the corresponding arithmetic cost is equal to $m^2 \nu(\mathbf{n}) + C m \nu(\mathbf{n})\log [\nu(\mathbf{n})] + o(m^2 \nu(\mathbf{n}))=m^2 \nu(\mathbf{n}) + o(m^2 \nu(\mathbf{n}))$
	arithmetic operations, if we make the natural assumption that $\log [\nu(\mathbf{n})] =o(m)$ and taking into account that $C$ is a constant;
	\item When solving one linear system with coefficient matrix \( P_{\widetilde{\mathbf{N}}} \), we perform $\frac{m^2}{2} +O(m)$ products of a constant times a vector of size
	$\nu(\mathbf{n})$, $\frac{m^2}{2} +O(m)$ additions between vectors of size $\nu(\mathbf{n})$, and $m+1$ solutions of linear systems with coefficient matrix $M_{\nu(\mathbf{n})}:=c I + D_{\mathbf{n}}(a)B_{\mathbf{n}}$ times a vector: hence the corresponding arithmetic cost is equal to $m^2 \nu(\mathbf{n}) + C' m \nu(\mathbf{n})\log [\nu(\mathbf{n})] + o(m^2 \nu(\mathbf{n}))=m^2 \nu(\mathbf{n}) + o(m^2 \nu(\mathbf{n}))$ arithmetic operations, if we make the natural assumption that $\log [\nu(\mathbf{n})] =o(m)$ and $C'$ depends linearly on the number of iterations we need preconditioning with $\widetilde P_{\nu(\mathbf{n})}:=c I + \hat a B_{\mathbf{n}}$ or with $\widetilde S_{\nu(\mathbf{n})}:=c I + \hat a S(B_{\mathbf{n}})$. As it can be seen from Table \ref{table:gmres2}, Table \ref{table:gmres3}, Table \ref{table:gmres4}, Table \ref{table:gmres5}, the number of iterations for both preconditioned GMRES is bounded by a constant independent of the matrix size $\nu(\mathbf{n})$ so that also $C'$ can be assumed to be constant.
\end{itemize}

We recall that $c I_{\mathbf{n}} + \hat a B_{\mathbf{n}}$ which is of Toeplitz lower triangular form and we have specialized techniques for lower triangular Toeplitz matrices of cost
$O(\nu(\mathbf{n})\log [\nu(\mathbf{n})])$ and this remains true for the circulant or $\omega$-circulant counterpart; see \cite{Bini,Diaz} and references therein. \ \\

\begin{table}[ht]
	\centering
	\caption{Iteration counts and CPU times of the GMRES method for solving the linear system with coefficient matrix \( \widetilde{M}_{\nu(\mathbf{n})} \) and $c=\gamma_0$, without preconditioner and with preconditioner \( \widetilde{P}_{\nu(\mathbf{n})} \), for different fractional orders and for \( a(x,y)=e^{x+y} \), with \( m=3 \).}
	\label{table:gmres2}
		\begin{adjustbox}{width=\textwidth}
	\begin{tabular}{ccccccccccccc}
		\toprule
		&
		\multicolumn{4}{c}{$\xi=0.2$} &
		\multicolumn{4}{c}{$\xi=0.5$} &
		\multicolumn{4}{c}{$\xi=0.8$}
		\\
		&
		\multicolumn{4}{c}{$\eta_{1}=0.2,\;\eta_{2}=0.8$} &
		\multicolumn{4}{c}{$\eta_{1}=0.5,\;\eta_{2}=0.5$} &
		\multicolumn{4}{c}{$\eta_{1}=0.8,\;\eta_{2}=0.2$}
		\\
		\cmidrule(lr){2-5}\cmidrule(lr){6-9}\cmidrule(lr){10-13}
		
		&
		\multicolumn{2}{c}{Iter} &
		\multicolumn{2}{c}{CPU} &
		\multicolumn{2}{c}{Iter} &
		\multicolumn{2}{c}{CPU} &
		\multicolumn{2}{c}{Iter} &
		\multicolumn{2}{c}{CPU}
		\\
		\cmidrule(lr){2-3}\cmidrule(lr){4-5}
		\cmidrule(lr){6-7}\cmidrule(lr){8-9}
		\cmidrule(lr){10-11}\cmidrule(lr){12-13}
		
		$n_{1}=n_{2}$ &
		- & $\widetilde P_{\nu(\mathbf{n})}$ &
		- & $\widetilde P_{\nu(\mathbf{n})}$ &
		- & $\widetilde P_{\nu(\mathbf{n})}$ &
		- & $\widetilde P_{\nu(\mathbf{n})}$ &
		- & $\widetilde P_{\nu(\mathbf{n})}$ &
		- & $\widetilde P_{\nu(\mathbf{n})}$
		\\
		\midrule
		
		$2^{4}$ & 31 & 22 & 0.00 & 0.00 & 27 & 22 & 0.02 & 0.01 & 31 & 22 & 0.01 & 0.01 \\
		$2^{5}$ & 61 & 25 & 0.06 & 0.07 & 38 & 25 & 0.04 & 0.09 & 59 & 25 & 0.05 & 0.08 \\
		$2^{6}$ & 116 & 26 & 1.67 & 1.28 & 56 & 26 & 0.56 & 1.17 & 113 & 26 & 1.27 & 1.17 \\
		$2^{7}$ & 214 & 27 & 46.36 & 21.34 & 83 & 27 & 20.54 & 21.03 & 207 & 27 & 35.68 & 18.57 \\
		\bottomrule
	\end{tabular}
	\end{adjustbox}
\end{table}
\begin{table}[ht]
	\centering
	\caption{Iteration counts and CPU times of the GMRES method for solving the linear system with coefficient matrix \( \widetilde{M}_{\nu(\mathbf{n})} \) and $c=\gamma_0$, without preconditioner and with preconditioner \( \widetilde{P}_{\nu(\mathbf{n})} \), for different fractional orders and for \( a(x,y)=1+{x+y} \), with \( m=3 \).}
	\label{table:gmres3}
		\begin{adjustbox}{width=\textwidth}
	\begin{tabular}{ccccccccccccc}
		\toprule
		&
		\multicolumn{4}{c}{$\xi=0.2$} &
		\multicolumn{4}{c}{$\xi=0.5$} &
		\multicolumn{4}{c}{$\xi=0.8$}
		\\
		&
		\multicolumn{4}{c}{$\eta_{1}=0.2,\;\eta_{2}=0.8$} &
		\multicolumn{4}{c}{$\eta_{1}=0.5,\;\eta_{2}=0.5$} &
		\multicolumn{4}{c}{$\eta_{1}=0.8,\;\eta_{2}=0.2$}
		\\
		\cmidrule(lr){2-5}\cmidrule(lr){6-9}\cmidrule(lr){10-13}
		
		&
		\multicolumn{2}{c}{Iter} &
		\multicolumn{2}{c}{CPU} &
		\multicolumn{2}{c}{Iter} &
		\multicolumn{2}{c}{CPU} &
		\multicolumn{2}{c}{Iter} &
		\multicolumn{2}{c}{CPU}
		\\
		\cmidrule(lr){2-3}\cmidrule(lr){4-5}
		\cmidrule(lr){6-7}\cmidrule(lr){8-9}
		\cmidrule(lr){10-11}\cmidrule(lr){12-13}
		
		$n_{1}=n_{2}$ &
		- & $\widetilde P_{\nu(\mathbf{n})}$ &
		- & $\widetilde P_{\nu(\mathbf{n})}$ &
		- & $\widetilde P_{\nu(\mathbf{n})}$ &
		- & $\widetilde P_{\nu(\mathbf{n})}$ &
		- & $\widetilde P_{\nu(\mathbf{n})}$ &
		- & $\widetilde P_{\nu(\mathbf{n})}$
		\\
		\midrule
		
		$2^{4}$ & 29 & 15 & 0.00 & 0.00 & 24 & 14 & 0.00 & 0.00 & 28 & 14 & 0.00 & 0.00 \\
		$2^{5}$ & 50 & 15 & 0.02 & 0.02 & 32 & 15 & 0.01 & 0.03 & 49 & 15 & 0.02 & 0.02 \\
		$2^{6}$ & 88 & 16 & 0.59 & 0.41 & 44 & 16 & 0.30 & 0.41 & 86 & 15 & 0.59 & 0.37 \\
		$2^{7}$ & 156 & 19 & 14.98 & 7.02 & 61 & 16 & 6.87 & 7.18 & 153 & 16 & 18.62 & 7.41 \\
		\bottomrule
	\end{tabular}
	\end{adjustbox}
\end{table}

\begin{table}[ht]
	\centering
	\caption{Iteration counts and CPU times of the GMRES method for solving the linear system with coefficient matrix \( \widetilde{M}_{\nu(\mathbf{n})} \) and $c=\gamma_0$, without preconditioner and with preconditioner \( \widetilde{S}_{\nu(\mathbf{n})} \), for different fractional orders and for \( a(x,y)=e^{x+y} \), with \( m=3 \).}
	\label{table:gmres4}
		\begin{adjustbox}{width=\textwidth}
	\begin{tabular}{ccccccccccccc}
		\toprule
		&
		\multicolumn{4}{c}{$\xi=0.2$} &
		\multicolumn{4}{c}{$\xi=0.5$} &
		\multicolumn{4}{c}{$\xi=0.8$}
		\\
		&
		\multicolumn{4}{c}{$\eta_{1}=0.2,\;\eta_{2}=0.8$} &
		\multicolumn{4}{c}{$\eta_{1}=0.5,\;\eta_{2}=0.5$} &
		\multicolumn{4}{c}{$\eta_{1}=0.8,\;\eta_{2}=0.2$}
		\\
		\cmidrule(lr){2-5}\cmidrule(lr){6-9}\cmidrule(lr){10-13}
		
		&
		\multicolumn{2}{c}{Iter} &
		\multicolumn{2}{c}{CPU} &
		\multicolumn{2}{c}{Iter} &
		\multicolumn{2}{c}{CPU} &
		\multicolumn{2}{c}{Iter} &
		\multicolumn{2}{c}{CPU}
		\\
		\cmidrule(lr){2-3}\cmidrule(lr){4-5}
		\cmidrule(lr){6-7}\cmidrule(lr){8-9}
		\cmidrule(lr){10-11}\cmidrule(lr){12-13}
		
		$n_{1}=n_{2}$ &
		- & $\widetilde S_{\nu(\mathbf{n})}$ &
		- & $\widetilde S_{\nu(\mathbf{n})}$ &
		- & $\widetilde S_{\nu(\mathbf{n})}$ &
		- & $\widetilde S_{\nu(\mathbf{n})}$ &
		- & $\widetilde S_{\nu(\mathbf{n})}$ &
		- & $\widetilde S_{\nu(\mathbf{n})}$
		\\
		\midrule
		
		$2^{4}$ & 31 & 25 & 0.00 & 0.01 & 27 & 21 & 0.02 & 0.02 & 31 & 24 & 0.01 & 0.01 \\
		$2^{5}$ & 61 & 27 & 0.06 & 0.03 & 38 & 24 & 0.04 & 0.03 & 59 & 26 & 0.05 & 0.03 \\
		$2^{6}$ & 116 & 29 & 1.67 & 0.42 & 56 & 25 & 0.56 & 0.37 & 113 & 27 & 1.27 & 0.32 \\
		$2^{7}$ & 214 & 29 & 46.36 & 6.64 & 83 & 25 & 20.56 & 4.78 & 207 & 27 & 35.68 & 4.69 \\
		\bottomrule
	\end{tabular}
	\end{adjustbox}
\end{table}

\begin{table}[ht]
	\centering
	\caption{Iteration counts and CPU times of the GMRES method for solving the linear system with coefficient matrix \( \widetilde{M}_{\nu(\mathbf{n})} \) and $c=\gamma_0$, without preconditioner and with preconditioner \( \widetilde{S}_{\nu(\mathbf{n})} \), for different fractional orders and for \( a(x,y)=1+{x+y} \), with \( m=3 \).}
	\label{table:gmres5}
		\begin{adjustbox}{width=\textwidth}
	\begin{tabular}{ccccccccccccc}
		\toprule
		&
		\multicolumn{4}{c}{$\xi=0.2$} &
		\multicolumn{4}{c}{$\xi=0.5$} &
		\multicolumn{4}{c}{$\xi=0.8$}
		\\
		&
		\multicolumn{4}{c}{$\eta_{1}=0.2,\;\eta_{2}=0.8$} &
		\multicolumn{4}{c}{$\eta_{1}=0.5,\;\eta_{2}=0.5$} &
		\multicolumn{4}{c}{$\eta_{1}=0.8,\;\eta_{2}=0.2$}
		\\
		\cmidrule(lr){2-5}\cmidrule(lr){6-9}\cmidrule(lr){10-13}
		
		&
		\multicolumn{2}{c}{Iter} &
		\multicolumn{2}{c}{CPU} &
		\multicolumn{2}{c}{Iter} &
		\multicolumn{2}{c}{CPU} &
		\multicolumn{2}{c}{Iter} &
		\multicolumn{2}{c}{CPU}
		\\
		\cmidrule(lr){2-3}\cmidrule(lr){4-5}
		\cmidrule(lr){6-7}\cmidrule(lr){8-9}
		\cmidrule(lr){10-11}\cmidrule(lr){12-13}
		
		$n_{1}=n_{2}$ &
		- & $\widetilde S_{\nu(\mathbf{n})}$ &
		- & $\widetilde S_{\nu(\mathbf{n})}$ &
		- & $\widetilde S_{\nu(\mathbf{n})}$ &
		- & $\widetilde S_{\nu(\mathbf{n})}$ &
		- & $\widetilde S_{\nu(\mathbf{n})}$ &
		- & $\widetilde S_{\nu(\mathbf{n})}$
		\\
		\midrule
		
		$2^{4}$ & 29 & 17 & 0.00 & 0.00 & 24 & 13 & 0.00 & 0.00 & 28 & 16 & 0.00 & 0.00 \\
		$2^{5}$ & 50 & 18 & 0.02 & 0.01 & 32 & 14 & 0.01 & 0.01 & 49 & 17 & 0.02 & 0.01 \\
		$2^{6}$ & 88 & 18 & 0.59 & 0.14 & 44 & 14 & 0.30 & 0.10 & 86 & 17 & 0.59 & 0.12 \\
		$2^{7}$ & 156 & 18 & 14.98 & 1.73 & 61 & 15 & 6.87 & 1.68 & 153 & 17 & 18.62 & 2.49 \\
		\bottomrule
	\end{tabular}
	\end{adjustbox}
\end{table}
\FloatBarrier

Finally, we introduce an alternative Strang-like preconditioner \( S_{\widetilde{\mathbf{N}}} \) for the linear system \eqref{eq:lin-syst2}.
\smallskip

\noindent A Strang-like preconditioner for our block structure is defined by replacing each lower triangular Toeplitz matrix in the various blocks with its natural circulant approximation and by replacing coefficient $a(x,y)$ with its integral average $\hat a$ so that $D_{\mathbf{n}}(a)$ is replaced by
$\hat aI_{\mathbf{n}}$ and $G_{\mathbf{n}}$ reduces to $$G_{\mathbf{n}} = \hat a B_{\mathbf{n}}, \qquad B_{\mathbf{n}}=\left(B_{\mathbf{n}}^{x}+B_{\mathbf{n}}^{y}\right),$$
that is
\begin{equation}\label{eq:strang-precond}
	\renewcommand{\arraystretch}{1.5}
	S_{\widetilde{\mathbf{N}}} := \left[ \begin{array}{c|c}
		S(U_m) \otimes I_{\mathbf{n}} + I_m\otimes \hat{a} S(B_{\mathbf{n}}) & \mathrm{O}_{\mathbf{n}m \times \mathbf{n}} \\
		\hline
		\mathbf{e}_m^T \otimes I_\mathbf{n} & \frac{\lambda}{\alpha_{m}} D_{\mathbf{n}}(a) S(B_{\mathbf{n}})
	\end{array} \right].
\end{equation}
Since circulant matrices are diagonalizable by the Fourier matrices, we write
\begin{equation*}
	S(U_m) = F_m \Lambda_m^{(1)} F_m^{*},
	\qquad
	S(B_{\mathbf{n}}) = F_{\mathbf{n}} \Lambda_{\mathbf{n}}^{(2)} F_{\mathbf{n}}^{*},
\end{equation*}
where $\Lambda_m^{(1)}$ and $\Lambda_\mathbf{n}^{(2)}$ are the diagonal matrices containing the eigenvalues and can be computed starting from the coefficients of $U_m$, $B_{\mathbf{n}}$ via the FFT.
\smallskip

Define the unitary transformation

\begin{equation*}
	\renewcommand{\arraystretch}{1.5}
	Q := \left[ \begin{array}{c|c}
		F_m\otimes F_{\mathbf{n}} & \mathrm{O}_{{\mathbf{n}}m \times {\mathbf{n}}} \\
		\hline
		\mathrm{O}_{{\mathbf{n}}\times{\mathbf{n}}m} &  F_{\mathbf{n}}
	\end{array} \right],
\end{equation*}
Then
\begin{equation}\label{eq:p11p}
	\renewcommand{\arraystretch}{1.5}
	Q^{*} S_{\widetilde{\mathbf{N}}} Q := \left[ \begin{array}{c|c}
		\Lambda_m^{(1)}\otimes I_{\mathbf{n}} + I_m\otimes  \hat{a} \Lambda_\mathbf{n}^{(2)}  & \mathrm{O}_{{\mathbf{n}}m \times {\mathbf{n}}} \\
		\hline
		\mathbf{e}_m^T F_m^{*} \otimes I_\mathbf{n} & \frac{\lambda}{\alpha_m} \hat{a} \Lambda_{\mathbf{n}}^{(2)}
	\end{array} \right].
\end{equation}

The preconditioner $S_{\widetilde{\mathbf{N}}}$ is cheaper than $P_{\widetilde{\mathbf{N}}}$, but it cannot ensure the clustering at $1$ due to GLT arguments since
\[
\{S_{\widetilde{\mathbf{N}}}^{-1}A_{\widetilde{\mathbf{N}}}\}_{\widetilde{\mathbf{N}}} \sim_{\glt,\sigma} \frac{h_\xi(\theta_1) + a(x,y) \left(\nu_{1}g_{\eta_{1}}(\theta_{2})+\nu_{2}g_{\eta_{2}}(\theta_{3})\right)}{h_\xi(\theta_1) + \hat a \left(\nu_{1}g_{\eta_{1}}(\theta_{2})+\nu_{2}g_{\eta_{2}}(\theta_{3})\right)}
\]
and hence it is less effective that $P_{\widetilde{\mathbf{N}}}$ in terms of iteration count, also because the asymptotic behavior of the eigenvalues is not favourable.

	\chapter{Structure-preserving preconditioning for variable-coefficient space-fractional diffusion equations and $\theta$-Method}\label{chapter:space-FDE}

In this chapter, based on the results presented in \cite{Khan2026}, we investigate the spectral properties and preconditioning of linear systems arising from the discretization of time-dependent space-FDEs with variable coefficient. Additional details, numerical experiments, and reformulations are included for completeness of the thesis presentation. Fractional diffusion models have become increasingly relevant in the mathematical modeling of anomalous transport processes, especially in situations where the classical diffusion equation fails to capture long-range interactions and nonlocal effects \cite{Tadjeran2006}. Such models are particularly valuable for describing superdiffusive behavior, where particles disperse more rapidly than predicted by classical diffusion theory.

The numerical discretization of these equations leads to large linear systems that are dense, nonsymmetric, and often severely ill-conditioned. To address these challenges, we develop a preconditioning strategy based on the structural properties of the discretized operator. A detailed spectral analysis of the associated matrix sequences is carried out using the GLT theory, which provides a unified framework for describing the asymptotic singular value and spectral distributions of both the original and preconditioned matrices \cite{garoni2017, garoni2018}. The main novelty is that the analysis fully covers the nonconstant variable coefficient case, both in one and two spatial dimensions. Numerical experiments confirm the effectiveness of the proposed preconditioners in accelerating the convergence of GMRES.

In this chapter, we consider the following class of evolutionary partial differential equations:
\begin{equation}\label{eq1}
	\begin{cases}
		\displaystyle \frac{\partial}{\partial t}u(x,y,t) = \mathcal{L}u(x,y,t) + f(x,y,t), & (x,y,t) \in \Omega \times (0,T], \\
		u(x,y,t) = 0, & (x,y) \in \partial\Omega,\; t \in [0,T], \\
		u(x,y,0) = \phi(x,y), & (x,y) \in \Omega,
	\end{cases}
\end{equation}
where \( \Omega = (L_1,R_1) \times (L_2,R_2) \subset \mathbb{R}^{2} \) is open, \( f(x,y,t) \) is a given source term, \( u(x,y,t) \) denotes the unknown solution, and \( \phi(x,y) \) is the prescribed initial condition.

More specifically, we consider the operator
\[
\mathcal{L} = a(x,y)\left(\mathcal{D}_{x}^{\alpha_1} + \mathcal{D}_{y}^{\alpha_2}\right),
\] 
where \( a(x,y) \) is a smooth variable diffusion coefficient, \( \mathcal{D}_{x}^{\alpha_1} \) denotes the left-sided Riemann--Liouville fractional derivative of order \( \alpha_1 \in (1,2) \) with respect to \(x\), and \( \mathcal{D}_{y}^{\alpha_2} \) denotes the left-sided Riemann--Liouville fractional derivative of order \( \alpha_2 \in (1,2) \) with respect to \(y\).

In the one-dimensional case, the model reduces to
\begin{equation}\label{eq2-1d}
	\begin{cases}
		\displaystyle \frac{\partial}{\partial t} u(x,t) = a(x)\, \mathcal{D}_{x}^{\alpha}u(x,t) + f(x,t), & (x,t) \in (L_1,R_1) \times (0,T], \\
		u(x,t) = 0, & x \in \partial(L_1,R_1),\; t \in [0,T], \\
		u(x,0) = \phi(x), & x \in (L_1,R_1).
	\end{cases}
\end{equation}

The chapter is organized as follows. Section~\ref{sec:analysis} contains the structural and spectral analysis of the relevant matrix sequences. In Sections~\ref{sec:prec} and \ref{sec:prec2D}, we introduce suitable preconditioners and study their spectral behavior. Finally, in Section~\ref{sec:num}, we present spectral visualizations and numerical experiments, and discuss the numerical results.

\section{Discretized problem and GLT analysis}\label{sec:analysis}
In this section, we present the discretization of the variable-coefficient space-fractional diffusion equation \eqref{eq1}, posed on the spatial domain $\Omega = {(L_1, R_1) \times (L_2, R_2) \subset \mathbb{R}^2}$. The equation is first reformulated on a discrete space-time grid, leading to a sequence of linear systems governed by structured coefficient matrices that encode the nonlocal behavior of the model. We then analyze the spectral properties of these matrices using the GLT framework. In particular, we compute the GLT symbol of the resulting matrix sequence, which characterizes the asymptotic distribution of eigenvalues and singular values.
\subsection{Discretization}

We first treat the one-dimensional space setting and then the two-dimensional space setting.

\subsubsection{One-dimensional}
Let $M$ and $N$ be positive integers representing the number of temporal steps and the number of interior spatial grid points, respectively. The spatial  and temporal meshes are defined as 
\begin{alignat*}{3}
	x_i &= L_1 + i \Delta x, && \qquad\Delta x =  \frac{(R_1 - L_1)}{(N + 1)}, && \qquad 0 \leq i \leq N+1, \\
	t_m &= m\Delta t, && \qquad\Delta t=\frac{T}{M}, && \qquad 0 \leq m \leq M.
\end{alignat*}
\\

\noindent To approximate the left-sided Riemann--Liouville fractional derivative, we adopt the shifted Grünwald finite difference scheme \cite{MeerschaertTadjeran2006}, which is widely recognized for its improved stability and accuracy compared to the standard Grünwald approach. The continuous fractional operator \( \mathcal{D}_x^\alpha u(x,t) \) is discretized at each interior point \( x_i \) using the following expression
{\begin{equation}\label{equ4}
		\mathcal{D}_x^\alpha u(x, t)\big|_{x = x_i, \, t = t_{m}}= \frac{1}{(\Delta x)^{\alpha}} \sum_{k=0}^{i+1} g_{i - k + 1}^{(\alpha)} u(x_{k},t_{m})+ \mathcal{O}(\Delta x),
\end{equation}}
where
\begin{equation}\label{equ5}
	{g_0^{(\alpha)} = 1, \quad g_{k+1}^{(\alpha)} = \left(1 - \frac{\alpha + 1}{k+1} \right) g_{k}^{(\alpha)}, \quad k \ge 0.}
\end{equation}
This recurrence can be rewritten in closed form as
\begin{equation}
	\label{eq:g-binomial}
	{g^{(\alpha)}_k = (-1)^{k} \binom{\alpha}{k}, \quad k \ge 0,}
\end{equation}
where
\[
\binom{\alpha}{k} = \frac{\alpha (\alpha-1) \cdots (\alpha-k+1)}{k!}
\]
denotes the generalized binomial coefficient.
\begin{proposition}\label{alpha coff}
	Let $\alpha \in (1,2)$ and $g_k^{(\alpha)}$ be defined as in \eqref{equ5}. Then we have
	{\[
		\begin{cases}
			g_0^{(\alpha)} = 1, \quad g_1^{(\alpha)} = -\alpha, \quad
			g_0^{(\alpha)} > g_2^{(\alpha)} > g_3^{(\alpha)} > \dots > 0, \\[4pt]
			\displaystyle \sum_{k=0}^{\infty} g_k^{(\alpha)} = 0,\quad
			\displaystyle \sum_{k=0}^{n} g_k^{(\alpha)} < 0, \quad n \ge 1.
		\end{cases}
		\]}
\end{proposition}
When combined with the variable diffusion coefficient \( a(x) \), the spatial discretization for the fractional diffusion term at time level   {\( m \)} becomes
\begin{equation}\label{equ6}
	{a(x) \, \mathcal{D}_x^\alpha u(x,t)\big|_{x = x_i, \, t = t_{m}}
		= \frac{a(x_{i})}{(\Delta x)^{\alpha}} \sum_{k=0}^{i+1} g_{i - k + 1}^{(\alpha)} u(x_{k},t_{m}) + \mathcal{O}(\Delta x).
}\end{equation}
{Finally,  by coupling  the \(\theta\)-method \cite{pan2016} to discretize the temporal derivative in  \eqref{eq2-1d} with  the approximation  in \eqref{equ6}}, we obtain the following time-stepping scheme
{\begin{equation}\label{equ7}
		\frac{u_i^{(m+1)} - u_i^{(m)}}{\Delta t} = \theta  \left( \frac{a(x_i)}{(\Delta x)^\alpha} \sum_{k=0}^{i+1} g_k^{(\alpha)} u_{i-k+1}^{(m+1)} \right) + (1 - \theta) \left( \frac{a(x_i)}{(\Delta x)^\alpha} \sum_{k=0}^{i+1} g_k^{(\alpha)} u_{i-k+1}^{(m)} \right) + f_i^{(m,\theta)},
\end{equation}}
where $0\leq\theta\leq 1$ is a weighting parameter, \( u_i^{(m)} \) approximates \( u(x_i, t_{m}) \),
and   {$f_i^{(m,\theta)} =  \theta f(x_i,t_{m+1}) + (1-\theta)f(x_i,t_{m}) $}.
The value of \( \theta \) determines the nature of the method: \( \theta = 0 \) corresponds to the Forward Euler method, \( \theta = 1 \) yields the Backward Euler method, and \( \theta = \frac{1}{2} \) gives the Crank--Nicolson scheme. \par
{By defining $u_N^{(m)}=\left [u_i^{(m)}\right ]^T_{i=1,\ldots,N}$, that is the numerical approximation of the solution \( u(x, t_m) \) at all points of the interior spatial grid   {\( x_i, {i=1, \ldots,N}\)}, and accordingly, the vector $f_N^{(m,\theta)}=\left[f_i^{(m,\theta)}\right]^T_{i=1,\ldots,N}$, and by multiplying both sides of \eqref{equ7} by \( \Delta t \) and simplifying, we obtain}
\begin{equation*}
	u_n^{(m+1)}-u_n^{(m)}=   {\theta G_{\alpha,N} \, u_N^{(m+1)} + (1 - \theta)  G_{\alpha,N} \, u_N^{(m)} + \Delta t  f_N^{(m,\theta)}, \quad m = 0, \dots, M-1,}
\end{equation*}
where
\begin{equation}\label{D*G}
	G_{\alpha,N}= \gamma {\mathcal D}_{N}(a) \overline{G}_{\alpha, N},
\end{equation}
with $\gamma = {\Delta t}/{(\Delta x)^{\alpha}}$,
${\mathcal D}_{N}(a) = \mathrm{diag} \left( a ( x_i) \right)_{i=1,\ldots,N}$  denoting the \( N \)-th diagonal sampling matrix and
\begin{equation} \label{eq:Gbar}
	\overline{G}_{\alpha, N} =
	\begin{bmatrix}
		g_1^{(\alpha)} & g_0^{(\alpha)}             & 0             & \cdots & 0 \\
		g_2^{(\alpha)} & g_1^{(\alpha)} & g_0^{(\alpha)}            & \ddots & \vdots \\
		\vdots & g_2^{(\alpha)} & g_1^{(\alpha)} & \ddots & 0 \\
		\vdots        & \ddots        & \ddots        & \ddots & g_0^{(\alpha)}  \\
		g_{N}^{(\alpha)} &\cdots & \cdots & g_2^{(\alpha)} & g_1^{(\alpha)}
	\end{bmatrix} \in \mathbb{R}^{N \times N},
\end{equation}
{denoting the Hessenberg Toeplitz matrix, that is
	\begin{equation}\label{eq8bis}
		\left (I_N - \theta  G_{\alpha,N} \right ) u_N^{(m+1)} = \left (I_N + (1 - \theta)  G_{\alpha,N} \right ) u_N^{(m)}
		+ \Delta t  f_N^{(m,\theta)}.
\end{equation}}
In terms of notations, notice that ${\mathcal D}_{N}(a)$ coincides with $D_N(a)$ as in (\ref{eq:da}) with $d=r=1$ if the domain of definition of the function $a$ is $[0,1]$. In general $[L_1,R_1]\neq [0,1]$ and hence ${\mathcal D}_{N}(a)=D_N(\hat a)$ with $\hat a(\hat x)=a(L_1 + (R_1-L_1)\hat x)$ and $x=L_1 + (R_1-L_1)\hat x$, $\hat x\in [0,1]$: this specification is useful for applying the GLT axioms when deducing the Weyl distributions of the related matrix sequences, as indicated in Remark \ref{rem:extending2}.

Equation \eqref{eq8bis} leads to the following   {all-at-once formulation} linear system
{\begin{align}
		\begin{bmatrix}
			I_{N} - \theta  G_{\alpha,N} & & & \\
			-(I_{N} + (1 - \theta)  G_{\alpha,N}) & I_{N} - \theta  G_{\alpha,N} & & \\
			& \ddots & \ddots & \\
			& & -(I_{N} + (1 - \theta)  G_{\alpha,N}) & I_{N} - \theta  G_{\alpha,N}
		\end{bmatrix}
		\begin{bmatrix}
			u_N^{(1)} \\
			u_N^{(2)} \\
			\vdots \\
			u_N^{(M)}
		\end{bmatrix} & \nonumber \\
		& \hskip -6cm =
		\begin{bmatrix}
			\Delta t \, f_N^{(1,\theta)} + (I_{N}  + (1 - \theta)  G_{\alpha,N}) u_N^{(0)}  \\
			\Delta t \, f_N^{(2,\theta)} \\
			\vdots \\
			\Delta t \, f_N^{(M,\theta)}
		\end{bmatrix}. \label{equ9}
\end{align}}
This system can be compactly expressed using Kronecker products in the following form
{\begin{equation}\label{equ10}
		\bar{A} \bar{\mathbf{u}} = \bar{\mathbf{f}},
	\end{equation}
}%
where
{\begin{align}
		\bar{A} &= (I_M \otimes (I_N - \theta  G_{\alpha,N}) - Z_{M} \otimes (I_N + (1 - \theta)  G_{\alpha,N})) \in \mathbb{R}^{NM \times NM}, \label{eq:Abar} \\
		%
		Z_{N} &=
		\begin{bmatrix}
			0      &        &        &        \\
			1      & 0      &        &        \\
			& \ddots & \ddots &        \\
			&  & 1      & 0
		\end{bmatrix}\in \mathbb{R}^{M \times M} \nonumber,
	\end{align}
}
and $\bar{\mathbf{u}}, \bar{\mathbf{f}} \in \mathbb{R}^{NM \times 1}$.
%
%
The coefficient matrix in \eqref{eq:Abar} can therefore be expressed compactly using Kronecker products as
{\[
	\bar{A} = H_{\theta,M} \otimes   G_{\alpha,N} + H_{M} \otimes I_N.
	\]}
where the matrices   {$H_{M}$} and $H_{\theta,M}$ are defined as follows
\[
{H_{M}} =
\begin{bmatrix}
	1           &             &             &        \\
	-1          & 1           &             &        \\
	& \ddots      & \ddots      &        \\
	&             & -1          & 1
\end{bmatrix}
\in \mathbb{R}^{M \times M}, \quad
{H_{\theta,M} = -
	\begin{bmatrix}
		\theta      &             &             &        \\
		1 - \theta  & \theta      &             &        \\
		& \ddots      & \ddots      &        \\
		&             & 1 - \theta  & \theta
	\end{bmatrix}
}\in \mathbb{R}^{M \times M}.
\]
To calculate the solution of \eqref{equ10}
in a parallel pattern form, we exchange the order of the Kronecker product in   {\( \bar{A}\)}. This leads to an equivalent system
{
	\begin{equation}\label{equ11}
		\widetilde{A}  \widetilde{\mathbf{u}}=\widetilde{\mathbf{f}},
\end{equation}}
where
\[
\widetilde{A} = G_{\alpha,N} \otimes H_{\theta,M} + I_N \otimes H_{M},
\] and   {\( \widetilde{\mathbf{f}} \)} is obtained by reshaping   {\( \bar{\mathbf{f}} \in \mathbb{R}^{NM} \)} into an \( N \times M \) matrix and then stacking its rows into a single column $MN\,\times \,1$ vector. Equivalently,   {\( \widetilde{\mathbf{f}} \)} is formed by stacking the columns of the transpose of the reshaped matrix.\\
By pre-multiplying both sides of \eqref{equ11} with \( (I_N \otimes   {H_{\theta,M})}^{-1} \), we obtain an equivalent formulation
{\begin{equation}\label{equ12}
		{A_{\alpha,N,M}} \mathbf{u} = \mathbf{f},
	\end{equation}
}where
{\[
	{A_{\alpha,N,M}} =  G_{\alpha,N} \otimes I_M + I_N \otimes H_{\theta,M}^{-1} H_{M} =  G_{\alpha,N} \otimes I_M + I_N \otimes Q_{\theta, M},
	\]
}
with the matrix \( Q_{\theta, M} = H_{\theta,M}^{-1} H_{M}\) defined as
\begin{equation}\label{Q mat}
	{Q_{\theta, M} = -
		\begin{bmatrix}
			\frac{1}{\theta} & & &  &  \\
			\frac{-1}{\theta^{2}} & \frac{1}{\theta} & &  &  \\
			\frac{-(\theta-1)}{\theta^3} & \frac{-1}{\theta^{2}} & \frac{1}{\theta} &  &  \\
			\vdots & \ddots & \ddots & \ddots &  \\
			\frac{-(\theta-1)^{M-2}}{\theta^{M}} & \cdots & \frac{-(\theta-1)}{\theta^3} & \frac{-1}{\theta^{2}} & \frac{1}{\theta}
		\end{bmatrix}
	}
	\in \mathbb{R}^{M \times M}.
\end{equation}
The transformed right-hand side is given by
{\(
	\mathbf{f} = (I_N\otimes H_{\theta,M})^{-1} \, \widetilde{\mathbf{f}}=(I_N\otimes H_{\theta,M}^{-1}) \, \widetilde{\mathbf{f}}.
	\)}
\subsubsection{  {Two-dimensional} }
Let $ N_1 $, $N_{2}$ and $ M $ be positive integers representing the number of interior spatial and temporal steps, respectively. 
The spatial  and temporal meshes are defined as 
\begin{alignat*}{3}
	x_i &= L_1 + i \Delta x, && \qquad\Delta x =  \frac{(R_1 - L_1)}{(N_{1} + 1)}, && \qquad 0 \leq i \leq N_{1}+1, \\
	y_j &= L_2 + j \Delta y, && \qquad\Delta y =  \frac{(R_2 - L_2)}{(N_{2} + 1)}, && \qquad 0 \leq j \leq N_{2}+1, \\
	t_m &= m\Delta t, && \qquad\Delta t=\frac{T}{M}, && \qquad 0 \leq m \leq M.
\end{alignat*}
\\

\noindent As in the one-dimensional setting, we adopt the shifted Grünwald finite difference scheme \cite{MeerschaertTadjeran2006} to approximate the left-sided Riemann-Liouville fractional derivatives.
The continuous fractional operator \( \mathcal{D}_{x}^{\alpha_1} u(x,y,t) \) is discretized at each interior point \( x_{i} \) using the expression below
\begin{equation}\label{equ44}
	\mathcal{D}_{x}^{\alpha_1}u(x,y,t) \big|_ {x= x_{i}, y = y_{j}, t = t_{m}} = \frac{1}{(\Delta x)^{\alpha_1}} \sum_{k=0}^{i+1} g_{i - k + 1}^{(\alpha_1)} u(x_{k},y_{j},t_{m})+ \mathcal{O}(\Delta x),
\end{equation}
with related properties of $g^{(\alpha_1)}_{k}$  already given in the one-dimensional setting in Proposition~\ref{alpha coff}.
Similarly, to discretize the \( \mathcal{D}_{y}^{\alpha_2} u(x,y,t) \) at each interior point \( y_{j} \), we deduce the following relation
\begin{equation}\label{equ66}
	\mathcal{D}_{y}^{\alpha_2} u(x, y,t) \big|_{x= x_{i}, y = y_{j}, t = t_{m}}= \frac{1}{(\Delta y)^{\alpha_2}} \sum_{l=0}^{j+1} g_{j - l + 1}^{(\alpha_2)} u(x_{i},y_{l},t_{m})+ \mathcal{O}(\Delta y),
\end{equation}
with related properties of $g^{(\alpha_2)}_{l}$  already given in the one-dimensional setting in Proposition~\ref{alpha coff}.
\par
Applying the spatial discretization \eqref{equ44}–\eqref{equ66} and the temporal discretization, we use the \(\theta\)-method \cite{pan2016} on \eqref{eq1}, which yields the following time-stepping scheme
\begin{align}\label{equ777}
	\frac{u_{i,j}^{(m+1)} - u_{i,j}^{(m)}}{\Delta t} & = \ \theta \ a(x_{i},y_{j}) \left( \frac{1}{(\Delta x)^{\alpha_1}} \sum_{k=0}^{i+1} g_k^{(\alpha_1)} u_{i-k+1,j}^{(m+1)} + \frac{1}{(\Delta y)^{\alpha_2}} \sum_{l=0}^{j+1} g_l^{(\alpha_2)} u_{i,j-l+1}^{(m+1)}\right) \nonumber\\
	& \hskip -2cm + (1 - \theta) \ a(x_{i},y_{j}) \left( \frac{1}{(\Delta x)^{\alpha_1}} \sum_{k=0}^{i+1} g_k^{(\alpha_1)} u_{i-k+1,j}^{(m)} + \frac{1}{(\Delta y)^{\alpha_2}} \sum_{l=0}^{j+1} g_l^{(\alpha_2)} u_{i,j-l+1}^{(m)} \right) + f_{i,j}^{m,\theta},
\end{align}
where
$u^{(m)}_{i,j}$ approximates $u(x_{i}, y_{j}, t_m)$ and
$f_{i,j}^{(m,\theta)} =  \theta f(x_{i},y_{j},t_{m+1}) + (1-\theta)f(x_{i},y_{j},t_{m})$.
By defining the vectors $$u_{\textbf{n}}^{(m)}=
\left [u_{i,j}^{(m)}\right ]^T_{i=1,\ldots,N_1, \ j=1,\ldots,N_2},\qquad f_{\textbf{n}}^{(m,\theta)}=\left[f_{i,j}^{(m,\theta)}\right]^T_{i=1,\ldots,N_1, \ j=1,\ldots,N_2},\;\;\; \textbf{n}=N_{1}N_{2}$$ according to the lexicographic order and by
multiplying both sides of \eqref{equ777} by \( \Delta t \) and simplifying, we obtain
\begin{equation}\label{equ8}
	u_{\textbf{n}}^{(m+1)}-u_{\textbf{n}}^{(m)}= \theta  G_{\boldsymbol{\alpha},\textbf{n}} \, u_{\textbf{n}}^{(m+1)} + (1 - \theta)   G_{\boldsymbol{\alpha},\textbf{n}} \, u_{\textbf{n}}^{(m)} + \Delta t  f_{\textbf{n}}^{(m,\theta)}, \quad m = 0, \dots, M-1,
\end{equation}
where
\begin{align}
	G_{\boldsymbol{\alpha},\textbf{n}} & = D_\textbf{n}(a) U_{\boldsymbol{\alpha},\textbf{n}}
	\in \mathbb{R}^{\textbf{n}\times\textbf{n}},\label{D*U} \\
	{U_{\boldsymbol{\alpha},\textbf{n}}} & {=   \gamma_1 \overline{G}_{\alpha_1,n_1}\otimes I_{n_2} + I_{n_1} \otimes \gamma_2 \overline{G}_{\alpha_2,n_2}
		\in \mathbb{R}^{\textbf{n}\times\textbf{n}}},\label{U mat 2D} \\
	{\mathcal D}_{\textbf{n}}(a) &= \underset{j=1,\dots,n_2}{\underset{i=1,\dots,n_1}{\operatorname{diag}}}\left( a\left(x_{i},y_{j}\right)\right)\in \mathbb{R}^{\textbf{n}\times\textbf{n}},
\end{align}
with $\overline{G}_{\alpha_1,N_1}$ and $\overline{G}_{\alpha_2,N_2}$ defined as in~\eqref{eq:Gbar}.
In terms of notations, notice that ${\mathcal D}_{\textbf{n}}(a)$ coincides with $D_{\textbf{n}}(a)$ as in (\ref{eq:da}) with $d=2, r=1$ if the domain of definition of the function $a$ is $[0,1]^2$. In general $[L_1,R_1]\times [L_2,R_2]\neq [0,1]^2$ and hence ${\mathcal D}_{\textbf{n}}(a)=D_{\textbf{n}}(\hat a)$ with $\hat a(\hat x,\hat y)=a(L_1 + (R_1-L_1)\hat x,L_2 + (R_2-L_2)\hat y)$ and $(x,y)=(L_1 + (R_1-L_1)\hat x,L_2 + (R_2-L_2)\hat y)$, $(\hat x,\hat y)\in [0,1]^2$: as in the one-dimensional setting, the latter setting is a key point for deducing the Weyl distributions of the related matrix sequences, using the GLT axioms; see Remark \ref{rem:extending2}.
\par
By combining the equations at all levels  into the single comprehensive linear system, known as the all-at-once formulation, and repeating verbatim the manipulations considered in the one-dimensional setting, which are not formally influenced by the two-dimensional space variable,  we obtain the equivalent formulation
\begin{equation}\label{equ122}
	{A_{\boldsymbol{\alpha},\mathbf{n},M}} \mathbf{u} = \mathbf{f},
\end{equation}
with matrix
\[ {A_{\boldsymbol{\alpha},\mathbf{n},M}} = G_{\boldsymbol{\alpha},\textbf{n}} \otimes I_M + I_{\textbf{n}}  \otimes Q_{\theta, M}, \]
where $Q_{\theta, M}$ is define as in \eqref{Q mat},
and transformed right-hand side
\[
\mathbf{f} = (I_{\textbf{n}} \otimes H_\theta)^{-1} \, \widetilde{\mathbf{f}} =(I_{\textbf{n}} \otimes H_\theta^{-1}) \, \widetilde{\mathbf{f}},
\]
where
\( \widetilde{\mathbf{f}} \) is obtained by reshaping the original right-hand side
\[ \bar{\mathbf{f}} =
\begin{bmatrix}
	\Delta t \, \mathbf{f}_{\textbf{n}}^{(1,\theta)} + (I_{\textbf{n}}  + (1 - \theta) G_{\boldsymbol{\alpha},\textbf{n}}  u_{\textbf{n}}^{(0)} \\
	\Delta t \, \mathbf{f}_{\textbf{n}}^{(2,\theta)} \\
	\vdots \\
	\Delta t \, \mathbf{f}_{\textbf{n}}^{(M,\theta)}
\end{bmatrix}
\in \mathbb{R}^{\textbf{n}M \times 1},
\]
into an \( \textbf{n} \times M \) matrix and then stacking its rows into a single column $M\textbf{n}\,\times \,1$ vector. Equivalently, \( \widetilde{\mathbf{f}} \) is formed by stacking the columns of the transpose of the reshaped matrix.

\subsection{GLT, singular value, and eigenvalue analysis}
In the current section, we first introduce a basic definition, and then proceed with the GLT analysis of the coefficient matrix sequences derived from the previous approximation methods. More precisely the section is divided into two parts: in the first we study the matrix sequences arising from the space discretizations, both in 1D and 2D. Then, in the second part, we provide the GLT analysis for the global matrix sequences, incorporating also the term coming from the time discretization.
\subsubsection{GLT analysis of the space matrix sequences}
\begin{definition}\label{def1}
	The Wiener class is the set of functions
	\[
	f(\boldsymbol{\theta})
	=
	\sum_{\mathbf{k}\in\mathbb{Z}^d} f_{\mathbf{k}}\,
	e^{\,i\,\langle \mathbf{k},\boldsymbol{\theta}\rangle},
	\qquad \boldsymbol{\theta}\in[-\pi,\pi]^d,
	\]
	such that
	\[
	\sum_{\mathbf{k}\in\mathbb{Z}^d} |f_{\mathbf{k}}| < \infty ,
	\]
	where \(i=\sqrt{-1}\) and
	\(\langle\mathbf{k},\boldsymbol{\theta}\rangle
	= k_1\theta_1+\cdots+k_d\theta_d\).
	The Wiener class forms a sub-algebra of the continuous \(2\pi\)-periodic
	functions defined on \([-\pi,\pi]^d\).
\end{definition}
In this view, we begin by reporting the symbol associated with the matrix sequence $\overline{G}_{\alpha, n}$ defined in equation~\eqref{D*G}.
\begin{proposition}\label{symb f_alpha}
	Let $\alpha \in (1,2)$. The symbol associated with $\overline{G}_{\alpha, N}$ belongs to the Wiener class, and
	\[
	f_\alpha(\xi_1) = \sum_{k=-1}^{\infty} g^{(\alpha)}_{k+1} e^{ik\xi_1}
	{=  e^{-i\xi_1}
		\left(1 + e^{i(\xi_1 + \pi)}\right)^{\alpha}, \quad \xi_1 \in [-\pi,\pi].
	}
	\]
\end{proposition}
\begin{proof}
	Let us observe that $\overline{G}_{\alpha, N} = [\, g^{(\alpha)}_{i-j+1} \,]_{i,j=1}^N$ with $g^{(\alpha)}_k = 0$ for $k<0$,  and let us define the function  $\sum_{k=-1}^{\infty} g^{(\alpha)}_{k+1} e^{ik\xi_1}$.
	When $\alpha \in (1,2)$, it is easy to see that $f_{\alpha}(\xi_1)$ lies in the Wiener class.
	In detail, from Proposition~\ref{alpha coff} we know that
	\[
	{g^{(\alpha)}_0 = 1 > 0,\quad g^{(\alpha)}_1 = -\alpha < 0,\quad g^{(\alpha)}_k > 0 \ \ \forall k \ge 2,\quad g^{(\alpha)}_k = 0 \ \ \forall k < 0.}
	\]
	Then
	\begin{align}
		\sum_{k=-1}^{\infty} \big| g_{k+1}^{(\alpha)} \big|
		&= \big| g_{0}^{(\alpha)} \big| + \big| g_{1}^{(\alpha)} \big| + \sum_{k=1}^{\infty} \big| g_{k+1}^{(\alpha)} \big| \nonumber \\
		&=g_{0}^{(\alpha)} +\left( - g_{1}^{(\alpha)}\right) + \sum_{k=1}^{\infty} g_{k+1}^{(\alpha)}
		\quad \text{(since $g_{0}^{(\alpha)}>0$, $g_{1}^{(\alpha)}<0$, $g_{k+1}^{(\alpha)}>0$ for $k\ge 1$)} \nonumber \\
		&= g_{0}^{(\alpha)}  - g_{1}^{(\alpha)} + \sum_{k=1}^{\infty} g_{k+1}^{(\alpha)}. \label{eq:1}
	\end{align}
	From Proposition~\eqref{alpha coff}, we have $\sum_{k=0}^{\infty} g_{k}^{(\alpha)} = 0.$
	Therefore,
	\begin{equation}\label{eq:2}
		{\sum_{k=1}^{\infty} g_{k+1}^{(\alpha)} = \sum_{k=2}^{\infty} g_{k}^{(\alpha)} = \sum_{k=0}^{\infty} g_{k}^{(\alpha)} - g_{0}^{(\alpha)} - g_{1}^{(\alpha)}
			=  - g_{0}^{(\alpha)} - g_{1}^{(\alpha)}
			= -1 + \alpha
		}
	\end{equation}
	%
	{Now, plugging \eqref{eq:2} in \eqref{eq:1}, we find
		\begin{equation}
			\sum_{k=-1}^{\infty} \big| g_{k+1}^{(\alpha)} \big|
			=  g_{0}^{(\alpha)} -g_{1}^{(\alpha)}  + ( -1 + \alpha )
			= 1 + \alpha + ( -1 + \alpha ) = 2\alpha.
	\end{equation}}
	Since this sum is finite, $f_{\alpha}(\xi_1)$ belongs to the Wiener class for $\alpha \in (1,2)$.
	\par
	To obtain an explicit formula for the symbol $f_{\alpha}(\xi_1)$, let us recall the definition of $g^{(\alpha)}_k$ given in~\eqref{eq:g-binomial} and rewrite $f_{\alpha}(\xi_1)$ in the following manner
	\[
	f_{\alpha}(\xi_1) = \sum_{k=-1}^{\infty} g^{(\alpha)}_{k+1} e^{ik\xi_1}
	= \sum_{k=0}^{\infty} g^{(\alpha)}_k e^{i(k-1)\xi_1}
	= \sum_{k=0}^{\infty}  {(-1)^{k}} \binom{\alpha}{k} e^{i(k-1)\xi_1}.
	\]
	We now write $(-1)^k$ as $e^{ik\pi}$ so that
	\[
	{f_{\alpha}(\xi_1) = \sum_{k=0}^{\infty} \binom{\alpha}{k} e^{i(k-1)\xi_1} e^{ik\pi}
		= e^{-i\xi_1} \sum_{k=0}^{\infty} \binom{\alpha}{k} e^{ik(\xi_1+\pi)}.
	}
	\]
	Applying the well-known binomial series
	\[
	(1+z)^{\alpha} = \sum_{k=0}^{\infty} \binom{\alpha}{k} z^k, \quad z \in \mathbb{C}, \ |z| \le 1,\ \alpha>0,
	\]
	with $z = e^{i(\xi_1+\pi)}$, we finally deduce a very informative expression i.e.
	\[
	{f_{\alpha}(\xi_1) = e^{-i\xi_1} \left( 1 + e^{i(\xi_1+\pi)} \right)^{\alpha}.
	}
	\]
\end{proof}
In the two-dimensional case,  {$\overline{G}_{\alpha_1,N_{1}}$} has the same symbol as obtained in Proposition~\ref{symb f_alpha}.
Similarly, the matrix  {$\overline{G}_{\alpha_2, N_{2}}$} has the same Toeplitz-like structure as $\overline{G}_{\alpha_1,N_{1}}$, with parameter $\alpha_2$ in place of $\alpha_1$.
Therefore, its associated symbol is given by
\[
{f_{\alpha_2}(\xi_{2}) =  {e^{-i \xi_{2}}} \left( 1 + e^{i(\xi_{2} + \pi)} \right)^{\alpha_2},
	\quad \xi_{2} \in [-\pi, \pi],
}
\]
which is the exact analogue of $f_{\alpha}(\xi_{1})$ obtained in Proposition~\ref{symb f_alpha}.
{
	\begin{proposition}\label{symb q_theta}
		Let $Q_{\theta, M} = H_{\theta, M}^{-1} H_M$.
		Then,
		\begin{equation}\label{Q symb}
			\{ Q_{\theta, M} \}_M \sim_{\mathrm{GLT}} q_{\theta} (\xi_{3})
		\end{equation}
		where
		$q_{\theta} (\xi_{3})=\frac{1 - e^{i\xi_3}}{\theta + (1-\theta) e^{i\xi_3}}$.
	\end{proposition}
	\begin{proof}
		In both the one and two-dimensional cases, both $H_M$ and $H_{\theta,M}$ are Toeplitz matrices generated by the symbols
		\[
		\rho_1(\xi_3) = 1 - e^{i\xi_3}, \qquad \rho_2(\xi_3) = \theta + (1-\theta) e^{i\xi_3},
		\]
		respectively, where $\xi_3 \in [-\pi, \pi]$. According to the GLT algebra and the standard properties of Toeplitz sequences \cite[Section 6.2]{garoni2017}, the matrix sequence $Q_{\theta, M} = H_{\theta, M}^{-1} H_M$ is associated with the GLT symbol
		\[
		q_{\theta} (\xi_{3}) = \frac{\rho_1(\xi_3)}{\rho_2(\xi_3)}.
		\]
		Therefore,
		\[
		\{ Q_{\theta, M} \}_M \sim_{\mathrm{GLT}} q_{\theta} (\xi_{3}).
		\]
\end{proof}}
%
\begin{theorem}\label{1D th:}
	If \( a : [0, 1] \to \mathbb{R} \) is continuous a.e., then
	\begin{itemize}
		\item[$a1)$]
		$\{  {{\mathcal D}_{N}(a)\overline{G}_{\alpha,N}} \}_N \sim_{\mathrm{GLT}} a(x)f_\alpha(\xi_{1})$,
		\item[$a2)$]
		$\{  {{\mathcal D}_{N}(a)\overline{G}_{\alpha,N}} \}_N \sim_{\sigma} a(x)f_\alpha(\xi_{1})$,
	\end{itemize}
	with $f_\alpha$ as in Proposition~\ref{symb f_alpha}.
	
\end{theorem}
\begin{proof}
	To prove~$a1)$, we first recall that the argument relies on the fundamental properties of the GLT class, in particular the second part in GLT \ref{glt:symbols} and the third part stated in GLT \ref{glt:algebra}. More precisely,  we have
	\[
	{\mathcal D}_N(a)=D_{N}(a) =\underset{i=1,2,\ldots,N }{ \text{diag}}
	a(x_{i})
	\]
	and hence
	\[
	\{{{\mathcal D}_N(a)=D_N(a)}\}_N \sim_{\mathrm{GLT}} a(x)
	\]
	by using GLT \ref{glt:symbols}. Furthermore, by Proposition~\ref{symb f_alpha} via axiom GLT \ref{glt:symbols}, we find
	\[
	\{ {\overline{G}_{\alpha,N}}\}_N \sim_{\mathrm{GLT}} f_\alpha(\xi_{1}).
	\]
	In conclusion by equation \eqref{D*G} and GLT \ref{glt:algebra}, we obtain $a1)$ i.e.
	\[
	\{  {{\mathcal D}_{N}(a) \overline{G}_{\alpha, N}\}_N}  \sim_{\mathrm{GLT}} a(x)f_\alpha(\xi_{1}),
	\]
	so that, by invoking GLT \ref{glt:distrib}, we infer $a2)$.
\end{proof}
And what about
\begin{itemize}
	\item[$a3)$] $\{  {{\mathcal D}_{N}(a)\overline{G}_{\alpha,N}} \}_N \sim_{\lambda} a(x)f_\alpha(\xi_{1})$?
\end{itemize}
Since \( G_{\alpha,N} \) is not Hermitian, its eigenvalue distribution cannot be immediately deduced from GLT \ref{glt:distrib}, even if other tools are available \cite{eig-nonH}.

\begin{figure}
	\centering
	\includegraphics[width=\textwidth]{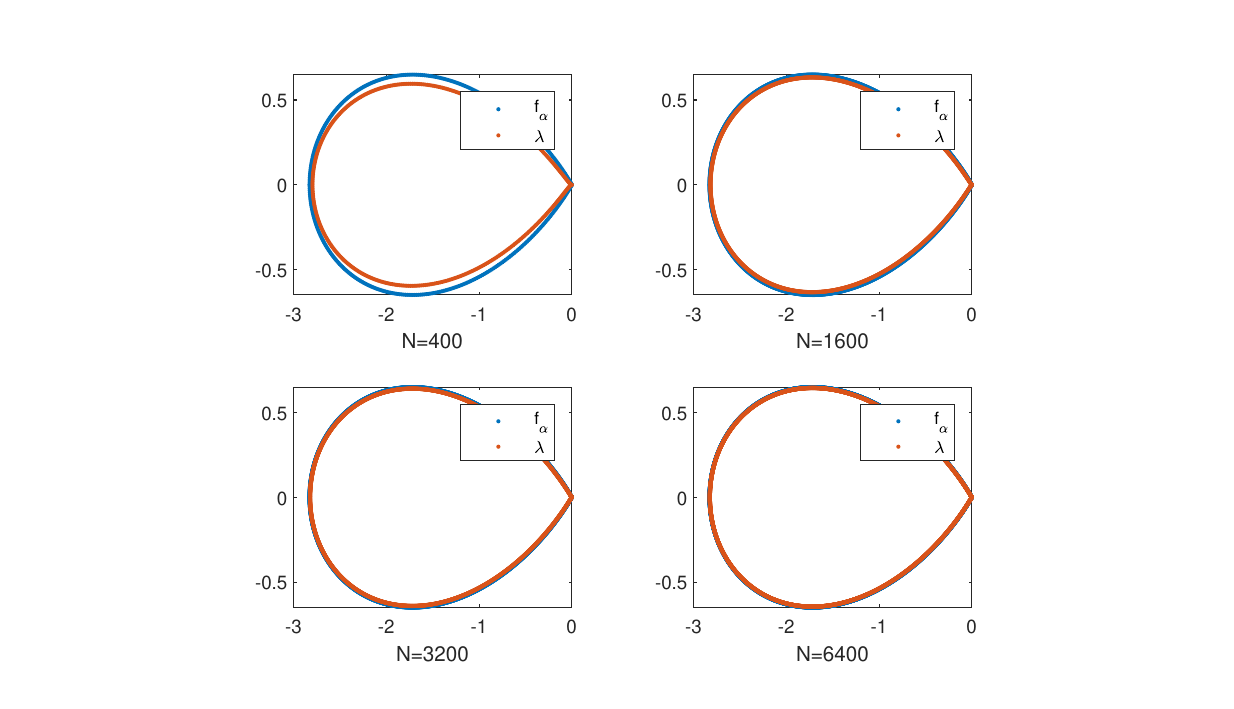}
	\caption{Comparison between the symbol $f_{\alpha}(\xi_{1})$ (blue stars) and the eigenvalues of ${\overline{G}}_{\alpha,N}$ (red stars) for different values of $N$ - case $\alpha = 1.5$. 
	}
	\label{fig:1d_const}
\end{figure}

\begin{figure}
	\centering
	\includegraphics[width=\textwidth]{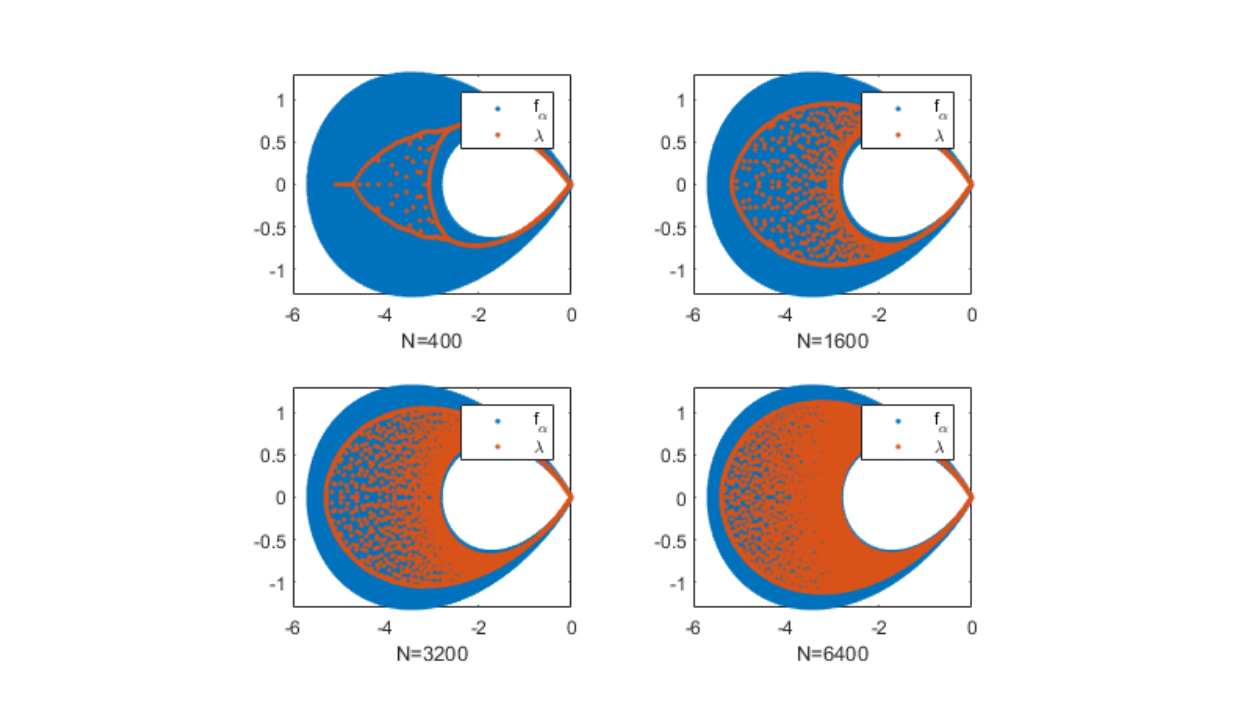}
	\caption{Comparison between the symbol $a(x)f_{\alpha}(\xi_{1})$ (blue stars) and the eigenvalues of $D_{N}(a) \,\overline{G}_{\alpha,N}$ (red stars) for different values of $N$ - case $a(x) = x^2+1$, $\alpha = 1.5$. 
	}
	\label{fig:1d_var}
\end{figure}

To illustrate the relationship between the eigenvalue distributions and the associated symbols, we present two representative examples for which Figures~\ref{fig:1d_const}–\ref{fig:1d_var} indicate that the question in $a3)$ has a positive answer, which is clear in the constant coefficient case. Indeed, although the convergence is slow, we believe that the answer is positive even in the variable coefficient case.

In the one-dimensional case, we first consider the constant-coefficient matrix
\(\overline{G}_{\alpha,N}\) and compare its eigenvalues with sampled values of the symbol
\( f_{\alpha}(\xi_{1}) \) (Figure~\ref{fig:1d_const}). The two graphs almost superpose and hence we observe a positive answer to the question in $a3)$ when $a(x)$ is a constant function. The fact holds also formally because of the results in \cite{eig-nonH}, where the authors study the canonical distribution of Toeplitz matrix sequences for a rich class of parametric generating functions which includes $f_\alpha$ for every $\alpha\in (1,2)$. Hence
\[
\{ \overline{G}_{\alpha,N} \}_N \sim_{\lambda} f_\alpha(\xi_{1}).
\]
We then examine the variable-coefficient case \( D_{N}(a) \,\overline{G}_{\alpha,N} \) with \( a(x) = x^{2} +1\), comparing the eigenvalues with sampled values of the symbol \( a(x) f_{\alpha}(\xi_{1}) \) (Figure~\ref{fig:1d_var}). Also in this more complicate setting, the numerics suggest that the answer to the question in $a3)$ is positive as reported
Figure~\ref{fig:1d_var}.
We also tried the more degenerate case of ${a(x)=x^2}$ and also in this case the agreement is remarkable. In fact, since $a(x)$ is continuous with range $[0,1]$, the range of $a(x)f_\alpha(\xi_{1})$ for fixed $\xi_1$ is exactly the segment connecting the complex zero with the value $f_\alpha(\xi_{1})$: as a consequence the total range of the GLT symbol $a(x)f_\alpha(\xi_{1})$ is exactly the convex hull of the range of $f_\alpha(\xi_{1})$. On the other hand, since \(1+x^2\) is continuous with range \([1,2]\), the
	eigenvalues fill the disconnected (highly nonconvex) region
	\[
	\bigcup_{y\in[1,2]} y\,\mathrm{range}(f_{\alpha})
	= \bigl\{\,y\,z : y\in[1,2],\ z\in\mathrm{range}(f_{\alpha})\,\bigr\},
	\]
	i.e.\ the union of the scaled copies of \(\mathrm{range}(f_{\alpha})\) by
	the real factors \(y\in[1,2]\), in agreement with the behaviour depicted in
	Figure~\ref{fig:1d_var}.
\vspace{0.07cm}

\noindent
Conversely, the theoretical proof is not available yet and it will be the subject of future investigations, starting from the key tools given in \cite{eig-nonH}.
\vspace{0.2cm}

\noindent The question \(a3)\) concerning the eigenvalue distribution of the non-Hermitian matrix sequence \(\{ \mathcal{D}_N(a)\overline{G}_{\alpha,N} \}_N\) deserves careful attention. Since the matrices are non-Hermitian, the GLT symbol obtained in \(a1)\) directly yields the singular value distribution but does not automatically imply the eigenvalue distribution. The numerical evidence presented in Figures~\ref{fig:1d_const}--\ref{fig:1d_var} suggests that the same symbol may describe the eigenvalue distribution in this setting, and this is indeed rigorously established in the constant coefficient case by the results of Bogoya, Böttcher, and Grudsky \cite{eig-nonH}. However, for the variable coefficient case, a rigorous proof remains open. Thus, while the numerical observations are compelling, the eigenvalue distribution should be regarded as a conjecture supported by numerical evidence rather than a proven result.
\begin{remark}\label{rem:extending1}
	The inclusion of the factor $\gamma=\Delta t / (\Delta x)^{\alpha}$ in the analysis is possible and indeed not difficult, so extending
	the scope of Theorem \ref{1D th:}. Assuming that $a : [0, 1] \to \mathbb{R}$ is continuous a.e. and that
	\[
	\lim_{\Delta t,\Delta x\rightarrow 0^+}\gamma(\Delta t,\Delta x)={\lim_{\Delta t,\Delta x\rightarrow 0^+}\frac{\Delta t}{(\Delta x)^\alpha}=}\ \gamma^*
	\]
	we have
	\begin{eqnarray}\label{gen1-1Dth}
		\{ G_{\alpha,N} \}_N & \sim_{\mathrm{GLT}} & \gamma^* a(x)f_\alpha(\xi_{1}), \\ \label{gen2-1Dth}
		\{ G_{\alpha,N} \}_N & \sim_{\sigma} & \gamma^* a(x)f_\alpha(\xi_{1}), \\ \label{gen3-1Dth}
		\{ {\overline{G}_{\alpha,N}}\}_N & \sim_{\lambda} & \gamma^* f_\alpha(\xi_{1}).
	\end{eqnarray}
	Alternatively, since we solve linear systems, one can get rid of the scaling factor very easily by writing in full generality
	\begin{eqnarray}\label{gen1-bis-1Dth}
		\{{\gamma}^{-1} G_{\alpha,N} \}_N & \sim_{\mathrm{GLT}} &  a(x)f_\alpha(\xi_{1}), \\ \label{gen2-bis-1Dth}
		\{{\gamma}^{-1} G_{\alpha,N} \}_N & \sim_{\sigma} &  a(x)f_\alpha(\xi_{1}), \\ \label{gen3-bis-1Dth}
		\{ {\gamma}^{-1}{\overline{G}_{\alpha,N}}\}_N & \sim_{\lambda} &  f_\alpha(\xi_{1}),
	\end{eqnarray}
	with the only hypothesis that $a : [0, 1] \to \mathbb{R}$ is continuous a.e. which is equivalent to the Riemann integrability of $a$.
	We remark that the relations \eqref{gen1-bis-1Dth}, \eqref{gen2-bis-1Dth}, \eqref{gen3-bis-1Dth} are meaningful also when the {limit $\gamma^*$} is equal to zero, while relations \eqref{gen1-1Dth}, \eqref{gen2-1Dth}, \eqref{gen3-1Dth} become somehow trivial and not informative {for $\gamma^*=0$}.
	\par
	We finally notice that the statements $a1)$, $a2)$, (\ref{gen1-1Dth}), (\ref{gen2-1Dth}), (\ref{gen1-bis-1Dth}), (\ref{gen2-bis-1Dth}) are all valid if $a$ is assumed complex-valued. However, the interest in the complex-valued case is at the moment quite limited, since in the modeling of the problem the function $a$ is required nonnegative, often strictly positive.
\end{remark}
\par
Now we consider a similar analysis in the  {two-level} space setting,  where we have to consider the whole matrix $U_{\alpha,\mathbf{n}}$, by joining the tools in Theorem \ref{1D th:} and Remark \ref{rem:extending1}.
\begin{theorem}\label{2D th:}
	{Let $\gamma_{r}=\Delta t / (\Delta x)^{\alpha_{r}}$, $r=1,2$, with $\alpha_1, \alpha_2 \in (1,2)$.}
	If $a : [0,1]^2 \to \mathbb{R}$ is continuous a.e. and
	\[
	\lim_{\Delta t,\Delta x\rightarrow 0^+}\gamma_r(\Delta t,\Delta x)= {\lim_{\Delta t,\Delta x\rightarrow 0^+} \frac{\Delta t}{(\Delta x)^{\alpha_r}} =}\  \gamma_r^*<\infty,\;\;\;\;r=1,2, 
	\]
	then
	\begin{itemize}
		\item [$b1)$] $\{G_{\boldsymbol{\alpha},\mathbf{n}}\}_\mathbf{n} \sim_{\mathrm{GLT}} a(x,y)\,\big(\gamma_1^* f_{\alpha_1}(\xi_1)
		+ \gamma_2^* f_{\alpha_2}(\xi_2)\big)$,
		\item [$b2)$] $\{G_{\boldsymbol{\alpha},\mathbf{n}}\}_\mathbf{n} \sim_{\sigma} a(x,y)\,\big(\gamma_1^* f_{\alpha_1}(\xi_1)
		+ \gamma_2^* f_{\alpha_2}(\xi_2)\big)$,
	\end{itemize}
	{where
		\[
		f_{\alpha_1}(\xi_1) = { e^{-i\xi_1}}\big(1+e^{i(\xi_1+\pi)}\big)^{\alpha_1}, \quad
		f_{\alpha_2}(\xi_2) =  {e^{-i\xi_2}}\big(1+e^{i(\xi_2+\pi)}\big)^{\alpha_2}, \quad
		(\xi_1,\xi_2) \in [-\pi,\pi]^2,
		\]
		and where necessarily $\gamma_2^*=0$ if $\alpha_1>\alpha_2$, and alternatively $\gamma_1^*=0$ if $\alpha_2>\alpha_1$.}
	In the special case where $\max\{\gamma_1^*,\gamma_2^*\}=0$, a normalization as in Remark \ref{rem:extending1} leads to
	\begin{itemize}
		\item [{$b1-1)$}] $\{{\gamma}^{-1}G_{\boldsymbol{\alpha},\mathbf{n}}\}_\mathbf{n} \sim_{\mathrm{GLT}} a(x,y)\,\big(\hat \gamma_1 f_{\alpha_1}(\xi_1)
		+ \hat \gamma_2 f_{\alpha_2}(\xi_2)\big)$,
		\item [{$b2-2)$}] $\{{\gamma}^{-1}G_{\boldsymbol{\alpha},\mathbf{n}}\}_\mathbf{n} \sim_{\sigma} a(x,y)\,\big(\hat \gamma_1 f_{\alpha_1}(\xi_1)
		+ \hat \gamma_2 f_{\alpha_2}(\xi_2)\big)$,
	\end{itemize}
	with $\gamma=\max\{\gamma_1,\gamma_2\}$. Here $\hat \gamma_1=1$ and $\hat \gamma_2=0$ if $\alpha_1>\alpha_2$, $\hat \gamma_2=1$ and $\hat \gamma_1=0$ if $\alpha_2>\alpha_1$, while the case $\alpha_1=\alpha_2$ leads to $\hat \gamma_1=1$ and $\hat \gamma_2\in (0,1)$ if $\gamma_1>\gamma_2$; $\hat \gamma_2=1$ and $\hat \gamma_1\in (0,1)$ if $\gamma_2>\gamma_1$.
\end{theorem}

\begin{figure}
	\centering
	\includegraphics[width=\textwidth]{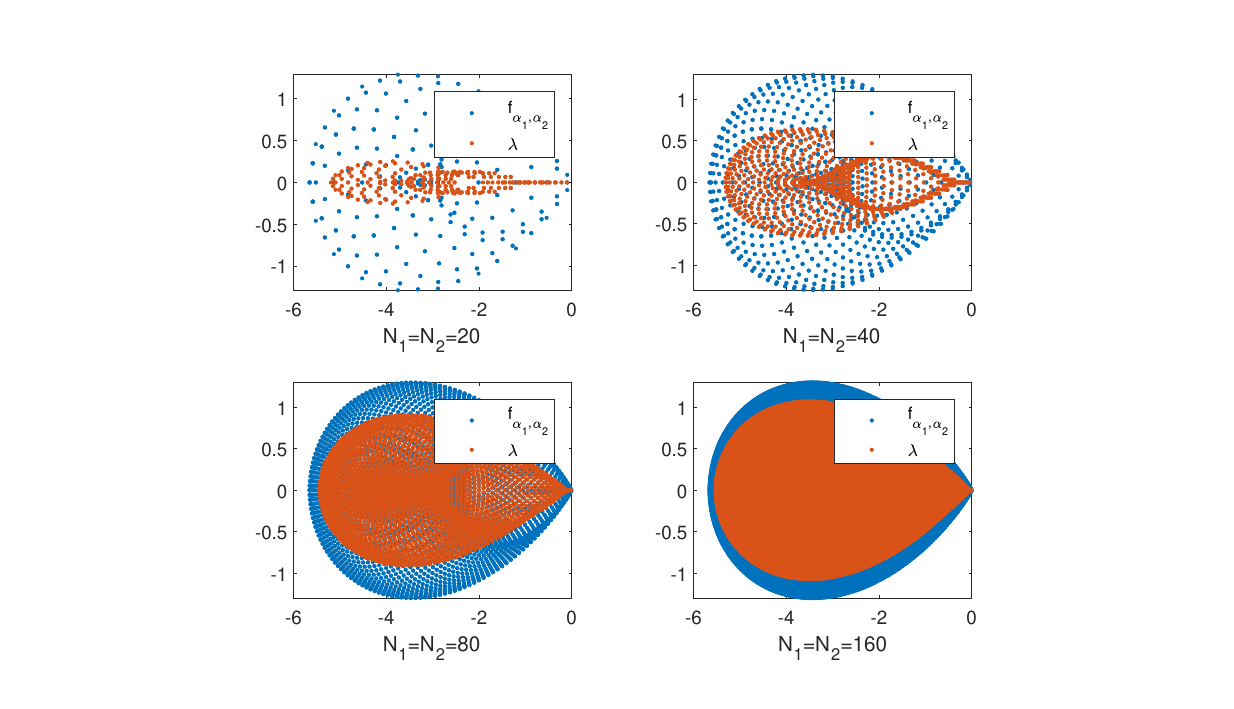}
	\caption{Comparison between the symbol $ {f_{\alpha_1}(\xi_{1})+f_{\alpha_2}(\xi_{2})}$ (blue stars) and the eigenvalues of
		{$\overline{G}_{\alpha_1,N_1}\otimes I_{N_2} + I_{N_1} \otimes \overline{G}_{\alpha_2,N_2}$} (red stars) for different values of $N_1=N_2$ - case $\alpha_1 = \alpha_2 = 1.5$, {$\gamma_1=\gamma_2=1$}.
	}
	\label{fig:2d_const}
\end{figure}

\begin{figure}
	\centering
	\includegraphics[width=\textwidth]{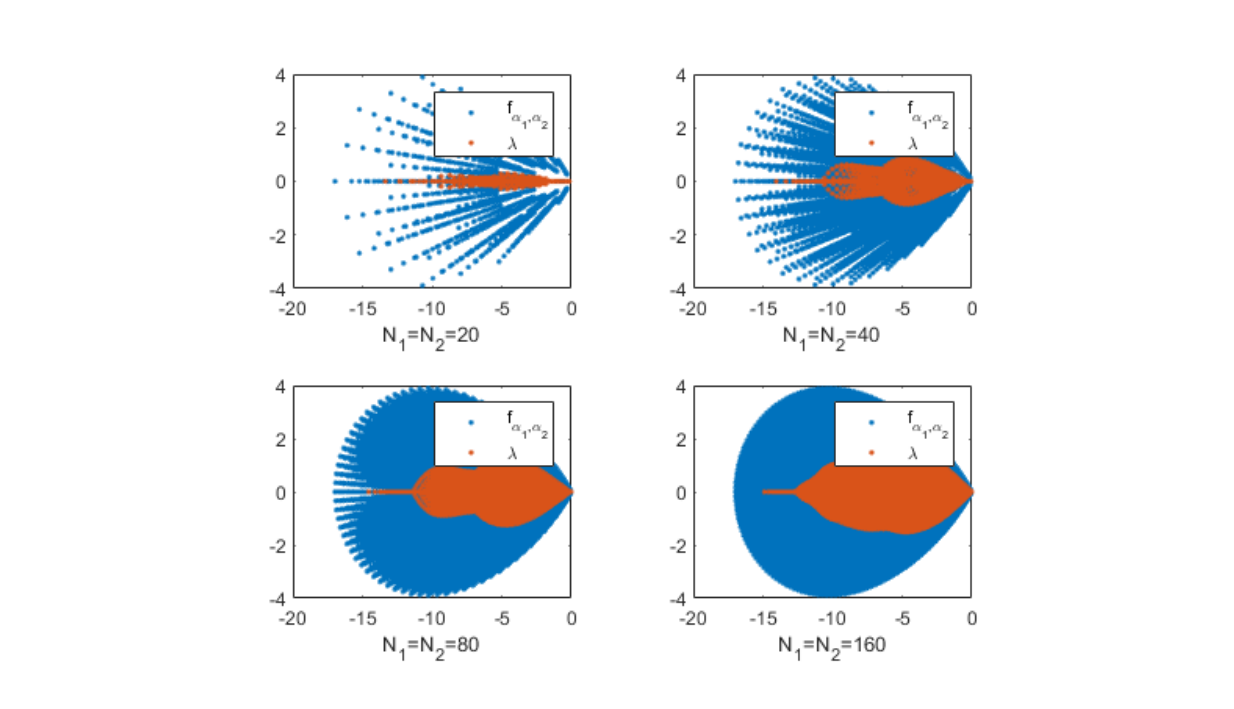}
	\caption{Comparison between the symbol $a(x_1,x_2)(f_{\alpha_1}(\xi_{1})+f_{\alpha_2}(\xi_{2}))$ (blue stars) and the eigenvalues of $D_{\mathbf{n}}(a)({\overline{G}_{\alpha_1,N_1}\otimes I_{N_2} + I_{N_1} \otimes \overline{G}_{\alpha_2,N_2}})$ (red stars) for different values of $N_1=N_2$ - case $a(x_1,x_2) = x_1^2+x_2^2+1$, $\alpha_1 = \alpha_2 = 1.5$, {$\gamma_1=\gamma_2=1$.}
		. 
	}
	\label{fig:2d_var}
\end{figure}

\begin{proof}
	It is sufficient to prove $b1)$, because $b2)$  follows from $b1)$ and GLT \ref{glt:distrib}. From the definition in~\eqref{D*U}, the matrix sequence $\{G_{\boldsymbol{\alpha},\mathbf{n}}\}_\mathbf{n}$ can be expressed as
	\[
	\{G_{\boldsymbol{\alpha},\mathbf{n}}\}_\mathbf{n} = \{D_{\mathbf{n}}(a) \, U_{\boldsymbol{\alpha},\mathbf{n}}\}_\mathbf{n},
	\]
	where $D_{\mathbf{n}}(a)$ is the diagonal sampling matrix and from~\eqref{U mat 2D} the matrix $U_{\boldsymbol{\alpha},\mathbf{n}}$ has the Kronecker sum form
	\[
	{U_{\boldsymbol{\alpha},\mathbf{n}} = \gamma_1 \overline{G}_{\alpha_1,N_1}\otimes I_{N_2} +
		I_{N_1} \otimes \gamma_2 \overline{G}_{\alpha_2,N_2}
	}
	\]
	%
	where $\overline{G}_{\alpha_1,N_1}$ and $\overline{G}_{\alpha_2,N_2}$  are one-level Toeplitz matrices generated by the shifted Grünwald coefficients in the $x$- and $y$-directions, respectively. Thus, by the second part of property GLT \ref{glt:symbols}, we deduce that
	\[
	\{D_{\mathbf{n}}\}_\mathbf{n}\sim_{\mathrm{GLT}} a(x,y).
	\]
	Using the first part of property GLT \ref{glt:symbols}, we easily infer that the matrix sequences satisfy the relations
	\[
	\{\overline{G}_{\alpha_1,N_1}\}_{N_1} \sim_{\mathrm{GLT}}  f_{\alpha_1}(\xi_1), \qquad
	\{\overline{G}_{\alpha_2,N_2}\}_{N_2} \sim_{\mathrm{GLT}}  f_{\alpha_2}(\xi_2).
	\]
	Invoking the $*$-algebra structure of the GLT class expressed by GLT \ref{glt:algebra}, we obtain that
	\[
	\{U_{\boldsymbol{\alpha},\mathbf{n}}\}_\mathbf{n} \sim_{\mathrm{GLT}} \gamma_1 f_{\alpha_1}(\xi_1) + \gamma_2 f_{\alpha_2}(\xi_2).
	\]
	Finally, by applying axiom GLT \ref{glt:algebra}, third item, to the sequences
	$\{D_{\mathbf{n}}\}_\mathbf{n}$ and $\{U_{\boldsymbol{\alpha},\mathbf{n}}\}_\mathbf{n}$,
	we conclude that
	\[
	\{G_{\boldsymbol{\alpha},\mathbf{n}}\}_\mathbf{n} \sim_{\mathrm{GLT}}
	a(x,y)\,\big( {f_{\alpha_1}(\xi_1) + f_{\alpha_2}(\xi_2)}\big).
	\]
	Then $b2)$ is an obvious consequence of $b1)$ and axiom GLT \ref{glt:distrib}.
	\par {Finally,  it is evident that $U_{\boldsymbol{\alpha},\mathbf{n}}$
		can be written as}
	\begin{align*}
		{   U_{\boldsymbol{\alpha},\mathbf{n}}} & = { \gamma_1 \overline{G}_{\alpha_1,N_1}\otimes I_{N_2} +
			I_{N_1} \otimes \gamma_2 \overline{G}_{\alpha_2,N_2}   }\\
		& = {\gamma_1 \left ( \overline{G}_{\alpha_1,N_1}\otimes I_{N_2} +
			(\Delta x)^{\alpha_1-\alpha_2} I_{N_1} \otimes \overline{G}_{\alpha_2,N_2}
			\right)}
	\end{align*}
	{so that the thesis follows if $\alpha_1 > \alpha_2$ and vice versa.}
\end{proof}
And what about
\begin{itemize}
	\item [$b3)$] $\{G_{\boldsymbol{\alpha},\mathbf{n}}\}_\mathbf{n} \sim_{\lambda} a(x,y)\,\big( {{\gamma_1^* f_{\alpha_1}(\xi_1) + \gamma_2^* f_{\alpha_2}(\xi_2)}}\big)$?
\end{itemize}
As in the one-dimensional case, the eigenvalue distribution is not an automatic GLT consequence for non-Hermitian matrices, but it is partially supported by the numerical evidence presented; see Figures~\ref{fig:2d_const}–\ref{fig:2d_var}.
In the two-dimensional case, we first consider the constant-coefficient matrix
{\( {U}_{\boldsymbol{\alpha},\mathbf{n}} \)} in equation~\eqref{U mat 2D} and compare its eigenvalues with sampled values of the symbol  {\( f_{\alpha_1}(\xi_{1}) + f_{\alpha_2}(\xi_{2}) \)} (Figure~\ref{fig:2d_const}).
We then examine the variable-coefficient case  {\( {D}_{\mathbf{n}}(a) {U}_{\boldsymbol{\alpha},\mathbf{n}} \)} in equation~\eqref{D*U}, with \( a(x,y) = x^{2} + y^{2}+1 \), comparing the eigenvalues with sampled values of the symbol
{\( a(x,y) \left( f_{\alpha_1}(\xi_{1}) + f_{\alpha_2}(\xi_{2}) \right) \)} (Figure~\ref{fig:2d_var}).
It is evident that the  {two-dimensional} setting is more intricate and we are not sure that the tools from \cite{eig-nonH} are useful as in the {one-dimensional} setting. The question remains open and it is left for future researches.
\begin{remark}\label{rem:extending2}
	Looking back at the original formulation in (\ref{eq2-1d}), we observe that the original domain
	{in the one-dimensional space case}
	is $[L_1,R_1]$ and in general $[L_1,R_1]\neq [0,1]$.
	For enlarging the scope of Theorem \ref{1D th:}, using the GLT theory, it is enough to consider ${\mathcal D}_{N}(a)=D_N(\hat a)$ with $\hat a(\hat x)=a(L_1 + (R_1-L_1)\hat x)$ and $x=L_1 + (R_1-L_1)\hat x$, $\hat x\in [0,1]$. In this way Theorem \ref{1D th:} and Remark \ref{rem:extending1} still hold with $\hat a$ instead of $a$ in the given statements.
	In the two-dimensional setting, looking at (\ref{eq1}), again we have $[L_1,R_1]\times [L_2,R_2]\neq [0,1]^2$ in general.
	However  ${\mathcal D}_{\textbf{n}}(a)=D_{\textbf{n}}(\hat a)$ with $\hat a(\hat x,\hat y)=a(L_1 + (R_1-L_1)\hat x,L_2 + (R_2-L_2)\hat y)$ and $(x,y)=(L_1 + (R_1-L_1)\hat x,L_2 + (R_2-L_2)\hat y)$, $(\hat x,\hat y)\in [0,1]^2$ and consequently Theorem \ref{2D th:} holds with the only change of  $\hat a$ instead of $a$.
	\par
	The flexibility of the GLT tools is of use in our context because also a general domain could be considered as far as it is Peano-Jordan measurable and in applications all the bounded domains have Riemann integrable characteristic function; see \cite[pp. 395–-399]{SerraCapizzano2003} for the first example and the idea of reduced GLT sequences in \cite[Section 3.1.4]{SerraCapizzano2006}, with the exhaustive theoretical development in \cite{Barbarino2022}.
	\par
	{We finally notice that the statements $b1)$, $b2)$ are valid if $a$ is assumed complex-valued. However, the interest in the complex-valued case is at the moment quite limited, since in the modeling of the problem the function $a$ is required nonnegative, often strictly positive.}
\end{remark}

\subsubsection{GLT analysis of the space-time matrix sequences}
We are ready to provide the GLT analysis for the space-time matrix sequences, both in 1D in space and 2D in space.
\begin{theorem}\label{1D orig mat}
	Let $a:[0,1]\to\mathbb{R}$ be continuous a.e. and let $\gamma=\Delta t / (\Delta x)^{\alpha}$, $\gamma^*$ as in Remark \ref{rem:extending1}. Then
	\begin{enumerate}
		\item[$c1)$] $\{ {A_{\alpha, N, M}}\}_{N, M} \sim_{\mathrm{GLT}}\gamma^* a(x)f_{\alpha}(\xi_1)+q_{\theta}(\xi_3)$,
		\item[$c2)$] $\{ {A_{\alpha, N, M}}\}_{N, M}\sim_{\sigma}       \gamma^* a(x)f_{\alpha}(\xi_1)+q_{\theta}(\xi_3)$,
	\end{enumerate}
	with ${A_{\alpha,N, M}}=G_{\alpha, N}\otimes I_M+I_N\otimes Q_{\theta,M}$, taking into account \eqref{equ11}-\eqref{equ12}.
\end{theorem}
\begin{proof}
	It is enough to prove $c1)$, since $c2)$ can be directly inferred from $c1$ and axiom \textbf{GLT~1}.
	From Theorem~\ref{1D th:} in the one-dimensional setting, we already know that the sequence of matrices
	$\{G_{\alpha,N}\}_{N}$ is distributed in the GLT sense as $\gamma^* a(x)\,f_{\alpha}(\xi_1)$.
	Moreover, by Proposition~\ref{symb q_theta}, the sequence $\{Q_{\theta,M}\}_{M}$ admits the GLT symbol
	\[
	q_{\theta}(\xi_3) \;=\; \frac{1-\mathrm{e}^{\mathrm{i}\xi_3}}{\theta+(1-\theta)\,\mathrm{e}^{\mathrm{i}\xi_3}},
	\]
	which represents the temporal discretization contribution arising from the $\theta$-method. Therefore, by combining these two GLT symbols through the Kronecker sum structure of $ {A_{\alpha,N, M}}$
	and invoking the $*$-algebra property of the GLT class (GLT \ref{glt:algebra}), we deduce that
	\[
	\{ {A_{\alpha,N, M}}\}_{N, M}=\{G_{\alpha,N}\otimes I_M+I_N\otimes Q_{\theta,M}\}_{N, M} \sim_{\mathrm{GLT}} \gamma^*a(x)f_{\alpha}(\xi_1)+q_{\theta}(\xi_3),
	\]
	which establishes the validity of $c1$. The statement $c2)$ then follows directly from $c1)$ and axiom GLT \ref{glt:distrib}.
\end{proof}
\begin{theorem}\label{2D orig mat}
	Let $a:[0,1]^2\to\mathbb{R}$ be continuous a.e. and let {$\gamma_r=\Delta t / (\Delta x)^{\alpha_r}$ and $\gamma^*_r$, $r=1,2$},
	as in Theorem \ref{2D th:}.
	{Furthermore, let
		\[
		{A_{\boldsymbol{\alpha},\mathbf{n},M}=G_{\boldsymbol{\alpha},\textbf{n}}\otimes I_M + I_{\mathbf{n}}\otimes Q_{\theta,M},}
		\]
		as defined in~\eqref{equ122}, with $f_{\alpha_1}$ and $f_{\alpha_2}$ taken from Proposition~\ref{symb f_alpha} and its $\alpha_2$–counterpart, and $q_{\theta}$ given by~\eqref{Q symb}.}
	Then
	\begin{enumerate}
		\item[$d1)$] $ {\{A_{\boldsymbol{\alpha},\mathbf{n},M}\}_{\mathbf{n},M}}\sim_{\mathrm{GLT}}a(x,y)(\gamma_1^*f_{\alpha_1}(\xi_1)+\gamma_2^*f_{\alpha_2}(\xi_2))+q_{\theta}(\xi_3),$
		\item[$d2)$] $ {\{A_{\boldsymbol{\alpha},\mathbf{n},M}\}_{\mathbf{n},M}}\sim_{\sigma}a(x,y)(\gamma_1^* f_{\alpha_1}(\xi_1)+ \gamma_2^* f_{\alpha_2}(\xi_2))+q_{\theta}(\xi_3),$
	\end{enumerate} 	%
	{where necessarily $\gamma_2^*=0$ if $\alpha_1>\alpha_2$ and alternatively $\gamma_1^*=0$ if $\alpha_2>\alpha_1$.}
\end{theorem}
\begin{proof}
	To establish the result, it is enough to verify $d1)$, since $d2)$ is an immediate consequence of $d1)$ and axiom GLT \ref{glt:distrib}.
	From Theorem~\ref{2D th:} we already know that
	\[
	{\{G_{\boldsymbol{\alpha},\textbf{n}}\}_{\mathbf{n}}}\sim_{\mathrm{GLT}}a(x,y)(\gamma_1^* f_{\alpha_1}(\xi_1)+ \gamma_2^* f_{\alpha_2}(\xi_2)),
	\]
	{where necessarily $\gamma_2^*=0$ if $\alpha_1>\alpha_2$ and alternatively $\gamma_1^*=0$ if $\alpha_2>\alpha_1$.}
	On the other hand, equation~\eqref{Q symb} ensures that $
	\{Q_{\theta,M}\}_{M}\sim_{\mathrm{GLT}}q_{\theta}(\xi_3).
	$
	By combining these two symbols according to the Kronecker-sum form of $A_{\boldsymbol{\alpha},\mathbf{n},M}$ and applying the $\ast$-algebra rules of the GLT class (GLT \ref{glt:algebra}), we arrive at
	\begin{align*}
		\{A_{\boldsymbol{\alpha},\mathbf{n},M}\}_{\mathbf{n},M}
		&=\{G_{\boldsymbol{\alpha},\textbf{n}}\otimes I_M+I_{\mathbf{n}}\otimes Q_{\theta,M}\}_{\textbf{n},M}\\
		&\sim_{\mathrm{GLT}}a(x,y)(\gamma_1^* f_{\alpha_1}(\xi_1)+ \gamma_2^* f_{\alpha_2}(\xi_2))+q_{\theta}(\xi_3),
	\end{align*}
	which confirms $d1)$. The property $d2)$ then follows from $d1)$ and axiom GLT \ref{glt:distrib}.
\end{proof}
And what about
\begin{itemize}
	\item[$c3)$] $\{ {A_{\alpha,N, M}}\}_{N, M}\sim_{\lambda} \gamma^* a(x)f_{\alpha}(\xi_1)+q_{\theta}(\xi_3)$?
	\item[$d3)$] $\{ {A_{\boldsymbol{\alpha},\mathbf{n},M}\}_{\mathbf{n},M}}\sim_{\lambda}a(x,y)(\gamma_1^* f_{\alpha_1}(\xi_1)+ \gamma_2^* f_{\alpha_2}(\xi_2))+q_{\theta}(\xi_3)$?
\end{itemize}
{Since the matrices $A_{\alpha,N, M}, {A_{\boldsymbol{\alpha},\mathbf{n},M}}$ are not Hermitian, their eigenvalue distribution is not guaranteed by the use of the GLT theory. Moreover the present case seems different from the previous ones reported in Figures \ref{fig:1d_const}, \ref{fig:1d_var}, \ref{fig:2d_const}, \ref{fig:2d_var}: indeed the numerical results seem to indicate that the statements $c3)$, $d3)$ are not true in the current setting, even when $a$ is chosen as a constant function.}
\vspace{0.3cm}

\noindent As noted in the discussion following Theorems~\ref{1D orig mat} and~\ref{2D orig mat}, the eigenvalue distribution of the space-time matrix sequences \(\{A_{\alpha,N, M}\}_{N, M}\) and \(\{A_{\alpha,\mathbf{n},M}\}_{\mathbf{n},M}\) is not guaranteed by the GLT theory, since these matrices are non-Hermitian. The numerical results reported in this chapter indicate that the eigenvalue distribution may not follow the GLT symbol in the current setting, even when \(a\) is chosen as a constant function. This contrasts with the spatial discretization case, where the numerical evidence was consistent with the symbol. The distinction is important: for non-Hermitian matrices, the GLT symbol provides a rigorous description of the singular value distribution, but the eigenvalue distribution requires additional analysis. In this work, we focus on preconditioning strategies that ensure weak clustering at one for the preconditioned matrix sequences, which is sufficient to guarantee fast convergence of the GMRES method without requiring a full characterization of the eigenvalue distribution of the original coefficient matrices.

\section{Preconditioning strategy: 1D space setting}\label{sec:prec}

Since the considered coefficient matrices are non-symmetric and exhibit a Toeplitz-like structure, Krylov subspace methods such as GMRES~\cite{Saad1986GMRES} are well-suited for solving the linear systems in \eqref{equ12}, benefiting from efficient fast matrix–vector products.
To further enhance the efficiency and robustness of these iterative solvers, the use of preconditioning techniques is an essential strategy \cite{Benzi2002Preconditioning}.
\par
Here, we propose a block circulant preconditioner for {$A_{N, M}$}, obtained by replacing the Toeplitz block $\overline{G}_{\alpha,N}$ with a circulant counterpart $C_{\alpha,N}$, while keeping $Q_{\theta,M}$ unchanged. {Since our constant term $\gamma = \Delta t / (\Delta x)^{\alpha}$ is sufficiently large}, we apply the circulant approximation only to $\overline{G}_{\alpha,N}$ and retain $Q_{\theta,M}$ as it is.
In order to preserve the circulant structure, instead of using the diagonal matrix $D_{N}(a)$, we replace it with a constant scalar multiplier given by the average of the diffusion coefficients
\[
d = \frac{1}{N} \sum_{i=1}^{N} d_i,
\]
where $d_i$ are the diagonal entries of $D_N(a)$. Of course the choice of the circulant counterpart is far from unique and several proposals are treated in three decades of relevant literature. From theoretical studies, we know that Strang-type approaches are better than Frobenius optimal approximations, especially for matching the small eigenvalues and especially when the generating functions are smooth. However, here we have lower Hessemberg structures and hence, together with the usual Strang approximation, we take the Strang-type matrix chosen taking into account the largest coefficients, i.e. those on the first column, except the last which is replaced by that on upper main diagonal.
{More precisely, we consider these two possible circulant approximations, where the first preconditioner is defined as
	\begin{equation}\label{precoS}
		{\widetilde{P}_{\alpha,N, M}}=\gamma  d \, \widetilde{C}_{\alpha,N} \otimes I_M + I_N \otimes Q_{\theta,M}
	\end{equation}
	where $\widetilde{C}_{\alpha,N}$ is the standard Strang circulant approximation, while the latter is defined as
	\begin{equation}\label{precog0}
		{\widehat{P}_{\alpha,N, M}}= \gamma d \, \widehat{C}_{\alpha,N} \otimes I_M + I_N \otimes Q_{\theta,M}
	\end{equation}
	where $\widehat{C}_{\alpha,N}$ is the circulant matrix generated by the vector $\left[g_1^{(\alpha)}, \ldots, g_{N-1}^{(\alpha)}, g_0^{(\alpha)}\right]$,
	that is by the first column of $\overline{G}_{\alpha,N}$ where the last entry $g_N^{(\alpha)}$ has been replaced by the coefficient $g_0^{(\alpha)}$, which is nonnegigible while $g_N^{(\alpha)}$ tends to zero as $N$ tends to infinity; see also the Riemann-Lebesgue lemma \cite{zyg-bible}.}
\par
It is well known that any circulant matrix $C_{N}$ can be diagonalized by the Fourier matrix, i.e.,
{\[
	C_{N} = F_N \Lambda_{N} F_N^{*},
	\]
}
where $F_N$ is the unitary matrix given by \eqref{fourier-matrix} and $\Lambda_{N}$ is a diagonal matrix containing the eigenvalues of $C_{N}$.
\begin{theorem}\label{thm:bd_P1}
	{Let {$P  \in\{\widetilde{P}_{\alpha,N,M},\widehat{P}_{\alpha,N,M}\}$} as defined  in~\eqref{precoS} and in~\eqref{precog0}, respectively}. Set $U = F_{N} \otimes I_{M}$.
	Then  {$P$} is unitarily  similar to a block diagonal matrix, namely
	\begin{equation}\label{permut}
		U^{*}PU = \operatorname{blkdiag}\bigl(\gamma d\lambda_{j}I_{M} + Q_{\theta ,M}\bigr)_{j=1}^{N},
	\end{equation}
	where $\lambda_{j}$, $j=1,\ldots,N$, denote the eigenvalues of the circulant block  {$\tilde{C}_{\alpha,N}$ or $\hat{C}_{\alpha,N}$, respectively}. Moreover, each diagonal block
	\begin{equation}\label{spectra}
		D_{j} = \gamma d\,\lambda_{j} I_{M} + Q_{\theta,M}, \ \ j=1,\ldots,N,
	\end{equation}
	is an invertible lower triangular Toeplitz matrix (ILTT).
\end{theorem}
\begin{proof}
	Substituting the diagonalization of  {the circulant matrix }into~\eqref{precoS}  {(or into~\eqref{precog0} equivalently)} gives
	\[
	P = \gamma d\,F_{N}\Lambda_{\alpha,N}F_{N}^{*}\otimes I_{M} + I_{N}\otimes Q_{\theta,M}.
	\]
	{Then, using the mixed-product property $(A\otimes B)(C\otimes D) = AC \otimes BD$, it holds that}
	\begin{align*}
		U^{*} P U &= (F_{N}^{*} \otimes I_{M})
		\left( \gamma d\,F_{N}\Lambda_{\alpha,N}F_{N}^{*} \otimes I_{M} + I_{N} \otimes Q_{\theta,M} \right)
		(F_{N} \otimes I_{M}) \\
		&=  \gamma d\,F_{N}^{*}F_{N}\Lambda_{\alpha,N}F_{N}^{*}F_{N} \otimes I_{M} + F_{N}^{*}F_{N} \otimes Q_{\theta,M}  \\
		&=  \gamma d\,\Lambda_{\alpha,N} \otimes I_{M} + I_{N} \otimes Q_{\theta,M}.
	\end{align*}
	%
	%
	Since $\Lambda_{\alpha,N}=\operatorname{diag}(\lambda_{1},\ldots,\lambda_{N})$, it follows that
	\[
	\Lambda_{\alpha,N}\otimes I_{M} = \operatorname{blkdiag}(\lambda_{1}I_{M},\ldots,\lambda_{N}I_{M}),
	\qquad
	I_{N}\otimes Q_{\theta,M} = \operatorname{blkdiag}(Q_{\theta,M},\ldots,Q_{\theta,M}).
	\]
	Hence their sum is block diagonal with blocks $\gamma  d\,\lambda_{j}I_{M}+Q_{\theta,M}$, $j=1,\ldots,N$. Finally, note that $Q_{\theta,M}$ is a lower triangular Toeplitz matrix, and adding a scalar multiple of the identity preserves both the Toeplitz and lower triangular structure. Thus, taking into account that $H_{\theta,M}$ is invertible, each block
	\[
	D_{j} = \gamma d\,\lambda_{j} I_{M} + Q_{\theta,M} = H_{\theta,M}^{-1}\left(\gamma d\,\lambda_{j}  H_{\theta,M} + H_{M}\right)
	\]
	is lower triangular Toeplitz. Now we can conclude that $D_{j}$ is nonsingular because $\gamma d\,\lambda_{j}  H_{\theta,M} + H_{M}$ is lower bidiagonal and its invertibility is decided by the diagonal elements, which are of the form $1-\gamma d\,\lambda_{j}\theta$ which can not vanish since the eigenvalues $\lambda_j$ have all nonzero imaginary part.
	With the latter the proof is concluded.
\end{proof}
\begin{remark}\label{rem:efficiency}
	The block diagonalization in Theorem~\ref{thm:bd_P1} shows that applying  {$P^{-1}$} to a given vector $v$ can be performed through the following three steps.
	\begin{enumerate}
		\item Step 1: compute $
		\dot{v} = (F_{N}^{*} \otimes I_{M}) v$.
		
		\item Step 2: compute
		$\ddot{v} = \operatorname{blkdiag}\big(D_{j}^{-1}\big)_{j=1}^{N}\,\dot{v},
		\qquad D_{j} = \gamma d\,\lambda_{j} I_{M} + Q_{\theta,M}$.
		
		\item Step 3: compute $
		\dddot{v} = (F_{N} \otimes I_{M}) \ddot{v}$.
	\end{enumerate}
	The resulting total cost is of the order $O(NM\log(N))$ which is optimal in our setting, since it is the same order of the cost of a matrix-vector product with the original coefficient matrix.

	From Theorem~\ref{thm:bd_P1}, both Steps~1 and~3 involve Fourier transforms in the spatial dimension and can be performed in $\mathcal{O}(NM \log N)$ operations using FFT.
	Step~2 requires solving $N$ independent systems with coefficient matrices $D_{j}$, which are ILTT matrices.
	Since the inverse of an ILTT matrix is also ILTT, its first column can be computed in $\mathcal{O}(M \log M)$ operations by the divide-and-conquer inversion algorithm~\cite{Commenges1984}.
	Therefore, Step~2 can be carried out in $\mathcal{O}(NM \log M)$ operations.
	In total, applying  {$P^{-1} v$} requires $
	\mathcal{O}(NM \log N + NM \log M)$,
	which is nearly linear in the problem size.
	Moreover, the $N$ systems in Step~2 are completely independent and can be solved in parallel, making the preconditioner highly efficient in practice.
	
	Of course, if the matrix $Q_{\theta,M}$ is not stored explicitly as in (\ref{Q mat}), then it can be viewed as  $H_{\theta,M}^{-1}H_{M}$ so that
	\[
	D_{j} = H_{\theta,M}^{-1}\left(\gamma d\,\lambda_{j}  H_{\theta,M} + H_{M}\right),
	\]
	which leads to a more efficient computation since  $H_{\theta,M}$ and $\gamma d\,\lambda_{j}  H_{\theta,M} + H_{M}$, $j=1,\ldots,N$, are all lower bidiagonal.
	With this choice, for solving $n$ linear systems with coefficient matrices $D_j$, $j=1,\ldots,N$, the related cost is reduced from $\mathcal{O}(NM \log M)$ to $\mathcal{O}(NM)$.
\end{remark}

We conclude the current section with a GLT/clustering analysis of the preconditioned matrix sequences.
\begin{theorem} \label{thm:inv(P)A clustering 1D}
	{Let {$P_{N,M}  \in\{\widetilde{P}_{\alpha,N,M},\widehat{P}_{\alpha,N,M}\}$}, as defined  in~\eqref{precoS} and in~\eqref{precog0}, respectively.
		Then,
		\[
		{\{P_{N,M}\}_{N,M}}\sim_{\mathrm{GLT}}\gamma^* \hat d f_{\alpha}(\xi_1)+q_{\theta}(\xi_3),
		\]
		with {$\hat d=\int_\Omega a(x)$} and $d=d_N$ converging to $\hat d$ as $N$ tends to infinity and
		\[
		{\{P_{N,M}^{-1} A_{\alpha,N,M}\}_{N,M}}\sim_{\mathrm{GLT}} \frac{\gamma^* a(x)f_{\alpha}(\xi_1)+q_{\theta}(\xi_3)}{\gamma^*\hat d f_{\alpha}(\xi_1)+q_{\theta}(\xi_3)}.
		\]
	}
\end{theorem}
\begin{proof}
	From the literature \cite{ChanNg1996,superopt-ESC}, it is known that $\{\overline{G}_{\alpha,N}-\widetilde{C}_{\alpha,N}\}_N\sim_\sigma 0$. The same clustering at zero is true for $\{\overline{G}_{\alpha,N}-\widehat{C}_{\alpha,N}\}_N$  since the approximation is even better (the only neglected coefficient is $g_N^{(\alpha)}$). In both cases the clustering is strong since the generating function is continuous and $2\pi$-periodic \cite{ChanNg1996,superopt-ESC}.
	Now by axiom GLT \ref{glt:symbols}, part 3, $\{\overline{G}_{\alpha,N}-{C}_{\alpha,N}\}_N\sim_\sigma 0$ if and only if $\{\overline{G}_{\alpha,N}-{C}_{\alpha,N}\}_N\sim_{\mathrm{GLT}} 0$, {${C}_{\alpha,N} \in \{\widetilde{C}_{\alpha,N}, \widehat{C}_{\alpha,N}\}$}. Since $\{\overline{G}_{\alpha,N}\}_N$ is also a GLT matrix sequence, by axiom GLT \ref{glt:algebra}, part 3, it is clear that both the approximating matrix sequences $\{\widetilde{C}_{\alpha,N}\}_N$,
	$\{\hat{C}_{\alpha,N}\}_N$ are both GLT matrix sequences with the same GLT symbol as $\{\overline{G}_{\alpha,N}\}_N$, as in Theorem \ref{1D th:} and Remark \ref{rem:extending1}.
	
	By following the very same steps as in the proof of Theorem \ref{1D th:}, Theorem \ref{1D orig mat} and of the reasoning in Remark \ref{rem:extending1}, we end up with the conclusion
	\[
	{\{P_{N,M}\}_{N,M}}\sim_{\mathrm{GLT}}\gamma^* \hat d f_{\alpha}(\xi_1)+q_{\theta}(\xi_3)
	\]
	with {$\hat d=\int_\Omega a(x)$} and $d=d_N$ converging to $\hat d$ as $N$ tends to infinity. As a consequence
	\[
	{\{P_{N,M}^{-1} A_{\alpha,N,M}\}_{N,M}}\sim_{\mathrm{GLT}} \frac{\gamma^* a(x)f_{\alpha}(\xi_1)+q_{\theta}(\xi_3)}{\gamma^*\hat d f_{\alpha}(\xi_1)+q_{\theta}(\xi_3)},
	\]
	by using axiom GLT \ref{glt:algebra}, part 4 and part 3.
	\par
	Notice that this GLT symbol is equal to $1$ in the constant coefficient case since $a=\hat d$, so ensuring a clustering at $1$ for the preconditioned matrix sequence.
\end{proof}

\section{Preconditioning strategy: 2D space setting}\label{sec:prec2D}

In the two-dimensional case, the spatial discretization yields a two-level Toeplitz structure, while the temporal discretization adds a further level so that a three-level matrix is obtained, in the all-at-once formulation. To construct an efficient preconditioner for this system, we approximate only the spatial part using a multilevel circulant matrix, while the temporal block $Q_{\theta,M}$ remains unchanged. The spatial operator {$G_{\boldsymbol{\alpha},\textbf{n}}$} is the diagonal sampling matrix {$D_\textbf{n}(a)$} multiplied with the two-level Toeplitz matrix
\[
{U_{\boldsymbol{\alpha},\textbf{n}}} {=    \gamma_1 \overline{G}_{\alpha_1,N_1}\otimes I_{N_2} + I_{N_1} \otimes \gamma_2 \overline{G}_{\alpha_2,N_2}}.
\]
To form the preconditioner, we replace {$U_{\boldsymbol{\alpha},\textbf{n}}$} with its two-level circulant approximation
\[
{C_{\boldsymbol{\alpha},\textbf{n}} = \gamma_1 C_{\alpha_1,N_1} \otimes I_{N_2} + I_{N_1} \otimes \gamma_2 C_{\alpha_2,N_2},
}\]
where {$C_{\alpha_r,N_r}$,  $r=1,2$} are the circulant matrices {already} considered in the one-dimensional case, {i.e.,
	$C_{\alpha_r,N_r} \in \{\widetilde{C}_{\alpha_r,N_r},\widehat{C}_{\alpha_r,N_r}\}$.
}
\par
{Since the constants
	$\gamma_1 = \Delta t/(\Delta x)^{\alpha_1}$ and
	$\gamma_2 = \Delta t/(\Delta y)^{\alpha_2}$
	are sufficiently large,} the circulant approximation is applied only to the spatial blocks, while the time stepping matrix $Q_{\theta,M}$ is left unchanged.
{As before, to} preserve the circulant structure, the diagonal sampling matrix {$D_\textbf{n}(a)$} is replaced by the spatial average of the diffusion coefficient
\[
d = \frac{1}{N_1 N_2} \sum_{{i_1}=1}^{N_1}\sum_{{i_2}=1}^{N_2} a{_{i_1,i_2}},
\]
and the preconditioner is defined as
\begin{equation}\label{eq:P2a}
	{\widetilde{P}_{\boldsymbol{\alpha},\textbf{n},M} = d \, \widetilde{C}_{\boldsymbol{\alpha},\textbf{n}} \otimes I_M + I_{\textbf{n}} \otimes Q_{\theta,M},
}\end{equation}
or
\begin{equation}\label{eq:P2b}
	{\widehat{P}_{\boldsymbol{\alpha},\textbf{n},M} = d \, \widehat{C}_{\boldsymbol{\alpha},\textbf{n}} \otimes I_M + I_{\textbf{n}} \otimes Q_{\theta,M},
}\end{equation}
{where, clearly, it holds
	\[
	C_{\boldsymbol{\alpha},\textbf{n}} = F_{\boldsymbol{\textbf{n}}}^* \  \Lambda_{\boldsymbol{\alpha},\textbf{n}} \ F_{\boldsymbol{\textbf{n}}}
	\]
	with $C_{\boldsymbol{\alpha},\textbf{n}} \in \{\widetilde{C}_{\boldsymbol{\alpha},\textbf{n}},\widehat{C}_{\boldsymbol{\alpha},\textbf{n}}\}$, $F_{\boldsymbol{\textbf{n}}} = F_{N_1} \otimes F_{N_2}$ and $\Lambda_{\boldsymbol{\alpha},\textbf{n}} = \gamma_1 \Lambda_{\alpha_1,N_1} \otimes I_{N_2} + I_{N_1} \otimes \gamma_2 \Lambda_{\alpha_2,N_2}$.
}%

Now we give only the statement of results regarding the invertibility of $P$, the computation cost related to the chosen preconditioner, and the  clustering results for the related preconditioned matrix sequence. The tools and the structure of the proofs follow verbatim those of Section \ref{sec:prec}.

\begin{theorem}
	Let {$P \in \{\widetilde{P}_{\boldsymbol{\alpha},\textbf{n},M},\widehat{P}_{\boldsymbol{\alpha},\textbf{n},M}\}$} be defined as in~\eqref{eq:P2a} and as in~\eqref{eq:P2b}, respectively. Set $U = F_{\textbf{n}}\otimes I_{M}$.
	Then $P$ is unitarily similar to a block diagonal matrix
	\[
	U^{*}PU = \operatorname{blkdiag}\bigl(d\lambda_{j}I_{M} + Q_{\theta ,M}\bigr)_{j=1}^{N_{1}N_{2}},
	\]
	where $\lambda_{j}$ are the eigenvalues of
	$C_{\boldsymbol{\alpha},\textbf{n}}$, $j=1,\ldots, N_{1}N_{2}$.
	Each diagonal block
	\[
	D_{j} = d\,\lambda_{j} I_{M} + Q_{\theta,M}, \ \ \ j=1,\ldots,N_{1}N_{2},
	\]
	is an invertible lower triangular Toeplitz matrix.
\end{theorem}
\begin{remark}
	The computation of a vector $P^{-1}v$ is carried out in three steps.
	\begin{enumerate}
		\item  Step 1: compute $\dot{v} = (F_{\boldsymbol{\textbf{n}}}^* \otimes I_M) v$.
		This requires $O(N_1 N_2 M \log (N_1 N_2))$ operations.
		\item Step 2: compute $\ddot{v} =
		\operatorname{blkdiag}\big(D_j^{-1}\big)_{j=1}^{N_1 N_2} \dot{v},$ \\
		where each system involves an ILTT matrix.
		The first column of each inverse can be computed in $O(M \log M)$ operations.
		The total cost is $O(N_1 N_2 M \log M)$.
		All systems are independent and can be solved in parallel.
		\item Step 3: compute
		$v = (F_{\boldsymbol{\textbf{n}}} \otimes I_M)\ddot{v}$,
		with cost $O(N_1 N_2M \log (N_1 N_2))$.
	\end{enumerate}
	Therefore, computing $P^{-1} v$ requires $O\!\left(N_1 N_2 M\log (N_1 N_2)
	+ N_1 N_2 M\log M\right),$
	which is nearly linear in the problem size.
	\par
	Furthermore, in line with the arguments in Remark \ref{rem:efficiency}, for solving $N_1 N_2$ linear systems with coefficient matrices $D_j$, $j=1,\ldots,N_1 N_2$, the related cost is reduced from $\mathcal{O}(N_1 N_2 M \log M)$ to $\mathcal{O}(N_1 N_2 M)$.
\end{remark}
\begin{theorem} \label{thm:inv(P)A clustering 2D}
	{Let {$P  \in\{\widetilde{P}_{\boldsymbol{\alpha},\textbf{n},M},
			\widehat{P}_{\boldsymbol{\alpha},\textbf{n},M}\}$}, as defined  in~\eqref{eq:P2a} and in~\eqref{eq:P2b}, respectively.
		Then,
		\[
		{\{P_{\textbf{n},M}\}_{\textbf{n},M}}\sim_{\mathrm{GLT}}
		\hat d (\gamma_1^*f_{\alpha_1}(\xi_1)+\gamma_2^*f_{\alpha_2}(\xi_2))
		+q_{\theta}(\xi_3),
		\]
		with {$\hat d=\int_\Omega a(x_1,x_2)$} and {$d=d_\textbf{n}$} converging to $\hat d$ as $\textbf{n}$ tends to infinity, and
		\[
		{\{P_{\textbf{n},M}^{-1} A_{\textbf{n},M}\}_{\textbf{n},M}}\sim_{\mathrm{GLT}} \frac{\gamma^* a(x_1,x_2)(\gamma_1^*f_{\alpha_1}(\xi_1)+\gamma_2^*f_{\alpha_2}(\xi_2))+q_{\theta}(\xi_3)}{\gamma^* \hat{d} (\gamma_1^*f_{\alpha_1}(\xi_1)+\gamma_2^*f_{\alpha_2}(\xi_2))+q_{\theta}(\xi_3)},
		\]
		where necessarily $\gamma_2^* = 0$ if $\alpha_1 > \alpha_2$ and alternatively $\gamma_1^* = 0$ if $\alpha_2 > \alpha_1$.
	}
\end{theorem}

\section{Numerical experiments}\label{sec:num}
We conclude our analysis by numerically evaluating the performance of the proposed preconditioning strategy.
The objective is to examine its capability to cluster the eigenvalues around~1, reduce the condition number,
and enhance the efficiency of the GMRES method in terms of iteration counts and CPU time.

All numerical experiments were conducted in \textsc{MATLAB}~R2022b on a laptop equipped with an 11th~Gen
Intel(R)~Core(TM)~i5-1155G7 CPU @~2.50\,GHz, 16\,GB RAM, running Windows~11~Pro
(version~23H2, build~22631.5189).
\subsection{Eigenvalue distribution}
\begin{figure}%
	\centering
	\includegraphics[width=\textwidth]{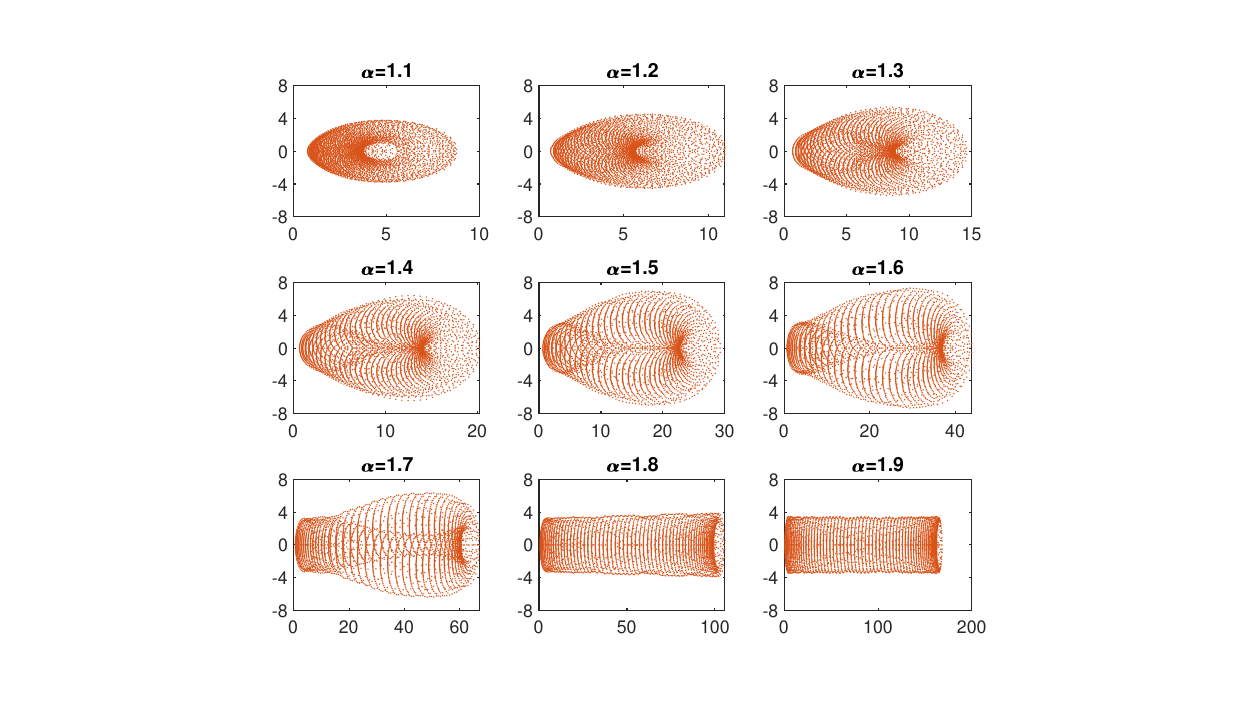}
	\caption{Eigenvalues of $A_{\alpha,N,M}$  - $N=M=2^6$, $a(x)=1$.} \label{Fig:Aes1}
\end{figure}

We first assess the performance of the preconditioners {$\widetilde{P}_{\alpha,N,M}$ and $\widehat{P}_{\alpha,N,M}$} by examining their influence on the eigenvalue distribution and condition number of the coefficient matrix.
Figure~\ref{Fig:Aes1} presents the eigenvalue plots in the complex plane for different values of the fractional order $\alpha$, with the problem size fixed at {$N=M = 2^6$}.
The eigenvalues of the unpreconditioned matrix are widely scattered along the positive real axis, indicating poor spectral clustering and, consequently, a high level of ill-conditioning. {The behaviour is worsening as the fractional diffusion parameter increases as expected.}
%

As shown in Table~\ref{tab:cond-numbers}, the unpreconditioned matrix exhibits rapidly growing 2-norm condition numbers, even for moderate values of $N$ and $M$. Though the condition numbers of the preconditioned matrix remain bounded and stable, it's not so lower.
Indeed, a more refined spectral analysis highlights that the condition number alone does not provide an exhaustive description of the eigenvalues distribution properties. \par
\begin{figure}%
	\centering
	\includegraphics[width=\textwidth]{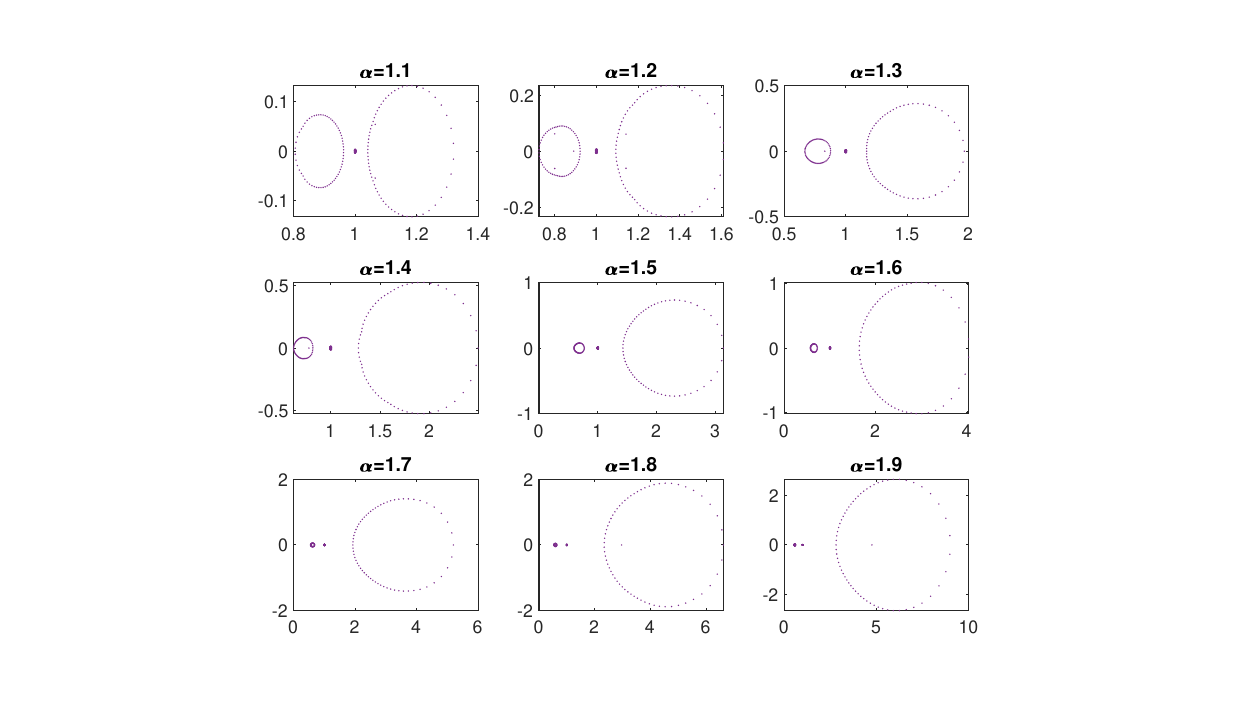}
	\includegraphics[width=\textwidth]{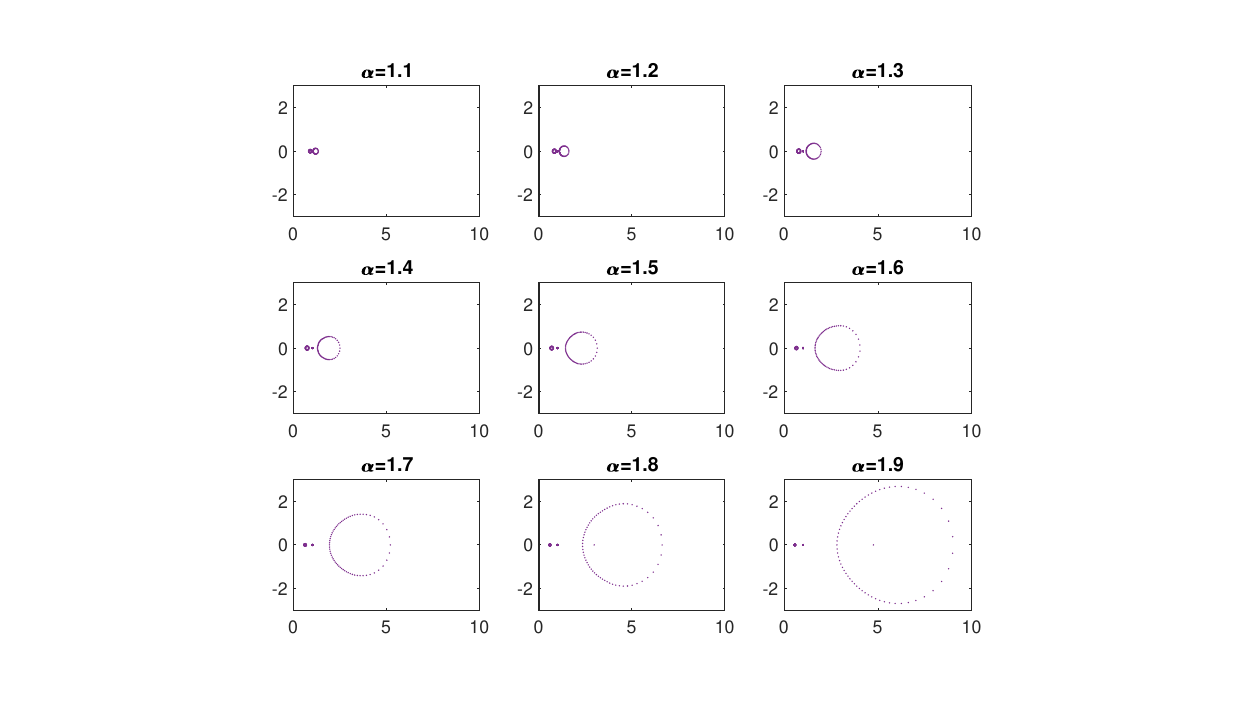}
	\caption{Eigenvalues of $\widetilde{P}^{-1}_{\alpha,N,M}A_{\alpha,N,M}$ - $N=M=2^6$, $a(x)=1$.} \label{Fig:PcircStranges1}
\end{figure}
\begin{figure}%
	\centering
	\includegraphics[width=\textwidth]{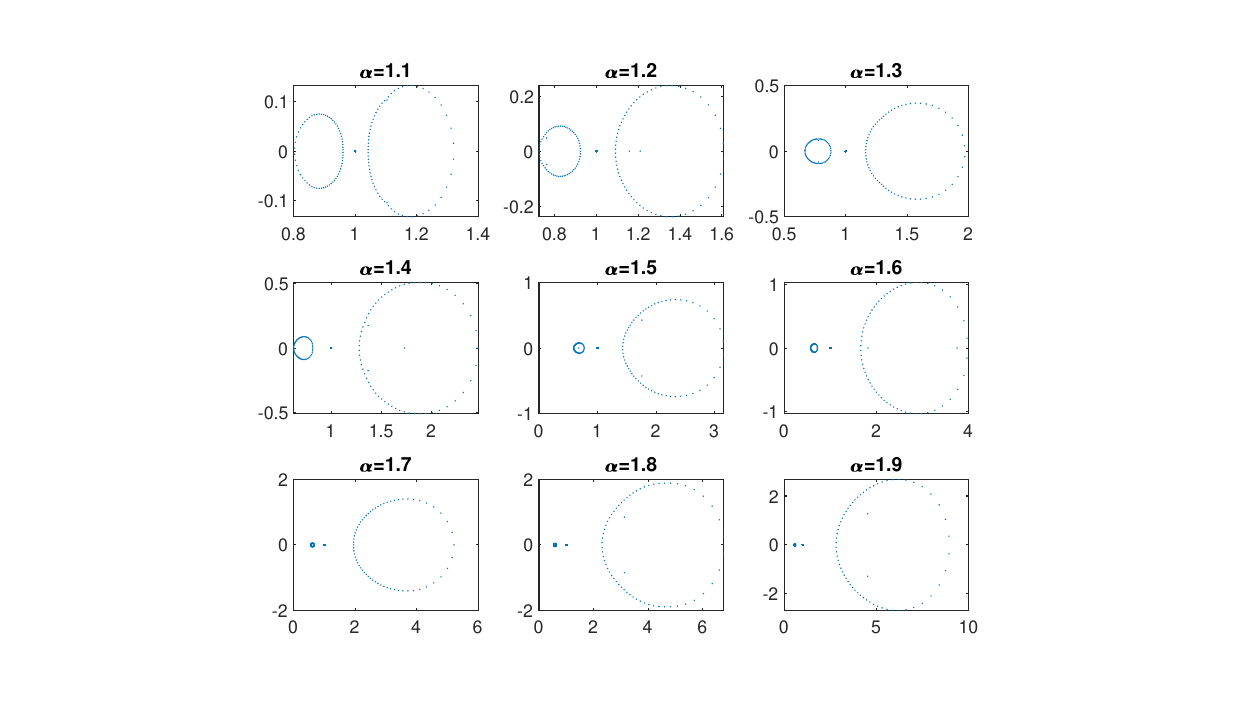}
	\includegraphics[width=\textwidth]{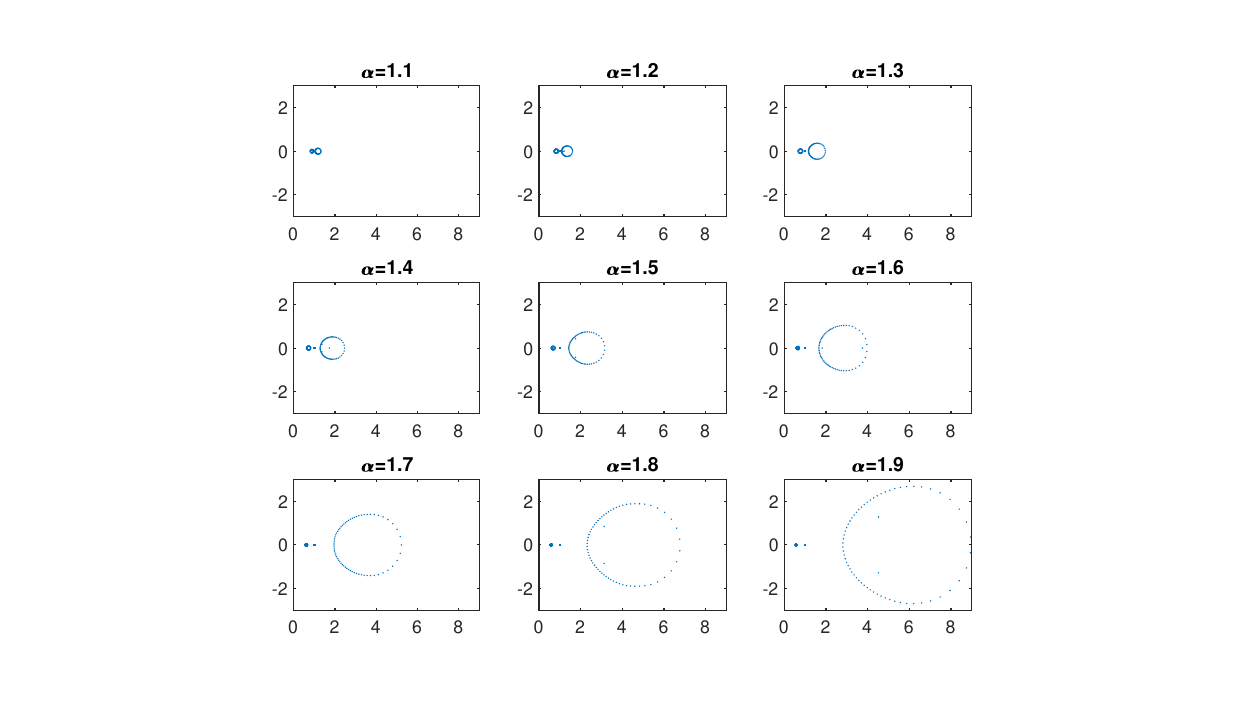}
	\caption{Eigenvalues of $\widehat{P}^{-1}_{\alpha,N,M}A_{\alpha,N,M}$  - $N=M=2^6$, $a(x)=1$.} \label{Fig:Pcircg0es1}
\end{figure}

In Figures~\ref{Fig:PcircStranges1}-\ref{Fig:Pcircg0es1} are reported the eigenvalues of preconditioned matrices {$\widetilde{P}_{\alpha,N,M}^{-1}A_{\alpha,N,M}$ and $\widehat{P}_{\alpha,N,M}^{-1}A_{\alpha,N,M}$}
for different values of the fractional order $\alpha$, with the problem size fixed at {$N=M = 2^6$.
	It evident as we can observe a pronounced clustered distributions around the unity.
	The same behaviour is observed in the non constant coefficient case: the related matrix sequences not cluster at $1$ but at the range of
	the nice function reported in Theorem \ref{thm:inv(P)A clustering 1D}.}%
\begin{table}
	\centering
	\renewcommand{\arraystretch}{1}
	\setlength{\tabcolsep}{3pt}
	\begin{adjustbox}{width=\textwidth}
		\begin{tabular}{|c c |ccc|ccc|ccc|}
			\hline
			\multicolumn {2}{|c}{} & \multicolumn {3}{|c}{$\alpha = 1.3$} & \multicolumn {3}{|c}{$\alpha = 1.5$} & \multicolumn {3}{|c|}{$\alpha = 1.8$} \\
			\hline
			$N$ & $M$ & $A$ & $\widetilde{P}^{-1}A$ & $\widehat{P}^{-1}A$ & $A$ & $\widetilde{P}^{-1}A$ & $\widehat{P}^{-1}A$ & $A$ & $\widetilde{P}^{-1}A$ & $\widehat{P}^{-1}A$ \\
			\hline
			$2^4$ &$2^4$ &2.013e+02 &1.727e+01 &1.953e+01 &1.634e+02 &2.662e+01 &3.238e+01 &1.483e+02 &6.765e+01 &8.167e+01 \\
			$2^4$ &$2^5$ &7.475e+02 &1.728e+01 &1.955e+01 &5.604e+02 &2.664e+01 &3.241e+01 &3.957e+02 &6.770e+01 &8.173e+01 \\
			$2^4$ &$2^6$ &2.937e+03 &1.729e+01 &1.955e+01 &2.160e+03 &2.664e+01 &3.241e+01 &1.417e+03 &6.771e+01 &8.174e+01 \\
			$2^4$ &$2^7$ &1.170e+04 &1.729e+01 &1.955e+01 &8.562e+03 &2.665e+01 &3.241e+01 &5.520e+03 &6.772e+01 &8.175e+01 \\
			\hline
			$2^5$ &$2^4$ &2.352e+02 &3.878e+01 &4.367e+01 &2.273e+02 &6.930e+01 &8.360e+01 &3.173e+02 &2.190e+02 &2.604e+02 \\
			$2^5$ &$2^5$ &7.854e+02 &3.882e+01 &4.372e+01 &6.177e+02 &6.937e+01 &8.367e+01 &5.477e+02 &2.191e+02 &2.606e+02 \\
			$2^5$ &$2^6$ &3.007e+03 &3.883e+01 &4.373e+01 &2.226e+03 &6.938e+01 &8.369e+01 &1.548e+03 &2.192e+02 &2.607e+02 \\
			$2^5$ &$2^7$ &1.190e+04 &3.883e+01 &4.373e+01 &8.684e+03 &6.939e+01 &8.370e+01 &5.650e+03 &2.192e+02 &2.607e+02 \\
			\hline
			$2^6$ &$2^4$ &3.221e+02 &9.102e+01 &1.023e+02 &4.187e+02 &1.878e+02 &2.256e+02 &9.089e+02 &7.328e+02 &8.655e+02 \\
			$2^6$ &$2^5$ &8.632e+02 &9.113e+01 &1.024e+02 &7.886e+02 &1.880e+02 &2.258e+02 &1.123e+03 &7.334e+02 &8.661e+02 \\
			$2^6$ &$2^6$ &3.096e+03 &9.115e+01 &1.024e+02 &2.379e+03 &1.880e+02 &2.258e+02 &2.057e+03 &7.335e+02 &8.663e+02 \\
			$2^6$ &$2^7$ &1.206e+04 &9.116e+01 &1.024e+02 &8.859e+03 &1.880e+02 &2.258e+02 &6.088e+03 &7.336e+02 &8.664e+02 \\
			\hline
			$2^7$ &$2^4$ &5.478e+02 &2.189e+02 &2.457e+02 &9.705e+02 &5.199e+02 &6.231e+02 &2.959e+03 &2.495e+03 &2.936e+03 \\
			$2^7$ &$2^5$ &1.064e+03 &2.191e+02 &2.460e+02 &1.317e+03 &5.204e+02 &6.237e+02 &3.166e+03 &2.497e+03 &2.938e+03 \\
			$2^7$ &$2^6$ &3.276e+03 &2.192e+02 &2.461e+02 &2.836e+03 &5.205e+02 &6.238e+02 &4.029e+03 &2.497e+03 &2.939e+03 \\
			$2^7$ &$2^7$ &1.227e+04 &2.192e+02 &2.461e+02 &9.270e+03 &5.205e+02 &6.239e+02 &7.811e+03 &2.497e+03 &2.939e+03 \\
			\hline
		\end{tabular}
	\end{adjustbox}
	\caption{Spectral condition numbers of unpreconditioned $A$ and preconditioned matrices $\widetilde{P}^{-1}A =\widetilde{P}^{-1}_{\alpha,N,M}A_{\alpha,N,M}$ and $\widehat{P}^{-1}A = \widehat{P}^{-1}_{\alpha,N,M}A_{\alpha,N,M}$ - case $a(x)=1$.}
	\label{tab:cond-numbers}
\end{table}
%
%
%
%
\begin{table}
	\centering
	\renewcommand{\arraystretch}{1}
	\setlength{\tabcolsep}{3pt}
	\begin{adjustbox}{width=\textwidth}
		\begin{tabular}{|c c |ccc|ccc|ccc|}
			\hline
			\multicolumn {2}{|c}{} & \multicolumn {3}{|c}{$\alpha = 1.3$} & \multicolumn {3}{|c}{$\alpha = 1.5$} & \multicolumn {3}{|c|}{$\alpha = 1.8$} \\
			\hline
			$N$ & $M$ & $A$ & $\widetilde{P}^{-1}A$ & $\widehat{P}^{-1}A$ & $A$ & $\widetilde{P}^{-1}A$ & $\widehat{P}^{-1}A$ & $A$ & $\widetilde{P}^{-1}A$ & $\widehat{P}^{-1}A$ \\
			\hline
			$2^4$ &$2^4$ &1.678e+02 &1.535e+01 &1.805e+01 &1.465e+02 &2.556e+01 &3.285e+01 &1.549e+02 &8.477e+01 &1.080e+02 \\
			$2^4$ &$2^5$ &5.939e+02 &1.536e+01 &1.806e+01 &4.582e+02 &2.558e+01 &3.288e+01 &3.493e+02 &8.484e+01 &1.081e+02 \\
			$2^4$ &$2^6$ &2.307e+03 &1.536e+01 &1.807e+01 &1.724e+03 &2.558e+01 &3.288e+01 &1.168e+03 &8.485e+01 &1.081e+02 \\
			$2^4$ &$2^7$ &9.162e+03 &1.536e+01 &1.807e+01 &6.796e+03 &2.558e+01 &3.288e+01 &4.472e+03 &8.486e+01 &1.081e+02 \\
			\hline
			$2^5$ &$2^4$ &2.168e+02 &3.207e+01 &3.756e+01 &2.417e+02 &6.240e+01 &7.951e+01 &4.020e+02 &2.592e+02 &3.245e+02 \\
			$2^5$ &$2^5$ &6.392e+02 &3.209e+01 &3.759e+01 &5.411e+02 &6.245e+01 &7.957e+01 &5.805e+02 &2.594e+02 &3.247e+02 \\
			$2^5$ &$2^6$ &2.365e+03 &3.210e+01 &3.760e+01 &1.803e+03 &6.246e+01 &7.959e+01 &1.365e+03 &2.594e+02 &3.248e+02 \\
			$2^5$ &$2^7$ &9.286e+03 &3.210e+01 &3.760e+01 &6.901e+03 &6.246e+01 &7.959e+01 &4.650e+03 &2.594e+02 &3.248e+02 \\
			\hline
			$2^6$ &$2^4$ &3.499e+02 &7.158e+01 &8.368e+01 &5.309e+02 &1.614e+02 &2.048e+02 &1.293e+03 &8.164e+02 &1.013e+03 \\
			$2^6$ &$2^5$ &7.559e+02 &7.164e+01 &8.375e+01 &8.100e+02 &1.615e+02 &2.049e+02 &1.461e+03 &8.171e+02 &1.014e+03 \\
			$2^6$ &$2^6$ &2.475e+03 &7.166e+01 &8.377e+01 &2.034e+03 &1.616e+02 &2.050e+02 &2.181e+03 &8.173e+02 &1.014e+03 \\
			$2^6$ &$2^7$ &9.428e+03 &7.166e+01 &8.378e+01 &7.123e+03 &1.616e+02 &2.050e+02 &5.343e+03 &8.173e+02 &1.014e+03 \\
			\hline
			$2^7$ &$2^4$ &6.949e+02 &1.670e+02 &1.950e+02 &1.372e+03 &4.339e+02 &5.492e+02 &4.452e+03 &2.639e+03 &3.262e+03 \\
			$2^7$ &$2^5$ &1.077e+03 &1.671e+02 &1.952e+02 &1.636e+03 &4.343e+02 &5.496e+02 &4.617e+03 &2.641e+03 &3.265e+03 \\
			$2^7$ &$2^6$ &2.750e+03 &1.672e+02 &1.952e+02 &2.781e+03 &4.344e+02 &5.498e+02 &5.293e+03 &2.642e+03 &3.265e+03 \\
			$2^7$ &$2^7$ &9.692e+03 &1.672e+02 &1.953e+02 &7.759e+03 &4.344e+02 &5.498e+02 &8.201e+03 &2.642e+03 &3.265e+03 \\
			\hline
		\end{tabular}
	\end{adjustbox}
	\caption{Spectral condition numbers of unpreconditioned $A$ and preconditioned matrices $\widetilde{P}^{-1}A =\widetilde{P}^{-1}_{\alpha,N,M}A_{\alpha,N,M}$ and $\widehat{P}^{-1}A = \widehat{P}^{-1}_{\alpha,N,M}A_{\alpha,N,M}$ - case $a(x)=x^2+1$.}
	\label{tab:cond-numbers-es2}
\end{table}
%
\begin{figure}%
	\centering
	\includegraphics[width=\textwidth]{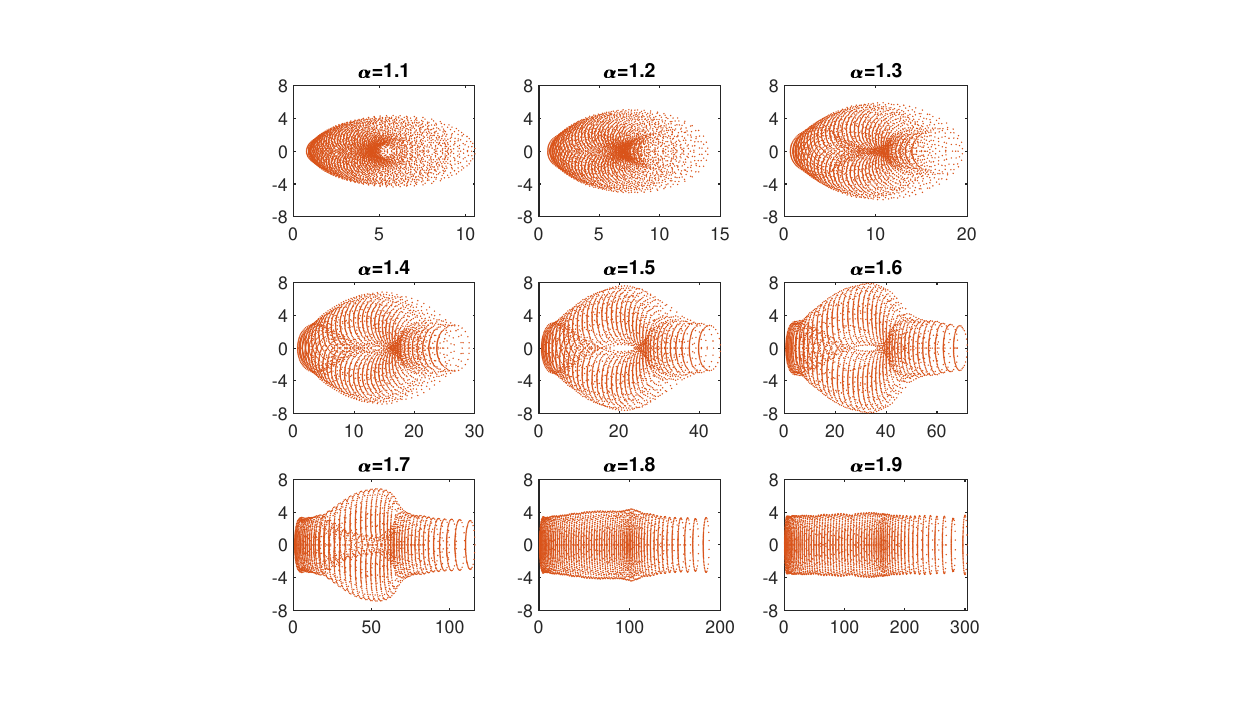}
	\caption{Eigenvalues of $A_{\alpha,N,M}$  - $N=M=2^6$, $a(x)=x^2+1$.} \label{Fig:Aes2}
\end{figure}
\begin{figure}%
	\centering
	\includegraphics[width=\textwidth]{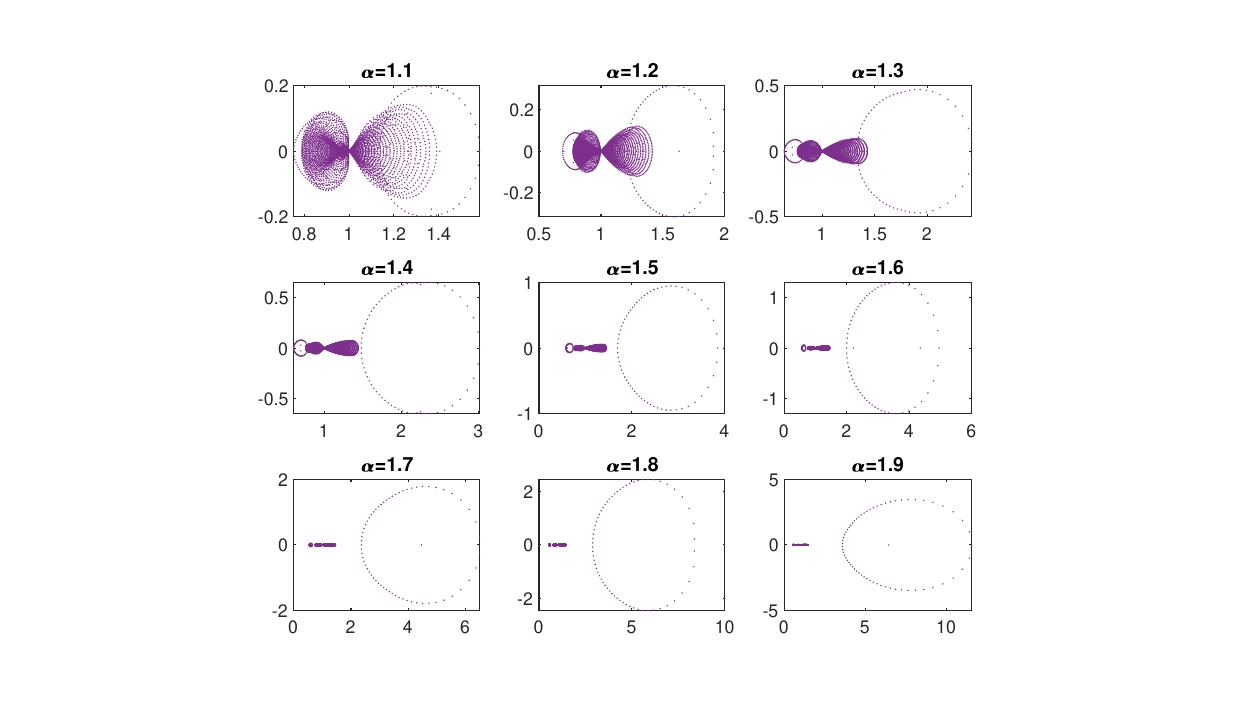}
	\includegraphics[width=\textwidth]{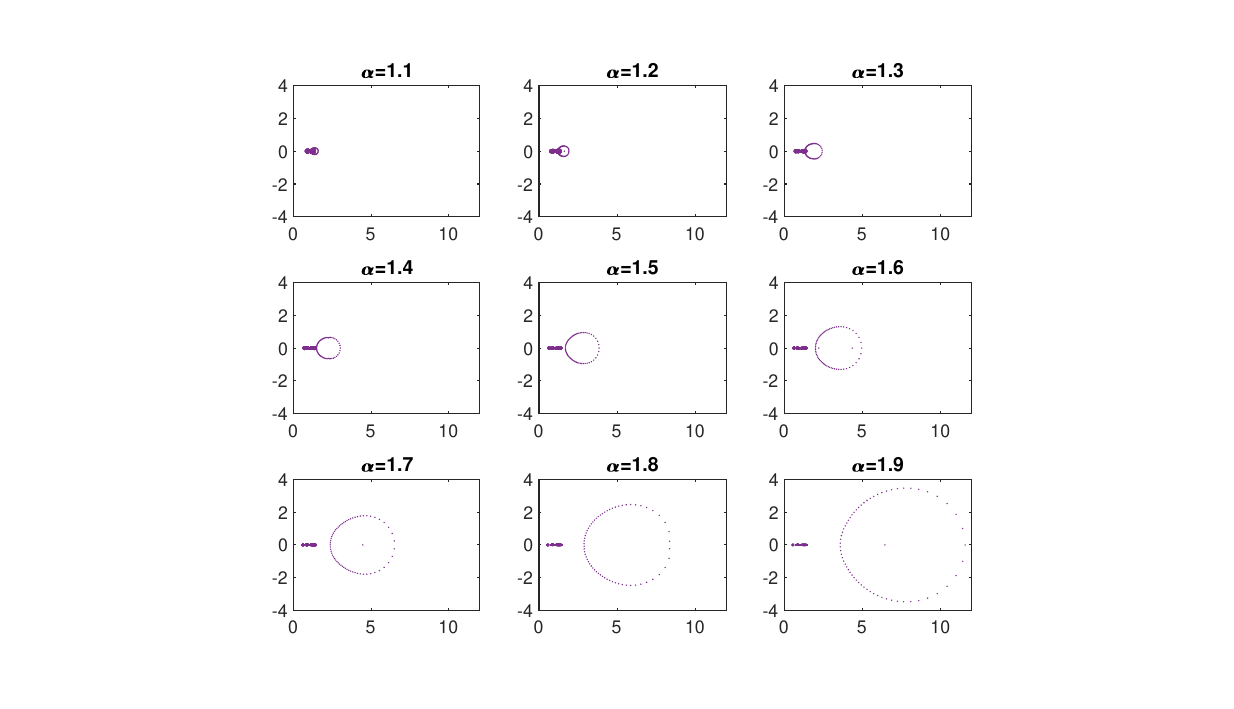}
	\caption{Eigenvalues of $\widetilde{P}^{-1}_{\alpha,N,M}A_{\alpha,N,M}$  - $N=M=2^6$, $a(x)=x^2+1$.} \label{Fig:PcircStranges2}
\end{figure}
\begin{figure}%
	\centering
	\includegraphics[width=\textwidth]{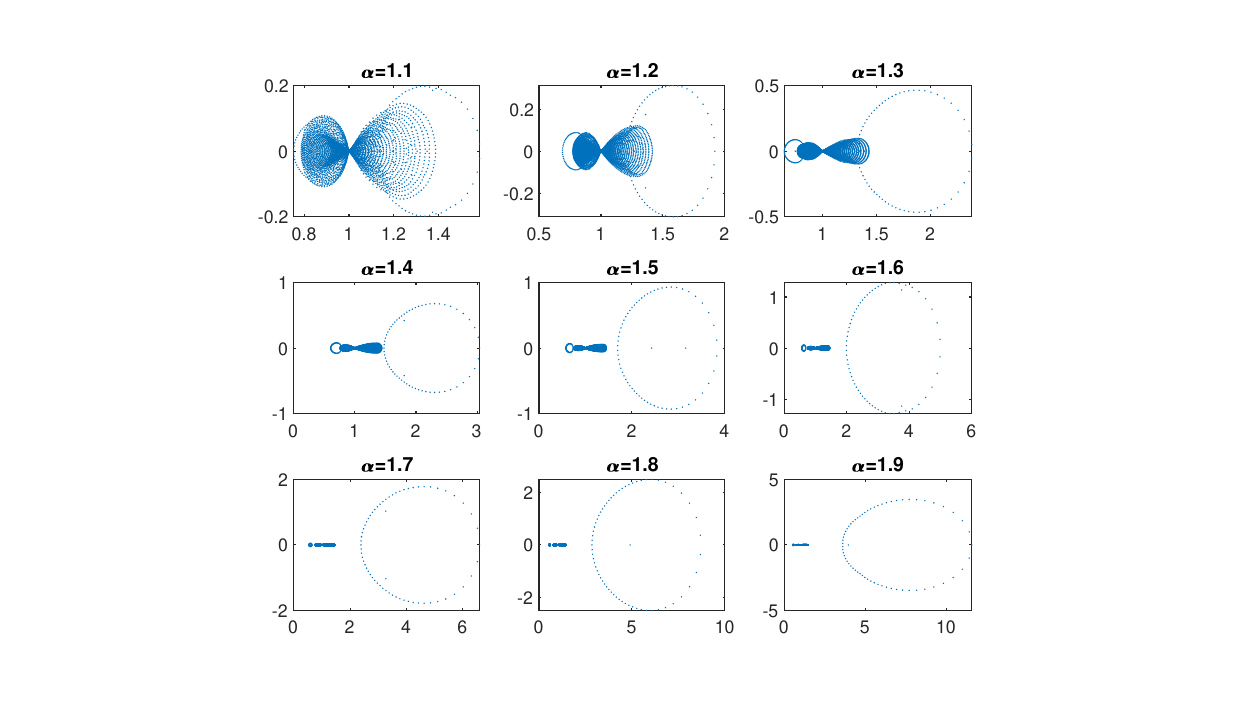}
	\includegraphics[width=\textwidth]{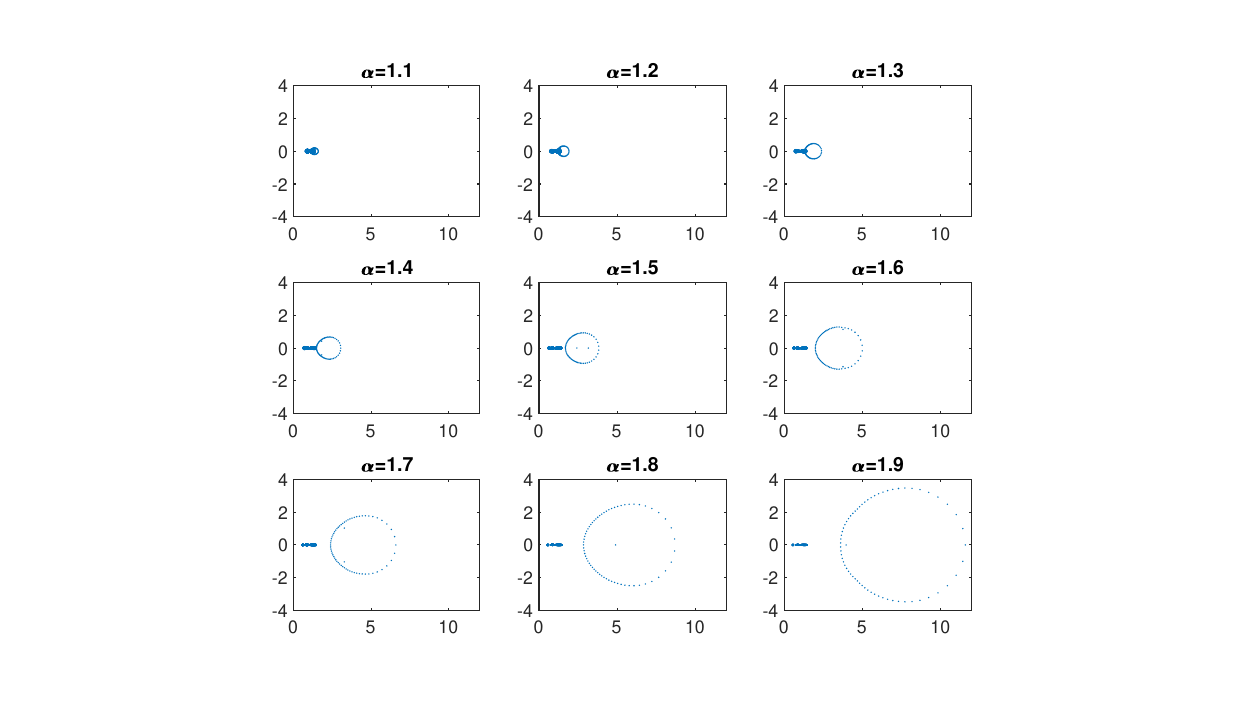}
	\caption{Eigenvalues of $\widehat{P}^{-1}_{\alpha,N,M}A_{\alpha,N,M}$ - $N=M=2^6$, $a(x)=x^2+1$.} \label{Fig:Pcircg0es2}
\end{figure}
%
\begin{figure}%
	\centering
	\includegraphics[width=\textwidth]{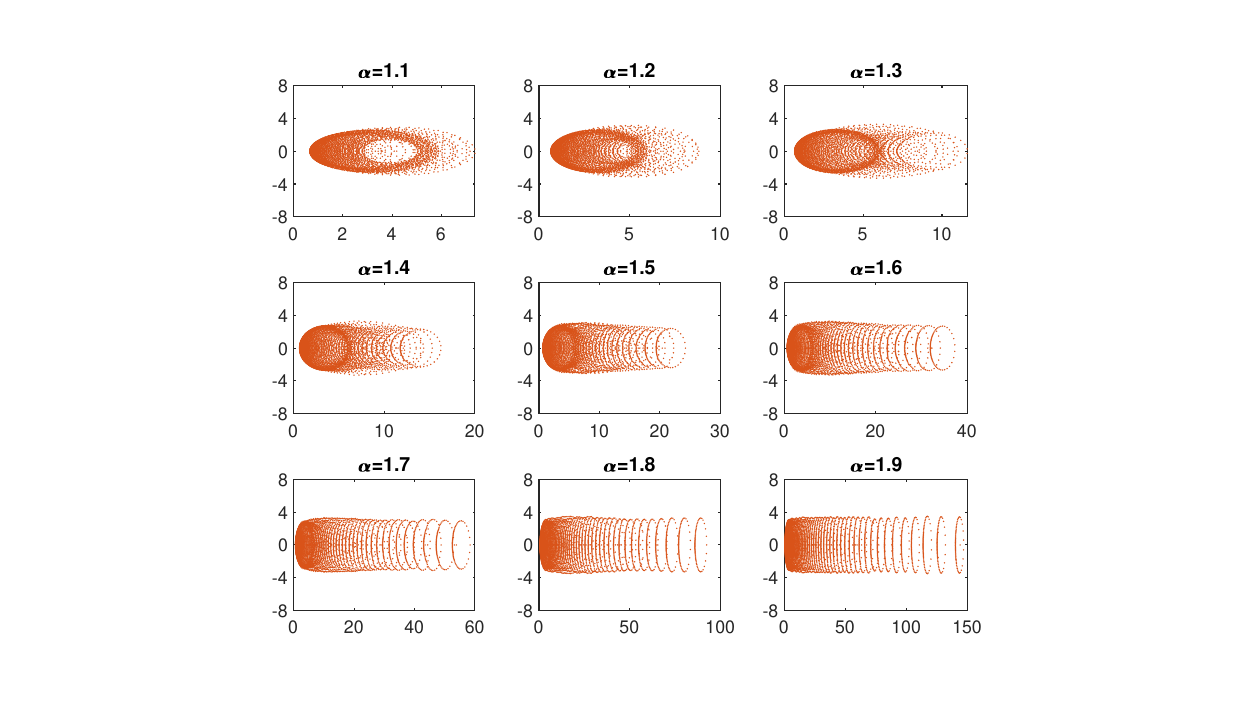}
	\caption{Eigenvalues of $A_{\alpha,N,M}$  - $N=M=2^6$, $a(x)=x^2$.} \label{Fig:Aes4}
\end{figure}
\begin{figure}%
	\centering
	\includegraphics[width=\textwidth]{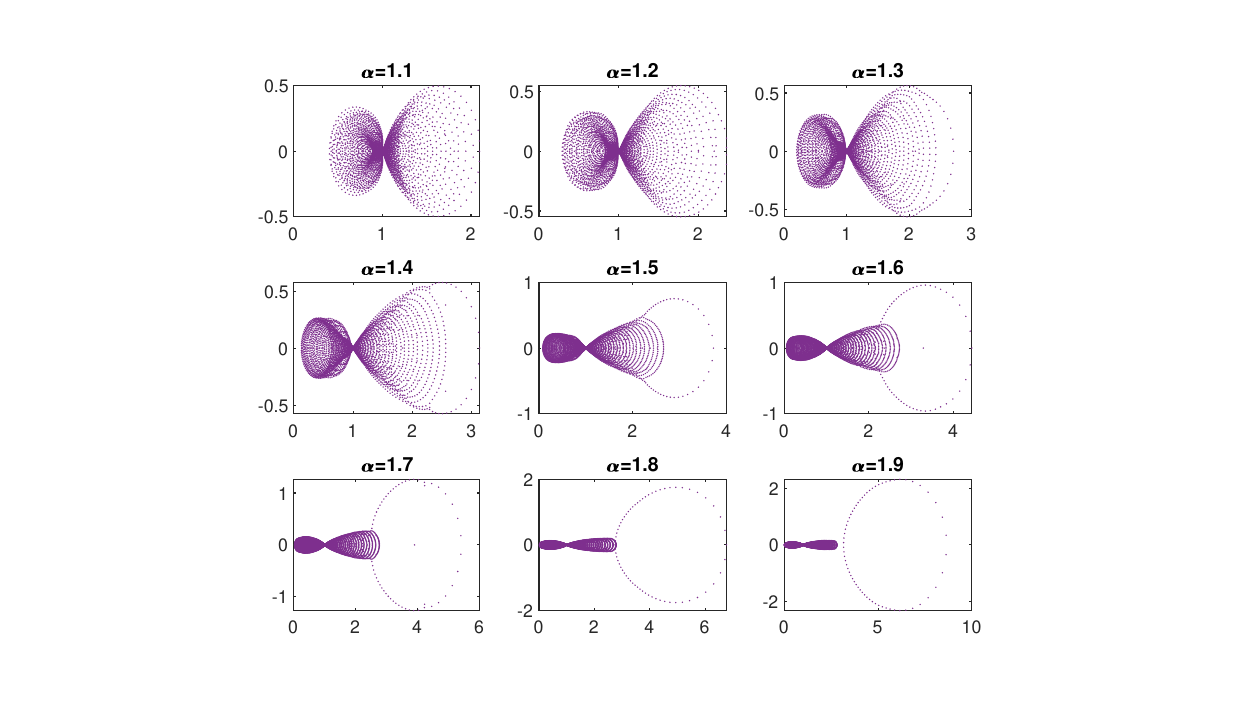}
	\includegraphics[width=\textwidth]{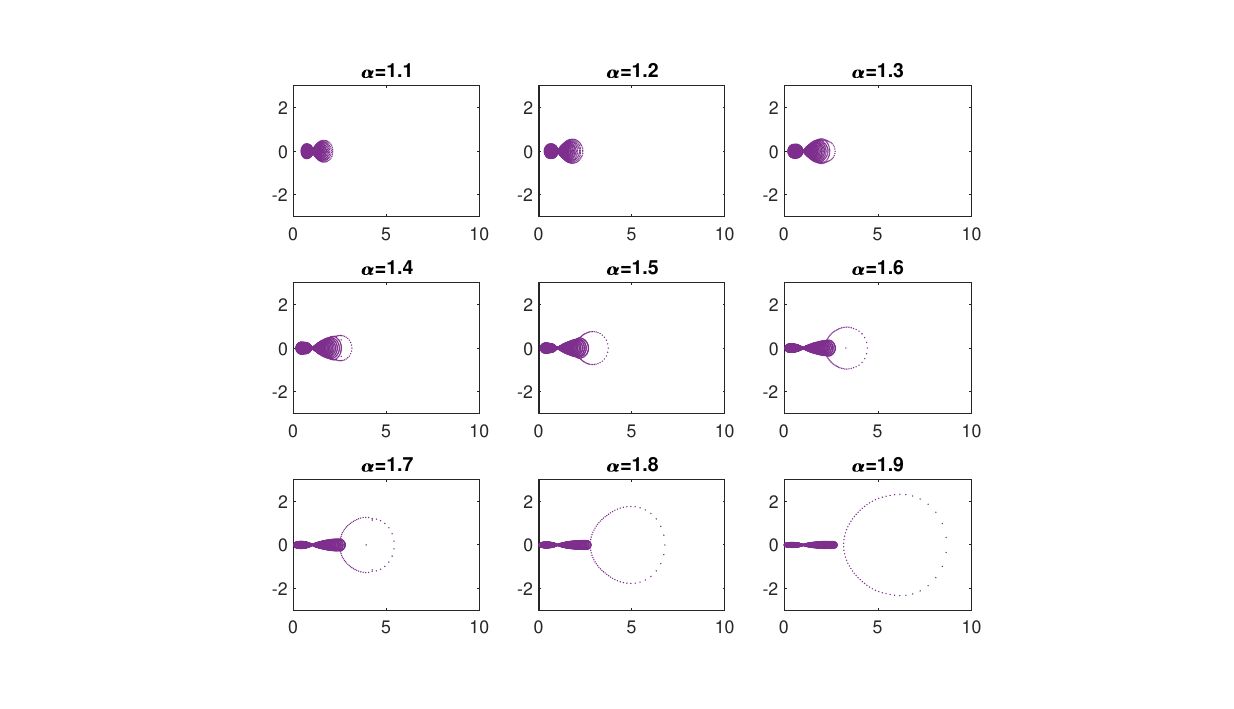}
	\caption{Eigenvalues of $\widetilde{P}^{-1}_{\alpha,N,M}A_{\alpha,N,M}$ - $N=M=2^6$, $a(x)=x^2$.} \label{Fig:PcircStranges4}
\end{figure}
\begin{figure}%
	\centering
	\includegraphics[width=\textwidth]{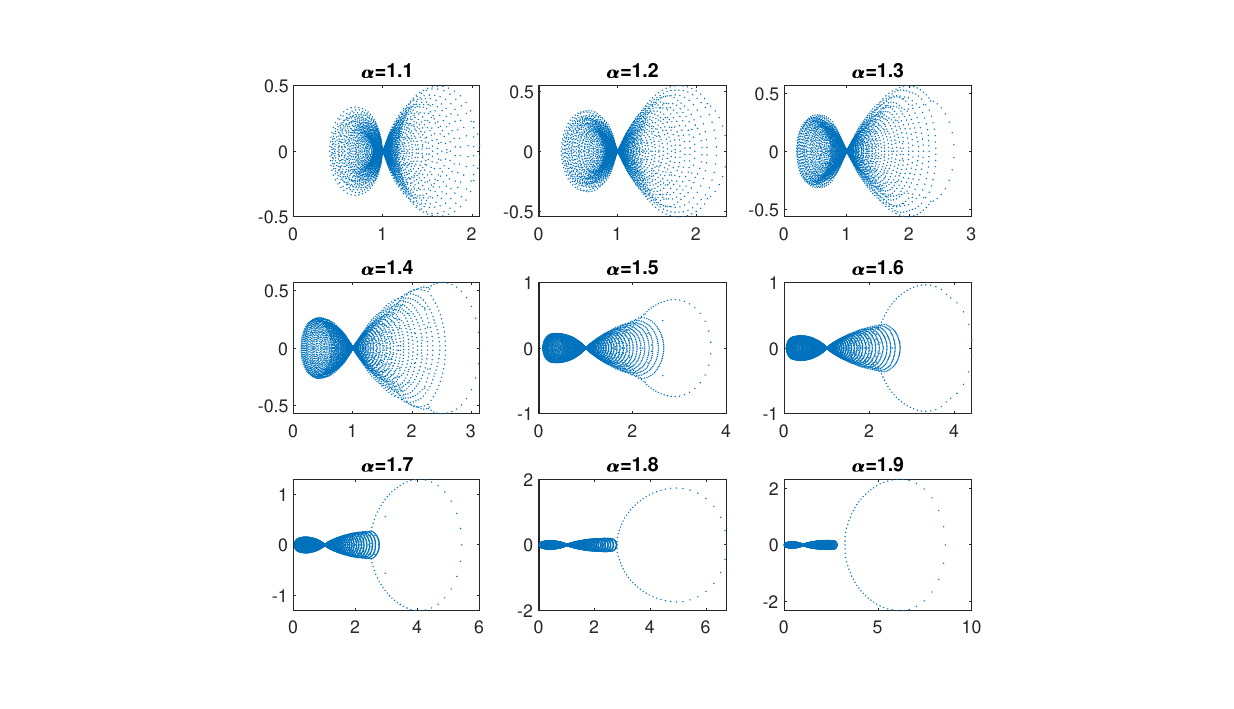}
	\includegraphics[width=\textwidth]{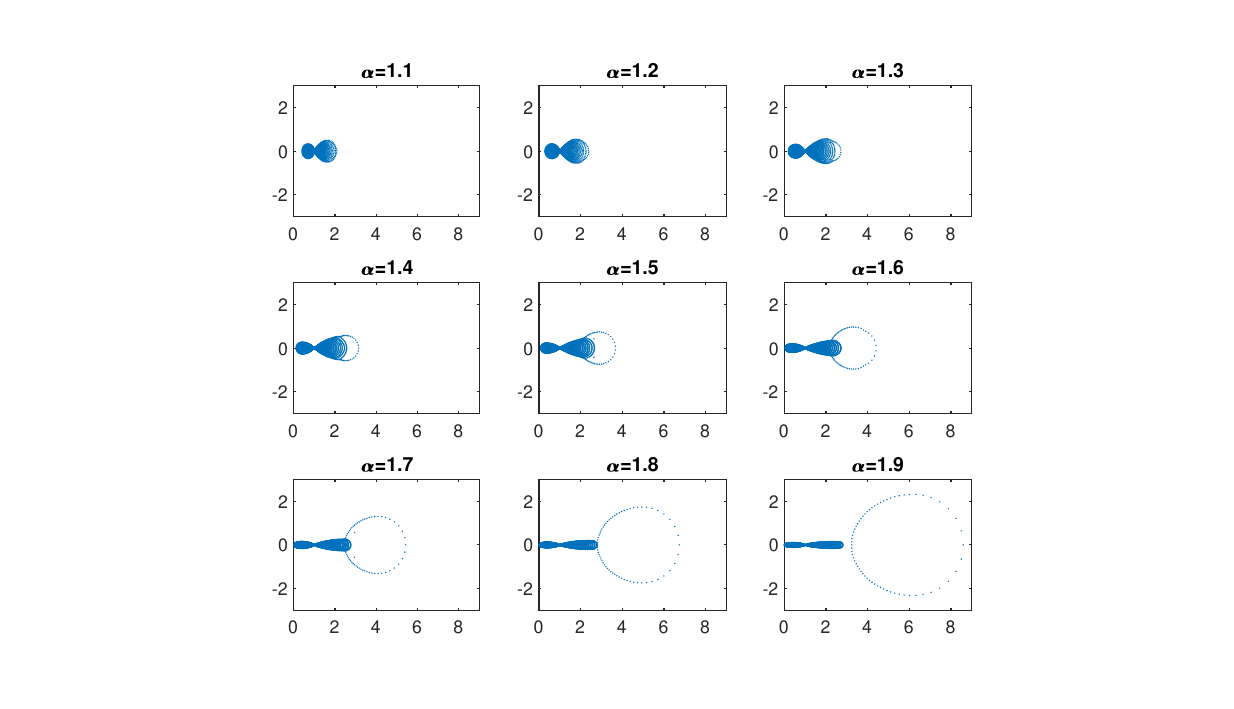}
	\caption{Eigenvalues of $\widehat{P}^{-1}_{\alpha,N,M}A_{\alpha,N,M}$  - $N=M=2^6$, $a(x)=x^2$.} \label{Fig:Pcircg0es4}
\end{figure}
\clearpage
{\subsection{GMRES performance in the $d$-dimensional space setting, $d=1,2$}}
We consider the one-dimensional variable-coefficient space-fractional diffusion equation (\ref{eq2-1d}):
\begin{equation*}
	\frac{\partial u(x,t)}{\partial t} = x^2\, D_{x}^{\alpha} u(x,t) + f(x,t), \quad x \in (0, 1),\, t \in (0, 1],
\end{equation*}
with homogeneous Dirichlet boundary conditions $u(0, t) = u(1, t) = 0$ and initial condition $u(x, 0) = 0$.

We choose the exact solution as
\begin{equation*}
	u(x, t) = t^2 x^2 (1-x)^2.
\end{equation*}
The corresponding source term is given by
\begin{align*}
	f(x, t) =\; & 2t\, x^2 (1-x)^2  - x^2\, t^2 \left[
	\frac{2}{\Gamma(3-\alpha)} x^{2-\alpha}
	- \frac{12}{\Gamma(4-\alpha)} x^{3-\alpha}
	+ \frac{24}{\Gamma(5-\alpha)} x^{4-\alpha}
	\right].
\end{align*}
Thus, in this section, we investigate the practical performance of the GMRES iterative solver for the all-at-once linear system arising from the one-dimensional discretization given in equation~\eqref{equ12}, with the time-stepping parameter set to $\theta = 0.5$ (Crank–Nicolson method). We compare the behavior of GMRES applied to both the original (non-preconditioned) system and the system preconditioned by a circulant structure-preserving preconditioners as defined in Section \ref{sec:prec}, to illustrate the impact of preconditioning on convergence and computational efficiency. First, the experiments are conducted for the one-dimensional constant-coefficient fractional diffusion problem ($a(x)=1$), and then for the one-dimensional variable-coefficient fractional diffusion problem
with coefficient $a(x) = x^2+1$, and $a(x) = x^2$.
The initial vector is set to zero, the tolerance is $10^{-8}$, and the maximum number of iterations equals the matrix dimension. Tables~\ref{tab:gmres-iter-es1}, \ref{tab:gmres-iter-es2}, and \ref{tab:gmres-iter-es4} summarize the numerical performance of the GMRES method both for the unpreconditioned system and the preconditioned ones.
Here, $N$ denotes the number of spatial grid points and $M$ the number of time levels and the average number of GMRES iterations per time step are reported in the columns.
Tables~\ref{tab:Timings_es1} and \ref{tab:Timings_es2}
report the averaged total computational time (in seconds) for solving the entire problem 10 times.
\par
As the discretization becomes finer, the unpreconditioned system  exhibits a substantial growth in both iteration counts and CPU times, reflecting the increasing ill-conditioning of the coefficient matrix.
The computational effort is mainly influenced by the spatial parameter $N$, as larger $N$ values lead to denser matrices and slower convergence, while variations in $M$ have a comparatively smaller effect.
\par
When the preconditioners are applied, GMRES shows a significant improvement in convergence.
The iteration numbers remain nearly constant for increasing $N$ and $M$, and the corresponding CPU times grow only mildly with the problem size.
This behavior demonstrates the mesh-independent convergence of the preconditioned system and confirms that the preconditioners effectively captures the spectral structure of the matrix, resulting in a well-clustered spectrum and reduced computational cost.
\begin{table}
	\centering
	\renewcommand{\arraystretch}{1}
	\setlength{\tabcolsep}{6pt}
	\begin{tabular}{|c c |ccc|ccc|ccc|}
		\hline
		\multicolumn {2}{|c}{} & \multicolumn {3}{|c}{$\alpha = 1.3$} & \multicolumn {3}{|c}{$\alpha = 1.5$} & \multicolumn {3}{|c|}{$\alpha = 1.8$} \\
		\hline
		$N$ & $M$ & $A$ & $\widetilde{P}$ &  $\widehat{P}$ & $-$ & $\widetilde{P}$ &  $\widehat{P}$ & $-$ & $\widetilde{P}$ &  $\widehat{P}$  \\
		\hline
		$2^4$ &$2^4$ &51 &13 &13 &59 &16 &16 &81 &18 &18 \\
		$2^4$ &$2^5$ &79 &13 &13 &89 &16 &16 &120 &20 &20 \\
		$2^4$ &$2^6$ &199 &13 &13 &200 &16 &16 &223 &20 &20 \\
		$2^4$ &$2^7$ &626 &13 &13 &609 &16 &16 &604 &20 &20 \\
		\hline
		$2^5$ &$2^4$ &82 &15 &14 &98 &17 &17 &140 &20 &20 \\
		$2^5$ &$2^5$ &112 &15 &15 &135 &18 &18 &200 &24 &24 \\
		$2^5$ &$2^6$ &229 &15 &15 &247 &18 &18 &327 &24 &24 \\
		$2^5$ &$2^7$ &717 &15 &15 &699 &18 &18 &751 &24 &24 \\
		\hline
		$2^6$ &$2^4$ &145 &16 &16 &175 &20 &19 &253 &22 &22 \\
		$2^6$ &$2^5$ &181 &16 &16 &232 &21 &20 &355 &28 &28 \\
		$2^6$ &$2^6$ &300 &16 &16 &365 &21 &20 &545 &29 &28 \\
		$2^6$ &$2^7$ &791 &16 &16 &841 &21 &20 &1102 &29 &28 \\
		\hline
		$2^7$ &$2^4$ &273 &17 &17 &329 &21 &21 &481 &22 &22 \\
		$2^7$ &$2^5$ &318 &17 &17 &427 &24 &23 &666 &33 &34 \\
		$2^7$ &$2^6$ &451 &18 &17 &617 &24 &23 &996 &34 &35 \\
		$2^7$ &$2^7$ &957 &18 &17 &1180 &24 &23 &1864 &34 &35 \\
		\hline
	\end{tabular}
	\caption{Iterations  required for solving the unpreconditioned $A$ and preconditioned linear system with $\widetilde{P} =\widetilde{P}_{\alpha,N,M}$ and $\widehat{P} = \widehat{P}_{\alpha,N,M}$ with the GMRES method - $a(x)=1$.}
	\label{tab:gmres-iter-es1}
\end{table}
\begin{table}
	\centering
	\renewcommand{\arraystretch}{1}
	\setlength{\tabcolsep}{3pt}
	\begin{tabular}{|c c |ccc|ccc|ccc|}
		\hline
		\multicolumn {2}{|c}{} & \multicolumn {3}{|c}{$\alpha = 1.3$} & \multicolumn {3}{|c}{$\alpha = 1.5$} & \multicolumn {3}{|c|}{$\alpha = 1.8$} \\
		\hline
		$N$ & $M$ & $A$ & $\widetilde{P}^{-1}A$ & $\widehat{P}^{-1}A$ & $A$ & $\widetilde{P}^{-1}A$ & $\widehat{P}^{-1}A$ & $A$ & $\widetilde{P}^{-1}A$ & $\widehat{P}^{-1}A$ \\
		\hline
		$2^4$ &$2^4$ &0.0035 &0.0132 &0.0058 &0.0023 &0.0069 &0.0066 &0.0005 &0.0081 &0.0061 \\
		$2^4$ &$2^5$ &0.0014 &0.0059 &0.0049 &0.0009 &0.0076 &0.0055 &0.0010 &0.0074 &0.0067 \\
		$2^4$ &$2^6$ &0.0043 &0.0085 &0.0077 &0.0043 &0.0090 &0.0077 &0.0051 &0.0113 &0.0092 \\
		$2^4$ &$2^7$ &0.0758 &0.0113 &0.0107 &0.0701 &0.0133 &0.0116 &0.0693 &0.0161 &0.0141 \\
		\hline
		$2^5$ &$2^4$ &0.0006 &0.0070 &0.0065 &0.0008 &0.0086 &0.0120 &0.0014 &0.0088 &0.0162 \\
		$2^5$ &$2^5$ &0.0015 &0.0098 &0.0081 &0.0024 &0.0112 &0.0099 &0.0058 &0.0124 &0.0118 \\
		$2^5$ &$2^6$ &0.0098 &0.0160 &0.0146 &0.0124 &0.0180 &0.0170 &0.0191 &0.0222 &0.0226 \\
		$2^5$ &$2^7$ &0.3930 &0.0250 &0.0255 &0.3764 &0.0297 &0.0290 &0.4587 &0.0393 &0.0501 \\
		\hline
		$2^6$ &$2^4$ &0.0023 &0.0371 &0.0141 &0.0032 &0.0180 &0.0179 &0.0079 &0.0221 &0.0266 \\
		$2^6$ &$2^5$ &0.0091 &0.0237 &0.0183 &0.0110 &0.0256 &0.0237 &0.0291 &0.0364 &0.0461 \\
		$2^6$ &$2^6$ &0.0917 &0.0457 &0.0403 &0.2056 &0.0571 &0.0542 &0.3715 &0.0830 &0.0958 \\
		$2^6$ &$2^7$ &1.7004 &0.0815 &0.0771 &1.8737 &0.0944 &0.0879 &3.8111 &0.2778 &0.2250 \\
		\hline
		$2^7$ &$2^4$ &0.0184 &0.0364 &0.0336 &0.0250 &0.0396 &0.0395 &0.0940 &0.0738 &0.0731 \\
		$2^7$ &$2^5$ &0.1173 &0.0551 &0.0916 &0.2095 &0.0733 &0.0825 &0.4578 &0.1577 &0.1877 \\
		$2^7$ &$2^6$ &0.6572 &0.1909 &0.1518 &1.4853 &0.1912 &0.3413 &3.3912 &0.2840 &0.3046 \\
		$2^7$ &$2^7$ &4.2326 &0.2572 &0.1810 &9.4866 &0.4533 &0.2786 &16.0755 &0.4172 &0.4750 \\
		\hline
	\end{tabular}
	\caption{Averaged timings for iterations in Table \ref{tab:gmres-iter-es1} repeated 10 times with  $A= A_{\alpha,N,M}$,
		$\widetilde{P}^{-1}A =\widetilde{P}^{-1}_{\alpha,N,M}A_{\alpha,N,M}$, and $\widehat{P}^{-1}A = \widehat{P}^{-1}_{\alpha,N,M}A_{\alpha,N,M}$ - case $a(x)=1$.}
	\label{tab:Timings_es1}
\end{table}
%
\begin{table}
	\centering
	\renewcommand{\arraystretch}{1}
	\setlength{\tabcolsep}{6pt}
	\begin{tabular}{|c c |ccc|ccc|ccc|}
		\hline
		\multicolumn {2}{|c}{} & \multicolumn {3}{|c}{$\alpha = 1.3$} & \multicolumn {3}{|c}{$\alpha = 1.5$} & \multicolumn {3}{|c|}{$\alpha = 1.8$} \\
		\hline
		$N$ & $M$ & $A$ & $\widetilde{P}$ &  $\widehat{P}$ & $-$ & $\widetilde{P}$ &  $\widehat{P}$ & $-$ & $\widetilde{P}$ &  $\widehat{P}$  \\
		\hline
		$2^4$ &$2^4$ &55 &15 &15 &67 &18 &18 &96 &21 &21 \\
		$2^4$ &$2^5$ &84 &16 &16 &98 &19 &19 &139 &24 &24 \\
		$2^4$ &$2^6$ &198 &16 &16 &201 &19 &19 &235 &24 &24 \\
		$2^4$ &$2^7$ &618 &16 &16 &604 &19 &19 &600 &24 &24 \\
		\hline
		$2^5$ &$2^4$ &91 &17 &18 &114 &21 &21 &171 &25 &25 \\
		$2^5$ &$2^5$ &127 &18 &18 &157 &23 &23 &244 &31 &31 \\
		$2^5$ &$2^6$ &234 &18 &18 &263 &23 &23 &381 &31 &31 \\
		$2^5$ &$2^7$ &701 &18 &18 &693 &23 &23 &793 &31 &31 \\
		\hline
		$2^6$ &$2^4$ &162 &19 &19 &207 &23 &23 &314 &26 &25 \\
		$2^6$ &$2^5$ &211 &20 &20 &273 &26 &26 &444 &36 &36 \\
		$2^6$ &$2^6$ &324 &20 &20 &411 &26 &26 &681 &38 &38 \\
		$2^6$ &$2^7$ &787 &20 &20 &875 &27 &26 &1282 &38 &38 \\
		\hline
		$2^7$ &$2^4$ &302 &21 &21 &394 &25 &25 &594 &26 &26 \\
		$2^7$ &$2^5$ &377 &22 &22 &506 &30 &30 &834 &40 &39 \\
		$2^7$ &$2^6$ &526 &22 &22 &730 &31 &30 &1287 &45 &45 \\
		$2^7$ &$2^7$ &1013 &22 &22 &1321 &31 &30 &2365 &46 &45 \\
		\hline
	\end{tabular}
	\caption{Iterations  required for solving the unpreconditioned $A$ and preconditioned linear system with $\widetilde{P} =\widetilde{P}_{\alpha,N,M}$ and $\widehat{P} = \widehat{P}_{\alpha,N,M}$ with the GMRES method - $a(x)=x^2+1$.}
	\label{tab:gmres-iter-es2}
\end{table}
\begin{table}
	\centering
	\renewcommand{\arraystretch}{1}
	\setlength{\tabcolsep}{3pt}
	\begin{tabular}{|c c |ccc|ccc|ccc|}
		\hline
		\multicolumn {2}{|c}{} & \multicolumn {3}{|c}{$\alpha = 1.3$} & \multicolumn {3}{|c}{$\alpha = 1.5$} & \multicolumn {3}{|c|}{$\alpha = 1.8$} \\
		\hline
		$N$ & $M$ & $A$ & $\widetilde{P}^{-1}A$ & $\widehat{P}^{-1}A$ & $A$ & $\widetilde{P}^{-1}A$ & $\widehat{P}^{-1}A$ & $A$ & $\widetilde{P}^{-1}A$ & $\widehat{P}^{-1}A$ \\
		\hline
		$2^4$ &$2^4$ &0.0008 &0.0059 &0.0045 &0.0003 &0.0063 &0.0049 &0.0079 &0.0181 &0.0077 \\
		$2^4$ &$2^5$ &0.0007 &0.0054 &0.0067 &0.0008 &0.0058 &0.0057 &0.0037 &0.0092 &0.0097 \\
		$2^4$ &$2^6$ &0.0044 &0.0087 &0.0079 &0.0039 &0.0087 &0.0086 &0.0062 &0.0134 &0.0107 \\
		$2^4$ &$2^7$ &0.0812 &0.0235 &0.0203 &0.0651 &0.0173 &0.0133 &0.0744 &0.0167 &0.0169 \\
		\hline
		$2^5$ &$2^4$ &0.0007 &0.0174 &0.0167 &0.0010 &0.0095 &0.0101 &0.0021 &0.0160 &0.0113 \\
		$2^5$ &$2^5$ &0.0048 &0.0192 &0.0138 &0.0026 &0.0123 &0.0116 &0.0061 &0.0174 &0.0157 \\
		$2^5$ &$2^6$ &0.0233 &0.0338 &0.0323 &0.0131 &0.0212 &0.0222 &0.0296 &0.0350 &0.0587 \\
		$2^5$ &$2^7$ &0.5191 &0.0555 &0.0515 &0.3705 &0.0354 &0.0364 &0.6981 &0.1120 &0.1018 \\
		\hline
		$2^6$ &$2^4$ &0.0053 &0.0346 &0.0385 &0.0042 &0.0199 &0.0223 &0.0163 &0.0395 &0.0371 \\
		$2^6$ &$2^5$ &0.0134 &0.0436 &0.0415 &0.0145 &0.0268 &0.0285 &0.0475 &0.0614 &0.0611 \\
		$2^6$ &$2^6$ &0.1129 &0.0620 &0.0863 &0.1322 &0.0492 &0.0691 &0.4915 &0.1285 &0.1319 \\
		$2^6$ &$2^7$ &1.5034 &0.1394 &0.1122 &1.8962 &0.1108 &0.1084 &4.3610 &0.2908 &0.2636 \\
		\hline
		$2^7$ &$2^4$ &0.0378 &0.0749 &0.0612 &0.0327 &0.0457 &0.0459 &0.1315 &0.0986 &0.0871 \\
		$2^7$ &$2^5$ &0.1386 &0.0995 &0.0609 &0.2669 &0.0844 &0.0802 &1.1671 &0.2150 &0.1563 \\
		$2^7$ &$2^6$ &0.9735 &0.2023 &0.1671 &1.3623 &0.1680 &0.2396 &5.9787 &0.4148 &0.7292 \\
		$2^7$ &$2^7$ &5.0880 &0.3373 &0.4153 &8.4437 &0.4889 &0.5794 &33.0959 &1.1614 &0.5574 \\
		\hline
	\end{tabular}
	\caption{Averaged timings for iterations in Table \ref{tab:gmres-iter-es2} repeated 10 times with  $A= A_{\alpha,N,M}$,
		$\widetilde{P}^{-1}A =\widetilde{P}^{-1}_{\alpha,N,M}A_{\alpha,N,M}$, and $\widehat{P}^{-1}A = \widehat{P}^{-1}_{\alpha,N,M}A_{\alpha,N,M}$ - case $a(x)=x^2+1$.}
	\label{tab:Timings_es2}
\end{table}
%
\begin{table}
	\centering
	\renewcommand{\arraystretch}{1}
	\setlength{\tabcolsep}{6pt}
	\begin{tabular}{|c c |ccc|ccc|ccc|}
		\hline
		\multicolumn {2}{|c}{} & \multicolumn {3}{|c}{$\alpha = 1.3$} & \multicolumn {3}{|c}{$\alpha = 1.5$} & \multicolumn {3}{|c|}{$\alpha = 1.8$} \\
		\hline
		$N$ & $M$ & $A$ & $\widetilde{P}$ &  $\widehat{P}$ & $-$ & $\widetilde{P}$ &  $\widehat{P}$ & $-$ & $\widetilde{P}$ &  $\widehat{P}$  \\
		\hline
		$2^4$ &$2^4$ &47 &26 &26 &59 &38 &38 &86 &65 &65 \\
		$2^4$ &$2^5$ &79 &27 &27 &87 &43 &43 &119 &82 &82 \\
		$2^4$ &$2^6$ &211 &28 &28 &210 &45 &45 &227 &96 &96 \\
		$2^4$ &$2^7$ &660 &28 &28 &656 &45 &45 &655 &99 &99 \\
		\hline
		$2^5$ &$2^4$ &71 &42 &42 &98 &66 &65 &154 &123 &123 \\
		$2^5$ &$2^5$ &100 &47 &47 &129 &84 &84 &207 &168 &167 \\
		$2^5$ &$2^6$ &232 &49 &49 &249 &101 &101 &343 &219 &218 \\
		$2^5$ &$2^7$ &767 &49 &49 &757 &103 &103 &780 &261 &260 \\
		\hline
		$2^6$ &$2^4$ &119 &66 &66 &170 &109 &109 &284 &216 &215 \\
		$2^6$ &$2^5$ &148 &86 &85 &210 &152 &151 &379 &271 &270 \\
		$2^6$ &$2^6$ &284 &91 &90 &349 &195 &193 &602 &353 &351 \\
		$2^6$ &$2^7$ &811 &93 &91 &843 &223 &223 &1161 &463 &461 \\
		\hline
		$2^7$ &$2^4$ &217 &99 &99 &317 &147 &146 &543 &198 &197 \\
		$2^7$ &$2^5$ &251 &122 &120 &378 &195 &193 &718 &256 &255 \\
		$2^7$ &$2^6$ &400 &138 &138 &571 &255 &250 &1124 &338 &336 \\
		$2^7$ &$2^7$ &921 &153 &146 &1128 &315 &311 &2062 &462 &458 \\
		\hline
	\end{tabular}
	\caption{Iterations  required for solving the unpreconditioned $A$ and preconditioned linear system with $\widetilde{P} =\widetilde{P}_{\alpha,N,M}$ and $\widehat{P} = \widehat{P}_{\alpha,N,M}$ with the GMRES method - $a(x)=x^2$.}
	\label{tab:gmres-iter-es4}
\end{table}
\clearpage
{\subsection*{GMRES performance in the two-dimensional setting}}
Finally, we report the results of a few numerical tests performed in the the two-dimensional setting, both in the constant-coefficient case and in the variable-coefficient setting, where the diffusion function $a(x,y)$ is chosen as the most natural bivariate extension of the ones previously considered.
Though, discretization parameters are a bit small, due to the global dimension of the all-at-once systems,  no relevant difference is observed in the convergence properties, which we already reported and commented in the one-dimensional setting; see Table \ref{tab:Test12d}, Table \ref{tab:Test22d}, and Table \ref{tab:Test42d} for the numerics in the two-dimensional space case.
\par
It is worth stressing that we do not report numerical experiments when $\alpha_1 \neq \alpha_2$ since the results in practice are of the same quality as in the case $\alpha_1=\alpha_2$ (and even slightly better indeed): the reason of such behavior is contained in Theorem \ref{2D th:}, Theorem \ref{2D orig mat}, Theorem \ref{thm:inv(P)A clustering 2D}, where it is evident that $\alpha_1 \neq \alpha_2$ implies a simplification in the spectral behavior, because the matrix sequence associated with the smaller parameter is spectrally negligible, when the matrix-sizes increase.

\begin{table}
	\centering
	\begin{tabular}{|c c c|ccc|ccc|ccc|}
		\hline
		\multicolumn {3}{|c}{} & \multicolumn {3}{|c}{$\alpha_1 =\alpha_2 = 1.3$} & \multicolumn {3}{|c}{$\alpha_1 =\alpha_2  = 1.5$} & \multicolumn {3}{|c|}{$\alpha_1 =\alpha_2  = 1.8$} \\
		\hline
		$N_1$ &$N_2$ & $M$ & $A$ & $\widetilde{P}^{-1}A$ & $\widehat{P}^{-1}A$ & $A$ & $\widetilde{P}^{-1}A$ & $\widehat{P}^{-1}A$ & $A$ & $\widetilde{P}^{-1}A$ & $\widehat{P}^{-1}A$ \\
		\hline
		$2^2$ &$2^2$ &$2^2$ &16 &10 &10 &19 &11 &11 &21 &12 &12 \\
		$2^2$ &$2^2$ &$2^3$ &22 &11 &11 &25 &12 &12 &27 &13 &14 \\
		$2^2$ &$2^2$ &$2^4$ &32 &11 &12 &35 &13 &13 &40 &15 &16 \\
		$2^2$ &$2^2$ &$2^5$ &65 &11 &12 &66 &13 &14 &68 &16 &16 \\
		\hline
		$2^3$ &$2^3$ &$2^2$ &27 &13 &13 &31 &14 &15 &37 &15 &15 \\
		$2^3$ &$2^3$ &$2^3$ &33 &13 &13 &37 &15 &15 &45 &17 &17 \\
		$2^3$ &$2^3$ &$2^4$ &45 &14 &14 &51 &16 &16 &63 &18 &18 \\
		$2^3$ &$2^3$ &$2^5$ &77 &14 &14 &81 &16 &17 &97 &20 &20 \\
		\hline
		$2^4$ &$2^4$ &$2^2$ &49 &15 &15 &55 &17 &17 &67 &18 &17 \\
		$2^4$ &$2^4$ &$2^3$ &56 &16 &16 &62 &19 &19 &79 &21 &21 \\
		$2^4$ &$2^4$ &$2^4$ &69 &16 &16 &79 &19 &19 &106 &23 &23 \\
		$2^4$ &$2^4$ &$2^5$ &103 &16 &16 &116 &20 &20 &155 &25 &25 \\
		\hline
		$2^5$ &$2^5$ &$2^2$ &91 &18 &18 &102 &20 &20 &124 &20 &20 \\
		$2^5$ &$2^5$ &$2^3$ &101 &18 &18 &112 &22 &22 &145 &24 &24 \\
		$2^5$ &$2^5$ &$2^4$ &118 &18 &18 &139 &23 &23 &189 &29 &29 \\
		$2^5$ &$2^5$ &$2^5$ &160 &19 &18 &187 &24 &24 &268 &31 &31 \\
		\hline
	\end{tabular}
	\caption{Iterations  required for solving the unpreconditioned $A$ and preconditioned linear system with $\widetilde{P} =\widetilde{P}_{\boldsymbol{\alpha},\textbf{n},M}$ and $\widehat{P} = \widehat{P}_{\boldsymbol{\alpha},\textbf{n},M}$ with the GMRES method - $a(x_1,x_2)=1$, $\theta = 0.5$.} \label{tab:Test12d}
\end{table}
\begin{table}
	\centering
	\begin{tabular}{|c c c|ccc|ccc|ccc|}
		\hline
		\multicolumn {3}{|c}{} & \multicolumn {3}{|c}{$\alpha_1 =\alpha_2 = 1.3$} & \multicolumn {3}{|c}{$\alpha_1 =\alpha_2  = 1.5$} & \multicolumn {3}{|c|}{$\alpha_1 =\alpha_2  = 1.8$} \\
		\hline
		$N_1$ &$N_2$ & $M$ & $A$ & $\widetilde{P}^{-1}A$ & $\widehat{P}^{-1}A$ & $A$ & $\widetilde{P}^{-1}A$ & $\widehat{P}^{-1}A$ & $A$ & $\widetilde{P}^{-1}A$ & $\widehat{P}^{-1}A$ \\
		\hline
		$2^2$ &$2^2$ &$2^2$ &20 &13 &13 &22 &14 &14 &24 &14 &14 \\
		$2^2$ &$2^2$ &$2^3$ &25 &14 &14 &28 &15 &15 &31 &16 &16 \\
		$2^2$ &$2^2$ &$2^4$ &36 &14 &15 &40 &16 &16 &45 &18 &19 \\
		$2^2$ &$2^2$ &$2^5$ &66 &14 &15 &67 &16 &17 &74 &20 &20 \\
		\hline
		$2^3$ &$2^3$ &$2^2$ &33 &15 &15 &38 &17 &17 &44 &17 &17 \\
		$2^3$ &$2^3$ &$2^3$ &40 &16 &16 &45 &18 &18 &54 &21 &20 \\
		$2^3$ &$2^3$ &$2^4$ &53 &17 &17 &62 &19 &19 &77 &23 &23 \\
		$2^3$ &$2^3$ &$2^5$ &82 &17 &17 &93 &20 &20 &118 &25 &25 \\
		\hline
		$2^4$ &$2^4$ &$2^2$ &59 &19 &18 &67 &21 &20 &84 &21 &21 \\
		$2^4$ &$2^4$ &$2^3$ &68 &19 &19 &78 &23 &22 &101 &25 &25 \\
		$2^4$ &$2^4$ &$2^4$ &86 &20 &20 &102 &23 &23 &139 &29 &29 \\
		$2^4$ &$2^4$ &$2^5$ &120 &21 &20 &146 &25 &24 &207 &32 &31 \\
		\hline
		$2^5$ &$2^5$ &$2^2$ &112 &22 &22 &127 &25 &24 &162 &27 &26 \\
		$2^5$ &$2^5$ &$2^3$ &125 &23 &23 &144 &27 &27 &193 &31 &30 \\
		$2^5$ &$2^5$ &$2^4$ &152 &23 &23 &183 &28 &28 &259 &36 &36 \\
		$2^5$ &$2^5$ &$2^5$ &198 &24 &24 &251 &30 &30 &381 &39 &39 \\
		\hline
	\end{tabular}
	\caption{Iterations  required for solving the unpreconditioned $A$ and preconditioned linear system with $\widetilde{P} =\widetilde{P}_{\boldsymbol{\alpha},\textbf{n},M}$ and $\widehat{P} = \widehat{P}_{\boldsymbol{\alpha},\textbf{n},M}$ with the GMRES method - $a(x_1,x_2)=x_1^2+x_2^2+1$, $\theta = 0.5$.} \label{tab:Test22d}
\end{table}
\begin{table}
	\centering
	\begin{tabular}{|c c c|ccc|ccc|ccc|}
		\hline
		\multicolumn {3}{|c}{} & \multicolumn {3}{|c}{$\alpha_1 =\alpha_2 = 1.3$} & \multicolumn {3}{|c}{$\alpha_1 =\alpha_2  = 1.5$} & \multicolumn {3}{|c|}{$\alpha_1 =\alpha_2  = 1.8$} \\
		\hline
		$N_1$ &$N_2$ & $M$ & $A$ & $\widetilde{P}^{-1}A$ & $\widehat{P}^{-1}A$ & $A$ & $\widetilde{P}^{-1}A$ & $\widehat{P}^{-1}A$ & $A$ & $\widetilde{P}^{-1}A$ & $\widehat{P}^{-1}A$ \\
		\hline
		$2^2$ &$2^2$ &$2^2$ &23 &18 &18 &27 &20 &21 &32 &23 &23 \\
		$2^2$ &$2^2$ &$2^3$ &27 &19 &18 &31 &22 &22 &38 &26 &26 \\
		$2^2$ &$2^2$ &$2^4$ &37 &19 &19 &42 &23 &23 &50 &29 &30 \\
		$2^2$ &$2^2$ &$2^5$ &72 &20 &20 &73 &25 &24 &79 &32 &32 \\
		\hline
		$2^3$ &$2^3$ &$2^2$ &41 &29 &29 &50 &36 &35 &63 &44 &44 \\
		$2^3$ &$2^3$ &$2^3$ &44 &29 &29 &54 &38 &37 &70 &52 &51 \\
		$2^3$ &$2^3$ &$2^4$ &54 &30 &30 &66 &42 &42 &88 &63 &63 \\
		$2^3$ &$2^3$ &$2^5$ &86 &32 &32 &96 &46 &46 &122 &74 &74 \\
		\hline
		$2^4$ &$2^4$ &$2^2$ &76 &41 &41 &93 &57 &56 &130 &79 &78 \\
		$2^4$ &$2^4$ &$2^3$ &77 &44 &43 &96 &61 &61 &141 &94 &93 \\
		$2^4$ &$2^4$ &$2^4$ &87 &48 &48 &111 &74 &74 &166 &125 &124 \\
		$2^4$ &$2^4$ &$2^5$ &116 &54 &53 &144 &91 &90 &217 &152 &150 \\
		\hline
		$2^5$ &$2^5$ &$2^2$ &149 &56 &55 &186 &86 &84 &253 &123 &122 \\
		$2^5$ &$2^5$ &$2^3$ &143 &61 &60 &187 &97 &95 &270 &158 &156 \\
		$2^5$ &$2^5$ &$2^4$ &152 &74 &73 &204 &123 &122 &315 &202 &196 \\
		$2^5$ &$2^5$ &$2^5$ &181 &93 &92 &247 &164 &161 &402 &249 &243 \\
		\hline
	\end{tabular}
	\caption{Iterations  required for solving the unpreconditioned $A$ and preconditioned linear system with $\widetilde{P} =\widetilde{P}_{\boldsymbol{\alpha},\textbf{n},M}$ and $\widehat{P} = \widehat{P}_{\boldsymbol{\alpha},\textbf{n},M}$ with the GMRES method - $a(x_1,x_2)=x_1^2+x_2^2$, $\theta = 0.5$.} \label{tab:Test42d}
\end{table}

	\chapter{Maximal results for geometric means of GLT sequences}\label{GM GLT2}


This chapter presents and discusses the main contributions of our recent work~\cite{GM2_GLT2025} on the spectral distribution of geometric means of HPD GLT matrix sequences. The spectral distribution of the geometric mean for two HPD GLT matrix sequences was studied in~\cite{ahmad2025matrix}. It was shown that if at least one of the GLT symbols \(\kappa\) and \(\xi\) is invertible almost everywhere, the geometric mean sequence remains GLT.
	In the scalar case \(r=1\), this hypothesis is equivalent to requiring that
	at least one of \(\kappa,\xi\) is nonzero almost everywhere; for matrix-valued
	symbols \(r>1\), the correct condition is invertibility a.e., i.e. the
	smallest eigenvalue of the symbol is positive on a set of full measure. The symbol is $(\kappa\xi)^{1/2}$ in the scalar or commuting case, and more generally $G(\kappa,\xi)=\kappa^{1/2}\left(\kappa^{-1/2}\xi\kappa^{-1/2}\right)^{1/2}\kappa^{1/2}$ in the noncommutative case. In this chapter, we address the open question of whether these assumptions can be relaxed. Motivated by the observation that invertibility is primarily required due to noncommutativity, we rigorously prove that when the GLT symbols $\kappa$ and $\xi$ commute, the invertibility condition can be completely removed in the case where $r=1$ and $d \geq 1$.
\smallskip

\noindent Furthermore, we investigate the limitations of this maximal result through both theoretical analysis and extensive numerical experiments. We show that, in the general block case ($r > 1$), if the GLT symbols $\kappa$ and $\xi$ are both non-invertible and do not commute, the geometric mean sequence $\{G(A_n, B_n)\}_n$ still forms a GLT matrix sequence with some symbol $\psi$; however, this symbol may not coincide with the classical geometric mean $G(\kappa, \xi)$, or, in some cases, $G(\kappa, \xi)$ may not be well defined. These findings demonstrate that our main results are indeed maximal: relaxing either the commutativity or invertibility assumptions can lead to fundamentally different or ill-posed spectral behavior.
\smallskip

\noindent A significant advance in this paper is the introduction and use of momentary GLT symbols (see \cite{bolten2022toeplitz,bolten2023note,new momentary}), which provide a finer description of the spectral distribution, especially for finite-size matrices and cases where the product $\kappa\xi$ is zero or the symbols are rank-deficient. Momentary symbols incorporate the effects of small-norm perturbations, allowing the observed richer spectral structures such as secondary eigenvalue distributions in numerical experiments to be explained within an extended GLT framework.
\smallskip

\noindent The work further investigates the extremal spectral behavior of the geometric mean, with particular focus on the asymptotic rate at which the minimal eigenvalues decay as the matrix size increases. Through detailed numerical experiments, we observe that the decay rate of the minimal eigenvalue is intimately connected to the order of vanishing (i.e., the multiplicity of zeros) of the associated GLT symbol.

\section{Geometric mean of GLT matrix sequences}\label{sec GM}
We begin by recalling the definition of the geometric mean of two HPD
matrices. Given two HPD matrices \(A,B\in\mathbb{C}^{n\times n}\), their
geometric mean is defined by
\begin{equation}\label{eq:gm-def}
	G(A,B)
	:= A^{1/2}\bigl(A^{-1/2}BA^{-1/2}\bigr)^{1/2}A^{1/2}=	G(B,A).
\end{equation}
This definition preserves symmetry, positive definiteness, and the key
variational properties of the two arguments \cite{Ando2004,bhatia2007}.
For matrix sequences \(\{A_n\}_n\) and \(\{B_n\}_n\) of HPD matrices, we write
\(\{G(A_n,B_n)\}_n\) for the corresponding sequence of geometric means.

The following result, due to \cite[Theorem 4]{ahmad2025matrix}, extends
\cite[Theorem 10.2]{garoni2017} to a broader setting.

\begin{theorem}[\!\!\cite{ahmad2025matrix}, Theorem 4]\label{th:two - r=1,d general}
	Let $r=1$ and $d\ge 1$.  Suppose \(\{A_{\bm{n}}\}_{\bm{n}} \sim_{\mathrm{GLT}} \kappa\) and \(\{B_{\bm{n}}\}_{\bm{n}} \sim_{\mathrm{GLT}} \xi\), where \(A_{\bm{n}}, B_{\bm{n}} \in \mathcal{P}_{\bm{n}}\) for every multi-index $\bm{n}$,{ with $\mathcal{P}_{\bm{n}}$ denoting the set of HPD matrices of multi-index $\bm{n}$}. Assume that at least one of between \(\kappa\) and \(\xi\) is nonzero almost everywhere. Then
	\begin{equation}\label{r=1, d general GLT}
		\{G(A_{\bm{n}}, B_{\bm{n}})\}_{\bm{n}} \sim_{\mathrm{GLT}}(\kappa \xi)^{1/2},
	\end{equation}
	and
	\begin{equation}\label{r=1, d general distributions}
		\{G(A_{\bm{n}}, B_{\bm{n}})\}_{\bm{n}} \sim_{\sigma, \lambda} (\kappa \xi)^{1/2}.
	\end{equation}
\end{theorem}

Using the topological a.c.s. notion, we now show that the thesis of the previous theorem holds without assuming that at least one between $\kappa$ and $\xi$ is invertible almost everywhere (a.e.), thus proving \cite[Conjecture 10.1]{garoni2017}. We make use of Axioms GLT \ref{glt:distrib}, GLT \ref{glt:symbols}, GLT \ref{glt:algebra}, GLT \ref{glt:acs}, GLT \ref{glt:hermitian}, GLT \ref{glt:function} as reported in Section~\ref{ast algebra}.

\begin{theorem}\label{theorem 1}
	Let $r=1$ and $d\ge 1$. 
	Assume \(\{A_{\bm{n}}\}_{\bm{n}} \sim_{\mathrm{GLT}} \kappa\) and \(\{B_{\bm{n}}\}_{\bm{n}}\sim_{\mathrm{GLT}} \xi\), where \(A_{\bm{n}}, B_{\bm{n}} \in \mathcal{P}_{\bm{n}}\) for every multi-index $\bm{n}$.  Then
	\[
	\{G(A_{\bm{n}}, B_{\bm{n}})\}_{\bm{n}} \sim_{\mathrm{GLT}} (\kappa \xi)^{1/2},
	\]
	\[
	\{G(A_{\bm{n}}, B_{\bm{n}})\}_{\bm{n}} \sim_{\sigma, \lambda} (\kappa \xi)^{1/2}.
	\]
\end{theorem}

\begin{proof}
	We have
	\begin{equation}\label{kappa equ}
		\{A_{\bm{n}}\}_{\bm{n}} \sim_{\mathrm{GLT}} \kappa,
	\end{equation}
	\begin{equation}\label{eta equ}
		\{B_{\bm{n}}\}_{\bm{n}} \sim_{\mathrm{GLT}} \xi,
	\end{equation}
	where $\kappa, \xi: [0,1]^d \times [-\pi,\pi]^d \to \mathbb{C}$ and $A_{\bm{n}}, B_{\bm{n}}$ are Hermitian positive definite matrices $\forall \bm{n}$. 
	Therefore $\kappa, \xi: [0,1]^d \times [-\pi,\pi]^d \to \mathbb{R}^+_0$ almost everywhere: furthermore, we assume that $ \mu_{2d}(\kappa \equiv 0) > 0$ and $\mu_{2d}(\xi \equiv 0) > 0$, so that the hypotheses of \cite[Theorem 10.2]{garoni2017} for $d=1$ and of Theorem \ref{th:two - r=1,d general} for $d>1$ are violated. Let $\varepsilon > 0$, let $\kappa_\varepsilon = \kappa + \varepsilon$ and let $A_{\bm{n}, \varepsilon} = A_{\bm{n}} + \varepsilon I_{N(\bm{n})},$ where $I_{N(\bm{n})}$ denotes the identity matrix of size $N(\bm{n})$. By the first part of GLT \ref{glt:symbols} (the identity is a special Toeplitz matrix with GLT symbol $1$), we have  $\{ \varepsilon I_{{\nu}(\bm{n})}=T_{\bm{n}}(\varepsilon) \}_{\bm{n}} \sim_{\mathrm{GLT}}\varepsilon,$ where $T_{\bm{n}}(\varepsilon)$ denotes the Toeplitz matrix generated by the constant function $\varepsilon$. Now, by exploiting linearity, i.e. the second item of GLT \ref{glt:algebra}, it follows that
	\[ 
	\{ A_{\bm{n}, \varepsilon} \}_{\bm{n}} \sim_{\mathrm{GLT}} \kappa_\varepsilon, \quad \kappa_\varepsilon \geq \varepsilon\, \text{ a.e.}, 
	\]
	since $\kappa, \xi \geq 0$ a.e. due to \eqref{kappa equ} and \eqref{eta equ}.
	Hence, we are again in the framework of Theorem \ref{th:two - r=1,d general}. Therefore, by  Theorem \ref{th:two - r=1,d general} we conclude that the sequence of geometric means satisfies the GLT relation
	\[ 
	\{G(A_{\bm{n},\varepsilon}, B_{\bm{n}}) \}_{\bm{n}} \sim_{\mathrm{GLT}} (\kappa_{\varepsilon} \xi)^{\frac{1}{2}}.
	\]
	Furthermore, we apply the full $*$-algebra framework of the GLT matrix sequences as in \cite[Section 2.9]{ahmad2025matrix}. More precisely, by third item of GLT \ref{glt:algebra} followed by GLT \ref{glt:function} with the function $f(z)=z^{\frac{1}{4}}$, we obtain
	\[ \{ (A_{\bm{n}, \varepsilon} B_{\bm{n}}^{2} A_{\bm{n}, \varepsilon})^{\frac{1}{4}} \}_{\bm{n}} \sim_{\mathrm{GLT}} (\kappa_\varepsilon^2 \xi^{2})^{\frac{1}{4}}=(\kappa_\varepsilon \xi)^{\frac{1}{2}}, \]
	which is the very same GLT symbol of $ \{G(A_{\bm{n},\varepsilon}, B_{\bm{n}}) \}_{\bm{n}} $. As a result, again by the second item of GLT \ref{glt:algebra}, the difference between these sequences satisfies the following asymptotic GLT relation
	\begin{equation}\label{rw inversion}
		\{G(A_{\bm{n}, \varepsilon}, B_{\bm{n}}) - (A_{\bm{n}, \varepsilon} B_{\bm{n}}^{2} A_{\bm{n}, \varepsilon})^{\frac{1}{4}} \}_{\bm{n}} \sim_{{\mathrm{GLT}} } 0.
	\end{equation}
	The previous relation is the key step since we have found a new GLT matrix sequence having the same GLT symbol as the geometric mean 
	matrix sequence $\{G(A_{\bm{n}, \varepsilon}, B_{\bm{n}})\}_{\bm{n}}$, but where no inversion is required.
Now given the structure of the performed operations, by virtue of Definition~\ref{def:acs}, the class 
	\[
	\{\{(A_{n,\epsilon_m}B_n^{2}A_{n,\epsilon_m})^{\frac14}\}_n\}_m,
	\qquad \epsilon_m=\tfrac{1}{m},
	\] 
	is obviously an a.c.s. for $\{ (A_{\bm{n}} B_{\bm{n}}^2 A_{\bm{n}})^{\frac{1}{4}} \}_{\bm{n}} $. As a consequence, by invoking (\ref{rw inversion}), also the class\(\{\{G(A_{n,\epsilon_m},B_n)\}_n\}_m\)
	is an a.c.s. for $\{ (A_{\bm{n}} B_{\bm{n}}^2 A_{\bm{n}})^{\frac{1}{4}} \}_{\bm{n}}, $
	with all the involved sequences being GLT matrix sequences, with symbols $ (\kappa_\varepsilon \xi)^{\frac{1}{2}}, (\kappa\xi)^{\frac{1}{2}}$, $\varepsilon>0$.
Therefore, since  $ \exists \lim_{m \to 0} (\kappa_{\varepsilon_{m}} \xi)^{\frac{1}{2}}= (\kappa\xi)^{\frac{1}{2}}$, by using the powerful Theorem~\ref{th: acs}, we deduce that $ \exists \lim_{m \to 0} \text{ (a.c.s.) } \{\{ G(A_{\bm{n}, \varepsilon_m}, B_{\bm{n}}) \}_{\bm{n}} \}_m= \{ (A_{\bm{n}} B_{\bm{n}}^2 A_{\bm{n}})^{\frac{1}{4}} \}_{\bm{n}} \sim_{\mathrm{GLT}} (\kappa\xi)^{\frac{1}{2}}$ and this limit coincides with $ \{ G(A_{\bm{n}}, B_{\bm{n}}) \}_{\bm{n}}$ in the a.c.s. topology, that is, 
	$\{G(A_{\bm{n}}, B_{\bm{n}}) - (A_{\bm{n}} B_{\bm{n}}^2 A_{\bm{n}})^{\frac{1}{4}} \}_{\bm{n}}\sim_{{\mathrm{GLT}}} 0.$
	Finally $\{(A_{\bm{n}} B_{\bm{n}}^2 A_{\bm{n}})^{\frac{1}{4}} \}_{\bm{n}}\sim_{{\mathrm{GLT}}} (\kappa\xi)^{\frac{1}{2}}$. Hence
	\[
	\{G(A_{\bm{n}}, B_{\bm{n}})\}_{\bm{n}} \sim_{\mathrm{GLT}} (\kappa \xi)^{1/2},
	\]
	so that, by GLT \ref{glt:distrib}, we infer
	\[
	\{G(A_{\bm{n}}, B_{\bm{n}})\}_{\bm{n}} \sim_{\sigma, \lambda} (\kappa \xi)^{1/2}.
	\]
\end{proof}
\begin{remark}[Intuition on the generalization]
	\label{gen vs barrier}
	One of the basic but key ingredients of the previous proof is that scalar-valued functions commute. Hence, it is reasonable to expect that the same proof also works in the case of $d$-level $r$-block GLT matrix sequences, when assuming that the symbols commute. We collect the result in Theorem \ref{theorem 2} whose proof follows verbatim that of Theorem \ref{theorem 1}, except for minimal changes. On the other hand, when both $\kappa$ and $\xi$ are not invertible almost everywhere (degenerate) and do not commute, the expression $G(\kappa,\xi)$ is not well defined. In fact, we could replace the inversion with the standard pseudoinversion, which is denoted as $X^+$ if $X$ is any complex matrix. However, we stress that the expression
	$\kappa^{\frac{1}{2}} ([\kappa^{\frac{1}{2}}]^+ \xi [\kappa^{\frac{1}{2}}]^+)^{\frac{1}{2}} \kappa^{\frac{1}{2}}$ is not the same as $\xi^{\frac{1}{2}} ([\xi^{\frac{1}{2}}]^+ \kappa [\xi^{\frac{1}{2}}]^+)^{\frac{1}{2}} \xi^{\frac{1}{2}}$ in general and this is a serious indication that the commutation between the GLT symbols is essential, when the symbols are both degenerate.
\end{remark}

We now recall the result proven in \cite{ahmad2025matrix}, concerning the general GLT setting with $d,r \ge 1$ and when at least one of the involved GLT symbols is invertible a.e.

\begin{theorem}[\!\!\cite{ahmad2025matrix}, Theorem 5]\label{th:two - r,d general}
	Let $r,d\ge 1$. 
	Suppose \(\{A_{\bm{n}}\}_{\bm{n}} \sim_{\mathrm{GLT}} \kappa\) and \(\{B_{\bm{n}}\}_{\bm{n}} \sim_{\mathrm{GLT}} \xi\), where \(A_{\bm{n}}, B_{\bm{n}} \in \mathcal{P}_{\bm{n}}\) for every multi-index $\bm{n}$. Assume that at least one of the minimal eigenvalue of \(\kappa\) and the minimal eigenvalue of \(\xi\) is nonzero almost everywhere. Then
	\begin{equation}\label{r,d general GLT}
		\{G(A_{\bm{n}}, B_{\bm{n}})\}_{\bm{n}} \sim_{\mathrm{GLT}}G(\kappa, \xi),
	\end{equation}
	and
	\begin{equation}\label{r,d general distributions}
		\{G(A_{\bm{n}}, B_{\bm{n}})\}_{\bm{n}} \sim_{\sigma, \lambda} G(\kappa,\xi).
	\end{equation}
	Furthermore $G(\kappa,\xi)=(\kappa \xi)^{1/2}$ whenever the GLT symbols  \(\kappa\) and \(\xi\) commute.
\end{theorem}

Assume that \(\{A_n\}_n \sim_{\mathrm{GLT}} \kappa\) and
\(\{B_n\}_n \sim_{\mathrm{GLT}} \xi\), where
\(A_n, B_n \in \mathcal{P}_{\mathbf{n}}\) for every multi-index
\(\mathbf{n}\), with \(\mathcal{P}_{\mathbf{n}}\) denoting the set of HPD
matrices of size \(rN(\mathbf{n})\), where
\(N(\mathbf{n}) = \prod_{r=1}^{d} n_r\).
\begin{theorem}\label{theorem 2}
	Let $r>1$ and $d\ge 1$. 
	Assume that  \(\{A_{\bm{n}}\}_{\bm{n}} \sim_{\mathrm{GLT}} \kappa\) and \(\{B_{\bm{n}}\}_{\bm{n}} \sim_{\mathrm{GLT}} \xi\), where \(A_{\bm{n}}, B_{\bm{n}} \in \mathcal{P}_{\bm{n}}\) for every multi-index $\bm{n}$, with $\mathcal{P}_{\bm{n}}$ denoting the set of HPD matrices of size $rN(\bm{n})$. Under the assumption that $\kappa$ and $\xi$ commute we infer 
	\[
	\{G(A_{\bm{n}}, B_{\bm{n}})\}_{\bm{n}} \sim_{\mathrm{GLT}} (\kappa \xi)^{1/2},
	\]
	\[
	\{G(A_{\bm{n}}, B_{\bm{n}})\}_{\bm{n}} \sim_{\sigma, \lambda} (\kappa \xi)^{1/2},
	\]
\end{theorem}
\begin{proof}
	Since $\{A_{\bm{n}}\}_{\bm{n}} \sim_{\mathrm{GLT}} \kappa$, $\{B_{\bm{n}}\}_{\bm{n}} \sim_{\mathrm{GLT}} \xi$, we deduce that $\kappa, \xi: [0,1]^d \times [-\pi,\pi]^d \to \mathbb{C}^{r\times r}$ are Hermitian nonnegative definite, because, by the assumptions, $A_{\bm{n}}, B_{\bm{n}}$ are HPD matrices $\forall \bm{n}$. 
	Therefore, the minimal eigenvalue of $\kappa$ and the minimal eigenvalue of $\xi$ are nonnegative almost everywhere. Here, we assume that $ \mu_{2d}(\lambda_{\min}(\kappa) \equiv 0) > 0$ and $\mu_{2d}(\lambda_{\min}(\xi) \equiv 0) > 0$, in such a way that the hypotheses of Theorem \ref{th:two - r,d general} are violated. 
	Let $\varepsilon > 0$, let $\kappa_\varepsilon = \kappa + \varepsilon I_r$, $I_r$ being the identity of size $r$, and let $A_{\bm{n}, \varepsilon} = A_{\bm{n}} + \varepsilon I_{rN(\bm{n})}$. Now $\{ \varepsilon I_{rN(\bm{n})}=T_{\bm{n}}(\varepsilon I_r) \}_{\bm{n}} \sim_{\mathrm{GLT}} \varepsilon I_r$ by the first part of GLT \ref{glt:symbols} (the identity $I_{rN(\bm{n})}$ is a special multilevel block Toeplitz matrix with GLT symbol $I_r$). Consequently, by exploiting linearity i.e. the second item of GLT \ref{glt:algebra}, it follows that
	\[ 
	\{ A_{\bm{n}, \varepsilon} \}_{\bm{n}} \sim_{\mathrm{GLT}} \kappa_\varepsilon, \quad \kappa_\varepsilon \geq \varepsilon I_r \text{ a.e.} 
	\]
	Since $\kappa, \xi$ are both nonnegative definite almost everywhere.
	Hence, by  Theorem \ref{th:two - r,d general}, we have $\{G(A_{\bm{n},\varepsilon}, B_{\bm{n}}) \}_{\bm{n}} \sim_{\mathrm{GLT}} (\kappa_{\varepsilon} \xi)^{\frac{1}{2}}$, 
	since the commutation between $\kappa$ and $\xi$ implies the commutation between $\kappa_\varepsilon = \kappa + \varepsilon I_r$ and $\xi$. 
	By exploiting the $*$-algebra features of the GLT matrix sequences, and specifically the third item of GLT \ref{glt:algebra} and GLT \ref{glt:function} with the function $f(z)=z^{\frac{1}{4}}$, we deduce
	$ \{ (A_{\bm{n}, \varepsilon} B_{\bm{n}}^{2} A_{\bm{n}, \varepsilon})^{\frac{1}{4}} \}_{\bm{n}} \sim_{\mathrm{GLT}} (\kappa_\varepsilon^2 \xi^{2})^{\frac{1}{4}}=(\kappa_\varepsilon \xi)^{\frac{1}{2}}$, 
	which is the very same GLT symbol of $ \{G(A_{\bm{n},\varepsilon}, B_{\bm{n}}) \}_{\bm{n}}$. As a result, again by the second item of GLT \ref{glt:algebra}, the difference between these sequences satisfies the following asymptotic relation
	\begin{equation}\label{reducing without inversion}
		\{G(A_{\bm{n}, \varepsilon}, B_{\bm{n}}) - (A_{\bm{n}, \varepsilon} B_{\bm{n}}^{2} A_{\bm{n}, \varepsilon})^{\frac{1}{4}} \}_{\bm{n}} \sim_{{\mathrm{GLT}} } 0.
	\end{equation}
	In other words, we have written a new GLT matrix sequence having the same symbol as the geometric mean 
	matrix sequence $\{G(A_{\bm{n}, \varepsilon}, B_{\bm{n}})\}_{\bm{n}}$, but where no inversion is required.
	By virtue of Definition~\ref{def:acs} , the class $ \{\{ (A_{\bm{n}, \varepsilon_m} B_{\bm{n}}^{2} A_{\bm{n}, \varepsilon_m})^{\frac{1}{4}} \}_{\bm{n}}\}_{m}$ is an a.c.s. for $\{ (A_{\bm{n}} B_{\bm{n}}^2 A_{\bm{n}})^{\frac{1}{4}} \}_{\bm{n}} $.\\ Therefore, by (\ref{reducing without inversion}), also the class $ \{\{ G(A_{\bm{n}, \varepsilon_m}, B_{\bm{n}}) \}_{\bm{n}}\}_m $
	is an a.c.s for $\{ (A_{\bm{n}} B_{\bm{n}}^2 A_{\bm{n}})^{\frac{1}{4}} \}_{\bm{n}}, $
	with all the involved sequences being GLT matrix sequences, with symbols $ (\kappa_\varepsilon \xi)^{\frac{1}{2}}, (\kappa\xi)^{\frac{1}{2}}$, $\varepsilon>0$. Since  $ \exists \lim_{\varepsilon \to 0} (\kappa_\varepsilon \xi)^{\frac{1}{2}}= (\kappa\xi)^{\frac{1}{2}}$, Theorem~\ref{th: acs}  implies $$ \exists \lim_{\varepsilon \to 0} \text{ (a.c.s.) } \{\{ G(A_{\bm{n}, \varepsilon}, B_{\bm{n}}) \}_{\bm{n}}, \varepsilon \}= \{ (A_{\bm{n}} B_{\bm{n}}^2 A_{\bm{n}})^{\frac{1}{4}} \}_{\bm{n}} \sim_{\mathrm{GLT}} (\kappa\xi)^{\frac{1}{2}}.$$ It is now clear that this limit coincides with $ \{ G(A_{\bm{n}}, B_{\bm{n}}) \}_{\bm{n}} $
	in the a.c.s. topology, that is, 
	$\{G(A_{\bm{n}}, B_{\bm{n}}) - (A_{\bm{n}} B_{\bm{n}}^2 A_{\bm{n}})^{\frac{1}{4}} \}_{\bm{n}}\sim_{{\mathrm{GLT}}} 0.$
	In this manner $\{(A_{\bm{n}} B_{\bm{n}}^2 A_{\bm{n}})^{\frac{1}{4}} \}_{\bm{n}} \sim_{{\mathrm{GLT}}} (\kappa\xi)^{\frac{1}{2}}$. Consequently
	$\{G(A_{\bm{n}}, B_{\bm{n}})\}_{\bm{n}} \sim_{\mathrm{GLT}} (\kappa \xi)^{1/2}$, so that, by GLT \ref{glt:distrib}, we finally obtain
	$\{G(A_{\bm{n}}, B_{\bm{n}})\}_{\bm{n}} \sim_{\sigma, \lambda} (\kappa \xi)^{1/2}$.
\end{proof}
\section{Numerical experiments}\label{Num_Exp}

As is well known, many localization, extremal, and distribution results hold when
\(d\)-level, \(r\)-block Toeplitz matrix sequences are considered and these
results are somehow summarized in specific analytic features of the generating
function of the corresponding matrix sequence. In turn, the generating function
is also the \(d\)-variate, \(r\times r\) matrix-valued GLT symbol of the
\(d\)-level, \(r\)-block Toeplitz matrix sequence.

Although the distribution results are also valid for general \(d\)-level,
\(r\)-block GLT matrix sequences, this is no longer true in general for the
extremal behavior and for the localization results, unless we add supplementary
assumptions, such as requiring the sequence to be generated by a
\emph{matrix-valued linear positive operator} (LPO), i.e.\ a linear operator
from the space of matrix-valued functions to the space of matrices that
preserves positive semidefiniteness; see
\cite{LPO-rev,SerraCapizzanoTilli-LPO}. We observe that the geometric mean can
be seen as a monotone operator with respect to its two variables, and the
monotonicity is implied when we consider an LPO, even if the converse is not
true. Hence, it is also interesting to verify which properties are maintained
by a geometric mean of two \(d\)-level, \(r\)-block GLT matrix sequences in
terms of its GLT symbol, when the latter exists.
\smallskip

\noindent According to the previous discussion, the remainder of the section considers numerical experiments in the following directions:

\begin{itemize}
	\item the validation of the distribution results in the commuting setting as in Theorem~\ref{theorem 1} and Theorem~\ref{theorem 2};
	\item connections of the previous results with the notion of Toeplitz and GLT momentary symbols and the extremal behavior compared to the GLT symbol; 
	\item evidence that Theorem~\ref{theorem 1} and Theorem~\ref{theorem 2} are maximal, by taking GLT matrix sequences with noncommuting symbols which are both not invertible almost everywhere. 
\end{itemize}
All numerical experiments are performed using MATLAB R2022b on a laptop equipped with an 11th Gen Intel(R) Core(TM) i5-1155G7 CPU running at 2.50 GHz, with 16 GB of RAM. The operating system was Windows 11 Pro (version 23H2, build 22631.5189).
\subsection{Validation of the distribution results} \label{ssec:num1}

\subsubsection{Example 1 }
Let $ d = 1 $ and consider the following two matrix-sequences $\{A_{n}\}_{n}$, $\{B_{n}\}_{n}$, with $
A_n = D_n(a) + \frac{1}{n^4} I_n$ and  $ B_n = T_n(3 + 2\cos(\theta)),$
where $T_n(\cdot) $ denotes the Toeplitz operator for $d=1 $, as introduced in Section~\ref{TM},  $D_n(a) $ represents the diagonal matrix generated by the continuous function
\[
a(x) =
\begin{cases}
	0, & x \in [0, \frac{1}{2}), \\
	1, & x \in [\frac{1}{2},1].
\end{cases}
\]
Because of Theorem~\ref{theorem 1}, the geometric mean of these two sequences, that is $\{ G(A_n, B_n) \}_n $, satisfies the GLT relation:
\[
\left\{ G\left(A_n, B_n\right) \right\}_n \sim_{\mathrm{GLT}} \sqrt{a(x)(3+2\cos(\theta))},
\]
where the corresponding GLT symbols are explicitly given by
$
\kappa = a(x)$ and $ \xi=3+2\cos(\theta).
$
\subsubsection*{Eigenvalue distribution}
\begin{figure}[H]
	\begin{minipage}{0.52\textwidth}
		\centering
		\includegraphics[width=\linewidth, height=6cm]{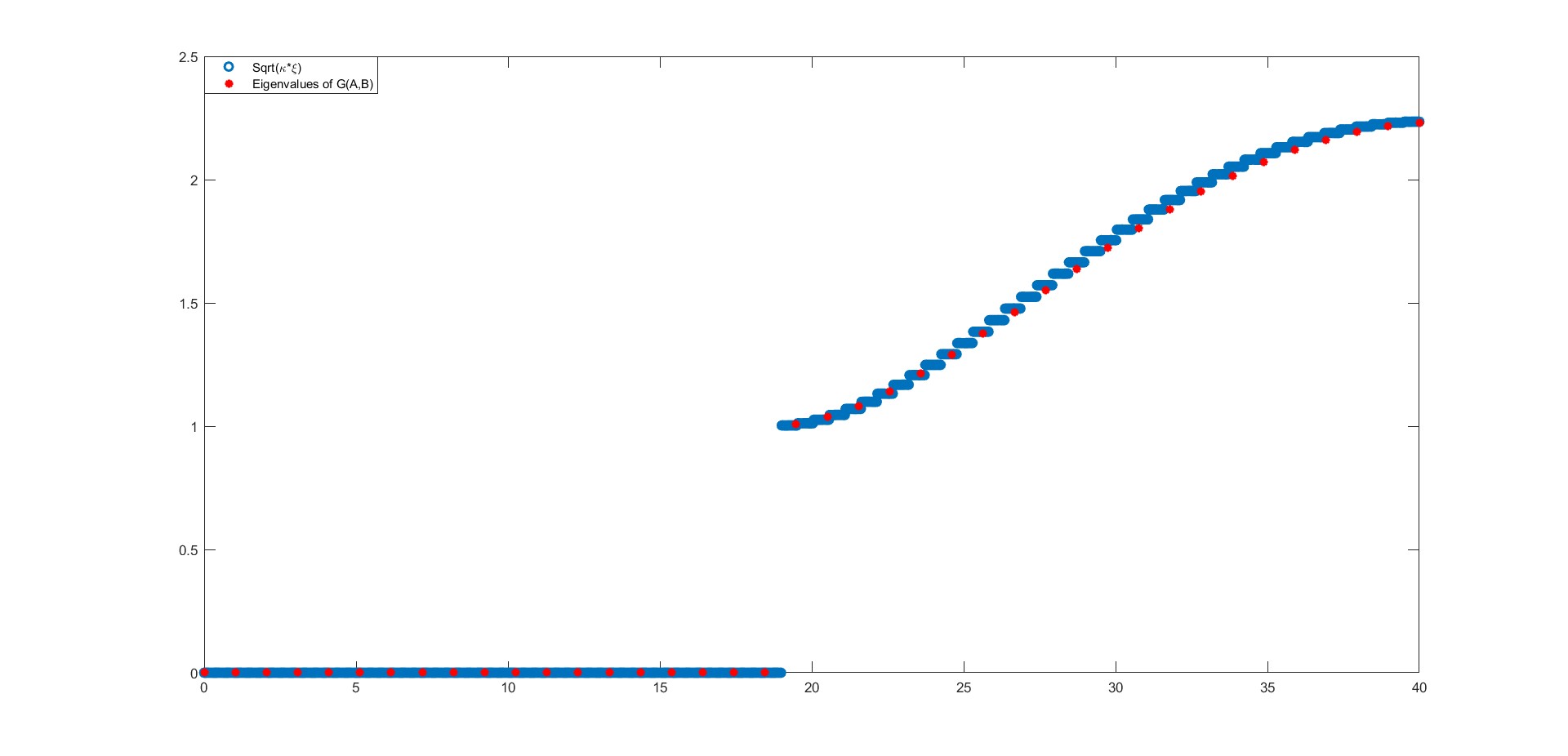} 
		\captionsetup{labelformat=empty}
		\caption*{}
	\end{minipage}\hfill
	\begin{minipage}{0.52\textwidth}
		\centering
		\includegraphics[width=\linewidth, height=6cm]{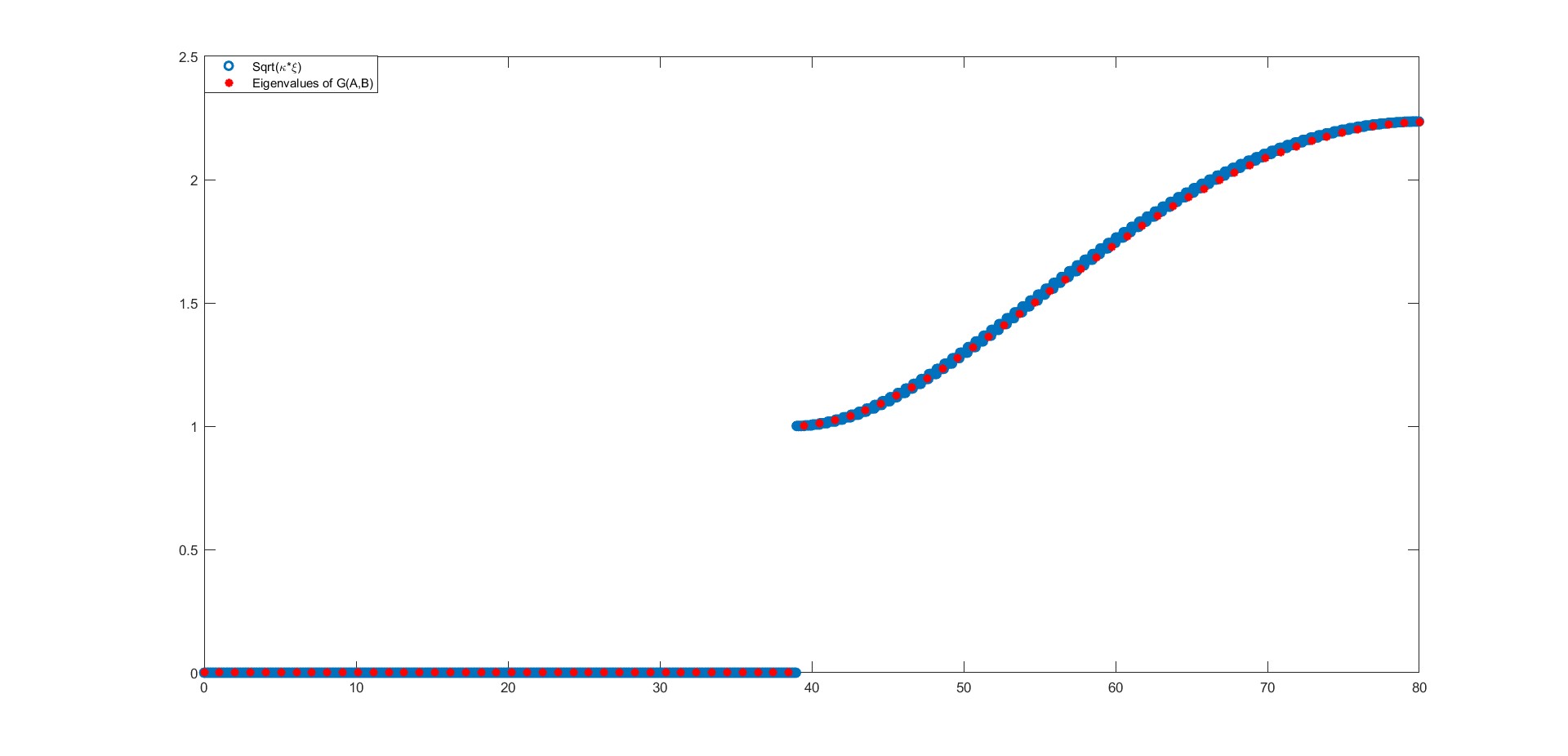} 
		\captionsetup{labelformat=empty}
		\caption*{}
	\end{minipage}
	\caption{{\textbf{ Example 1.}
			Comparison between eigenvalues of \(G(A_n,B_n)\) (red stars; \(n=40,80\)) versus the symbol \( (\kappa\xi)^{\frac{1}{2}} \) (blue circles)}}
	\label{fig:exp1}
\end{figure}

We numerically analyze the spectral behavior of the geometric mean $G(A_n, B_n) $ from \textbf{Example 1}. The eigenvalues of this geometric mean are computed for various increasing dimensions $n $ and compared with the uniformly sampled points from the GLT symbol $\sqrt{a(x)(3+2\cos(\theta))} $. In Figure \ref{fig:exp1}, numerical results strongly indicate that the GLT symbol accurately characterizes the eigenvalue distribution of the geometric mean. As $ n $ grows, the eigenvalue distribution closely matches the GLT symbol, giving evidence of the theoretical findings in Theorem~\ref{theorem 1}.

\subsection{Distribution results and momentary symbols} \label{ssec:num2}

\subsubsection{Example 2 }
Let $d = 1 $ and define the matrix sequences $\{A_{n}\}_{n}$, $\{B_{n}\}_{n}$, with $
A_n = D_n(a) + \frac{1}{n^4} I_n,$  $ C_n = T_n(3 + 2\cos(\theta))
$ and $
B_n = \left(D_n(1 - a) + \frac{1}{n^4} I_n \right) C_n \left(D_n(1 - a) + \frac{1}{n^4} I_n \right),
$ where \( D_n(a) \) is the diagonal matrix generated by the piecewise continuous function \( a(x) \) given in \textbf{Example 1}. Then, the geometric mean sequence explicitly satisfies the GLT relation:
\begin{align*}
	\left\{ G\left(A_n, B_n\right) \right\}_n &\sim_{\mathrm{GLT}} \sqrt{a(x)(1 - a(x))(3 + 2\cos\theta)} \\
	&= a(x)(1 - a(x))\sqrt{3 + 2\cos\theta} = 0
\end{align*}
where
\[
\kappa = a(x), \quad \xi=(1-a(x))(3+2\cos(\theta)).
\]
The zero-valued symbol arises naturally due to the definition of the piecewise function \( a(x) \):
\begin{itemize}
	\item For \( x \in [0, \frac{1}{2}) \), we have \( a(x) = 0 \). Thus, the term \( a(x)(1 - a(x)) \) becomes zero.
	\item For \( x \in [\frac{1}{2},1] \), we have \( a(x) = 1 \). In this interval, the factor \( 1 - a(x) \) becomes zero, making \( a(x)(1 - a(x)) = 0 \).
\end{itemize}
\subsubsection*{Eigenvalue distribution}
\begin{figure}[H]
	\begin{minipage}{0.52\textwidth}
		\centering
		\includegraphics[width=\linewidth, height=6cm]{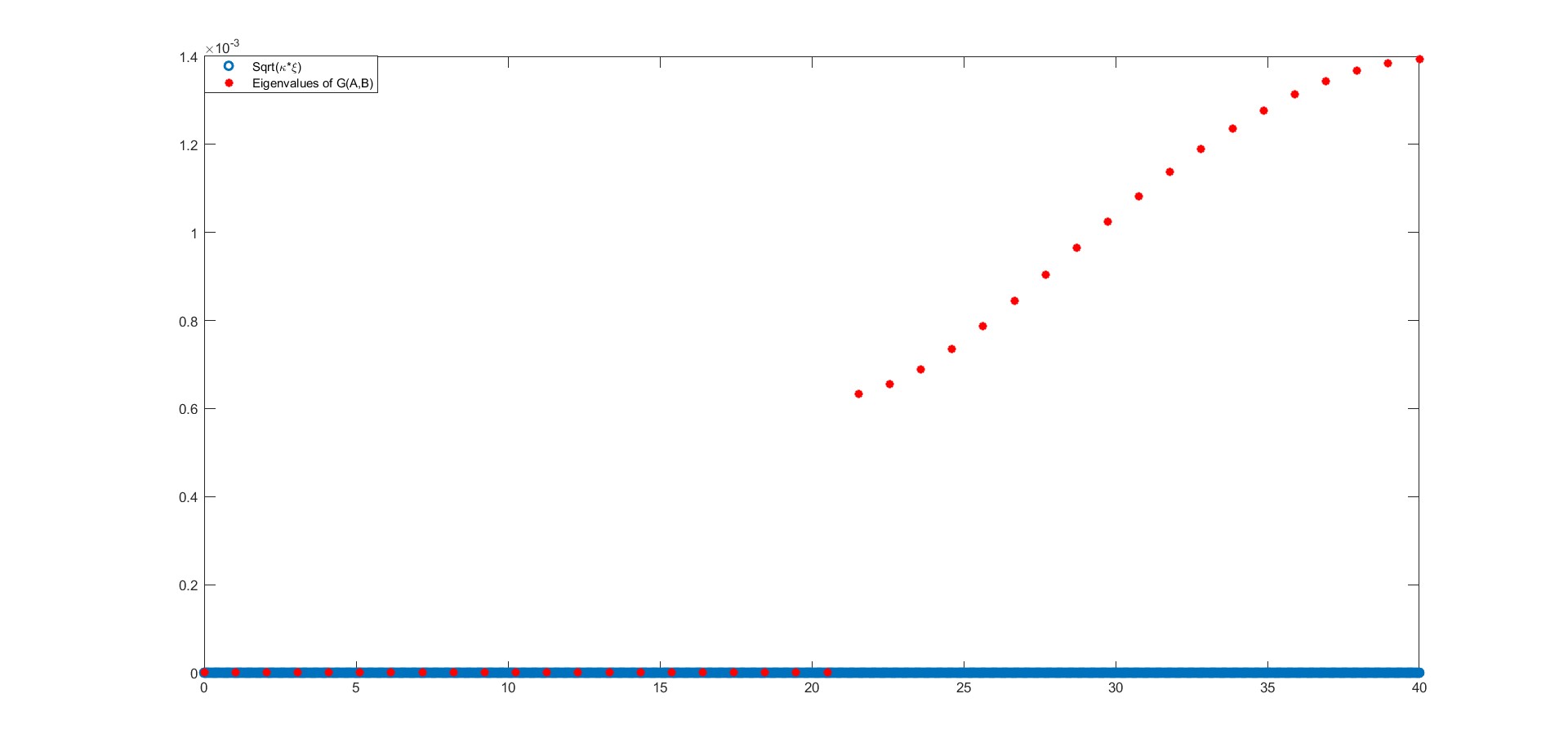} 
		\captionsetup{labelformat=empty}
		\caption*{}
	\end{minipage}\hfill
	\begin{minipage}{0.52\textwidth}
		\centering
		\includegraphics[width=\linewidth, height=6cm]{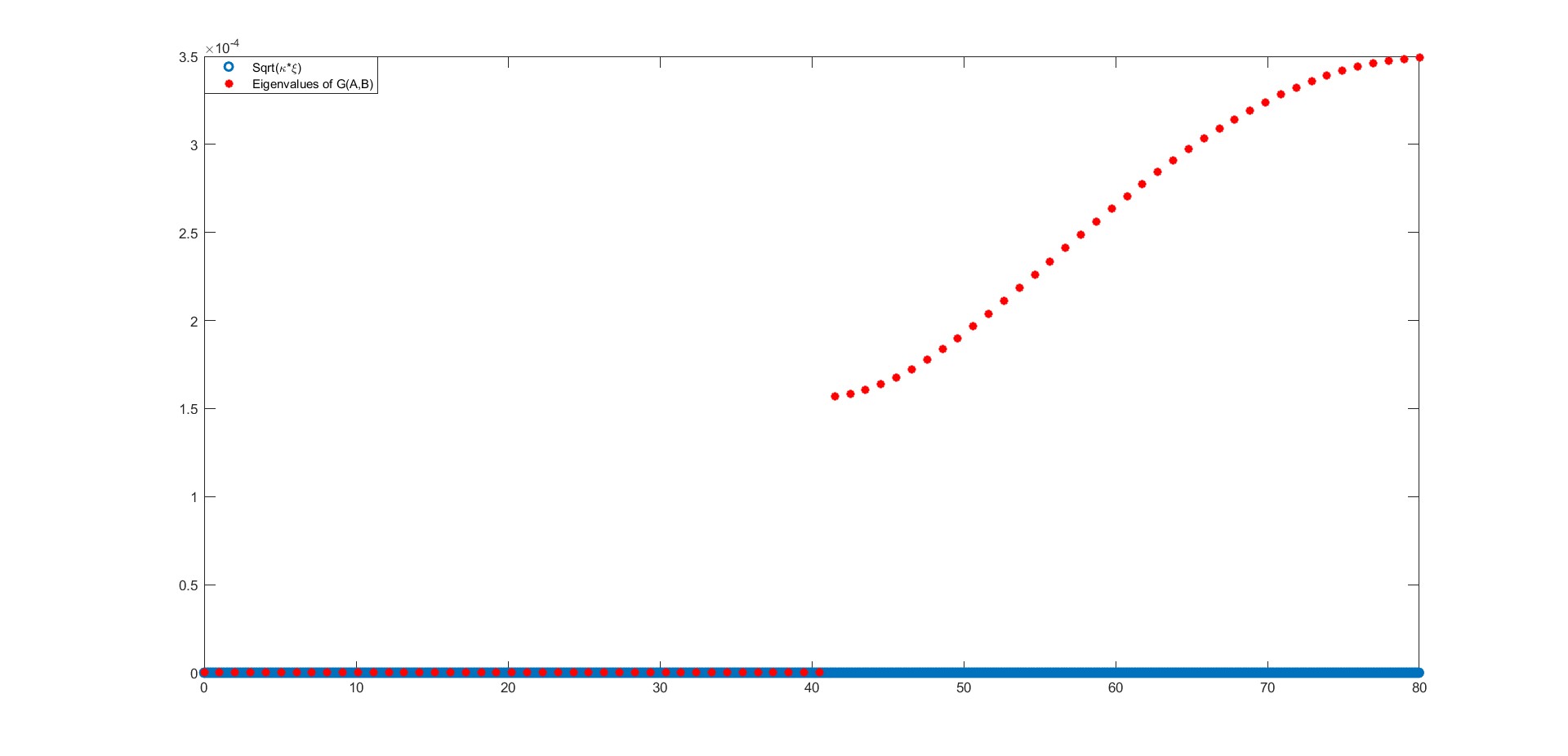} 
		\captionsetup{labelformat=empty}
		\caption*{}
	\end{minipage}
	\caption{{\textbf{ Example 2.}
			Comparison between eigenvalues of \(G(A_n,B_n)\) (red stars; \(n=40,80\)) versus the symbol \( (\kappa\xi)^{\frac{1}{2}} \) (blue circles)}}
	\label{fig:exp2}
\end{figure}

The eigenvalue distribution of \( \{G(A_n, B_n)\}_{n} \) is numerically analyzed for increasing matrix dimensions \( n \). Since the GLT symbol for this example is zero, by Theorem~\ref{theorem 1} we expect the eigenvalues to concentrate along the zero line. Interestingly, numerical results reveal an additional spectral structure beyond Theorem~\ref{theorem 1}.
{For smaller values of $ n $, eigenvalues align closely with zero, as predicted by the GLT symbol. However, some eigenvalues appear above this level, forming a secondary symbol. This phenomenon is attributed to {momentary symbols}, as discussed in the introduction (see \cite{bolten2022toeplitz,bolten2023note,new momentary}). As \( n \) increases, these elevated eigenvalues shift upwards reaching approximately \( 10^{-3} \) for \( n = 40 \) (shown in the Figure \ref{fig:exp2}) and \( 10^{-4} \) for \( n = 80 \).}

\subsection{Minimal eigenvalues and conditioning} \label{ssec:num3}

In this section, we analyze the extremal spectral behavior and conditioning of the geometric mean sequence \( \{ G(A_n, B_n) \}_n \), focusing on its dependence on the analytical properties of the corresponding GLT symbol. This approach follows prior studies on extremal eigenvalues in structured matrix settings, particularly in Toeplitz matrices \cite{extreme3,extreme1,extreme2,extreme4} and block Toeplitz matrices \cite{SerraCapizzano1999a,SerraCapizzano1999b}, as well as variable coefficient differential operators, including multilevel cases with \( d > 1 \) \cite{extreme diff1,extreme diff2,extr2-glt}. In the case of variable coefficient differential operators it is worth noticing that extremal spectral results do not stem from the GLT theory but from a combination of GLT tools and properties which are typical of linear positive operators 
\cite{SerraCapizzanoTilli-LPO,LPO-rev}.\\
Here, we restrict our attention to the unilevel scalar setting with \( d = r = 1 \), considering \textbf{Example 1} and \textbf{Example 2}.
\subsubsection{Example 1: minimal eigenvalue}
\begin{itemize}
	\item $X_{n}=G(A_{n},B_{n})$
	\item Take $n_{j}=40\cdot2^{j}, \quad\quad j=0,1,2,3, $
	\item Compute $\tau_{j}= \lambda_{\min}(X_{n_{j}}), \quad \quad j=0,1,2,3, $
	\item Compute $\alpha_{j}=\log_{2}(\frac{\tau_{j}}{\tau_{j+1}}), \quad \quad j=0,1,2. $ 
\end{itemize}
\vspace{9pt}

\begin{table}[H]
	\centering
	\caption{Numerical behaviour of the minimal \\eigenvalue.}
	\label{tab:tau_alpha2}
	\begin{tabular}{|c
			|S[table-format=1.4e+2]
			|S[table-format=1.4]|}
		\hline
		{$n$} & {$\tau_j$} & {$\alpha_j$}\\
		\hline
		40  & 6.3265e-04 & 2.0132 \\
		80  & 1.5673e-04 & 2.0064 \\
		160 & 3.9093e-05 & 2.0074 \\
		320 & 9.7675e-06 & \multicolumn{1}{c|}{}     \\
		\hline
	\end{tabular}
\end{table}
As it can be seen, the quantity $\alpha_j$ stabilizes around $2$ as $n$ increases. This is in perfect agreement with the fact that the minimal eigenvalue of $A_n$ converges to zero as $n^{-4}$, while the minimal eigenvalue of $B_n$ converges monotonically from above to $\min \xi=1$. The key point is that the minimal eigenvalue of the geometric mean behaves asymptotically as the geometric mean of the minimal eigenvalues of $A_n$ and $B_n$, respectively. Similar remarks can be made in the subsequent case with respect to \textbf{Example 2}.

The observed numerical evidences are not implied by the theoretical derivations and this is an interesting fact that deserves to be investigated theoretically in the future. In particular, we would like to prove formally that that $\alpha_j$ converges to $2$ for the first example, Table~\ref{tab:tau_alpha2}, and it converges to 4 in the second example, Table~\ref{tab:tau_alpha3}, even from the preliminary
numerical results the convergence does not look monotonic.

\subsubsection{Example 2: minimal eigenvalue}
\begin{itemize}
	\item $X_{n}=G(A_{n},B_{n})$
	\item Take $n_{j}=40\cdot2^{j}, \quad\quad j=0,1,2,3, $
	\item Compute $\tau_{j}= \lambda_{\min}(X_{n_{j}}), \quad \quad j=0,1,2,3, $
	\item Compute $\alpha_{j}=\log_{2}(\frac{\tau_{j}}{\tau_{j+1}}), \quad \quad j=0,1,2. $ 
\end{itemize}
\vspace{9pt}
\begin{table}[H]
	\centering
	\caption{Numerical behaviour of the minimal\\ eigenvalue.}
	\label{tab:tau_alpha3}
	\begin{tabular}{|c
			|S[table-format=1.4e+2]
			|S[table-format=1.4]|}
		\hline
		{$n$} & {$\tau_j$} & {$\alpha_j$} \\
		\hline
		40  & 3.9177e-07 & 4.0003 \\
		80  & 2.4480e-08 & 4.0040 \\
		160 & 1.5259e-09 & 4.0009 \\
		320 & 9.5367e-11 & \multicolumn{1}{c|}{}      \\
		\hline
	\end{tabular}
\end{table}


\subsection{Numerical study: noncommuting, rank-deficient symbols}
\label{sec:non_commuting_rank_def}

In this section, through numerical test, we investigate the eigenvalue distribution of the geometric means $\{G(A_{n}, B_{n})\}_{n}$ with
\[
G(A_n,B_n)\;=\;
A_n^{1/2}\!\bigl(A_n^{-1/2}B_nA_n^{-1/2}\bigr)^{1/2}\!A_n^{1/2},
\]
where $\{A_n\}_n$ and $\{B_n\}_n$ are HPD GLT matrix sequences with $r\times r$ matrix-valued symbols $\kappa,\xi$.
The GLT matrix sequences are specifically chosen with GLT symbols so that:

\begin{enumerate}
	\item they are rank-deficient on a subset of positive Lebesgue measure; and
	\item they do not commute pointwise on a set of positive measure.
\end{enumerate}

These hypotheses clearly set the problem outside the scope of commuting or almost everywhere invertible cases treated in Theorems~\ref{theorem 1} and in Theorem~\ref{theorem 2}. \\
Our tests provide evidence for two phenomena:
\begin{itemize}
	\item[\textbf{(i)}]
	The support of the asymptotic eigenvalue distribution of $\{G(A_n,B_n)\}_{n}$ coincides with the intersection set $\operatorname{ess\,Ran}(\kappa)$ $\cap$ $\operatorname{ess\,Ran}(\xi)$ see Section~\ref{subsubsec:GM_conjecture}.
	\item[\textbf{(ii)}] The eigenvalues appear to converge in distribution to a candidate symbol which is rank-deficient and, in general, distinct from the classical geometric mean $G(\kappa,\xi)$.
\end{itemize}
Clearly, the relation \( G(A_n, B_n) \sim_{\mathrm{GLT}} G(\kappa, \xi) \) cannot hold under hypotheses 1.-2., and the experiments show the substantial differences in this case compared to the commuting or invertible scenarios. If this case can be studied within the GLT framework, it requires a more technical application of GLT and a.c.s. theory, along with a deeper functional calculus of matrix means.

\subsubsection{Setup of the numerical experiments}
\label{sec:NumSetup}
Below, we provide the setup of the numerical tets. \\
We construct four pairs of unilevel GLT matrix sequences $\{A_n\}_n, \{B_n\}_n,$ whose symbols
simultaneously satisfy the hypothesis stated at the beginning of
Section~\ref{sec:non_commuting_rank_def} (rank deficiency and
non-commutation). These examples are organized according to the property of rank loss:

\begin{itemize}\setlength\itemsep{4pt}
	\item \textbf{Case 1} - each symbol is full rank on a set of positive
	measure;
	\item \textbf{Case 2} - each symbol is rank-deficient almost everywhere.
\end{itemize}

\medskip\noindent
\textbf{Case 1, Example 1.}

Define
\[
f(\theta)\; =\; \begin{cases}
	0, & \theta\in[-\pi,0],\\[4pt]
	\theta, & \theta\in(0,\pi],
\end{cases}
\qquad
g(\theta) = f(-\theta),
\]
and set
\[
F(\theta)\;=\;f(\theta)\otimes
\begin{bmatrix}2&1\\[2pt]1&2\end{bmatrix},\qquad
G(\theta)\;=\;g(\theta)\otimes
\begin{bmatrix}3&1\\[2pt]1&1\end{bmatrix}.
\]
With the $2n\times2n$ Toeplitz matrices $T_n(F)$ and $T_n(G)$, define
\begin{equation}\label{eq:case1_ex1}
	A_n=T_n(F)+\frac{1}{n^{3}}I_{2n},\qquad
	B_n=T_n(G)+\frac{1}{n^{3}}I_{2n}.
\end{equation}

\medskip\noindent
\textbf{Case 1, Example 2.}

Let $\chi_{[-a,a]}$ be the characteristic function of $[-a,a]$, $0<a<\pi$, and set
\[
f(\theta)\,=\,\chi_{[-1/2,\,1/2]}(\theta),\qquad
g(\theta)\,=\,\chi_{[-1/4,\,1/4]}(\theta),
\]
\[
A\;=\;\begin{bmatrix}2&1\\[2pt]1&2\end{bmatrix},\qquad
B\;=\;\begin{bmatrix}3&1\\[2pt]1&1\end{bmatrix},
\]
\[
F(\theta)\;=\;f(\theta)\otimes A,\qquad
G(\theta)\;=\;g(\theta)\otimes B.
\]
Define
\begin{equation}\label{eq:case1_ex2}
	A_n=T_n(F)+\frac{1}{n^{3}}I_{2n},\qquad
	B_n=T_n(G)+\frac{1}{n^{3}}I_{2n}.
\end{equation}

\medskip\noindent
\textbf{Case 2, Example 1.}

With $f(\theta)=2-\cos\theta$, $g(\theta)=3+\cos\theta$ and
rank one blocks
\[
A\;=\;\begin{bmatrix}1&1\\[2pt]1&1\end{bmatrix},\qquad
B\;=\;\begin{bmatrix}1&2\\[2pt]2&4\end{bmatrix},
\]
let
\(F(\theta)\;=\;f(\theta)\otimes A\) and
\(G(\theta)\;=\;g(\theta)\otimes B\).
Set
\[
A_n=T_n(F)+\frac{1}{n^{2}}I_{2n},\qquad
B_n=T_n(G)+\frac{1}{n^{2}}D_n(b),\qquad
b(x)=(1+x) I_2.
\]
Here, both symbols are rank $1$ almost everywhere.

\medskip\noindent
\textbf{Case 2, Example 2.}

Define piecewise-linear scalar functions in $[0,1]$:
\[
a(x)\;=\;\begin{cases}
	1-2x, & 0\le x\le\tfrac12,\\[2pt]
	0, & \tfrac12<x\le1,
\end{cases}
\qquad
b(x)\;=\;\begin{cases}
	0, & 0\le x\le\tfrac13\text{ or } \tfrac23\le x\le1,\\[2pt]
	x-\tfrac13, & \tfrac13<x\le\tfrac12,\\[2pt]
	\tfrac23-x, & \tfrac12<x<\tfrac23,
\end{cases}
\]
and let $D_n(a)$, $D_n(b)$ be the corresponding sampling matrices.
With $f(\theta)=2+\cos\theta$, $g(\theta)=3+\cos\theta$ and
\[
A\;=\;\begin{bmatrix}2 & 0 & 1\\[2pt]
	0 & 2 & 1\\[2pt]
	1 & 1 & 1\end{bmatrix},\qquad
B\;=\;\begin{bmatrix}2 & 1 & 0\\[2pt]
	1 & 1 & 1\\[2pt]
	0 & 1 & 2\end{bmatrix},
\]
define
\(F(\theta)\;=\;f(\theta)\otimes A\) and
\(G(\theta)\;=\;g(\theta)\otimes B\), and set
\begin{align}\label{eq:case2_ex2}
	A_n &\;=\; \bigl(D_n^{1/2}(a)\,\otimes\, I_3\bigr)\,
	T_n(F)\, \bigl(D_n^{1/2}(a)\,\otimes\, I_3\bigr)
	+\frac{1}{5n}\,I_{3n}, \\
	B_n&\;=\; \bigl(D_n^{1/2}(b)\,\otimes\, I_3\bigr)\,
	T_n(G)\,
	\bigl(D_n^{1/2}(b)\,\otimes\, I_3\bigr)
	+\frac{1}{5n}\,I_{3n}.
\end{align}

\medskip
In every example we denote by $\kappa$ and $\xi$ the GLT symbols of
$\{A_n\}_n$ and $\{B_n\}_n$, respectively.

\subsubsection{A candidate symbol for the geometric mean}
\label{subsubsec:GM_conjecture}

Let $\{A_n\}_n$ and $\{B_n\}_n$ be GLT matrix sequences whose matrix-valued symbols
\[
\kappa,\;\xi : [0,1]^d\times[-\pi,\pi]^d \;\longrightarrow\; \mathbb{C}^{\,r\times r},
\]
are positive-semidefinite. \\
For any $\varepsilon>0$, define the strictly positive symbols
\(
\kappa_\varepsilon=\kappa+\varepsilon I_r,\;
\xi_\varepsilon=\xi+\varepsilon I_r.
\)
Since the matrix geometric mean is monotone continuous in each argument, the limit below exists point-wise {and it is unique} (\cite[p.\,3]{Ando2004}):
\begin{definition}[Candidate symbol]\label{def:candidate}
	For almost every $(x,\theta)$ set
	\begin{equation}\label{eq:cand_symbol}
		\widetilde G(\kappa,\xi)(x,\theta)
		:=\lim_{\varepsilon\rightarrow 0}
		G\bigl(\kappa_\varepsilon(x,\theta),\,
		\xi_\varepsilon (x,\theta)\bigr).
	\end{equation}
\end{definition}
The matrix-valued symbol $\widetilde G(\kappa,\xi)$ is \emph{essentially} positive semidefinite (see \cite{Ando2004}) and satisfies the essential support identity
\[
\operatorname{ess\,Ran}\,\widetilde G(\kappa,\xi)\;=\;
\operatorname{ess\,Ran}\!\bigl[(x,\theta)\mapsto
\operatorname{ran}\kappa(x,\theta)\;\cap\;\operatorname{ran}\xi(x,\theta)\bigr].
\]

Here $\operatorname{ran}M$ denotes the column space of an individual matrix $M$, whereas
$\operatorname{ess\,Ran}(.)$ is the essential range of a Grassmannian valued measurable map. The essential formulation of the subspace condition can also be formulated as the pre-image of the essential numerical range in the style of [\cite{nomi-mazza}, Definition\, 1].

Taking into account the Definition \ref{def:candidate} of the candidate symbol, we propose the following conjecture:
\begin{conjecture}[GLT closure under geometric mean]
	\label{conj:GM_glt}
	Let $\{A_n\}_n$ and $\{B_n\}_n$ be HPD
	GLT matrix sequences with
	\[
	\{A_n\}_n\sim_{\mathrm{GLT}}\kappa,
	\qquad
	\{B_n\}_n\sim_{\mathrm{GLT}}\xi.
	\]
	Then the geometric mean matrix sequence is again GLT and
	\[
	\{\,G(A_n,B_n)\,\}_n \;\sim_{\operatorname{GLT}}\; \widetilde G(\kappa,\xi),
	\]
	where $\widetilde G(\kappa,\xi)$ is the candidate symbol
	defined in~\eqref{eq:cand_symbol}.
\end{conjecture}

The heuristic justification of the conjecture is the following. For each fixed $\varepsilon>0$ the shifted sequences
$\{A_n+\varepsilon I\}_n$ and $\{B_n+\varepsilon I\}_n$ are GLT sequences symbols
$\kappa_\varepsilon$ and $\xi_\varepsilon$ and their geometric mean sequences satisfy
\(
\{\,G(A_n+\varepsilon I,B_n+\varepsilon I)\,\}_n
\sim_{\operatorname{GLT}}
G(\kappa_\varepsilon,\xi_\varepsilon).
\)
{Properties of monotone functions and the matrix geometric mean axioms described in \cite{Ando2004}}, combined with closeness of the GLT $\ast$-algebra (GLT \ref{glt:acs}) under a.c.s convergence, suggest that letting $\varepsilon\rightarrow0$ the limit of $ \{\,G(A_n+\varepsilon I,B_n+\varepsilon I)\,\}_n$ is well defined in the {a.c.s topology and this limit is $ \{\,G(A_n,B_n)\,\}_n$, up to zero-distributed perturbations }, leading to the GLT algebraic part of Conjecture~\ref{conj:GM_glt}, a full proof when
$\kappa$ and $\xi$ do not commute point-wise remains open. Meanwhile, all
numerical evidence below supports the conjecture on the spectral distribution part of our hypothesis. 
\subsubsection*{Candidate symbols predicted by Conjecture~\ref{conj:GM_glt}}

For the four test pairs introduced in Section \ref{sec:NumSetup}, Conjecture~\ref{conj:GM_glt} predicts
the following GLT symbols for $\{G(A_n,B_n)\}_n$.
Throughout
\(
\mathbf 1=(1,1,1)^{\mathsf T},\;
J=\mathbf 1\,\mathbf 1^{\mathsf T}\in\mathbb{C}^{3\times3}.
\)

\paragraph{Case 1, Example 1.}
\[
\widetilde G(\kappa,\xi)(x,\theta)=0_{2\times2}.
\]

\paragraph{Case 1, Example 2.}
\[
\widetilde G(\kappa,\xi)(x,\theta)=
\chi_{[-\frac14,\frac14]}(\theta)\;
\otimes C,
\]
\[
C \; =\;
\frac{1}{6^{1/4}\sqrt{\,2+\sqrt6\,}}
\begin{bmatrix}
	2\sqrt2+3\sqrt3 & \sqrt2+\sqrt3\\[4pt]
	\sqrt2+\sqrt3 & 2\sqrt2+\sqrt3
\end{bmatrix}.
\]

\paragraph{Case 2, Example 1.}
\[
\widetilde G(\kappa,\xi)(x,\theta)=0_{2\times2}.
\]

\paragraph{Case 2, Example 2.}
\[
\widetilde G(\kappa,\xi)(x,\theta)=
\sqrt{\,\bigl(2+\cos\theta\bigr)\,\bigl(3+\cos\theta\bigr)\,h(x)\,k(x)}\;
\otimes J.
\]

\subsubsection{Numerical verification of the Conjecture}

Let
\[
G_n := G(A_n,B_n)\in\mathbb{C}^{d_n\times d_n},
\qquad d_n=\text{size}(G_n).
\]
We test whether the empirical eigenvalue distribution of $\{G_n\}_n$
converges to the distribution induced by the
candidate symbol
$\widetilde G(\kappa,\xi)(x,\theta)$ defined in
Section~\ref{subsubsec:GM_conjecture}.
The comparison follows the standard rearrangement strategy.

\begin{itemize}\setlength\itemsep{6pt}
	
	\item[\textbf{1.}] \textbf{Sampling the symbol.}
	Evaluate $\widetilde G(\kappa,\xi)$ on a tensor grid
	$\{x_j\}_{j=1}^{M_x}\!\subset[0,1]$,
	$\{\theta_i\}_{i=1}^{M_\theta}\!\subset[-\pi,\pi]$ with a total number of
	$M_x \times M_\theta\simeq2000$ samples.
	For each point $(x_j,\theta_i)$ compute the
	$r$ eigenvalues of the function; then merge and sort in non-decreasing order all the computed values.
	This yields the empirical quantile function of the eigenvalues.
	
	\item[\textbf{2.}] \textbf{Eigenvalue computation.}
	Matrices $A_n,B_n$ are generated for
	$n\in\{40,80,160,320\}$, and eigenvalues are ordered in non-decreasing order, to obtain the quantile spectral distribution function of $\{G_n\}_n$. \\
	We also record
	$\lambda_{\min}(G_n)$, $\lambda_{\max}(G_n)$ and the condition
	number~$\mathrm{cond}_2(G_n)$ for later stability analysis.
	
	\item[\textbf{3.}] \textbf{Quantile plot matching.}
	For each $n$ we plot the sorted eigenvalues
	$\bigl(\lambda_i(G_n)\bigr)$ against the corresponding rearranged eigenvalues of the sample values of $\widetilde G$.
	A superposition of the two plots for growing values of $n$ gives an experimental indication of the asymptotic convergence of the empirical spectral distribution of $\{G_n\}_n$ to that of the symbol.
	
	\item[\textbf{4.}] \textbf{Support analysis.}
	To estimate the measure of the zero set
	$\{(x,\theta):\widetilde G=0\}$, we count the fraction of eigenvalues below the threshold $0.1$  compared to the matrix size; the value is an estimate of the complement of the predicted measure of the support as $n$ grows.
\end{itemize}
All tests for this section are run on an Intel Core Ultra 7 155H CPU (22 threads, 30.9 GiB RAM) under Ubuntu 24.10, Linux 6.11.0-29, with Python~3.12.7, NumPy~1.26.4 and SciPy~1.13.1 linked against OpenBLAS.
\medskip
\begin{remark}[Measure-theoretic meaning of the rearrangement test]
	\label{rearrangement}
	As anticipated implicitly, the comparison of sorted eigenvalues versus the sorted eigenvalues of the sampled symbol  has a precise measure-theoretic foundation based on the quantile approximation theory \cite{Bogoya2016}. For completeness, we report the construction in \cite{Bogoya2016} for the unilevel setting $d=1$; the multilevel case is analogous.
 	Let $r$ denote the size of the matrix
 $\widetilde G(\kappa,\xi)(x,\theta)$ and put
 $D=[0,1]\times[-\pi,\pi]$ (so $\mu(D)=2\pi$).	Define the probability space by the triple
		\[
		(\Omega ,\mathcal{F},\mathbb{P}),\qquad
		\Omega = D\times \{1,\ldots ,r\} ,\qquad
		\mathcal{F} = \mathcal{B}(D)\otimes 2^{\{1,\ldots ,r\}}.
		\]
		Here, \(\mathcal{B}(D)\) denotes the Borel \(\sigma\)-algebra on \(D\),
		\(2^{\{1,\ldots ,r\}}\) denotes the power set of \(\{1,\ldots ,r\}\),
		and the product measure is
		\[
		\mathbb{P}(A\times B)
		= \frac{\mu(A)}{2\pi}\,\frac{|B|}{r},
		\qquad A\in\mathcal{B}(D),\ B\subseteq\{1,\ldots ,r\}.
		\]
		With this choice,
		\[
		\mathbb{P}(\Omega)
		= \frac{\mu(D)}{2\pi}\cdot\frac{r}{r}
		= 1,
		\]
		so that \((\Omega,\mathcal{F},\mathbb{P})\) is indeed a probability space.
\vspace*{0.15cm}\\
	Introduce the random variable
	\[
	X(x,\theta,i):=\lambda_i\!\bigl(\widetilde G(\kappa,\xi)(x,\theta)\bigr),
	\qquad (x,\theta,i)\in\Omega,
	\]
	
	where the eigenvalues are ordered non-decreasingly.
	The quantile function of $X$ is precisely the \emph{non-decreasing rearrangement} $\widetilde G^{\dagger}$ of the matrix-valued symbol.\\
	Because spectral distributions are defined only up to
	measure-preserving rearrangements, we deduce that
	\[
	\{G_n\}_n\sim_{\lambda} \widetilde G
	\;\Longrightarrow\;
	\{G_n\}_n\sim_{\lambda} \widetilde G^{\dagger},
	\]
	so verifying convergence to the quantile $\widetilde G^{\dagger}$ is equivalent and numerically simpler than verifying convergence to the symbol
	$\widetilde G(\kappa,\xi)$ itself. \\
	Furthermore, whenever $\widetilde G^{\dagger}$ is continuous at $t\in(0,1)$, the sorted
	eigenvalues satisfy
	\(
	\lambda_{\lceil t d_n\rceil}(G_n)\to \widetilde G^{\dagger}(t)
	\)
	as $n\to\infty$, and the empirical proportion of eigenvalues below a small threshold
	$t_0$ (we take $t_0=0.1$) converges to
	\(
	\mathbb P[X\le t_0],
	\)
	therefore it can be used as an estimate of the measure of the rank-deficient subset
	$\{(x,\theta):\widetilde G(\kappa,\xi)(x,\theta)=0\}$.
\end{remark}

\subsubsection{Numerical results}
\label{sec:numerical-results}
\begin{figure}[H]
	\centering
	\begin{subfigure}[t]{0.49\textwidth}
		\centering
		\includegraphics[width=\textwidth]{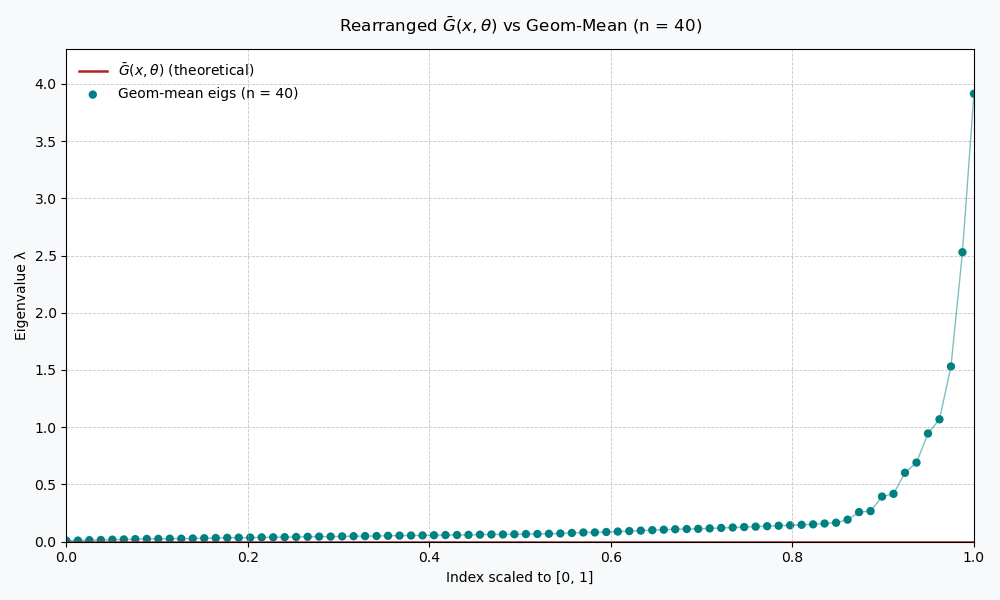}
	\end{subfigure}
	\hfill
	\begin{subfigure}[t]{0.49\textwidth}
		\centering
		\includegraphics[width=\textwidth]{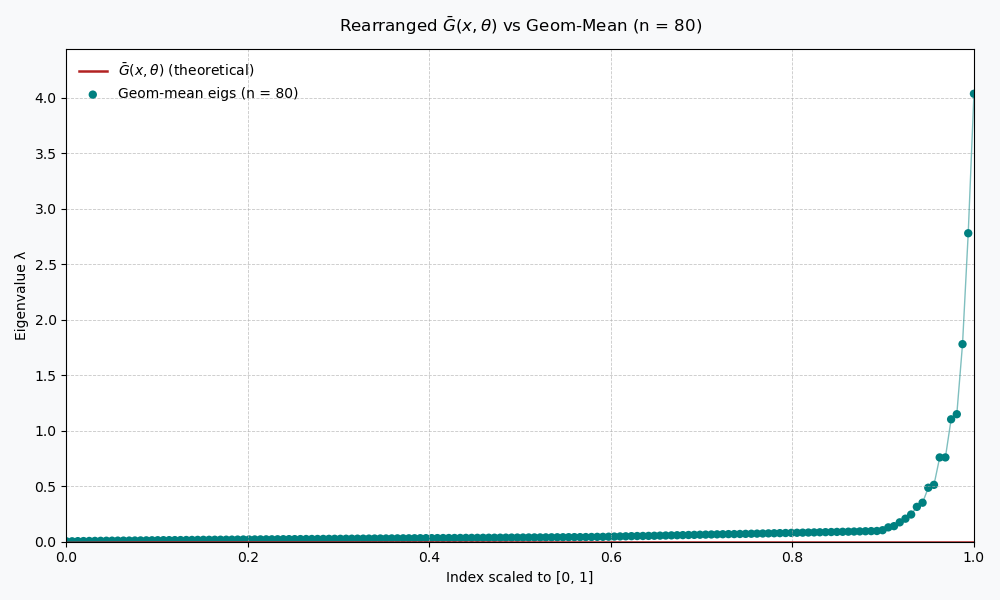}
	\end{subfigure}
	\\[-1em] 
	
	\begin{subfigure}[t]{0.49\textwidth}
		\centering
		\includegraphics[width=\textwidth]{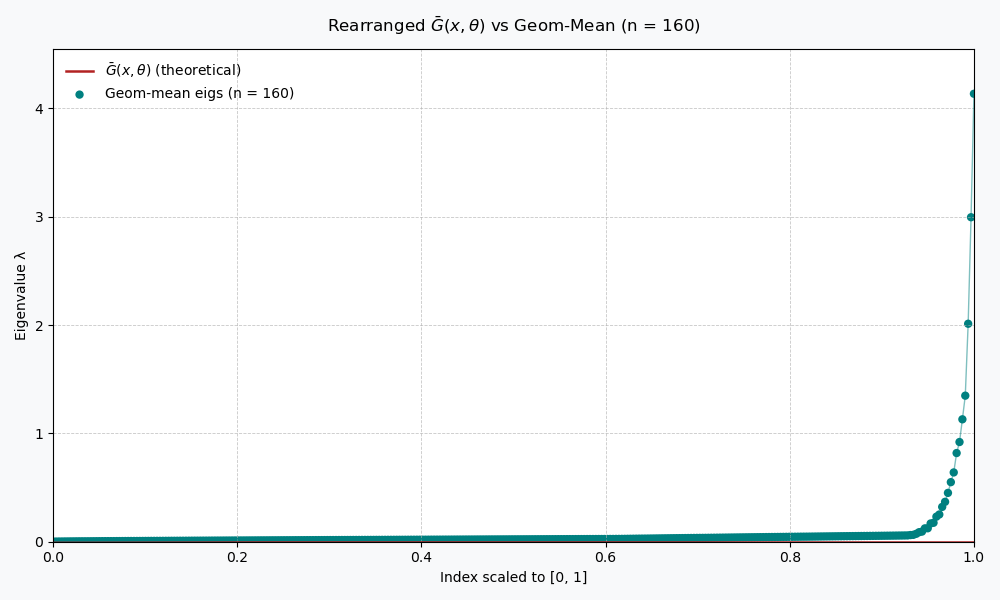}
	\end{subfigure}
	\hfill
	\begin{subfigure}[t]{0.49\textwidth}
		\centering
		\includegraphics[width=\textwidth]{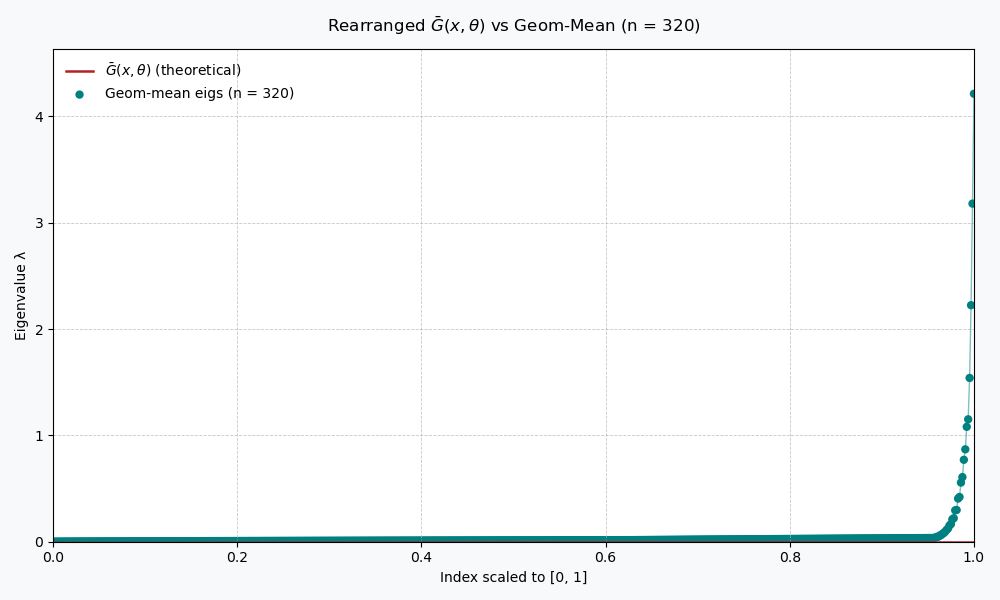}
	\end{subfigure}
	\caption{{\textbf{Case 1 - Example 1.}
			Sorted eigenvalues of \(G(A_n,B_n)\) (colored markers; \(n=40,80,160,320\))
			versus the rearranged distribution of the symbol
			\(\widetilde G(\kappa,\xi)\) (solid red line)}}.
	\label{fig:case1_ex1}
\end{figure}

\begin{figure}[H]
	\centering
	\begin{subfigure}[t]{0.49\textwidth}
		\centering
		\includegraphics[width=\textwidth]{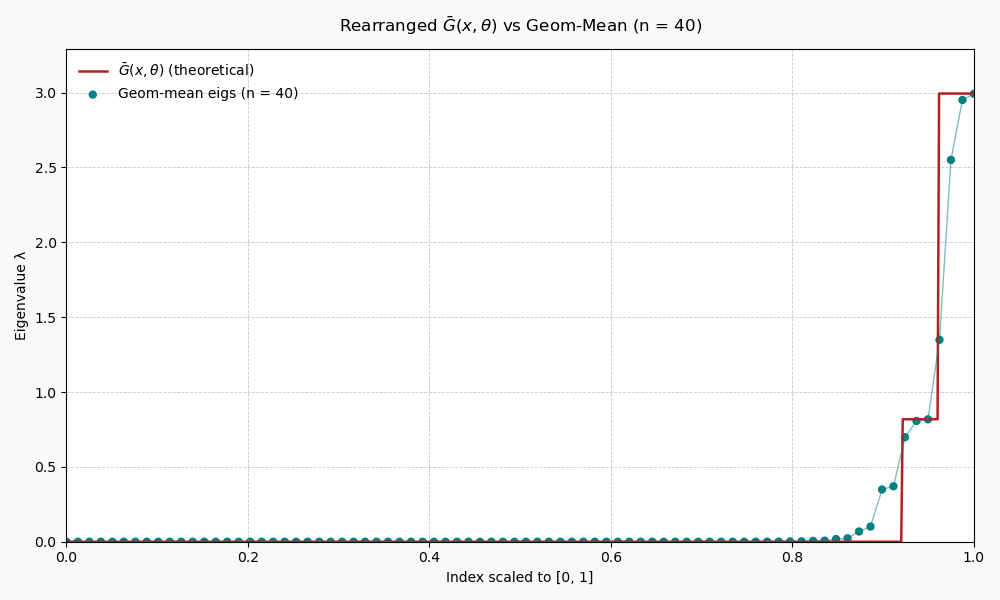}
	\end{subfigure}
	\hfill
	\begin{subfigure}[t]{0.49\textwidth}
		\centering
		\includegraphics[width=\textwidth]{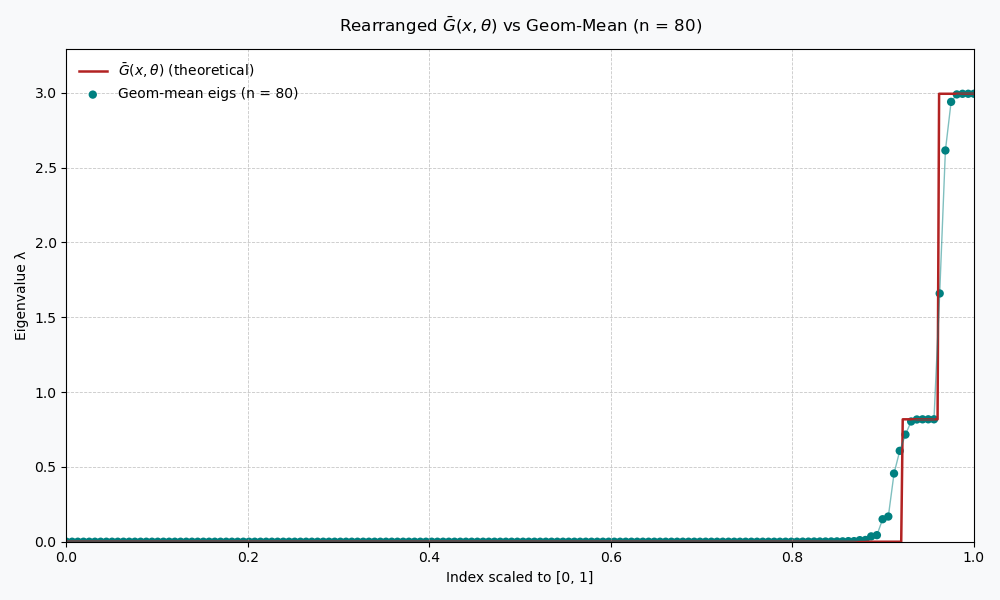}
	\end{subfigure}
	
	\begin{subfigure}[t]{0.49\textwidth}
		\centering
		\includegraphics[width=\textwidth]{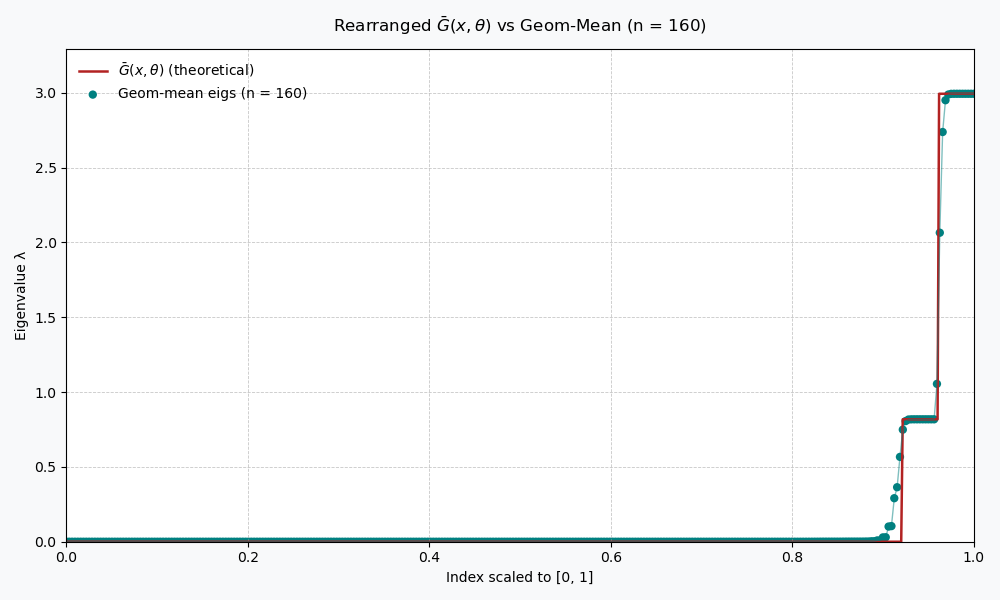}
	\end{subfigure}
	\hfill
	\begin{subfigure}[t]{0.49\textwidth}
		\centering
		\includegraphics[width=\textwidth]{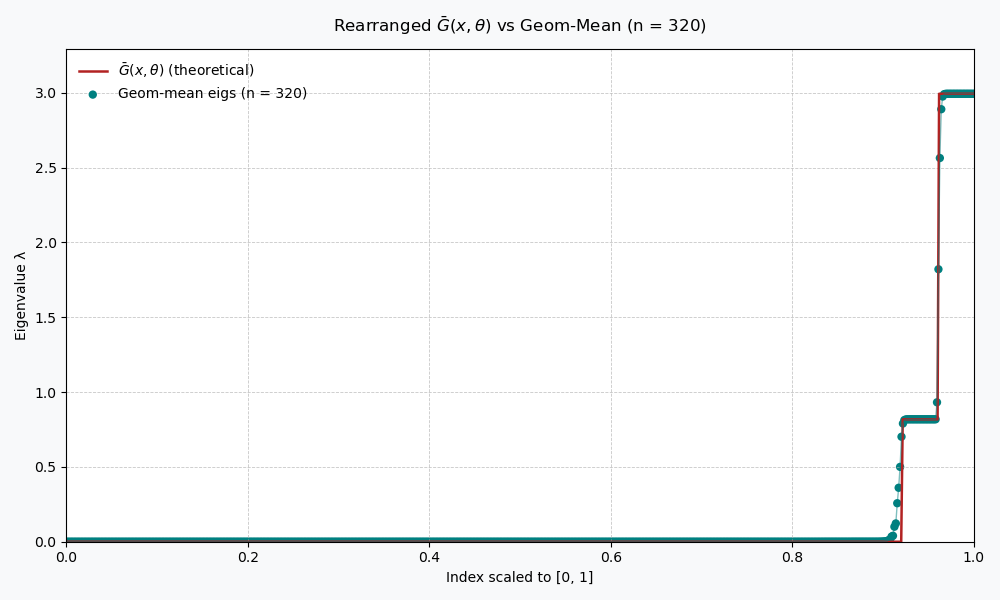}
	\end{subfigure}
	
	\caption{{\textbf{Case 1 - Example 2.}
			Comparison of the ordered eigenvalues of \(G(A_n,B_n)\)
			(colored markers) with the rearranged values of
			\(\widetilde G(\kappa,\xi)\) (solid red line) for
			\(n=40,80,160,320\).}}
	\label{fig:case1_ex2}
\end{figure}

In this section, we compare the sorted eigenvalues of each geometric mean matrix
\(G(A_n,B_n)\) with the rearranged eigenvalue distribution predicted by the
candidate symbol \(\widetilde{G}(\kappa,\xi)\). Plots are provided for each block size.
\paragraph{Case 1.}
Figures~\ref{fig:case1_ex1}--\ref{fig:case1_ex2} display the ordered
eigenvalues of \(G(A_n,B_n)\) (colored markers) together with the ordered
samples of \(\widetilde{G}(\kappa,\xi)\) (solid red line).

\smallskip
\emph{Example 1.}
Because the supports of \(\kappa\) and \(\xi\) are essentially disjoint,
\(\widetilde{G}(\kappa,\xi)\equiv 0\).
Most of the eigenvalues match the zero line, showing convergence to the zero-distribution. Only a small number of positive outliers remain above the line. For small values of $n$, the eigenvalue plot appears to be governed by a GLT {momentary} symbol that, when rearranged, follows a power law that rapidly collapses to the zero function as \(n \to \infty\).

\smallskip
\emph{Example 2.}
Here, the supports of \(\kappa\) and \(\xi\) intersect on \(\theta \in [-0.25,0.25]\). As a result, \(\widetilde{G}(\kappa,\xi)\) takes two values: zero outside the overlap and a positive constant matrix inside it. The eigenvalues reflect exactly this behavior: they form a long plateau at zero, followed by two shorter clusters at the positive eigenvalues determined by \(\widetilde{G}(\kappa,\xi)\). As we will see, the length of this zero plateau matches the respective Lebesgue measure of the complement of the support of \(\widetilde G\).

\FloatBarrier

\begin{figure}[H]
	\centering
	\begin{subfigure}[t]{0.49\textwidth}
		\centering
		\includegraphics[width=\textwidth]{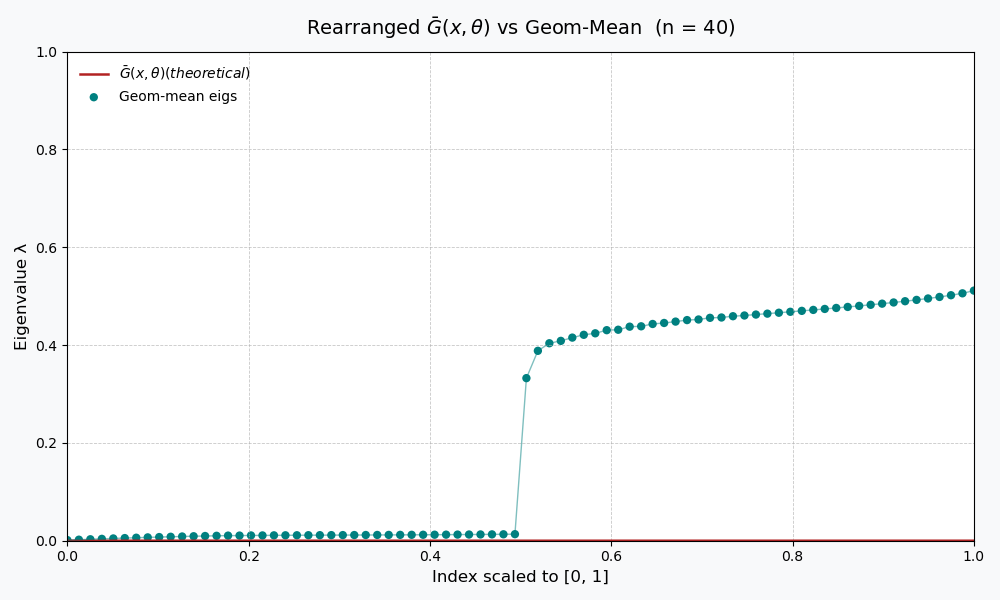}
	\end{subfigure}
	\hfill
	\begin{subfigure}[t]{0.49\textwidth}
		\centering
		\includegraphics[width=\textwidth]{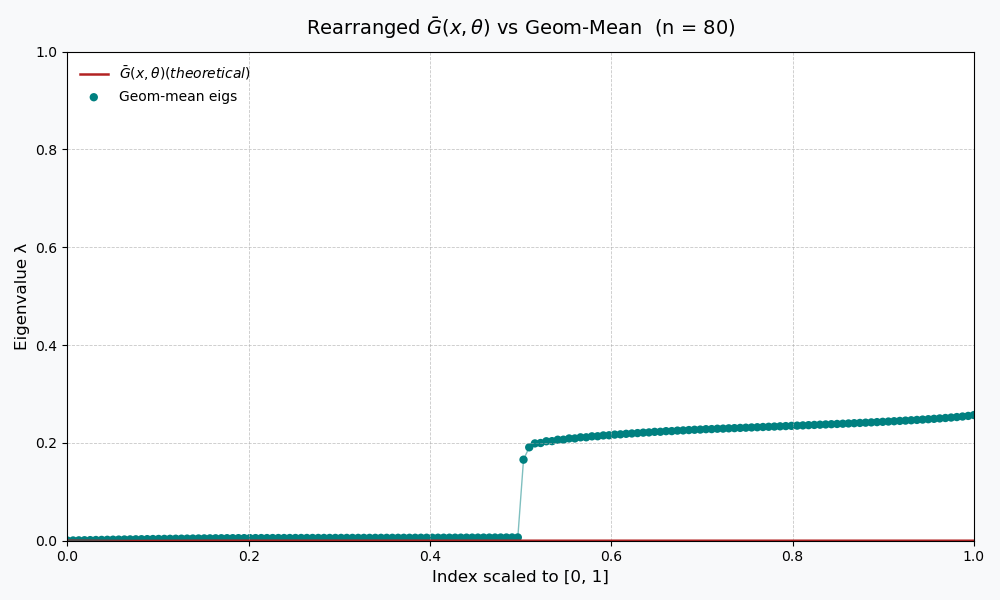}
	\end{subfigure}
	\\[-1em]
	
	\begin{subfigure}[t]{0.49\textwidth}
		\centering
		\includegraphics[width=\textwidth]{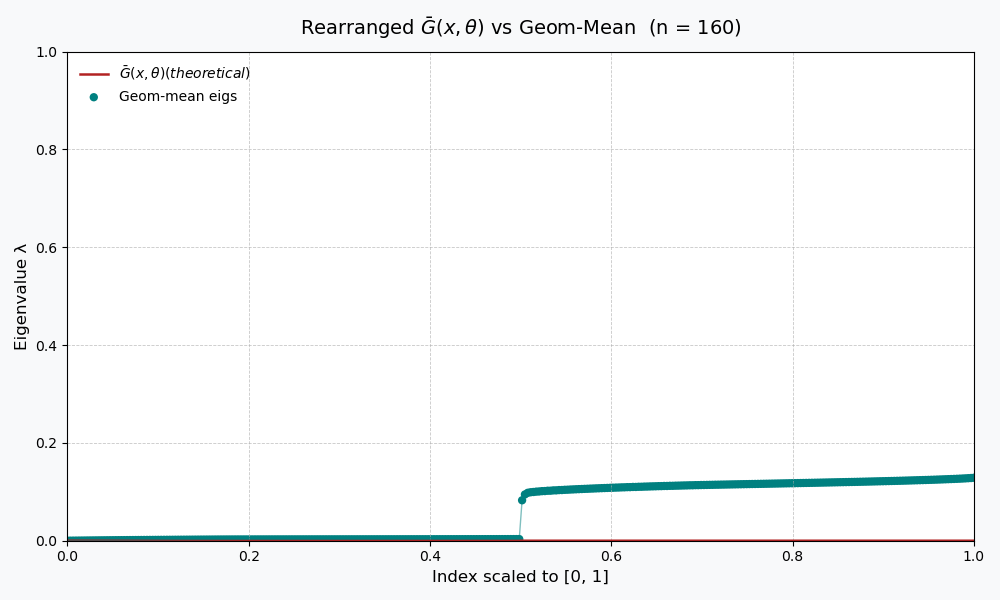}
	\end{subfigure}
	\hfill
	\begin{subfigure}[t]{0.49\textwidth}
		\centering
		\includegraphics[width=\textwidth]{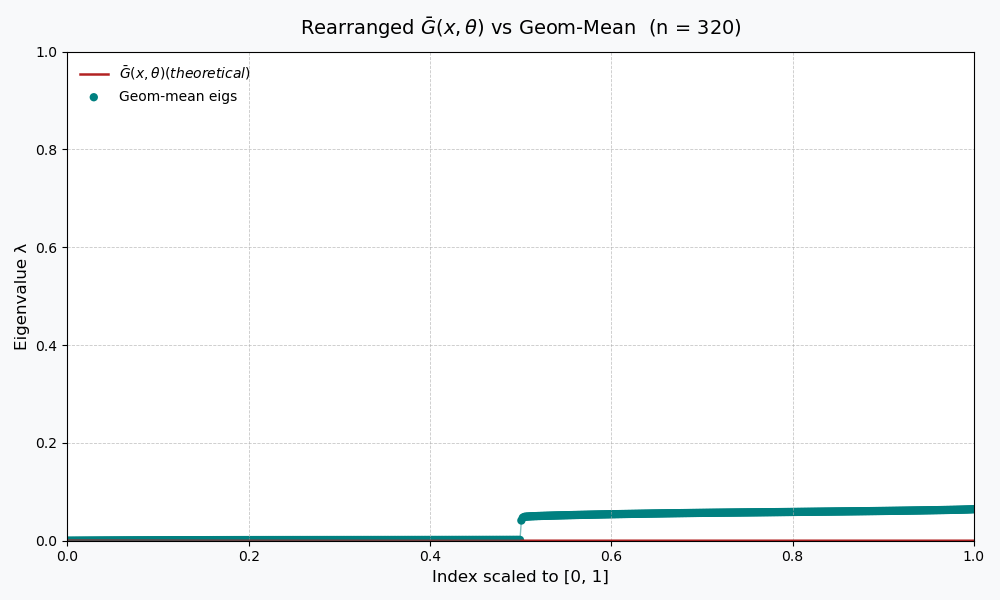}
	\end{subfigure}
	\caption{{\textbf{Case 2, Example 1.} 
			Sorted eigenvalues of \(G(A_n,B_n)\) (colored markers) versus the rearranged eigenvalue distribution predicted by \(G(\kappa,\xi)\) (solid red line), for
			\(n=40,80,160,320\).}}
	\label{fig:case2_ex1}
\end{figure}

\begin{figure}[H]
	\centering
	\begin{subfigure}[t]{0.49\textwidth}
		\centering
		\includegraphics[width=\textwidth]{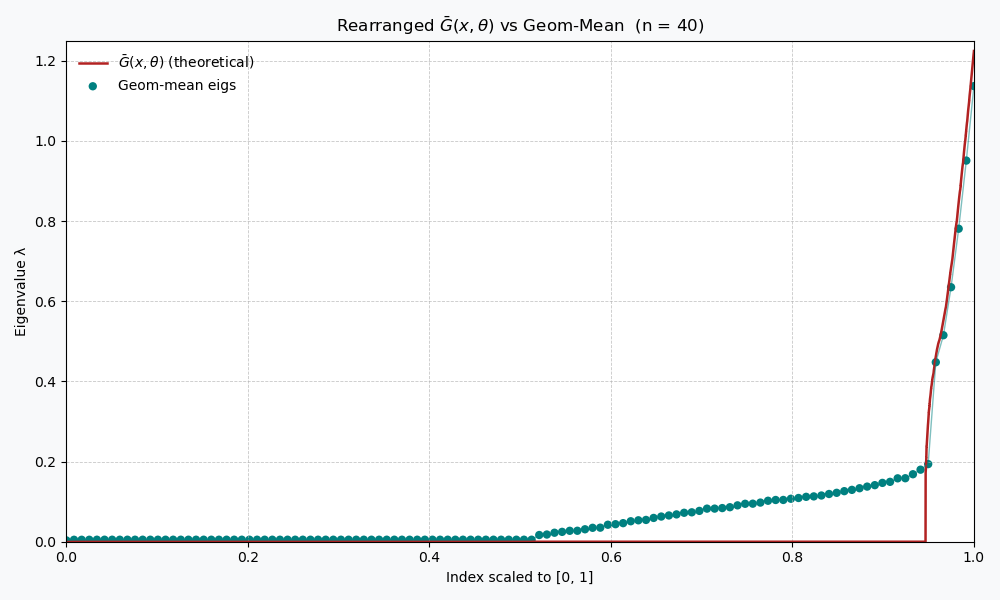}
	\end{subfigure}
	\hfill
	\begin{subfigure}[t]{0.49\textwidth}
		\centering
		\includegraphics[width=\textwidth]{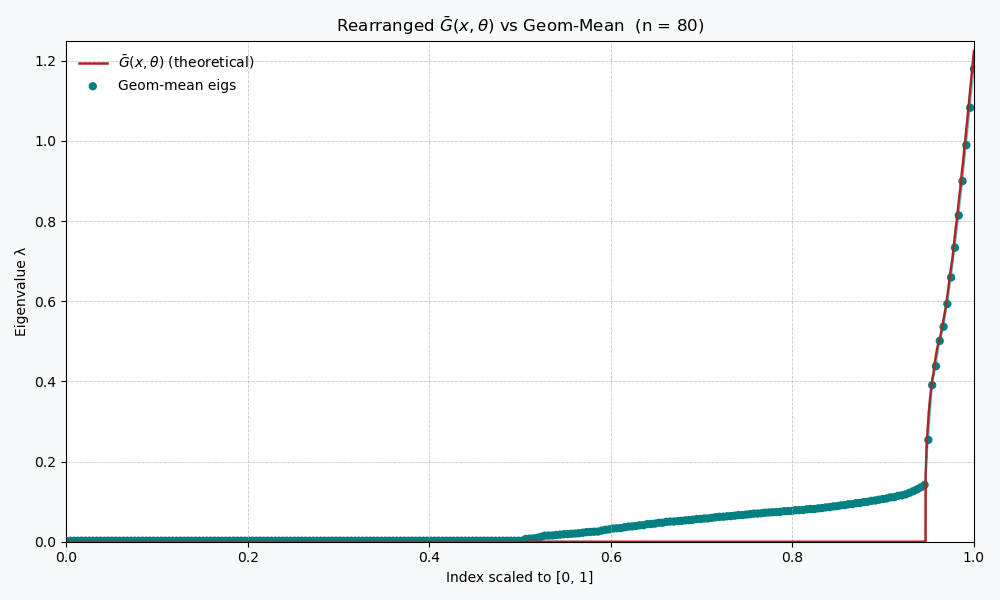}
	\end{subfigure}
	\\[-1em]
	
	\begin{subfigure}[t]{0.49\textwidth}
		\centering
		\includegraphics[width=\textwidth]{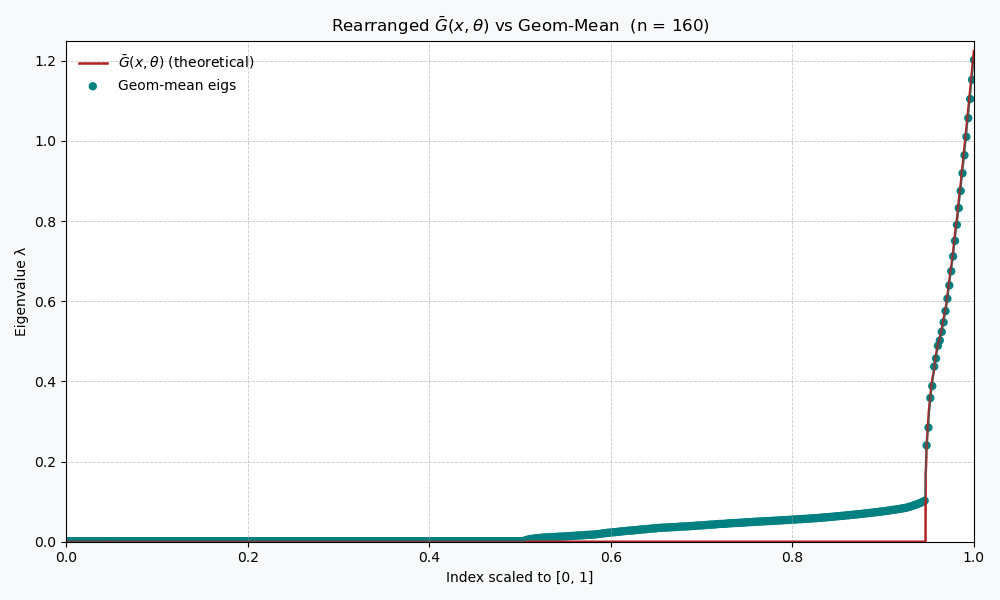}
	\end{subfigure}
	\hfill
	\begin{subfigure}[t]{0.49\textwidth}
		\centering
		\includegraphics[width=\textwidth]{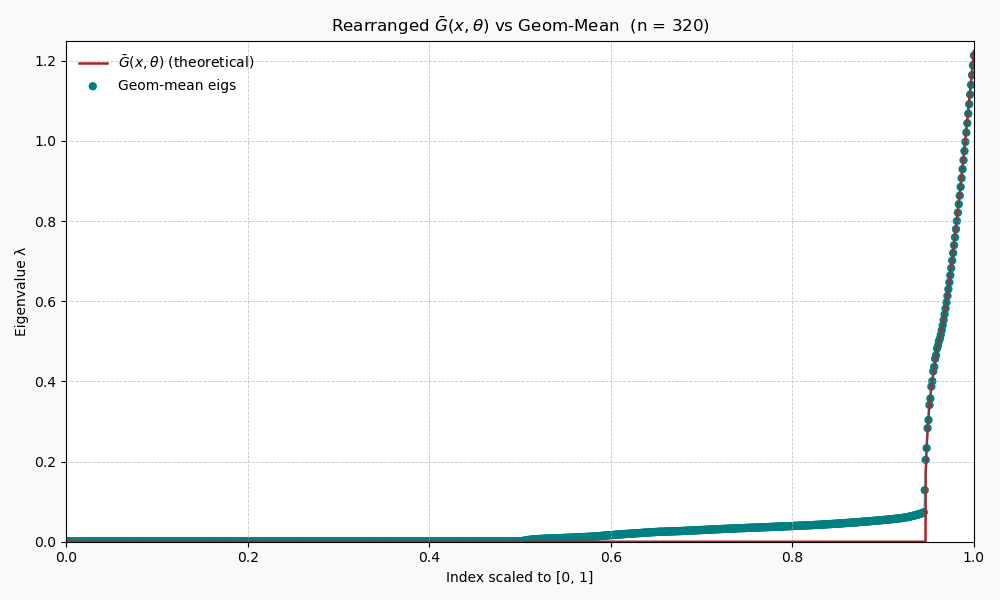}
	\end{subfigure}
	
	\caption{{\textbf{Case 2, Example 2.} Sorted eigenvalues of \(G(A_n,B_n)\) (colored markers) versus the rearranged eigenvalue distribution from \(G(\kappa,\xi)\) (solid red line), for
			\(n=40,80,160,320\).}}
	\label{fig:case2_ex2}
\end{figure}

\paragraph{Case 2.}
Figures~\ref{fig:case2_ex1}--\ref{fig:case2_ex2} repeat the experiment for
symbols that are rank-deficient on sets of full measure.

\smallskip
\emph{Example 1.}
Because the supports are disjoint, we have \(\widetilde{G}(\kappa,\xi)\equiv 0\).
The eigenvalues converge to the zero symbol. Despite the momentary symbol may seem to cause a rougher perturbation, we remark that, unlike in \textbf{Case 1 - Example 1}, this perturbation does not result in proper outliers since all perturbed eigenvalues converge to zero. In practice, we observe a better convergence.

\smallskip
\emph{Example 2.}
Here, the intersection set has rank one, leading to exactly one positive eigenvalue
in \(\widetilde{G}(\kappa,\xi)\).
The plot shows convergence despite the rank deficiency of the support. Because the norm of the diagonal perturbation introduced in \eqref{eq:case2_ex2} decays slowly, a momentary GLT perturbation \cite{new momentary,bolten2023note} is still observed, even in this non-zero case.

\FloatBarrier
\paragraph{Extremal eigenvalues.}
For each block size \(n\) we report
\[
\lambda_{\min}(G_n), \qquad
\lambda_{\max}(G_n), \qquad
\mathrm{cond}_2(G_n)=\frac{\lambda_{\max}(G_n)}{\lambda_{\min}(G_n)},
\]
together with the essential infimum and supremum of the conjectured
symbol, \(\widetilde G_{\min}\) and \(\widetilde G_{\max}\).
In this way, we can observe how rapidly the minimum eigenvalue decays and, in most cases, how the maximum eigenvalues converge to the essential supremum of the spectral symbol.
This provides a clear measure of the conditioning behaviour as \(n\) grows.

\medskip\noindent

\begin{table}[H]
	\centering
	\caption{Extremal eigenvalues and condition numbers for
		\textbf{Case 1 - Example 1} (symbol zero almost everywhere).}
	\label{tab:case1_ex1_range}
	\begin{tabular}{|c
			|S[table-format=1.4]
			|S[table-format=1.4]
			|S[table-format=1.4e+2]
			|S[table-format=1.8]
			|S[table-format=1.4e+2]|}
		\hline
		{$n$} & {\(\widetilde G_{\min}\)}
		& {\(\widetilde G_{\max}\)}
		& {Min.\ eig.}
		& {Max.\ eig.}
		& {\(\mathrm{cond}_2(G_n)\)}\\
		\hline
		40 & 0.0000 & 0.0000 & 8.7142e-03 & 3.91278029 & 4.490e+02 \\
		80 & 0.0000 & 0.0000 & 3.6849e-03 & 4.03656685 & 1.095e+03 \\
		160 & 0.0000 & 0.0000 & 1.4929e-03 & 4.13399129 & 2.769e+03 \\
		320 & 0.0000 & 0.0000 & 5.8335e-04 & 4.21189805 & 7.220e+03 \\
		\hline
	\end{tabular}
\end{table}
\paragraph{Case 1 - Example 2.}
Here \(\widetilde G_{\max}=2.9939\) is strictly positive. In this example, we observe a fast convergence of \(\lambda_{\max}(G_n)\) and an exponential decay of
\(\lambda_{\min}(G_n)\), which clearly causes the condition number \(\mathrm{cond}_2(G_n)\) to escalate, reaching
\(\mathcal{O}(10^{8})\) at \(n=320\).

\begin{table}[H]
	\centering
	\caption{Extremal eigenvalues and condition numbers for
		\textbf{Case 1 - Example 2} (symbol positive on an overlap set).}
	\label{tab:case1_ex2_range}
	\begin{tabular}{|c
			|S[table-format=1.4]
			|S[table-format=1.8]
			|S[table-format=1.4e+2]
			|S[table-format=1.8]
			|S[table-format=1.4e+2]|}
		\hline
		{$n$} & {\(\widetilde G_{\min}\)}
		& {\(\widetilde G_{\max}\)}
		& {Min.\ eig.}
		& {Max.\ eig.}
		& {\(\mathrm{cond}_2(G_n)\)}\\
		\hline
		40 & 0.0000 & 2.99393066 & 1.5625e-05 & 2.99257415 & 1.915e+05 \\
		80 & 0.0000 & 2.99393066 & 1.9531e-06 & 2.99393280 & 1.533e+06 \\
		160 & 0.0000 & 2.99393066 & 2.4414e-07 & 2.99393094 & 1.226e+07 \\
		320 & 0.0000 & 2.99393066 & 3.0518e-08 & 2.99393070 & 9.811e+07 \\
		\hline
	\end{tabular}
\end{table}
\paragraph{Case 2 - Example 1.}
The minimum eigenvalue appear to decay to 0 at a quadratic rate, while the largest eigenvalue converges to zero linearly, as determined by the momentary symbol. As a consequence, the growth of
\(\mathrm{cond}_2(G_n)\) is moderate.
\medskip
\smallskip

\begin{table}[H]
	\centering
	\caption{Extremal eigenvalues and condition numbers for
		\textbf{Case 2 - Example 1} (candidate symbol
		\(\widetilde G\equiv0\)).}
	\label{tab:case2_ex1_range}
	\begin{tabular}{|c
			|S[table-format=1.2]
			|S[table-format=1.2]
			|S[table-format=1.2e+2]
			|S[table-format=1.8]
			|S[table-format=1.2e+2]|}
		\hline
		{$n$} & {\(\widetilde G_{\min}\)}
		& {\(\widetilde G_{\max}\)}
		& {Min.\ eig.}
		& {Max.\ eig.}
		& {\(\mathrm{cond}_2(G_n)\)}\\
		\hline
		40 & 0.00 & 0.00 & 1.2615e-03 & 0.51134401 & 4.0536e+02 \\
		80 & 0.00 & 0.00 & 3.2103e-04 & 0.25691974 & 8.0029e+02 \\
		160 & 0.00 & 0.00 & 8.0993e-05 & 0.12877171 & 1.5899e+03 \\
		320 & 0.00 & 0.00 & 2.0342e-05 & 0.06447364 & 3.1695e+03 \\
		\hline
	\end{tabular}
\end{table}
\medskip
\smallskip

\paragraph{Case 2 - Example 2.}
Here \(\widetilde G_{\max}\simeq1.22\).
\(\lambda_{\max}(G_n)\) converges to this value and is included in the range of $\widetilde G$, whereas
\(\lambda_{\min}(G_n)\) decays linearly. This better conditioning compared to \textbf{Example~2} of the first case is reasonable, since the positive perturbations we introduced to the Toeplitz matrices of this example (Equation \ref{eq:case2_ex2} decay in norm more slowly than the ones in \textbf{Case 1 - Example 1}).

\medskip

\begin{table}[H]
	\centering
	\caption{Extremal eigenvalues and condition numbers for
		\textbf{Case 2 - Example 2} (symbol positive on an overlap set).}
	\label{tab:case2_ex2_range}
	\begin{tabular}{|c
			|S[table-format=1.2]
			|S[table-format=1.8]
			|S[table-format=1.2e+2]
			|S[table-format=1.8]
			|S[table-format=1.2e+2]|}
		\hline
		{$n$} & {\(\widetilde G_{\min}\)}
		& {\(\widetilde G_{\max}\)}
		& {Min.\ eig.}
		& {Max.\ eig.}
		& {\(\mathrm{cond}_2(G_n)\)}\\
		\hline
		40 & 0.00 & 1.22418242 & 3.24e-03 & 1.13658304 & 3.50e+02 \\
		80 & 0.00 & 1.22418242 & 1.46e-03 & 1.17911634 & 8.10e+02 \\
		160 & 0.00 & 1.22418242 & 5.65e-04 & 1.20156562 & 2.13e+03 \\
		320 & 0.00 & 1.22418242 & 8.19e-05 & 1.21306699 & 1.48e+04 \\
		\hline
	\end{tabular}
\end{table}

\paragraph{Zero-related statistics.}
For each block size \(n\), we report the fraction of eigenvalues whose
modulus does not exceed the fixed cutoff \(0.1\); see
Tables~\ref{tab:case1_ex1_zero}-\ref{tab:case2_ex2_zero}.
Each table lists the empirical proportion
(\emph{Prop.\(\le0.1\)}), the theoretical measure of the zero set of the
rearranged symbol (\emph{Target}), and the absolute difference (\emph{Error}).

\smallskip
\paragraph{Case 1 - Example 1.}
The conjectured symbol is the zero function, so the target measure is~1. The empirical value increases monotonically with~\(n\), remaining  within a small error.
Because the number of outliers is at most \(o(n)\), we expect that the proportion will converge to 1 as $n \to \infty $.
\paragraph{Case 1 - Example 2.}
Here, the zero set of the rearranged spectral symbol has a measure
\[
1-\frac{1}{4\pi}\approx 0.9204.
\]%
The measured proportions appear to approach this value with an
\(\mathcal O(n^{-1})\) rate of decay, with the error falling below \(10^{-3}\) by \(n=160\).

\medskip
\paragraph{Case 2 - Example 1.}
Here the expected spectral symbol is zero, but the momentary symbol perturbation,
as observed in Figure~\ref{fig:case2_ex1}, produces an apparent
plateau \(0.5\) up to \(n=320\), where the proportion jumps to~1 and the measured error vanishes.
Since the momentary perturbation does not create proper outliers, as observed in \textbf{Case 1 - Example 2}, a faster convergence can reasonably be expected in this case.

\paragraph{Case 2 - Example 2.}
In this case, the complement of the support of the rearranged distribution has measure
\[
\frac{17}{18}\approx 0.9444.
\]
The empirical measure of the zero cluster converges to this value, with the error appearing to follow the same power law observed in \textbf{Case 2 - Example 1}.

\medskip
\begin{table}[H]
	\centering
\noindent	\caption{Case 1, Example 1: proportion of eigenvalues below \(0.1\); target \(1\).}
	\label{tab:case1_ex1_zero}
	\begin{tabular}{|c|S|S|S|}
		\hline
		{$n$} & {Prop.\(\le0.1\)} & {Target} & {Error} \\
		\hline
		40 & 0.6375 & 1.0000 & 0.3625 \\
		80 & 0.8938 & 1.0000 & 0.1062 \\
		160 & 0.9438 & 1.0000 & 0.0562 \\
		320 & 0.9688 & 1.0000 & 0.0312 \\
		\hline
	\end{tabular}
\end{table}
\begin{table}[H]
	\centering
	\caption{Case 1, Example 2: proportion of eigenvalues below \(0.1\);
		target \(1-\frac{1}{4\pi}\approx0.92042\).}
	\label{tab:case1_ex2_zero}
	\begin{tabular}{|c|S|S|S|}
		\hline
		{$n$} & {Prop.\(\le0.1\)} & {Target} & {Error} \\
		\hline
		40 & 0.8750 & 0.9204 & 0.0454 \\
		80 & 0.8938 & 0.9204 & 0.0266 \\
		160 & 0.9031 & 0.9204 & 0.0173 \\
		320 & 0.9109 & 0.9204 & 0.0095 \\
		\hline
	\end{tabular}
\end{table}
\begin{table}[H]
	\centering
\noindent \caption{Case 2, Example 1: proportion of eigenvalues below \(0.1\); target \(1\).}
	\label{tab:case2_ex1_zero}
	\begin{tabular}{|c|S|S|S|}
		\hline
		{$n$} & {Prop.\(\le0.1\)} & {Target} & {Error} \\
		\hline
		40 & 0.5000 & 1.0000 & 0.5000 \\
		80 & 0.5000 & 1.0000 & 0.5000 \\
		160 & 0.5156 & 1.0000 & 0.4844 \\
		320 & 1.0000 & 1.0000 & 0.0000 \\
		\hline
	\end{tabular}
\end{table}
\begin{table}[H]
	\centering
	\caption{Case 2, Example 2: proportion of eigenvalues below \(0.1\);
		target \(\tfrac{17}{18}\approx0.9444\).}
	\label{tab:case2_ex2_zero}
	\begin{tabular}{|c|S|S|S|}
		\hline
		{$n$} & {Prop.\(\le0.1\)} & {Target} & {Error} \\
		\hline
		40 & 0.7667 & 0.9444 & 0.1777 \\
		80 & 0.8833 & 0.9444 & 0.0611 \\
		160 & 0.9438 & 0.9444 & 0.0006 \\
		320 & 0.9448 & 0.9444 & 0.0004 \\
		\hline
	\end{tabular}
\end{table}

	\chapter*{Conclusions and future work} 
\vspace{0.2cm}

Throughout the dissertation, we analyzed various structured matrix sequences arising from the discretization of multi-dimensional fractional partial differential equations. After introducing the analytical foundations of fractional operators, the spectral framework of matrix sequences, the main classes of structured matrices, and the GLT theory in Part I to provide a foundation for the numerical applications, in Part II we discretized several types of equations and addressed both direct and inverse problems. The associated linear systems are generally computationally challenging due to the dense and ill-conditioned nature of the coefficient matrices. However, we showed that the rich framework of GLT theory can effectively support the development of fast solvers, focusing in particular on preconditioning strategies. This final chapter summarizes our findings, highlights key similarities and differences across the studied scenarios, and suggests possible directions for future research.

\subsection*{Summary of research contributions}

In Chapter \ref{chapter:TFDE}, we studied ISPs for a two-dimensional generalized diffusion equation involving a Caputo-type time-fractional derivative with a general memory kernel and nonlocal boundary conditions. Using the spectral decomposition of the associated elliptic operator, we derived explicit series representations of the solutions in terms of bi-orthogonal systems of eigenfunctions. We established rigorous existence and uniqueness results for regular solutions and demonstrated that both ISPs are ill-posed in the Hadamard sense. Several special cases, including tempered, multi-term, mono-term, and classical diffusion models, were obtained as particular choices of the memory kernel, illustrating the unifying nature of the generalized formulation. An interesting direction for future research would be to combine the present theoretical study with the computational features of structured matrices and the asymptotic spectral theory for matrix sequences, in the spirit of the works \cite{comput1,comput2,comput3}. The design of fast regularization algorithms for the ISPs examined in this chapter presents a significant challenge.
\vspace{0.1cm}

In Chapter \ref{chapter:IFL}, we addressed an ISP for a space-time FDE involving a tempered fractional derivative in time and a Caputo fractional derivative in space. The ill-posed ISP was regularized using the quasi-boundary value method. The finite difference discretization led to a large-scale two-by-two block linear system with structured coefficients, for which we developed an efficient preconditioning strategy based on the block triangular structure. A detailed spectral analysis of the associated matrix sequences was conducted using the GLT framework, which provided a rigorous description of their asymptotic spectral behavior despite the presence of the variable coefficient that spoils the pure Toeplitz structure. Numerical experiments confirmed the stability, accuracy, and computational efficiency of the proposed approach in reconstructing the source term from noisy data. For future work, we plan to construct a more accurate preconditioner that accounts for the remaining outliers and possibly leads to a strong eigenvalue cluster at 1. In fact, even in presence of negative findings concerning matrix algebra preconditioners for multilevel Toeplitz systems, there exist cases especially in the fractional setting where optimality can still be reached \cite{against barrier2,against barrier1}. As a further step, a theoretical and practical comparison of our technique with direct-like solvers stemming either from block Gaussian elimination or from Sherman-Morrison-Woodbury methods would be valuable.
\vspace{0.1cm}

In Chapter \ref{chapter:RFDE}, we considered an ISP for a time-space FDE with variable coefficients, extending the previous work from the uni-dimensional to the two-dimensional setting. A comprehensive spectral analysis of the coefficient matrix sequences was conducted using GLT theory, which guided the construction of an effective preconditioner. We introduced a block triangular preconditioner and proved that the preconditioned matrix sequence was weakly clustered at one in the GLT, singular value and eigenvalue sense, ensuring fast convergence of the preconditioned GMRES method. Numerical experiments demonstrated the accuracy of the source reconstruction and the efficiency of the proposed preconditioning strategy across various fractional orders and problem sizes. The convergence behavior depends mainly on the spatial discretization, with limited sensitivity to temporal refinement, which is reasonable considering that the number of outliers in the cluster at one of the preconditioned matrix depends only on the spatial dimension. Future work will address extensions to more complex geometries, high-order finite elements or isogeometric analysis with intermediate regularity, and higher dimensions. Additionally, taking inspiration from the variable-coefficient approach in \cite{Toe-diff} and from the \(\omega\)-circulant approach \cite{Diaz}, other options for the choice of the preconditioners should be considered.
\vspace{0.1cm}

In Chapter \ref{chapter:space-FDE}, we studied the spectral properties and preconditioning of linear systems arising from the finite difference discretization of time-dependent space-FDEs with variable coefficients. The model involved a one-sided Riemann-Liouville fractional derivative multiplied by a variable coefficient, discretized using the shifted Grünwald formula in space and the \(\theta\)-method in time. The resulting all-at-once systems exhibited a \((d+1)\)-level Toeplitz-like matrix structure. Using GLT theory, we provided a detailed analysis of the spectral distribution of both the original and preconditioned matrix sequences and proposed structure-preserving preconditioners based on circulant approximations. The main novelty was that the analysis fully covered the nonconstant variable coefficient case, both in one and two spatial dimensions, extending several existing results in the literature which were often restricted to more simplified settings. Numerical experiments confirmed the theoretical predictions and demonstrated significant reductions in GMRES iteration counts and computational cost. Several open problems remain. The extension of the eigenvalue distribution results for non-Hermitian Toeplitz matrices to the variable coefficient case seems not trivial and is not covered by the GLT theory due to the intrinsic non-Hermitian nature of the underlying matrix sequences. Additionally, estimates of the number of outliers and of the asymptotic behavior of the outliers in the preconditioned matrix sequences would have precise practical impact in estimating a priori the number of needed iterations of the proposed preconditioned GMRES procedures for reaching a preassigned accuracy.
\vspace{0.1cm}

In Chapter \ref{GM GLT2}, we addressed a fundamental open question concerning the spectral distribution of geometric means of Hermitian positive definite GLT matrix sequences. We proved that when the GLT symbols commute, the invertibility condition on the symbols could be completely removed, resolving a long-standing conjecture. This result was established both in the scalar setting for any dimension and in the general block setting under the assumption that the symbols commute. We further demonstrated that the result is maximal: when the symbols do not commute and are simultaneously degenerate, the classical relation may fail and the geometric mean sequence may not coincide with the classical geometric mean of the symbols. Numerical experiments revealed new spectral phenomena explained by the theory of momentary GLT symbols, including secondary eigenvalue distributions and extremal spectral behavior connected to the order of vanishing of the associated GLT symbol. The numerical evidences open the door to further theoretical studies. Finally, the study of the spectral distribution when we consider the geometric mean of more than two matrix sequences is completely open and requires normwise error estimates for using the relevant tools.

\subsection*{Future research directions}

The same theoretical apparatus, including GLT theory, approximating classes of sequences, and non-Hermitian perturbations of Hermitian matrix sequences, could be applied to even more complex cases. For instance, high-order finite element approximations in \(d\) spatial dimensions lead to \(d\)-level Toeplitz matrices with matrix-valued generating functions, corresponding to the GLT symbol \cite{Beckermann2007}. Using similar approaches, we could also consider discontinuous Galerkin methods or isogeometric analysis. These cases lead to matrix-valued symbols in \(2d\) independent variables, as in the general GLT theory already available \cite{Barbarino2022,barbarino2020blockmulti,barbarino2020block1d,garoni2017,garoni2018}. Finally, as a future task in this research direction, we will consider the design of optimized code in order to match the theoretical complexity derived in the present work for the preconditioning strategies proposed in the research chapters. The implementation of efficient algorithms that fully exploit the structure of the preconditioners, combined with careful optimization and parallelization, would further enhance the practical applicability of our approach.
	
	\addcontentsline{toc}{chapter}{Conclusions and future work}
	
	\addcontentsline{toc}{chapter}{References}
\end{document}